%% file: main.tex
\documentclass[3p]{elsarticle}

\makeatletter
\def\ps@pprintTitle{%
    \let\@oddhead\@empty
    \let\@evenhead\@empty
    \let\@evenfoot\@oddfoot
    }
\makeatother
\input{package}

\newcommand{\cmark}{\text{\ding{51}}}
\newcommand{\xmark}{\text{\ding{55}}}

\usepackage{scalerel}
\newcommand\vfrac[2]{\ThisStyle{%
  \setbox0=\hbox{$\SavedStyle#1#2$}%
  \setbox2=\hbox{$\SavedStyle X$}%
  \ifdim\ht0>\ht2\setlength{\ht0}{\ht2}\fi%
  #1\mathord{\stretchto{\raisebox{2.3\LMpt}{$\SavedStyle/$}}{\ht0}}#2}}
  
\newcommand\Gaussian{\mathcal{N}}
\newcommand\PriorCov{\mathcal{C}}
\newcommand\transpose{\top}
\newcommand\identity{\text{\normalfont\ignorespaces Id}}
\newcommand\Forward{\boldsymbol{\mathcal{G}} }
\newcommand\FullRidgeForward{\widetilde{{\Forward}}}
\newcommand{\FullRidgeForwardOpt}{{\FullRidgeForward}_{\text{\normalfont\ignorespaces opt}}}
\newcommand{\RidgeForward}{\boldsymbol{\mathsf{g}}}
\newcommand{\RidgeForwardOpt}{{\boldsymbol{\mathsf{g}}}_{\text{\normalfont\ignorespaces opt}}}

\newcommand{\RidgeMisfitOpt}{{\RidgeMisfit}_{\text{\normalfont\ignorespaces opt}}}

\newcommand{\RidgeForwardOptApprox}{\boldsymbol{g}}

\newcommand\FullParam{m}

\newcommand\ParamSpace{\scrM}
\newcommand\Obs{\by}
\newcommand\noise{\bn}
\newcommand\NoiseCov{\Gamma_{n}}

\newcommand\ObsDim{{d_{\Obs}}}

\newcommand\RedCoordParam{\boldsymbol{x}}
\newcommand\refCoord{\boldsymbol{z}}

\newcommand\RedCoordDim{{d_{r}}}
\newcommand\RedCoordParamSpace{\mathbb{R}^{{\RedCoordDim}}}

\newcommand{\ObsSpace}{\mathbb{R}^{\ObsDim}}

\newcommand\FullPrior{\mu}
\newcommand\FullPost{\mu^{\Obs}}

\newcommand\FullPostOpt{\widetilde{\mu}^{\Obs}_{\text{\normalfont\ignorespaces opt}}}

\newcommand\Misfit{\Phi^{\Obs}}

\newcommand\RidgeMisfit{\widetilde{\Phi}^{\Obs}}

\newcommand{\Tr}{\text{\normalfont\ignorespaces Tr}}

\newcommand\Projector{\mathcal{P}}

\newcommand\RedCoordPrior{\pi}

\newcommand\RedCoordPostOpt{\widetilde{\pi}^{\Obs}_{\text{\normalfont\ignorespaces opt}}}

\newcommand\Encoder{\mathcal{E}_{r}}
\newcommand\Decoder{\mathcal{D}_{r}}

\newcommand\FullMap{\mathcal{T}}

\newcommand\Map{\boldsymbol{\mathsf{T}}}
\newcommand\MapClass{\mathscr{T}}
\newcommand\MapParams{{\boldsymbol\theta}}

\newcommand\normalization{Z}

\newcommand\KL{\mathcal{D}_{\text{\normalfont\ignorespaces KL}}}

\newcommand\ParamCM{{H_{\PriorCov}}}
\newcommand\ObsCM{{H_{\NoiseCov}}}

\newcommand\ProbSpace{\mathscr{P}}

\newcommand{\iid}{\stackrel{\text{\normalfont\ignorespaces i.i.d.}}{\sim}}
\newcommand\HS{\text{\normalfont\ignorespaces HS}}
\newcommand\TrClass{B_1^{+}}
\newcommand\meas{\text{\normalfont\ignorespaces d}}

\newcommand{\norm}[1]{\left\lVert #1 \right\rVert}
\newcommand\StateSpace{\mathscr{U}}
\newcommand\SolOp{\mathcal{F}}
\newcommand\ObsOp{\boldsymbol{\mathcal{O}}}

\newcolumntype{C}{ >{\centering\arraybackslash} m{0.5\textwidth} }
\newcolumntype{D}{ >{\centering\arraybackslash}m{0.05\textwidth}}
\newcolumntype{E}{ >{\centering\arraybackslash}m{0.20\textwidth}}
\newcolumntype{F}{ >{\centering\arraybackslash}m{0.23\textwidth}}
\newcolumntype{G}{ >{\centering\arraybackslash}m{0.16\textwidth}}
\newcolumntype{H}{ >{\centering\arraybackslash}m{0.08\textwidth}}
\newcolumntype{I}{ >{\centering\arraybackslash}m{0.1\textwidth}}

\definecolor{darkgreen}{rgb}{.15,.55,0}

\newcommand\wye[1][]{%
    \tikz\draw[thick, line cap=round,x=1ex,y=1ex,#1]
    (0,0) -- ++(90:1)
    (0,0) -- ++(-30:1)
    (0,0) -- ++(-150:1);
}

\crefname{assumption}{Assumption}{Assumptions}
\crefname{algorithm}{Algorithm}{Algorithms}
\crefname{definition}{Definition}{Definitions}
\crefname{remark}{Remark}{Remarks}
\crefname{equation}{}{}
\crefname{appendix}{}{}
\crefname{lemma}{Lemma}{Lemmas}
\crefname{proposition}{Proposition}{Propositions}
\crefname{corollary}{Corollary}{Corollaries}
\crefname{theorem}{Theorem}{Theorems}
\crefname{table}{Table}{Tables}
\crefname{section}{Section}{Sections}
\crefname{appendix}{}{}

\begin{document}
\begin{frontmatter}
\title{LazyDINO: Fast, Scalable, and Efficiently Amortized Bayesian Inversion via Structure-Exploiting and Surrogate-Driven Measure Transport}
\author[2]{Lianghao Cao\corref{eq}}
\ead{lianghao@caltech.edu}
\author[3]{Joshua Chen\corref{eq}}
\ead{joshuawchen@utexas.edu}
\author[1]{Michael Brennan}
\ead{mcbrenn@mit.edu}
\author[3]{Thomas O'Leary-Roseberry}
\ead{tom.olearyroseberry@utexas.edu}
\author[1]{Youssef Marzouk}
\ead{ymarz@mit.edu}
\author[3,4]{Omar Ghattas}
\ead{omar@oden.utexas.edu}

\address[1]{Center for Computational Science and Engineering, Massachusetts Institute of Technology, Cambridge, MA 02139, USA}

\address[2]{Department of Computing and Mathematical Sciences, California Institute of Technology
Pasadena, CA 91125, USA}

\address[3]{Oden Institute for Computational Engineering and Sciences, The University of Texas at Austin, 201 E. 24th Street, C0200, Austin, TX 78712, USA
}

\address[4]{Walker Department of Mechanical Engineering, The University of Texas at Austin, 204 E. Dean Keaton Street, Austin, TX 78712, USA
}

\cortext[eq]{Co-first author, both are corresponding authors.}

\begin{abstract}
\input{abstract}
\end{abstract}

\begin{keyword}
Bayesian inverse problem, variational inference, measure transport, surrogate model, dimension reduction, derivative-informed operator learning
\end{keyword}

\end{frontmatter}

\input{introduction}

\input{bips}

\input{transport}

\input{dino_approximated_likelihood}

\input{lazydino}

\input{results}

\input{conclusion}

\input{acknowledgment}

\addcontentsline{toc}{section}{References}
\bibliographystyle{model1-num-names}
\biboptions{sort,numbers,comma,compress}                 
\bibliography{main.bib}

\appendix

\input{appendix_glossary}
\input{appendix_inf_dimensions}
\input{appendix_proof}

\input{lazydino_detailed_appendix}
\input{appendix_pde_bayes}

\input{appendix_additional_numerics}

\input{map-diagnostics}
\input{appendix_numerical_results}

\end{document}

%% file: package.tex
\usepackage{tabularx}
\usepackage{amsmath}
\usepackage{amssymb}
\usepackage{amsthm}
\usepackage{mathrsfs}
\usepackage{bm}
\usepackage[dvipsnames]{xcolor}
\usepackage[final]{hyperref}
\usepackage[ruled,vlined]{algorithm2e}
\usepackage{enumitem}
\usepackage{array}
\usepackage{multirow}
\usepackage{makecell}
\usepackage[title]{appendix}
\usepackage{empheq}
\usepackage{cleveref}
\usepackage{cases}
\usepackage{textcomp}
\usepackage{gensymb}
\usepackage{subcaption}
\hypersetup{
	colorlinks=true,       
	linkcolor=blue,        
	citecolor=blue,        
	filecolor=magenta,     
	urlcolor=blue         
}
\usepackage{float}
\usepackage{multirow}

\usepackage[colorinlistoftodos,prependcaption,textsize=tiny]{todonotes}
\usepackage{regexpatch}

\usepackage{tikz}
\usetikzlibrary{calc,patterns,decorations.pathreplacing,calligraphy}
\usetikzlibrary{shapes,positioning}
\tikzstyle{bag} = [align=center]

\usepackage{graphicx} 
\usepackage{amsmath,amsfonts, mathtools}
\usepackage{tabularx}
\usepackage{amssymb}
\usepackage{amsthm}
\usepackage{mathrsfs}
\usepackage{bm}
\usepackage{graphicx}
\usepackage{float}
\usepackage{enumitem}
\usepackage{array}
\usepackage{multirow}
\usepackage{makecell}
\usepackage{empheq}
\usepackage{cases}
\usepackage{textcomp}
\usepackage{gensymb}
\usepackage{subcaption}
\usepackage{hyperref}
\usepackage{multicol}
\hypersetup{
	colorlinks=true,       
	linkcolor=blue,        
	citecolor=blue,        
	filecolor=magenta,     
	urlcolor=blue    
}
\usepackage{cleveref}
\usepackage{lineno}
\usepackage{hhline}
\usepackage{dsfont}

\include{command_lc}

\usepackage{soul}
\usepackage{booktabs}
\usepackage{diagbox}

\usepackage{algpseudocode}
\usepackage{amsmath}
\usepackage{amssymb}

\usepackage{mathtools}
\usepackage{pifont}
\usepackage{hhline}

\usepackage{changepage}
\usepackage{placeins}

%% file: command_lc.tex
\def\Dkl{\mathcal{D}_{\textrm{KL}}}

\newcommand{\R}{\mathbb{R}}

\newcommand{\grad}{\nabla}

\newtheorem{remark}{Remark}
\newcommand{\bdmc}[1]{\boldsymbol{\mathcal{#1}}}

\newcommand{\dd}{\,\textrm{d}}

\newcommand{\bzero}{\bm{0}}

\newcommand{\be}{\bm{e}}

\newcommand{\bg}{\bm{g}}

\newcommand{\bn}{\bm{n}}

\newcommand{\bs}{\bm{s}}

\newcommand{\bu}{\bm{u}}

\newcommand{\bw}{\bm{w}}
\newcommand{\bx}{\bm{x}}
\newcommand{\by}{\bm{y}}

\newcommand{\bA}{\bm{A}}
\newcommand{\bB}{\bm{B}}

\newcommand{\bF}{\bm{F}}

\newcommand{\bJ}{\bm{J}}

\newcommand{\bP}{\bm{P}}

\newcommand{\bV}{\bm{V}}

\newcommand{\calC}{{\mathcal{C}}}
\newcommand{\calD}{{\mathcal{D}}}
\newcommand{\calE}{{\mathcal{E}}}

\newcommand{\calI}{{\mathcal{I}}}

\newcommand{\calK}{{\mathcal{K}}}

\newcommand{\calO}{{\mathcal{O}}}
\newcommand{\calP}{{\mathcal{P}}}

\newcommand{\calR}{{\mathcal{R}}}

\newcommand{\calT}{{\mathcal{T}}}

\newcommand{\calW}{{\mathcal{W}}}

\newcommand{\scrA}{{\mathscr{A}}}

\newcommand{\scrM}{{\mathscr{M}}}

\newcommand{\scrP}{{\mathscr{P}}}

\newcommand{\scrT}{{\mathscr{T}}}
\newcommand{\scrU}{{\mathscr{U}}}

\newcommand{\scrX}{{\mathscr{X}}}

\crefname{algocf}{Algorithm}{Algorithms}
\Crefname{algocf}{Algorithm}{Algorithms}
\crefname{appendix}{}{}
\Crefname{appendix}{}{}
\crefname{figure}{Figure}{Figures}
\Crefname{figure}{Figure}{Figures}

\theoremstyle{plain}
\newtheorem{theorem}{Theorem}[section]

\newtheorem{corollary}[theorem]{Corollary}
\newtheorem{proposition}[theorem]{Proposition}

\newtheorem{assumption}{Assumption}[section]

\theoremstyle{definition}

%% file: abstract.tex
We present \texttt{LazyDINO}, a transport map variational inference method for fast, scalable, and efficiently amortized solutions of high-dimensional nonlinear Bayesian inverse problems with expensive parameter-to-observable (PtO) maps. Our method consists of an offline phase in which we construct a derivative-informed neural surrogate of the PtO map using joint samples of the PtO map and its Jacobian. During the online phase, when given observational data, we seek rapid posterior approximation using surrogate-driven training of a lazy map [Brennan et al., \textit{NeurIPS}, (2020)], i.e., a structure-exploiting transport map with low-dimensional nonlinearity. The trained lazy map then produces approximate posterior samples or density evaluations. Our surrogate construction is optimized for amortized Bayesian inversion using lazy map variational inference. We show that (i) the derivative-based reduced basis architecture [O'Leary-Roseberry et al., \textit{Comput. Methods Appl. Mech. Eng.}, 388 (2022)] minimizes the upper bound on the expected error in surrogate posterior approximation, and (ii) the derivative-informed training formulation [O'Leary-Roseberry et al., \textit{J.\ Comput.\ Phys.}, 496 (2024)] minimizes the expected error due to surrogate-driven transport map optimization. Our numerical results demonstrate that \texttt{LazyDINO} is highly efficient in cost amortization for Bayesian inversion. We observe one to two orders of magnitude reduction of offline cost for accurate posterior approximation, compared to simulation-based amortized inference via conditional transport and conventional surrogate-driven transport. In particular, \texttt{LazyDINO} outperforms Laplace approximation consistently using fewer than 1000 offline samples, while other amortized inference methods struggle and sometimes fail at 16,000 offline samples.

%% file: introduction.tex
\section{Introduction}

We investigate the solution of nonlinear \emph{Bayesian inverse problems} (BIPs), i.e., inferring uncertain parameters of computational models from sparse, noisy, and indirect observational data. Let $m\in\scrM$ denote the unknown model parameter and assume the observational data vector $\Obs\in\R^{d_{\Obs}}$ is given by:
\begin{equation*}
    \Obs = \Forward(\FullParam) + \noise, \quad \noise \sim\Gaussian(0,\NoiseCov),
\end{equation*}
where $\Forward:\scrM\to\R^{d_{\Obs}}$ is the parameter-to-observable (PtO) map and $\noise\in\R^{d_{\Obs}}$ is an unknown noise vector. Given a parameter prior distribution $\mu$, we seek to characterize the posterior distribution $\mu^{\Obs}$ defined via Bayes' rule 
\begin{equation*}
\underbrace{\meas\FullPost(m)}_{\displaystyle\mathclap{\text{Posterior}}}\propto \underbrace{\exp\left(-\frac{1}{2}\left\|\NoiseCov^{-1/2}\left(\Forward(\FullParam) - \Obs\right)\right\|^2\right)}_{\displaystyle\mathclap{\text{Likelihood}}}\underbrace{\meas\mu(m)}_{\displaystyle\mathclap{\text{Prior}}}.
\end{equation*}

We are particularly interested in continuum models for physical systems, e.g., parametric partial differential equations (PDEs), where the parameter $m\in\ParamSpace$ can have arbitrarily high dimensions, such as spatially varying parameter fields, and the PtO map $\Forward$ is defined implicitly through the solution of the governing equations \cite{Stuart2010,  buithanh2013computational, petra2014computational,Ghattas_Willcox_2021}. This type of BIP is challenging due to (i) the high computational cost of likelihood evaluations due to model solutions, (ii) the difficulty of characterizing high-dimensional posterior distributions due to the curse of dimensionality, and (iii) non-Gaussianity of the posterior distribution. These challenges are acute limitations when one seeks fast solutions of BIPs for a range of observational data, as in real-time uncertainty quantification for predictive digital twins \cite{kapteyn2021probabilistic} and optimal experimental design \cite{huan2024optimal}. Solving BIPs in this setting requires methods with \textit{amortized computational cost}---that is, most of the expensive computation is performed offline, i.e., before acquiring the data, and posterior characterization incurs a comparatively negligible cost once the data is available. These challenges demand methodological innovations beyond conventional approaches such as Markov chain Monte Carlo (MCMC). In this work, we integrate recent advances in dimension reduction, neural operator learning, and measure transport to derive a fast, scalable, and efficiently amortized method for BIPs that is well-suited to modern computing frameworks.

\subsection{Variational inference using lazy maps}

We consider using transport map variational inference (TMVI) to approximate the posterior $\mu^{\Obs}$. This method seeks to construct a parameterized transport map $\FullMap_{\MapParams}:\scrM\to\scrM$ between a reference distribution, which we take to be the prior $\mu$, and the target Bayesian posterior $\FullPrior^{\Obs}$. The map parameters can be found by minimizing the reverse KL divergence (rKL):
\begin{equation}\label{eq:rkl_optimization}   \min_{\MapParams}\KL({\FullMap_{\MapParams}}_{\sharp}\mu\|\FullPost),
\end{equation}
where $(\cdot)_{\sharp}$ denotes the pushforward of a probability distribution. Once the transport map is constructed, it allows for fast on-demand approximate posterior sampling through map evaluations  $\FullMap_{\MapParams}(m^{(j)})$ of the reference samples $m^{(j)}\iid \mu$. However, it can be difficult to represent expressive transport maps in high dimensions. For example, triangular maps  \cite{baptista2023representation} on $\R^n$ must describe $n$-variate functions and thus suffer from the curse of
dimensionality. Kernel-based methods, such as Stein variational inference \cite{liu2019stein,chen2019projected, detommaso2018stein}, lose expressiveness in high dimensions.
Flow-based methods \cite{rezende2015variational, papamakarios2019normalizing} often increase expressiveness by adding layers, which is typically performed ad hoc and require tuning.

Lazy maps \cite{brennan2020greedyinferencestructureexploitinglazy} are a class of transport maps that alleviate the curse of dimensionality by restricting the nonlinearity of the map to a relatively low-dimensional parameter latent space. Let $\Encoder:\scrM\to\R^{d_{r}}$ with $d_r\ll\text{dim}(\scrM)$ be a linear encoder such that $\text{Im}(\Encoder)=\R^{d_{r}}$ defines the parameter latent space. A lazy map has the following form: 
\begin{equation}\label{eq:lazy_map}
    \FullMap_{\MapParams} \coloneqq (\identity_{\scrM} - \Encoder\circ\Decoder) + \Decoder\circ \Map_{\MapParams} \circ \Encoder,
\end{equation}
where $\identity$ denotes the identity map, $\Decoder:\R^{d_{r}}\to\scrM$ is a linear decoder, and $\Map_{\MapParams}:\R^{d_{r}}\to\R^{d_r}$ is a parametrized latent space transport map. The lazy map approximates the posterior in the latent space using TMVI while the prior fills the complementary space of $\text{Im}(\Decoder)$. When the prior is Gaussian, and the encoder is a whitening transformation, the rKL minimization problem becomes
\begin{equation}\label{eq:intro_rkl_min}
\min_{\MapParams}\mathbb{E}_{m\sim\mu}\left[\frac{1}{2}\left\|\NoiseCov^{-1/2}\left(\left(\Forward\circ\FullMap_{\MapParams}\right)(m) - \Obs\right)\right\|^2 + \frac{1}{2}\norm{ \Map_{\MapParams}(\Encoder m)}^2 -\log|\det\grad \Map_{\MapParams}(\Encoder m)| \right].
\end{equation}

A key component of lazy maps is finding a parameter subspace that captures the discrepancy between the prior and posterior distribution. This parameter subspace is often known as the likelihood-informed subspace \cite{Cui_2014} or active subspace \cite{constantine2015active}, which is known to exist for a large class of inverse problems and can be found via solving eigenvalue problems based on score functions; see, e.g., \cite{Bui-ThanhGhattas12a, Bui-ThanhGhattas12, Bui-ThanhGhattas13a, ChenGhattas2019, ChenVillaGhattas2017, FlathWilcoxAkcelikEtAl11,
IsaacPetraStadlerEtAl15,spantini2015optimal}. By exploiting the structure of the BIPs, TMVI using lazy maps typically achieves high-quality posterior approximation more efficiently than TMVI without parameter reduction or with alternative reduction techniques.

Another fundamental challenge of TMVI, including when a lazy map is used, lies in the high cost of transport map training, which requires solving the stochastic and model-constrained rKL minimization problem in \cref{eq:intro_rkl_min}. 
Numerous evaluations of the PtO map $\Forward$ and the actions of its Jacobian $D\Forward$ are required within each optimization iteration. These evaluations involve repeated solutions of the governing equations of the computational models and their forward or adjoint sensitivities, which can be prohibitively expensive when these equations are, e.g., large-scale nonlinear PDEs. This cost barrier becomes further exacerbated when multiple posteriors need to be approximated for different instances of observational data.

\subsection{Derivative-informed surrogate for amortized lazy map variational inference}

In this work, we remove the computational bottleneck of model solutions in lazy map training by constructing a fast-to-evaluate ridge function surrogate of the PtO map using a neural network latent representation $\RidgeForward_{\boldsymbol{w}}:\R^{d_r}\to\ObsSpace$:
\begin{equation*}
    \Forward(m)\approx \boldsymbol{V}\RidgeForward_{\boldsymbol{w}}(\Encoder m),
\end{equation*}
where $\bw$ is the weight of the neural network and $\bV$ is a (reduced) basis on the data space. This surrogate architecture that uses the same parameter reduction technique as lazy maps is the derivative-informed neural network (\texttt{DIPNet}) in \cite{o2022derivative}, which belongs to a larger class of reduced basis neural operator \cite{Hesthaven2018, kovachki2023neural}. Once the surrogate is constructed, we perform TMVI in the parameter latent space using the surrogate-driven rKL objective for the given observational data $\Obs$:
\begin{equation*}
\min_{\MapParams}\mathbb{E}_{\refCoord\sim\mathcal{N}(0, \identity_{\R^{d_r}})}\left[\frac{1}{2}\left\|(\RidgeForward_{\boldsymbol{w}}\circ\Map_{\MapParams})(\refCoord) - \bV^*\Obs\right\|^2 + \frac{1}{2}\norm{ \Map_{\MapParams}(\refCoord)}^2 -\log|\det\grad \Map_{\MapParams}(\refCoord)| \right].
\end{equation*}
The forward and Jacobian evaluation costs of the neural network $\RidgeForward_{\boldsymbol{w}}$ are significantly lower than those of model and sensitivity solutions. As a result, the high cost of the PtO map and its Jacobian evaluations for optimizing lazy maps are amortized.

Conventionally, the neural network is trained offline using samples of the PtO map via mean squared error minimization \cite{kovachki2024operator, o2022derivative}:
\begin{equation*}
     \min_{\bw} \mathbb{E}_{m\sim\mu}\left[\left\|\bV^*\Forward(m) - \RidgeForward_{\bw}(\Encoder m)\right\|^2\right].
\end{equation*}
When a limited number of PtO map samples is used to estimate the expectation, the trained surrogate may be inadequate for lazy map training as the surrogate Jacobian error is not directly controlled. This, in turn, leads to inaccurate gradients of the rKL objective and substantial gaps between the objective values at the true optimum and the surrogate-approximated optimum (which we refer to herein as optimality gaps). Similarly, conventionally trained neural operator surrogates struggle to accelerate other gradient-based optimization in high or infinite dimensions; see, e.g., optimal design \cite{luo2023efficient} and geometric MCMC \cite{cao2024efficient}.

In this work, we follow \cite{o2024derivative, cao2024efficient} and train \texttt{DIPNet} PtO map surrogates using a derivative-informed learning method, which exerts surrogate error control in the Sobolev space of the PtO map latent representation using joint samples of the PtO map $\Forward$ and its Jacobian $D\Forward$:
\begin{equation*}
    \min_{\bw} \mathbb{E}_{m\sim\mu}\left[\left\|\bV^*\Forward(m) - \RidgeForward_{\bw}(\Encoder m)\right\|^2 + \norm{\bV^*D\Forward(m)\Decoder - \nabla\RidgeForward_{\bw}(\Encoder m)}_{F}^2\right].
\end{equation*}
We show that derivative-informed learning of \texttt{DIPNet} surrogate is equivalent to minimizing an upper bound on posterior approximation error as well as the optimality gap of surrogate-driven lazy map training for amortized Bayesian inversion. We refer to this surrogate construction with optimized architecture and training as reduced basis derivative-informed neural operator (\texttt{RB-DINO}) \cite{o2024derivative}.

\subsection{Solving high-dimensional nonlinear Bayesian inverse problems using \texttt{LazyDINO}}

Combining TMVI using lazy maps and \texttt{RB-DINO} surrogate construction creates a competitive method for amortized solutions of high-dimensional model-constrained BIPs, referred to as \texttt{LazyDINO}. The method is composed of offline and online phases.
\begin{itemize}[label={},leftmargin=0pt]
    \item \textbf{Offline phase.} We first generate samples of the PtO map and its Jacobian by solving the governing equations of the computational model and its forward or adjoint sensitivity. These samples are then used to construct a \texttt{RB-DINO} surrogate of the PtO map.
    \item \textbf{Online phase.} After collecting observational data, we seek rapid posterior approximation via \texttt{RB-DINO} surrogate-driven training of a latent space transport map. The trained transport map can be used to produce approximate posterior samples. This process can be repeated for different observational data, effectively amortizing the construction cost of the \texttt{RB-DINO} surrogate.
\end{itemize}
We provide extensive numerical studies that compare \texttt{LazyDINO} against a range of TMVI methods, including the Laplace approximation, simulation-based amortized inference (SBAI) via conditional transport, TMVI using lazy maps, and lazy maps combined with conventional surrogate construction (\texttt{LazyNO}); see \cref{subsubsec:related_work_amortized} and \cref{tab:transposed_baselines}. We test these methods on two BIPs, each with four different instances of observation data: (i) inferring the diffusivity field in a nonlinear reaction--diffusion PDE and (ii) inferring the heterogeneous material property of a hyperelastic material thin film.
We devise extensive posterior approximation tests, including moment discrepancies, probability density-based metrics, and various posterior visualizations (e.g., marginals, mean, MAP estimates, and point-wise marginal variance).

The main contributions of the \texttt{LazyDINO} method are summarized below.

\begin{enumerate}[label=(\textbf{C}{{\arabic*}})]
\item The \texttt{RB-DINO} surrogate construction is optimized for amortized lazy map variational inference.
\begin{itemize}[leftmargin=0pt]
    \item \textit{Surrogate architecture}.
    In \cref{theorem:kl_bound}, we derive upper bounds for the expected posterior approximation error when a neural ridge function surrogate replaces the PtO map in the likelihood. This result, which is a straightforward extension from those in \cite{cui2021data, cao2024efficient}, bounds the forward KL (fKL) averaged over the marginal observational data distribution by the sum of a parameter reduction error and a latent representation error. Minimizing this error upper bound gives rise to the \texttt{DIPNet} architecture \cite{o2022derivative}, where the parameter encoder is found by derivative-informed dimension reduction \cite{zahm2020gradient, constantine2015active} using samples of the PtO map Jacobian.
    \item \textit{Derivative-informed learning.} We show that the expected gradient error (\cref{theorem:kl_gradient}), and the expected optimality gap (\cref{corollary:optimality_gap}) in surrogate-driven lazy map optimization can be controlled by a weighted Sobolev norm of the surrogate approximation error. This error measure is consistent with the objective function in derivative-informed operator learning \cite{o2024derivative, cao2024efficient} that uses joint samples of the PtO map and its Jacobian for surrogate training. In other words, derivative-informed learning of reduced basis neural networks (\texttt{RB-DINO}) minimizes the expected error in the stochastic optimization of lazy maps due to the surrogate representation.
\end{itemize}
\item \texttt{LazyDINO} enables fast, scalable, and efficiently amortized solutions of high-dimensional Bayesian inverse problems.
\begin{itemize}[leftmargin=0pt]
    \item \textit{Scabalility.} The surrogate and transport map training in \texttt{LazyDINO} are independent of the parameter dimension as we co-design their latent representations in the same relatively low-dimensional parameter subspace that captures prior-to-posterior updates.
    \item \textit{Fast online inference.} Using a neural network surrogate rKL objective for transport map training, \texttt{LazyDINO} circumvents the computational bottleneck of model solutions and fully exploits GPU-based accelerations to rapidly approximate posteriors. We demonstrate that the optimize-then-sample approach of \texttt{LazyDINO} leads to faster online sampling than the typical inversion-to-sample approach of SBAI.
    \item \textit{Superior cost--accuracy trade-off in amortized Bayesian inversion.} The \texttt{RB-DINO} surrogate construction is highly efficient in cost amortization for solving BIPs, i.e., it achieves high posterior approximation error at low offline training cost. In our numerical example, we observed one to two orders of magnitude of cost reduction in offline computation for achieving similar accuracy in posterior approximation compared to \texttt{LazyNO} and SBAI. Moreover, \texttt{LazyDINO} consistently outperforms Laplace approximation at a small training sample regime ($<1,000$). In contrast, \texttt{LazyNO} and SBAI struggle to outperform Laplace approximation and, in some cases, failed at $16,000$ training samples.
\end{itemize}
\end{enumerate}

\subsection{Related works}

In the following subsections, we discuss related work in dimension reduction, surrogate modeling, and variational inference for BIPs.

\subsubsection{Baseline: The Laplace Approximation}
The Laplace approximation (LA) constructs a Gaussian approximation of the posterior, leading to efficient sampling and density evaluations \cite{laplace1986memoir}. This makes the LA a sensible baseline for settings that require fast approximate posterior sampling. The LA construction requires (1) MAP point estimation, followed by (2) covariance estimation via solving a generalized eigenvalue problem for the Hessian of the negative log-posterior at the MAP point. These Hessians often have low effective numerical rank, allowing for efficient implementations in practice \cite{,buithanh2013computational,Ghattas_Willcox_2021, 10.5555/2388996.2389000}. Details on LA are included in \cref{LA_background}.

\subsubsection{Dimension reduction for Bayesian inverse problems}
A common likelihood-independent dimension reduction technique is the Karhunen--Lo\`eve expansion, which represents the parameter in a finite (small) number of prior covariance eigenfunctions, see e.g.,  see~\cite{marzouk_kl_truncation}. 
Derivative-based dimension reduction techniques identify the parameter subspace that the likelihood is most sensitive to, in prior or posterior expectation,
and thereby provide more targeted dimension reduction and greater efficiency \cite{Cui_2014,
zahm2020gradient, zahm2022certified,bigoni2022nonlinear,ChenGhattas2019}. This subspace, often referred to as the likelihood-informed subspace, is related to the Fisher information and is shown to be optimal with respect to the KL divergence in~\cite{zahm2022certified}. In \cite{cui2021data}, the authors show that a subspace computed by averaging Fisher information over the prior distribution is optimal on average over the marginal distribution of observational data. 

\subsubsection{Surrogate models for Bayesian inverse problems}

Substantial work has been done on using surrogate models, e.g., polynomial approximation \cite{marzouk_kl_truncation, marzourk2007stochastic,marzouk2009stochastic,farcas2020multilevel} and model-order reduction \cite{galbally2010nonlinear, lieberman2010parameter, cui2015data}, to accelerate solutions of BIPs. Surrogate models are often used within multi-fidelity posterior sampling algorithms \cite{Peherstorfer2018, lykkegaard2023multilevel, cao2023residual}.

This work focuses on neural network surrogates due to their high flexibility, rich approximation properties, and scalability. Since our algorithmic framework can be applied to infinite-dimensional BIPs, we note the connection to neural operator surrogates \cite{kovachki2023neural} that map between function spaces with architectures and training agnostic to the discretization of these spaces. Notable architectures include reduced basis neural networks using linear \cite{Hesthaven2018,bhattacharya2021model,fresca2022pod,o2022derivative,o2022learning} or nonlinear \cite{lu2021learning, seidman2023variational} dimension reduction, and neural network integral kernels such as Fourier neural operator and its variants \cite{li2021fourier,cao2023lno, lanthaler2023nonlocal, li2020neural}. Notable training formulations include conventional supervised learning using input-output samples and physics-informed learning \cite{li2024physics, wang2021learning} with additional loss functions related to the residual of the implicit equation (e.g., PDE residuals).

Our surrogate architecture is based on the aforementioned derivative-based dimension reduction strategy, i.e., \texttt{DIPNet} in \cite{o2022derivative,o2022learning}. We advocate for the derivative-informed operator learning method for surrogate training, i.e., derivative-informed neural operator (DINO) \cite{o2024derivative}. This method has been successfully applied to surrogate-driven solutions of PDE-constrained optimization under uncertainties \cite{luo2023efficient}, Bayesian optimal experimental design \cite{go2023accelerating,go2024sequential}, and infinite-dimensional BIPs \cite{cao2024efficient}. Recent work \cite{qiu2024derivative} also explored training the deep operator network (DeepONet) architecture using this method.

\subsubsection{Transport map parameterizations and amortized inference} \label{subsubsec:related_work_amortized}

The \texttt{LazyDINO} algorithm performs TMVI for a reduced dimensional inference problem in a specifically chosen latent space. It assumes no particular map parameterization; rather it wraps around any provided transport map class. We briefly review popular transport map parameterizations and provide references for further reading.
Normalizing flows (see~\cite{rezende2015variational,tabak2013family,papamakarios2019normalizing,kobyzev2020normalizing}) form a broad class of methods that construct transport maps through compositions of neural networks with specific parameterizations. Autoregressive flows~\cite{dinh2016density,Kingma2016,papamakarios2017masked,huang2018neural,de2019block, Jaini2019}, a popular subclass, compose autoregressive (triangular) maps to allow efficient computation of Jacobian determinants~\cite{daniels2010monotone,wehenkel2019unconstrained}.
Several works seek an approximation to the Knothe-Rosenblatt (KR) rearrangement~\cite{Knothe1957ContributionsTT, 10.1214/aoms/1177729394}, a diffeomorphic triangular map that exists between any two distributions that are absolutely continuous with respect to a common measure, using orthonormal basis expansions (e.g., sparse polynomials) 
\cite{el2012bayesian,marzouk2016introduction, spantini2018inference, baptista2020adaptive,Zech2022Transport2,Westermann2023Transport} or neural networks \cite{Zech2022Transport1}.
One distinguishing feature of \texttt{LazyDINO} is its use of PtO map Jacobian evaluations during training. Other recent works also incorporate derivative information, such as by adding a Fisher divergence term to the training objective~\cite{zeghal2022neuralposteriorestimationdifferentiable, Brehmer_2020}, thus exploiting the differentiability of the log-likelihood.
Finally, we note the recent rise in inference methods that amortize transport-based posterior approximation \cite{durkan2018sequential, papamakarios2019sequential,greenberg2019automatic, papamakarios2018fastepsilonfreeinferencesimulation, baptista2024bayesian}. Simulation-based amortized
inference (SBAI)~\cite{Ganguly_2023} approaches parameterize transport maps to treat the conditioning variable (i.e., the observation) as a functional input and generate samples from the corresponding posterior. We compare \texttt{LazyDINO} with SBAI in our numerical examples.

\subsection{Notation}\label{notation}

\begin{itemize}

    \item We use bold symbols to denote finite-dimensional vectors, e.g., $\bx \in \mathbb{R}^{d_{\bx}}$ where $d_{\bx}$ is the dimension. We denote the 2-norm on finite-dimensional vector spaces as $\| \cdot \|$. We denote the Frobenius matrix norm as $\|\cdot \|_F$. We use math script to denote separable Hilbert spaces that have high or infinite dimensions, e.g., $\mathscr{X}$, with inner product $\langle \cdot, \cdot \rangle_{\mathscr{X}}$, norm $\norm{\cdot}_{\mathscr{X}}$ and element $x \in \mathscr{X}$. 
    \item We denote the Banach space of bounded linear operators from $\mathscr{X}_1$ to $\mathscr{X}_2$ as $B(\mathscr{X}_1,\mathscr{X}_2)$. We denote its subset of Hilbert--Schmidt (HS) operators as $\HS(\scrX_1, \scrX_2)$. We define $B(\mathscr{X})\vcentcolon= B(\mathscr{X},\mathscr{X})$ and similar for HS operators. We denote the set of positive, self-adjoint, and trace class operators on $\scrX$ as $\TrClass(\scrX)$. When $\scrX$ is a finite-dimensional vector space, $\TrClass(\scrX)$ consists of symmetric positive definite matrices.
    
    \item We denote the inner-product and norm weighted by a positive and self-adjoint operator $\mathcal{A}:\scrX\to\scrX$ as $\langle x_1, x_2\rangle_{\mathcal{A}}\coloneqq \langle \mathcal{A}^{1/2}x_1, \mathcal{A}^{1/2}x_2\rangle_{\mathscr{X}}$ and $\norm{x}_{\mathcal{A}}\coloneqq\sqrt{\langle\mathcal{A}^{1/2}x, \mathcal{A}^{1/2}x\rangle_{\mathscr{X}}}$  and $\mathcal{A}^{1/2}$ denotes the self-adjoint square root of $\mathcal{A}$. We note that the operator square root is not required in numerical implementation.

    \item The set of probability distributions defined using Borel $\sigma$-algebra on $\scrX$ is denoted by $\mathscr{P}(\mathscr{X})$. The density between two probability distributions (i.e., Radon--Nikodym derivative) $\mu_1,\mu_2\in\scrP(\scrX)$ evaluated at $x\in\scrX$ is given by $(\meas\mu_1/\meas\mu_2)(x)$. For probability distributions on finite-dimensional vector spaces, we do not distinguish between a distribution $\mu$ and its probability density function $\pi(\bx) = (\meas\mu/\meas\mu_L)(\bx)$, where $\mu_L$ is the Lebesgue measure. We use $\bx\sim\pi$ and $\bx\iid\pi$ to denote a $\pi$-distributed random variable and independent and identically distributed samples from $\pi$, respectively.
    \item We denote the diffeomorphism group of $\mathscr{X}$ as $\text{Diff}^1(\mathscr{X}):=\big\{\FullMap: \mathscr{X} \to \mathscr{X} \;| \; \FullMap$ is an automorphism, and  $\FullMap, \FullMap^{-1} \in C^1(\mathscr{X})\big\}$. For $\FullMap \in\text{Diff}^1(\scrX)$, we denote by $\FullMap_{\sharp}\mu$ and $\FullMap^{\sharp}\mu$ the pushforward and pullback of probability distributions in the sense that $\FullMap_{\sharp}\mu = \mu\circ\FullMap^{-1}$ and $\FullMap^{\sharp}\mu = \mu\circ\FullMap$.
 
\end{itemize}

\subsection{Outline of the paper}

The remainder of this work proceeds as follows: 
In~\cref{section:inf_dim_bayes_measure_transport_subspace}, we introduce the lazy map variational inference method for solving high-dimensional Bayesian inversion. 
Then, in ~\cref{section:dino_approximated_likelihood}, we introduce the optimized surrogate construction for amortized Bayesian inverse for lazy map variational inference.
In~\cref{section:lazy_dino_amortize}, we describe the \texttt{LazyDINO} algorithm in detail, including documentation of all offline and online procedures and its role in enabling amortized inference. 
In~\cref{section:numerical_study_setup}, we define the setup for numerical experiments, including the two infinite-dimensional PDE-constrained BIPs and all metrics utilized to measure the posterior approximation errors.
In~\cref{sec:numerical_results}, we present the numerical results and discuss the relative performance of the methods. 
Finally, we give concluding remarks in ~\cref{sec:conclusion}. We have additional results and discussions in the various appendices.

%% file: bips.tex
\section{Solving Bayesian inverse problems using lazy map variational inference}\label{section:inf_dim_bayes_measure_transport_subspace}
This section introduces our framework for solving high-dimensional nonlinear BIPs using lazy map variational inference (LMVI). In \cref{nonlinear_gaussian_prior_gaussian_noise}, we define the setting for BIPs considered in this work. Then, we describe posterior approximation via ridge functions and the resulting inference problems in a latent parameter space of lower dimensions in \cref{subsec:subsapce_decomposition,subsec:ridge_function,subsec:latent_space}. Lastly, we introduce LMVI that approximates the target posterior in the latent space in \cref{subsec:lazy_map,subsec:optimizing_lazy_map}.

\subsection{Nonlinear Bayesian inverse problems}\label{nonlinear_gaussian_prior_gaussian_noise}
We denote the unknown parameter of interest $\FullParam \in \ParamSpace$, where $\ParamSpace$ is a separable Hilbert space with inner product $\langle\cdot, \cdot\rangle_{\ParamSpace}$. Let $\Obs\in\ObsSpace$ denote observational data, and $\Forward:\ParamSpace\to \R^\ObsDim$ denote a nonlinear parameter-to-observable (PtO) map. We begin with the following standard assumptions.
\begin{assumption}[Gaussian prior distribution]\label{gaussian_prior}
We consider a Gaussian prior distribution $\FullPrior = \Gaussian(0,\PriorCov)\in\ProbSpace(\ParamSpace)$ with $\PriorCov\in\TrClass(\scrM)$.
\end{assumption}

\begin{assumption}[Additive Gaussian noise]\label{additive_gaussian}
    We assume the observed data has the following distribution:
\begin{equation}\label{eq:data_model}
\Obs\sim\Gaussian(\Forward(\FullParam),\NoiseCov),
\end{equation}
where $\NoiseCov\in \TrClass(\ObsSpace)$.
\end{assumption}
We consider the inverse problem to recover the parameter $\FullParam$ from an observed data vector $\Obs$. We are interested in characterizing the posterior distribution satisfying Bayes' rule, which we denote by $\FullPost \in \ProbSpace(\ParamSpace)$:
\begin{equation}\label{eq:generic_full_bip}
    \frac{\meas \FullPost}{ \meas \FullPrior}(\FullParam) = \frac{1}{\normalization^{\Obs}}\exp(-\Misfit(\FullParam)),\quad\ \Misfit(\FullParam) \coloneqq \frac{1}{2}\|\Forward(\FullParam)- \Obs\|_{\NoiseCov^{-1}}^2.
\end{equation}
Here, we define the \emph{potential} $\Misfit:\ParamSpace \to \R$ (i.e., the negative log-likelihood), and the normalization constant ${\normalization^{\by}}\coloneqq\mathbb{E}_{m\sim\mu}\left[\exp(-\Misfit(m))\right]$.

\begin{remark}
    We address two potential concerns regarding our BIP setting. Firstly, to infer parameters with non-Gaussian priors, one can perform inference in transformed coordinates distributed according to a Gaussian prior and subsequently sample from the posterior via the inverse transform; see, e.g., \cite{doi:10.1137/S1064827503424505, SRAJ2016205}. Secondly, even though we only consider the BIP arising from a single observation, this work straightforwardly extends to a collection of observations, e.g., $\Obs_1, \Obs_2,\dots\iid \mathcal{N}(\Forward(m), \NoiseCov)$, in which case the potential in \cref{eq:generic_full_bip} becomes the sum of the negative log-likelihood of each observation. This extension is non-trivial for many SBAI methods.
\end{remark}

We define the following Cameron--Martin spaces, separable Hilbert spaces with prior and noise precision--weighted inner products,
\begin{equation*}
    \ParamCM \coloneqq \left\{m\in\ParamSpace:\norm{m}_{\PriorCov^{-1}}<\infty\right\},\quad
    \ObsCM \coloneqq \left\{\by\in\ObsSpace:\norm{\Obs}_{\NoiseCov^{-1}}<\infty\right\},
\end{equation*}
where $\ParamCM$ and $\ObsCM$ are equipped with inner products $\langle\cdot,\cdot\rangle_{\PriorCov^{-1}}$ and $\langle\cdot,\cdot\rangle_{\NoiseCov^{-1}}$, respectively. Note that $\ParamCM$ is continuously embedded in $\ParamSpace$ and is isomorphic to $\ParamSpace$ with respect to the identity map only when $\ParamSpace$ is finite-dimensional.

In this work, we assume the PtO map is $H^1_{\mu}$-differentiable in the following sense.
\begin{assumption}[$H^1_\mu$-differentiable PtO map]
We assume the PtO map to live in the Sobolev space with Gaussian measure $H^1_\mu(\ParamSpace; \ObsCM)\coloneqq\{\Forward:\ParamSpace\to \ObsSpace: \norm{\Forward}_{H^1_{\mu}(\ParamSpace; \ObsCM)} <\infty\}$ where
\begin{equation}\label{assume_h1}
    \|\Forward\|^2_{H^1_{\mu}(\ParamSpace, \ObsCM)} \coloneqq \mathbb{E}_{\FullParam \sim \mu} \left[\|\Forward(\FullParam)\|^2_{\NoiseCov^{-1}} + \|D_H\Forward(\FullParam)\|_{\mathrm{HS}(\ParamCM,\ObsCM)}^2 \right],
\end{equation}
and $D_H\Forward: \ParamSpace\to\HS(\ParamCM,\ObsCM)$ is the stochastic derivative or the Malliavin derivative of $\Forward$ that satisfies
    \begin{equation}\label{eq:stochastic_derivative}
        \lim_{t\to 0}\norm{t^{-1}\left(\Forward(m+t\delta m) - \Forward(m)\right) - D_{H} \bdmc{G}(m)\delta m}_{\NoiseCov^{-1}} = 0\quad \forall \delta m\in \ParamCM\quad \mu\text{-a.e.}
    \end{equation}
\end{assumption}
If the Fr\'echet derivative $D\Forward: \ParamSpace\to \HS(\ParamSpace, \ObsSpace)$ exists, we have the following equality $\mu$-a.e.:
\begin{equation}\label{eq:gateaux_versus_stochastic_gateaux}
    D_H\Forward(m) = D\Forward(m)|_{\ParamCM},\quad D_H\Forward(m)^* = \PriorCov D\Forward(m)^*\NoiseCov^{-1},
\end{equation}
where $|_{\ParamCM}$ denotes the restriction of the function domain from $\ParamSpace$ to $\ParamCM$. We do not distinguish between $D\Forward(m)$ and $D_H\Forward(m)$ when it is clear that the domain is $\ParamCM$.

\subsection{Subspace decomposition of Bayesian inverse problems}\label{subsec:subsapce_decomposition}

Let $\Projector\in B(\ParamSpace)$ be a rank-$\RedCoordDim$ linear projection. We have the corresponding unique decomposition of $\ParamSpace$ into the image (i.e., the range space) and kernel (i.e., the null space) of $\Projector$:
\begin{equation*}
    \ParamSpace = \text{Im}(\Projector) \oplus \text{Ker}(\Projector),
    \quad
    \FullParam = \underbrace{\Projector \FullParam}_{\displaystyle \FullParam_r} + \underbrace{(\identity_{\ParamSpace}-\Projector)\FullParam}_{\displaystyle \FullParam_{\perp}}\quad\forall \FullParam\in\ParamSpace,
\end{equation*}
where $\oplus$ denotes the direct sum as defined above. We denote the prior and posterior marginals in $\text{Im}(\Projector)$ by the pushforward $\FullPrior_r\coloneqq \Projector_{\sharp}\mu$ and $\FullPost_r\coloneqq \Projector_{\sharp}\FullPost$, respectively. A probability distribution on $\ParamSpace$ can be decomposed into its marginal probability in $\text{Im}(\Projector)$ and its conditional probability in $\text{Ker}(\Projector)$ in the following sense. For any measurable subset $\scrA\subseteq\scrM$ and its decomposition $\scrA=\scrA_r\oplus\scrA_{\perp}$, where $\scrA_r\subseteq\text{Im}(\Projector)$ and $\scrA_{\perp}\subseteq\text{Ker}(\Projector)$, the prior and posterior probability concentrations on $\scrA$ are given by
\begin{gather}\label{eq:prior_disintegration}
    \FullPrior(\scrA)  = \int_{\scrA_r} \mu_{{\perp}|r} (\scrA_{\perp} |\FullParam_r)\meas\mu_r(\FullParam_r), 
    \quad  
    \mu^{\Obs}(\scrA) =  \int_{\scrA_r} \mu^{\Obs}_{{\perp}| r} (\scrA_{\perp} |\FullParam_r) \meas \mu_r^{\Obs}(\FullParam_r),
\end{gather}
where $\mu_{{\perp}| r}(\cdot| m_r), \mu^{\Obs}_{{\perp}|r}(\cdot| m_r)\in \ProbSpace(\text{Ker}(\Projector))$ are the prior and posterior conditionals in $\text{Ker}(\Projector)$. In particular, the prior marginal $\mu_r$, hereafter referred to as the \emph{subspace prior}, and conditional $\mu_{{\perp}| r} (\cdot|\FullParam_r)$ has closed forms given by:
\begin{equation}\label{eq:close_form_prior_disintegration}
        \mu_r = \Gaussian\left(0, \Projector\PriorCov\Projector^*\right),
        \quad 
        \mu_{{\perp}| r} (\cdot|\FullParam_r) = \Gaussian\left(\PriorCov\Projector^*(\Projector\PriorCov\Projector^*)^{-1}\FullParam_r - \FullParam_r, \PriorCov\Projector^* - \Projector\PriorCov\right),
\end{equation}
where $\Projector^*$ is the Hermitian adjoint of $\Projector$. These forms can be simplified for specific choices of $\Projector$, which is discussed in 
\cref{subsec:latent_space}.

\subsection{Posterior approximation using ridge functions}\label{subsec:ridge_function}

We proceed under the assumption that the projection $\Projector$ has been chosen such that the data $\Obs$ are uninformative of the parameter in $\text{Ker}(\Projector)$, i.e., the difference between the prior and the posterior is small in $\text{Ker}(\Projector)$. The process for choosing $\Projector$ will be delineated in \cref{section:dino_approximated_likelihood}. Under this assumption, we consider a ridge function approximation of the PtO map:
\begin{equation} \label{eq:ridge_definition}
    \FullRidgeForward: \text{Im}(\Projector)\to\ObsSpace,\quad \FullRidgeForward\circ \calP\approx \Forward.
\end{equation}
An example of such a ridge function is the conditional expectation of the PtO map, where the projected parameter input is lifted into the full space $\ParamSpace$ by filling $\text{Ker}(\Projector)$ with the prior conditional:
\begin{equation} \label{eq:conditional_expectation_ridge}
    \FullRidgeForward_{\text{opt}}(\Projector m) \coloneqq \mathbb{E}_{m_{\perp}\sim \mu_{\perp|r}(\cdot|m_r)}\left[\Forward(\Projector m + m_{\perp})\right].
\end{equation}
For a given projection, this ridge function is optimal with respect to the Bochner norm on $L^2_{\mu}(\scrM;\ObsCM)$ \cite{zahm2020gradient, cao2024efficient},
\begin{equation*}
    \mathbb{E}_{m\sim\mu}\left[\norm{\bdmc{G}(m)-\FullRidgeForward_{\text{opt}}\left(\Projector m\right)}^2_{\NoiseCov^{-1}}\right] = \inf_{\FullRidgeForward:\text{Im}(\Projector)\to\ObsSpace}\mathbb{E}_{m\sim\mu}\left[\norm{\Forward(m)-\FullRidgeForward(\Projector m)}_{\NoiseCov^{-1}}^2\right].
\end{equation*}
We refer to $\FullRidgeForward_{\text{opt}}\circ\Projector$ as the \emph{optimal ridge function}.

A ridge function approximation of the PtO map induces an approximate posterior $\widetilde{\FullPrior}^{\Obs}\in\mathscr{P}(\ParamSpace)$ given by
\begin{equation}\label{eq:full_ridge_posterior} 
    \RidgeMisfit(\Projector m)\coloneqq\frac{1}{2}\norm{\FullRidgeForward( \Projector m) - \Obs}^2_{\NoiseCov^{-1}},\quad
    \frac{\meas \widetilde{\FullPrior}^{\Obs}}{ \meas \mu}(m) = \frac{1}{\widetilde{\normalization}^{\Obs}}\exp\left(-\RidgeMisfit(\Projector m)\right),
\end{equation}
where $\RidgeMisfit\circ\Projector \approx\Misfit$. Since the ridge function does not act in $\text{Ker}(\Projector)$, the approximate posterior conditional $\widehat{\mu}^{\Obs}_{\perp|r}$ is proportional to the prior conditional $\mu_{\perp|r}$, and the following holds by Bayes' rule 
\begin{equation}\label{eq:subspace_ridge_posterior}
    \frac{\meas\widetilde{\FullPrior}^{\Obs}_{r}}{\meas\FullPrior_r}(\FullParam_r) = \frac{1}{\widetilde{\normalization}^{\Obs}_r}\exp(-\RidgeMisfit(\FullParam_r)),\quad \widetilde{\mu}^{\Obs}_{\perp|r} = \frac{\widetilde{\normalization}_r^{\Obs}}{\widetilde{\normalization}^{\Obs}} \mu_{\perp|r},
\end{equation}
where the \emph{subspace posterior}, $\widetilde{\FullPrior}^{\Obs}_r\in\mathscr{P}(\text{Im}(\Projector))$, is a marginal of $\widetilde{\FullPrior}^{\Obs}$. 

The quality of the ridge function $\FullRidgeForward\circ\Projector$ can be understood through statistical distances between the posterior $\mu$ and the approximate posterior $\widetilde{\mu}^{\Obs}$; these results are covered in \cref{error_analysis_section}. 

\subsection{Latent Bayesian inverse problems induced by ridge functions}\label{subsec:latent_space}
Let $\Psi_{r}=\{\psi_j\in\ParamSpace\}_{j=1}^{\RedCoordDim}$ denote a basis for $\text{Im}(\Projector)$, i.e., $\text{span}(\Psi_{r}) = \text{Im}(\Projector)$. We refer to $\Psi_{r}$ as a \textit{reduced basis} on $\ParamSpace$. The reduced basis defines an encoder $\Encoder$ and decoder $\Decoder$ pair:
\begin{equation}\label{eq:generic_encoder_decoder_definition}
    \begin{cases}
        \Encoder:\ParamSpace\ni \sum_{j=1}^{\RedCoordDim} \RedCoordParam_j\psi_j + m_\perp \mapsto \RedCoordParam\in\RedCoordParamSpace,\\
        \Decoder: \RedCoordParamSpace\ni\RedCoordParam\mapsto \sum_{j=1}^{\RedCoordDim} \RedCoordParam_j\psi_j \in \text{Im}(\Projector), 
    \end{cases}\quad \begin{cases}
        \Projector = \Decoder\circ \Encoder,\\
        \identity_{d_{r}} = \Encoder\circ\Decoder,
    \end{cases}
\end{equation}
where $\identity_{d_{r}}$ is the identity matrix in $\RedCoordParamSpace$. We refer to $\RedCoordParamSpace=\Encoder(\ParamSpace)$ as the \emph{latent parameter vector space}. The \emph{latent prior} $\pi\in\ProbSpace(\RedCoordParamSpace)$ and \emph{latent posterior} $\widetilde{\pi}^{\Obs}\in\ProbSpace(\RedCoordParamSpace)$ are defined as the pushforward of the prior marginal $\FullPrior_r$ and the subspace posterior $\widetilde{\FullPrior}_r^{\Obs}$, respectively, by the encoder $\Encoder$, and they satisfy Bayes' rule of the latent parameters:
\begin{equation}\label{eq:reduced_inference}
\begin{cases}
    \RedCoordPrior \coloneqq {\Encoder}_{\sharp} \FullPrior_{r}= \Gaussian(0, \Encoder\PriorCov\Encoder^*)\\
    \widetilde{\RedCoordPrior}^{\Obs} \coloneqq  {\Encoder}_{\sharp} \widetilde{\FullPrior}^{\Obs}_r
\end{cases},\quad
    \widetilde{\RedCoordPrior}^{\Obs}(\RedCoordParam)=\frac{1}{\widetilde{\normalization}^{\Obs}_r} \exp\left(-\RidgeMisfit(\Decoder\RedCoordParam)\right)\RedCoordPrior(\RedCoordParam).
\end{equation}
Given latent posterior samples $\RedCoordParam^{(j)}\iid\widetilde{\RedCoordPrior}^{\Obs}$ and prior conditional samples $m_{\perp}^{(j)}\iid\mu_{\perp|r}(\cdot|\Decoder\RedCoordParam^{(j)})$, we obtain approximate posterior samples $\FullParam^{(j)}\iid\widetilde{\FullPrior}^{\Obs}$, where $\FullParam^{(j)} = \Decoder\RedCoordParam^{(j)} + m_{\perp}^{(j)}$.

While there are many choices of reduced basis, we consider a class of $\ParamCM$--orthonormal reduced basis given as follows:
\begin{equation}\label{eq:as_encoder_decoder_def}
        \left\langle\psi_j, \psi_k\right\rangle_{\PriorCov^{-1}} = \delta_{jk},\quad j,k = 1,\dots,\RedCoordDim,\quad
        \Encoder m = \sum_{j=1}^{\RedCoordDim}\left\langle m, \psi_j\right\rangle_{\PriorCov^{-1}}\be_j,
\end{equation}
where $\delta_{jk}$ is the Kronecker delta, and $\be_j$ is the unit vector in the latent space $\RedCoordParamSpace$. Through the definition of the Hermitian adjoint on $\ParamSpace$, we have $\Encoder^* = \PriorCov^{-1}\Decoder$ and $\Decoder^* = \Encoder\PriorCov$, which implies $\Encoder\PriorCov\Encoder^* = \identity_{\R^{d_{r}}}$ and $\Projector^*=\PriorCov^{-1}\Projector\PriorCov$ due to \cref{eq:generic_encoder_decoder_definition}. Consequently, \cref{eq:reduced_inference,eq:close_form_prior_disintegration} yield
\begin{align}
    \pi &= \Gaussian(\boldsymbol{0}, \identity_{\R^{d_{r}}}), &&&(\text{whitened latent prior})\label{eq:whitened_subsapce_prior}\\
    \FullPrior_{{\perp}| r} (\cdot|m_r) &\equiv \FullPrior_{\perp},  &&&(\text{independence of marginals})\label{eq:marginal_independence}
\end{align}
where $\FullPrior_{\perp}=(\identity_{\ParamSpace}-\Projector)_{\sharp}\mu$ is the prior marginal in $\text{Ker}(\Projector)$. As a result of \cref{eq:marginal_independence}, sampling of the conditional $m_{\perp}^{(j)}\iid\FullPrior_{{\perp}| r}(\cdot|m_r^{(j)}) $ can be accomplished via sampling the full prior:
\begin{equation}\label{eq:sampling_fullpost}
    m_{\perp}^{(j)} = m_\text{pr}^{(j)} - \Projector m_{\text{pr}}^{(j)},\quad m_\text{pr}^{(j)}\sim \mu.
\end{equation}

%% file: transport.tex
\subsection{Lazy map variational inference}\label{subsec:lazy_map}
We consider TMVI that seeks a diffeomorphic deterministic coupling between a target and a reference distribution \cite{marzouk2016introduction}. For our BIP in \eqref{eq:generic_full_bip},  we aim to find a transport map $\FullMap\in\text{Diff}^1(\ParamSpace)$ that couples our reference, the prior $\mu\in\ProbSpace(\ParamSpace)$, to the target, the posterior $\mu^{\Obs}\in\ProbSpace(\ParamSpace)$:
\begin{equation*}
    \FullMap_{\sharp}\mu = \mu^{\Obs},\quad \FullMap^{\sharp} \mu^{\Obs}= \mu.
\end{equation*}
Given $\FullMap$, sampling from the posterior $m_{\text{post}}^{(j)}\iid \mu^{\Obs}$ is accomplished by evaluating $m_{\text{post}}^{(j)} = \FullMap(m_{\text{pr}}^{(j)})$, where $m_{\text{pr}}^{(j)}\iid\mu$. Since $\FullMap$ is typically unavailable in closed form for nonlinear BIPs, we consider classes of transport maps parametrized by weights $\MapParams\in\R^{d_{\theta}}$ such that ${\FullMap_{\MapParams}}_{\sharp}\mu$ and $\mu$ are mutually absolutely continuous. These weights are found via the solution of a stochastic optimization problem, whose goal is to find ${\FullMap_{\MapParams}}_{\sharp}\mu \approx \mu^{\Obs}$. The reverse Kullback-Leibler (rKL) divergence of the approximating distribution from the target posterior $\mu^{\Obs}$
\begin{equation*}
    \KL({\FullMap_{\MapParams}}_\sharp \mu || \mu^{\Obs}) =  \KL(\mu || {\FullMap_{\MapParams}}^\sharp \mu^{\Obs}) = \mathbb{E}_{m\sim\mu}\left[\text{log}\left(\frac{\meas\mu}{\meas(\mu^{\Obs}\circ\FullMap_{\MapParams})}(m)\right)\right]
\end{equation*}
is often used to measure the error of transport map posterior approximation and thereby employed as the objective function for optimization \cite{marzouk2016introduction, blei2017variational}. This objective is equivalent to the \emph{evidence lower bound} objective function.

To make TMVI tractable when $\ParamSpace$ has high or infinite dimensions, we use \emph{lazy maps}, proposed in \cite{brennan2020greedyinferencestructureexploitinglazy}, which leverages the subspace decomposition as in \cref{subsec:latent_space}:
\begin{align}\label{eq:lazymap}
    \Map_{\MapParams}\coloneqq\Map(\cdot,\MapParams)\in \scrT\subset\text{Diff}^1(\RedCoordParamSpace),\quad
    \FullMap_{\MapParams} \coloneqq \overbrace{(\identity_{\ParamSpace} -\Projector)}^{\mathclap{\text{identity in } \text{Ker}(\Projector)}} + \underbrace{\Decoder\circ \Map_{\MapParams}\circ\Encoder}_{\mathclap{\text{nonlinear transport in } \text{Im}(\Projector)}},
&&&(\text{Lazy Map})
\end{align}
where a latent space nonlinear transport is used to represent the coupling of the prior and the posterior in $\text{Im}(\Projector)$ while the prior is preserved in $\text{Ker}(\Projector)$. Here we consider a parametrized class $\mathscr{T}\subset\text{Diff}^1(\R^{d_r})$ with weights $\MapParams\subseteq\R^{d_{\MapParams}}$ that allows ${\Map_{\MapParams}}_{\sharp}\pi$ and $\pi$ to be mutually absolutely continuous and assume the diffeomorphic property is achieved by constraining the map to satisfy $\det\grad_{\refCoord}\Map_{\MapParams}(\refCoord)>0$ a.e.; see~\cite{marzouk2016introduction}.

The following proposition shows an equivalence in rKL between a lazy map defined on the parameter space $\ParamSpace$, and a transport map defined on the latent space $\RedCoordParamSpace$ through the optimal ridge function.
\cref{eq:reduced_inference}
\begin{proposition}\label{theorem:lazy_map_kl}
Given a linear projection $\Projector$ defined using a $\ParamCM$-orthonormal reduced basis and a latent space transport $\Map_{\MapParams}\in\scrT$, we have
    \begin{align}
  \KL({\FullMap_{\MapParams}}_{\sharp}\mu||\mu^{\Obs}) &= \KL({\Map_{\MapParams}}_{\sharp}\pi||\RedCoordPostOpt) + C_1\nonumber\\
        & = \mathbb{E}_{\refCoord\sim\pi}\left[\left(\RidgeMisfitOpt\circ\Decoder\circ \Map_{\MapParams}\right)(\refCoord) + \frac{1}{2}\norm{ \Map_{\MapParams}(\refCoord)}^2 -\log\det\grad_{\refCoord} \Map_{\MapParams}(\refCoord) \right] + C_2,\label{eq:reverse_KL}
    \end{align}
    where $\FullMap_{\MapParams}$ is the lazy map in \cref{eq:lazymap}$, \pi=\mathcal{N}(\boldsymbol{0},\identity_{\RedCoordParamSpace})$ is the whitened latent prior in \cref{eq:whitened_subsapce_prior}, $\RidgeMisfitOpt$ and $\RedCoordPostOpt$ are the approximate potential and latent posterior induced by $\FullRidgeForwardOpt\circ\Projector$ in \eqref{eq:conditional_expectation_ridge}, and $C_1$ and $C_2$ are constants that do not depend on $\MapParams$.
\end{proposition}
The proof of \cref{theorem:lazy_map_kl} is provided in \cref{inf_map}. Due to this result, we may formulate LMVI as
the following optimization problem with an equivalent latent representation
\begin{subequations}\label{eq:lazy_map_variational_inference}
\begin{align}
    &\min_{\MapParams\in\R^{\MapParams}} \Dkl({\FullMap_{\MapParams}}_{\sharp}\mu||\mu^{\Obs}), &&& (\text{Lazy map variational inference})\\
    &\min_{\MapParams\in\R^{\MapParams}} \Dkl({\Map_{\MapParams}}_{\sharp}\pi||\RedCoordPostOpt). &&& (\text{Equivalent latent representation})
\end{align}
\end{subequations}
The latent representation can be treated as a TMVI problem that seeks $\Map\in\text{Diff}^1(\RedCoordParamSpace)$ defined as follows.
\begin{align*}
    &{\Map}_{\sharp}\RedCoordPrior=\RedCoordPostOpt,
        \quad\Map_\sharp \RedCoordPrior (\RedCoordParam) = (\RedCoordPrior \circ \Map^{-1})(\RedCoordParam)  |\det\nabla \Map^{-1}(\RedCoordParam)|, &&& (\text{Latent space pushforward})\\
	&{\Map}^{\sharp}\RedCoordPostOpt=\RedCoordPrior,\quad\Map^\sharp \RedCoordPostOpt(\refCoord) =(\RedCoordPostOpt \circ \Map)(\refCoord)  |\det\nabla \Map(\refCoord)|. &&& (\text{Latent space pullback})
\end{align*}
Given such a transport map, sampling $\bx^{(j)}\iid\RedCoordPostOpt$ is accomplished by evaluating $\RedCoordParam^{(j)}=\Map(\refCoord^{(j)})$, where $\refCoord^{(j)} \iid\RedCoordPrior$. In turn, approximate posterior, or \emph{pushforward}, samples $\FullParam^{(j)}\iid \FullPostOpt$ can be drawn simply by lifting $\RedCoordParam^{(j)}$ into $\ParamSpace$ following \cref{eq:sampling_fullpost}.

\begin{table}[!htbp]
    \caption{A summary of BIP problems discussed in Section 2.}
    \begin{center}
    {\renewcommand{\arraystretch}{2.2}
    \begin{tabular}{ |c | c | c | c | c | c| }\hline
    BIP Name & \multicolumn{2}{c|}{Approximation} & \makecell{Parameter\\ space} & \multicolumn{2}{c|}{Prior and Posterior}  \\\hline\hline
    Original & \multicolumn{2}{c|}{None} & $\scrM$ &$\FullPrior, \FullPrior^{\Obs}$ & \cref{eq:generic_full_bip}\\\hline
    Subspace& $\Forward\approx\FullRidgeForward\circ\Projector$ & \cref{eq:ridge_definition} & $\text{Im}(\Projector)$ & $\FullPrior_r, \widetilde{\FullPrior}_r^{\Obs}$ & \cref{eq:subspace_ridge_posterior} \\\hline
    Latent & \makecell{$\Forward\approx\FullRidgeForward\circ\Projector$\\$\Projector=\Decoder\circ\Encoder$} & \makecell{\cref{eq:ridge_definition}\\\cref{eq:generic_encoder_decoder_definition}} & $\RedCoordParamSpace$ &$\RedCoordPrior,\widetilde{\RedCoordPrior}^{\Obs}$ & \cref{eq:reduced_inference}  \\\hline
    \makecell{Lazy map\\latent representation} & \makecell{$\Forward\approx\FullRidgeForwardOpt\circ\Projector$\\$\Projector=\Decoder\circ\Encoder$} & \makecell{\cref{eq:conditional_expectation_ridge}\\\cref{eq:generic_encoder_decoder_definition}} & $\RedCoordParamSpace$ & $\RedCoordPrior,\RedCoordPostOpt$ & \cref{eq:reduced_inference}  \\\hline
    \end{tabular}
    }
    \end{center}
\end{table}

\subsection{Stochastic optimization of lazy map and challenges}\label{subsec:optimizing_lazy_map}
For a given transport map parametrization $\Map_{\MapParams}\in \MapClass$, the map parameter vector $\MapParams$ are typically found via gradient-based stochastic optimization, which in turn requires evaluating Monte Carlo (MC) estimates of the gradient of the objective with respect to $\MapParams$. The shifted rKL objective, denoted as $\mathcal{L}^{\Obs}:\R^{d_{\MapParams}}\to\R$, can be expressed as an expectation of a single-sample MC estimator, $\mathcal{L}_1^{\Obs}:\ParamSpace\times\R^{d_{\MapParams}}\to\R$, defined as follows: 
\begin{align}
     \mathcal{L}^{\Obs}(\theta) \coloneqq& \mathbb{E}_{\FullParam\sim\FullPrior} \left[\mathcal{L}_1^{\Obs}(m,\MapParams) \right]:=\Dkl({\FullMap_{\MapParams}}_{\sharp}\mu||\mu^{\Obs}) - C_2,\\
    \mathcal{L}_1^{\Obs}(\Decoder\refCoord^{(j)} + m_{\perp}^{(j)}, \theta) =& \Misfit\left(\left(\Decoder\circ \Map_{\theta}\right)(\refCoord) + m_{\perp}^{(j)}\right) + \frac{1}{2}\norm{ \Map_{\theta}(\refCoord^{(j)})}^2 -\log\det\grad_{\refCoord} \Map_{\theta}(\refCoord^{(j)}),
    \label{eq:mtvi_loss}
\end{align}
where $\refCoord^{(j)}\iid\pi$ and $m_{\perp}^{(j)}\iid\mu_{\perp}$. 
The single-sample MC gradient estimator with respect to the map parameters $\MapParams$ takes the following form
\begin{equation}\label{eq:mtvi_gradient}
    \begin{aligned}
    \nabla_{\MapParams} \mathcal{L}^{\Obs}_1(\Decoder\refCoord^{(j)} + m_{\perp}^{(j)},\theta) &= \nabla_{\MapParams} \Map_{\MapParams}(\refCoord^{(j)})^{\transpose}(\Encoder\circ D_H\Misfit)\left((\Decoder\circ\Map)(\refCoord^{(j)}) + m_{\perp}^{(j)}\right)\\
    &\quad  + \nabla_{\MapParams}\Map_{\MapParams}(\refCoord^{(j)})^{\transpose}\Map_{\MapParams}(\refCoord^{(j)}) -
    \nabla_{\MapParams}(\log\det\grad_{\refCoord} \Map_{\theta})(\refCoord^{(j)}),
\end{aligned}
\end{equation}
where $D_H\Misfit(m)\coloneqq D_H\Forward(m)^*(\Forward(m) - \Obs)$ is the prior-preconditioned gradient of the potential; see \cref{eq:gateaux_versus_stochastic_gateaux}.

The MC gradient estimator of the rKL objective is then
\begin{equation*}
    \nabla_{\MapParams}\mathcal{L}^{\Obs}(\MapParams) \approx  \nabla_{\MapParams}\widehat{\mathcal{L}}^{\Obs}(\MapParams) = \frac{1}{N_\text{MC}} \sum_{j=1}^{N_\text{MC}} \nabla_{\MapParams}\mathcal{L}^{\Obs}_1(m^{(j)},\MapParams),\quad m^{(j)}\iid\mu.
\end{equation*}

\begin{remark}
In the numerical results, starting in~\cref{sec:numerical_results}, the lazy map optimization is implemented with an alternative form of the rKL objective where samples from $\mu_{\perp}$ are estimated as $\mathbb{E}[\mu_{\perp}] = 0$:
    \begin{equation*}
    \mathcal{L}^{\Obs}(\theta) \coloneqq \mathbb{E}_{\refCoord\sim\pi} \left[\mathcal{L}_1^{\Obs}(\Decoder\refCoord,\MapParams) \right].
    \end{equation*}
This leads to a different TMVI problem induced by the ridge function $\Forward\circ\Projector$ instead of $\FullRidgeForwardOpt\circ\Projector$ in \cref{theorem:lazy_map_kl}. We empirically found that it performs better under a limited computational budget, likely due to a superior bias--variance trade-off in TMVI.
\end{remark}

The computation cost of evaluating the second and third terms in \eqref{eq:mtvi_gradient} only depends on the parametrization of the transport map. Notably triangular transport maps \cite{marzouk2016introduction} and parameterization built as compositions of triangular maps (e.g., inverse autoregressive flows~\cite{kingma2017improvingvariationalinferenceinverse}), are structured such these terms are efficiently computable. The first term requires evaluating the PtO map and its prior-preconditioned gradient. When the PtO map is expensive to evaluate, as is the case with large-scale PDE-governed problems, optimizing for an accurate lazy map is prohibitively expensive. 

%% file: dino_approximated_likelihood.tex
\section{Optimized surrogate construction for lazy map variational inference}\label{section:dino_approximated_likelihood}

This section discusses the construction of a fast-to-evaluate neural network ridge function surrogate $\FullRidgeForward\circ\Projector\approx\Forward$ that leads to a small and controllable expected error in surrogate-driven LMVI, where the expectation is taken over the marginal distribution of data vectors with density $\gamma\in\ProbSpace(\R^\ObsDim)$ where $\gamma(\by)\propto Z^{\by}$ as in \cref{eq:data_model}. Our strategy for constructing this surrogate is given as follows.
\begin{enumerate}
\item Minimizing an upper bound on $\mathbb{E}_{\Obs\sim\gamma}[\Dkl( \FullPost||\widetilde{\FullPrior}^{\Obs})]$, the expected forward KL divergence (fKL) from the posterior $\FullPost$ to the approximate posterior defined by the surrogate, $\widetilde{\mu}^{\Obs}$.
\item Minimizing an upper bound on the expected optimality gap $\mathbb{E}_{\Obs\sim\gamma}\left[\sqrt{\mathcal{L}^{\Obs}(\widetilde{\MapParams}^{\by,\dagger}) - \mathcal{L}^{\Obs}(\MapParams^{\by,\dagger})}\right]$ for surrogate-driven LMVI, where $\MapParams^{\by,\dagger}$ is the true minimizer of the rKL objective and $\widetilde{\MapParams}^{\by,\dagger}$ is the minimizer found via the ridge function surrogate.
\end{enumerate}
In this section, we show the resulting ridge function surrogate is \texttt{DIPNet} \cite{o2022derivative} trained using the derivative-informed learning method \cite{o2024derivative} and it leads to a latent space surrogate rKL objective.

\subsection{Error analysis for surrogate-driven lazy map variational inference}\label{error_analysis_section}
Recall the definition of $\FullRidgeForwardOpt$ in \eqref{eq:conditional_expectation_ridge}. We define the following finite-dimensional latent representations as follows: 
\begin{subequations}
\begin{align}\label{eq:ridge_latent_representation}
    \RidgeForward \quad\,\coloneqq \boldsymbol{V}^*\circ\FullRidgeForward \circ\Decoder:\quad\;\R^{\RedCoordDim}\to \ObsSpace,&&&
    \FullRidgeForward \quad\,=   \boldsymbol{V}\circ\RidgeForward\circ \Encoder: \quad\:\text{Im}(\Projector)\to\ObsCM,&&&&&&&&&&&\\
     \RidgeForwardOpt \coloneqq \boldsymbol{V}^*\circ\FullRidgeForwardOpt \circ\Decoder:\R^{\RedCoordDim}\to \ObsSpace,&&&
    \FullRidgeForwardOpt =   \boldsymbol{V}\circ\RidgeForwardOpt\circ \Encoder: \text{Im}(\Projector)\to\ObsCM.&&&&&&&&&&&
\end{align}
\end{subequations}
Here $\boldsymbol{V}\in \HS(\R^{\ObsDim}, \ObsCM)$ is a full-rank matrix with columns consists of $\ObsCM$-orthonormal basis and $\boldsymbol{V}^*=\boldsymbol{V}^{\transpose}\NoiseCov^{-1}$ is its Hermitian adjoint that satisfies $\boldsymbol{V}^*\boldsymbol{V}=\identity_{\ObsDim}$. Note that $\boldsymbol{V}^*$ is a whitening transformation on the data space.

The following theorem provides an upper bound on the expected fKL between the true posterior and ridge function approximated posterior taken over the marginal distribution of data.

\begin{theorem}[Posterior approximation through a ridge function surrogate] \label{theorem:kl_bound}
Given $\Forward\in H^1_{\mu}(\ParamSpace; \ObsCM)$ and a projector $\mathcal{P}\in B(\ParamSpace)$ defined via an $\ParamCM$-orthonormal reduced basis as in \cref{eq:as_encoder_decoder_def}, we have the following inequality for the approximate posterior $\widetilde{\FullPrior}^{\Obs}$ in \cref{eq:subspace_ridge_posterior} defined via any ridge function $\FullRidgeForward\circ\Projector\in L^2_{\mu}(\ParamSpace; \ObsCM)$:
\begin{equation*}
\begin{aligned}
    \mathbb{E}_{\Obs\sim\gamma} \left[ \KL( \FullPost||\widetilde{\FullPrior}^{\Obs})\right] &\leq  \underbrace{\Tr_{\ParamCM}\left(\left(\identity_{\ParamCM}-\Projector\right)\mathcal{H}_A\left(\identity_{\ParamCM}-\Projector\right)\right)}_{\mathclap{\text{\normalfont\ignorespaces Parameter reduction error}}} + \underbrace{\mathbb{E}_{\refCoord\sim\pi}\left[\norm{\RidgeForwardOpt(\refCoord)-\RidgeForward(\refCoord)}^2\right]}_{\mathclap{\text{\normalfont\ignorespaces Latent representation error}}},
\end{aligned}
\end{equation*}
where
$\Tr_{\ParamCM}:\TrClass(\ParamCM)\to\R$ returns the trace of operators on $\ParamCM$, and $\mathcal{H}_A\in\TrClass(\ParamCM)$ is the expected prior--preconditioned Gauss--Newton Hessian of the potential:
\begin{equation}\label{eq:expected_gauss-newton_hessian}
    \mathcal{H}_A\coloneqq \mathbb{E}_{m \sim \mu} \left[ D_H \Forward(m)^*D_H \Forward(m) \right].
\end{equation}
Here $D_H\Forward(m)^*\in\HS(\ObsCM, \ParamCM)$ denotes the Hermitian adjoint of the stochastic derivative $D_H\Forward(m)$ in \cref{eq:stochastic_derivative}.
\end{theorem}
The bound decomposes the expected error into terms involving the parameter reduction error that depends on the choice of $\Projector$ and the discrepancy between the surrogate latent representation $\RidgeForward$ and the optimal latent representation $\RidgeForwardOpt$ of $\FullRidgeForwardOpt\circ\Projector$. The proof of \cref{theorem:kl_bound} can be found in \cref{app:kl_bound_proof} and follows from results in \cite{cui2021data,cao2024efficient, zahm2020gradient}.

For surrogate-driven LMVI, understanding the expected discrepancy between the posteriors and its transport targets (i.e., the surrogate approximated posteriors $\widetilde{\mu}^{\Obs}$) is insufficient, as the transport map is constructed through the process of gradient-based stochastic optimization and the accuracy of the rKL gradient approximation is also important. The following theorem establishes error upper bounds for the surrogate objective gradient.

\begin{theorem}[Surrogate approximation of the rKL objective gradient]\label{theorem:kl_gradient} Given $\Forward\in H^1_{\mu}(\ParamSpace;\ObsCM)$, a linear projector $\mathcal{P}\in B(\ParamSpace)$ defined using an $\ParamCM$-orthonormal reduced basis as in \cref{eq:as_encoder_decoder_def}. Assume we have a latent space transport $\Map_{\MapParams}\in\scrT$ with an essentially bounded density between ${\Map_{\MapParams}}_{\sharp}\pi$ and $\pi$ and an essentially-bounded Jacobian with respect to $\MapParams$. We have the following error upper bound for the approximate gradient of rKL objective $\widetilde{\mathcal{L}}^{\Obs}(\MapParams)$ given by a ridge function,
\begin{equation}
    \mathbb{E}_{\Obs\sim\gamma}\left[\|\nabla_{\MapParams} \mathcal{L}^{\Obs}(\MapParams) - \nabla_\MapParams \widetilde{\mathcal{L}}^{\Obs}(\MapParams)\|\right] \lesssim\left(\mathbb{E}_{\refCoord\sim\pi} \Big[  \|\RidgeForwardOpt(\refCoord) - \RidgeForward(\refCoord) \|^2 + \|\nabla\RidgeForwardOpt(\refCoord) - \nabla\RidgeForward (\refCoord) \|^2_{F}        \Big]\right)^{1/2},
\end{equation}\label{eq:gradient_err_bound}
where $\lesssim$ denotes bounded up to a multiplicative constant.
\end{theorem}
The proof of \cref{theorem:kl_gradient} is presented in \cref{app:kl_gradient_proof}. Our result states that the expected gradient error is controlled by the latent representation error measured in a $\pi$-weighted Sobolev norm on $H^1_{\pi}(\R^{d_r};\R^{d_y})$, which additionally contain the expected error in the Jacobian compared to the error measure using $\pi$-weighted Bochner norm on $L^2_{\pi}(\R^{d_r};\R^{d_y})$ in \cref{theorem:kl_bound}. Notably, the two error measures are generally not equivalent, and the Sobolev norm is stronger than the Bochner norm. This result reflects the fact that the gradient of the rKL objective involves the Jacobian of the PtO map, and the surrogate Jacobian accuracy affects the optimization of lazy maps. To further explore the consequences of surrogate Jacobian misfit, we consider the following corollary on the expected optimality gap for surrogate-driven LMVI under a stronger set of assumptions.

\begin{corollary}[Optimality gap for surrogate-driven LMVI]\label{corollary:optimality_gap}
Suppose the assumptions in \cref{theorem:kl_gradient} holds. Let $\MapParams^{\Obs,\dagger}$ and $\widetilde{\MapParams}^{\Obs,\dagger}$ denote $\gamma$-measurable functions that return the second order stationary points of $\mathcal{L}^{\Obs}$ and $\widetilde{\mathcal{L}}^{\Obs}$ $\gamma$-a.e., respectively. Let $B_r(\boldsymbol{x})$ denote a ball of radius $r$ centered at $\boldsymbol{x}$. We assume that $R^{\Obs}$ and $\lambda^{\Obs}$ are $\gamma$-essentially bounded from below by some positive constants such that (i) $\widetilde{\MapParams}^{\Obs,\dagger}\in B_{R^{\Obs}}(\MapParams^{\Obs,\dagger})$ $\gamma$-a.e., and (ii) $\nabla^2_\MapParams \mathcal{L}^{\Obs}(\MapParams) \succeq \lambda^{\Obs} \identity_{\R^{d_r}}$ for all $\MapParams \in B_{R^{\Obs}}(\MapParams^{\Obs,\dagger})$ $\gamma$-a.e.

We have the following upper bound on the optimality gap: 
\begin{equation}\label{loss_bound}
    \mathbb{E}_{\Obs\sim\gamma}\left[\sqrt{\mathcal{L}^{\Obs}(\widetilde{\MapParams}^{\Obs,\dagger}) - \mathcal{L}^{\Obs}(\MapParams^{\Obs,\dagger})}\right] \lesssim \left(\mathbb{E}_{\refCoord\sim\pi}\Big[  \|\RidgeForwardOpt(\refCoord) - \RidgeForward (\refCoord) \|^2 + \|\nabla\RidgeForwardOpt(\refCoord) - \nabla\RidgeForward(\refCoord) \|^2_{F} \Big]\right)^{1/2}.
\end{equation}
\end{corollary}
This result states that if the rKL objective is locally strongly convex near the true rKL minimizers and the minimizers found by the surrogate lands within those locally convex regions, we can bound the expected optimality gap by the latent representation error measured by the weighted Sobolev norm. The assumptions in \cref{corollary:optimality_gap} are commonly employed in the optimization and machine learning literature, either directly via local strong convexity or through the Polyak-\L{}ojasiewicz inequality, see for example \cite{bottou2018optimization}. These two results together demonstrate the need to control the surrogate Jacobian error in the context of LMVI.

Motivated by these results, we delineate the procedure for constructing a surrogate model for LMVI in the following subsections.

\subsection{Minimize the parameter reduction error: derivative-informed subspace}\label{sec:AS_dim_reduction}

We seek an $\ParamCM$-orthonormal reduced basis $\{\psi_j\}_{j=1}^{\RedCoordDim}$ such that $\text{span}(\{\psi_j\}_{j=1}^{\RedCoordDim})=\text{Im}(\Projector)$ and the parameter reduction error term in \cref{theorem:kl_bound} is minimized. This can be accomplished by finding the parameter subspace that the PtO map is most sensitive to in expectation. 
This subspace is often referred to as the derivative-informed subspace or active subspace \cite{zahm2020gradient}, and it can be computed from the dominant $\RedCoordDim$ eigenbases arising from the following eigenvalue problem in $\ParamCM$,

\begin{align}\label{eq:derivative_evp}
    \mathcal{H}_A\psi_j = \lambda_j\psi_j,\quad 
    \langle \psi_j,\psi_k\rangle_{\mathcal{C}^{-1}} = \delta_{jk},\quad \lambda_1\geq\lambda_2\geq\dots \geq 0
\end{align}
where $\mathcal{H}_A$ is the prior-preconditioned Gauss--Newton Hessian in \cref{eq:expected_gauss-newton_hessian}.
Under a stronger Fr\'echet differentiability assumption, the eigenvalue problem in $\ParamCM$ is equivalent to a more common form of a generalized eigenvalue problem in $\ParamSpace$ due to \cref{eq:gateaux_versus_stochastic_gateaux}:
\begin{align}\label{eq:derivative_gevp}
     \mathbb{E}_{m \sim \mu} \left[ D \Forward(m)^*\NoiseCov^{-1}D\Forward(m) \right]\psi_j = \lambda_j\PriorCov^{-1}\psi_j,\quad
    \langle \psi_k, \psi_j\rangle_{\mathcal{C}^{-1}} = \delta_{jk},\quad\lambda_1\geq\lambda_2\geq\dots\geq 0.
\end{align}  
The eigenvalue problem in \cref{eq:derivative_gevp} can be found in \cite{Isaac_2015, cui2021data,zahm2020gradient,zahm2022certified, bigoni2022nonlinear}. The minimum value of the parameter reduction error is
\begin{equation}\label{eq:as_ev_tail_sum}
    \min_{\substack{\Projector\in\{\text{rank}-\RedCoordDim\text{ linear} \\ \text{projection on $\ParamSpace$\}}}} \Tr_{\ParamCM}\left(\left(\identity_{\ParamCM}-\Projector\right)\mathcal{H}_A\left(\identity_{\ParamCM}-\Projector\right)\right) =  \sum_{j>\RedCoordDim}\lambda_j.
\end{equation}
This derivative-based reduced basis leads to an expected parameter reduction error proportional to the eigenvalue tail sum in \cref{eq:derivative_evp} corresponding to the discarded eigenbases. Existing bounds for the truncated Karhunen–Lo\'eve expansion of the parameter are strictly higher than \cref{eq:as_ev_tail_sum}; see \cite{zahm2020gradient,cao2024efficient}.

\subsection{Minimize the latent representation error: Conventional operator learning}\label{ridge_approx_subsection}
We first consider a neural operator ridge function using a neural network latent representation $\RidgeForward_{\text{NN}}:\R^{\RedCoordDim}\times\R^{d_{\boldsymbol{w}}}\to\ObsSpace$:
\begin{subequations}\label{eq:neural_ridge_function}
\begin{align}
        \RidgeForward_{\boldsymbol{w}}(\RedCoordParam) &\coloneqq\RidgeForward_{\text{NN}}(\RedCoordParam, \bw), &&& (\text{Neural latent representation})\\
        \FullRidgeForward_{\boldsymbol{w}}(\Projector m) &\coloneqq\boldsymbol{V}\RidgeForward_{\text{NN}}(\Encoder m, \boldsymbol{w}), &&& (\text{Neural operator ridge function})
\end{align}
\end{subequations}
where $\boldsymbol{w} \in \mathbb{R}^{d_{\boldsymbol{w}}}$ consists of trainable neural network weights. Neural network surrogates architecture using the derivative-informed subspace \cref{eq:derivative_gevp} are known as \texttt{DIPNet} \cite{o2022derivative}.

Motivated by \cref{theorem:kl_bound}, it would be sensible to find the neural network weights by minimizing the latent representation error, which also minimizes the upper bound on the expected surrogate posterior approximation error:
\begin{align}
    \min_{\boldsymbol{w}\in\R^{d_{\boldsymbol{w}}}} \mathbb{E}_{\refCoord\sim \pi} \left[\left\|\RidgeForwardOpt(\refCoord) - \RidgeForward_{\boldsymbol{w}}(\refCoord)\right\|^2\right].\label{eq:neural_latent_represnetation_learning}
\end{align}
However, estimating this objective function requires a nested MC method due to the $\text{Ker}(\Projector)$ marginalization in $\RidgeForwardOpt$; see definitions in \cref{eq:ridge_definition,eq:ridge_latent_representation}. Specifically, $\refCoord^{(j)}\sim\pi$, $1\leq j\leq N_{\text{out}}$, are used to estimate the objective, and $\FullParam_{\perp}^{(j,k)}\sim\mu_{\perp}$, $1\leq k\leq N_{\text{in}}$, are used to estimate the output of the optimal latent representation at each $\refCoord^{(j)}$. The nested MC sample generation requires $N_{\text{out}}\times N_{\text{in}}$ PtO map evaluations. However, the inner MC is unnecessary when $\Projector$ is chosen as in \cref{sec:AS_dim_reduction}, since the PtO map is insensitive to changes in $\text{Ker}(\Projector)$; see, e.g., ~\cite[Corollary 7.5]{zahm2022certified}.
Therefore, we consider the conventional operator learning method with error measure using the norm on the Bochner space $L^2_\mu(\ParamSpace,\ObsCM)$, i.e., a prior-weighted mean squared error objective:
\begin{align} \label{eq:l2_training}
    &\min_{\boldsymbol{w\in\R^{d_{\boldsymbol{w}}}}} \mathbb{E}_{m\sim\mu} \left[\left\|\Forward(m) - \FullRidgeForward_{\boldsymbol{w}}(\Projector m)\right\|^2_{\NoiseCov^{-1}}\right], &&& (\text{Conventional $L^2_{\mu}$ operator learning})\\
    &\min_{\boldsymbol{w\in\R^{d_{\boldsymbol{w}}}}} \mathbb{E}_{(\refCoord,m_{\perp})\sim\pi\otimes\mu_{\perp}}\bigg[\Big\| \underbrace{ \boldsymbol{V}^*\Forward(\Decoder \refCoord + m_{\perp})}_{\displaystyle\mathclap{\approx \RidgeForwardOpt(\refCoord)} } - \RidgeForward_{\boldsymbol{w}}(\refCoord)\Big\|^2\bigg]. &&& (\text{Equivalent latent representation})\label{eq:l2_learning_latent} 
\end{align}
The equivalent latent representation of the operator learning objective reveals that this objective can be derived from the neural latent representation error \cref{eq:neural_latent_represnetation_learning} using a single sample ($m_{\perp} \sim\mu_{\perp}$) estimate of the marginalization in $\RidgeForwardOpt$. For notational convenience, we will use $\bg^{(j)}$ to denote the i.i.d.\ \emph{whitened PtO samples} used to estimate the conventional $L^2_{\mu}$ operator learning objective:
\begin{align} \label{eq:latent_parameter}
    \boldsymbol{g}^{(j)} \coloneqq \boldsymbol{V}^*\Forward(m^{(j)})\in \mathbb{R}^{\ObsDim},\quad m^{(j)}\iid\mu.&&& (\text{whitened PtO sample})
\end{align}

We refer to \texttt{DIPNet} surrogates trained using the conventional $L^2_{\mu}$ operator learning method in \cref{eq:l2_training} as \texttt{RB-NO} (reduced basis neural operator) in contrast to surrogates construction introduced in the following subsection.

\subsection{Minimizing the expected optimality gap: Derivative-informed operator learning}\label{H1mu}

While the conventional $L^2_\mu$ learning problem presented in \cref{ridge_approx_subsection} is suitable for constructing a neural operator ridge function, \cref{theorem:kl_gradient} shows that controlling latent representation error in $H^1_{\pi}$ controls both the expected gradient error, as well as the expected optimality gap between the exact and surrogate variational inference objective functions. To this end, we consider minimizing the latent representation error measured by the $\pi$-weighted Sobolev norm on $H^1_{\pi}(\R^{\RedCoordDim}; \ObsSpace)$:
\begin{align}\label{eq:h1_learning_latent}
    \min_{\boldsymbol{w}\in\R^{d_{\boldsymbol{w}}}} \mathbb{E}_{\refCoord\sim \pi} \left[\left\|\RidgeForwardOpt(\refCoord) - \RidgeForward_{\boldsymbol{w}}(\refCoord) \right\|^2 + \left\|\grad\RidgeForwardOpt(\refCoord) - \grad_{\refCoord}\RidgeForward_{\boldsymbol{w}}(\refCoord) \right\|_F^2\right].
\end{align}
However, as discussed in the previous subsection, estimating the objective function above also requires a nested MC method. To circumvent this issue, we adopt an derivative-informed $H^1_{\mu}$ operator learning objective following \cite{o2024derivative,cao2024efficient}:
\begin{align} \label{eq:h1_objective}
\begin{split}
    &\min_{\boldsymbol{w}\in\R^{d_{\boldsymbol{w}}}} \mathbb{E}_{m\sim\mu}\Big[\left\|\Forward(m) - \FullRidgeForward_{\boldsymbol{w}}(\Projector m)\right\|^2_{\NoiseCov^{-1}}\\
    &\qquad + \left\|D \Forward(m) - D (\FullRidgeForward_{\boldsymbol{w}}\circ\Projector)(m)\right\|^2_{\HS(\ParamCM,\ObsCM)} \Big]
\end{split}
     &&& (\text{Derivative-informed operator learning})\\
\begin{split}\label{eq:h1_training_equations}
    &\min_{\boldsymbol{w}\in\R^{d_{\boldsymbol{w}}}} \mathbb{E}_{(\refCoord,m_{\perp})\sim\pi\otimes\mu_{\perp}}\Big[\left\|\boldsymbol{V}^*\Forward(\Decoder\refCoord + m_{\perp}) - \RidgeForward_{\boldsymbol{w}}(\refCoord)\right\|^2 \\
    &\qquad+ \big\| \underbrace{\boldsymbol{V}^*\circ D \Forward(\Decoder \refCoord + m_{\perp})\circ \Decoder}_{\mathclap{\approx \nabla\RidgeForwardOpt(\refCoord)} } - \grad_{\refCoord}\RidgeForward_{\boldsymbol{w}}(\refCoord)\big\|^2_{F} \Big]
\end{split}
     &&& (\text{Equivalent latent representation})
\end{align}
The equivalent latent representation reveals that the derivative-informed learning objective can be derived from \cref{eq:h1_learning_latent} using a single-sample ($m_{\perp} \sim\mu_{\perp}$) estimate for the marginalization in both $\RidgeForwardOpt$ and $\nabla_{\refCoord}\RidgeForwardOpt$; see~\cref{gradient_conditional_expectation_commute} for a discussion on the marginalization in $\nabla_{\refCoord}\RidgeForwardOpt$. We refer to \texttt{DIPNet} surrogates trained using the derivative-informed $H^1_{\mu}$ operator learning method as \texttt{RB-DINO} (reduced basis derivative-informed neural operator).

We emphasize that one only needs samples of the latent representation of the derivative for $H^1_{\mu}$ operator learning compared to $L^2_{\mu}$ operator learning, 
\begin{align} \label{eq:reduced_jacobian}
    \boldsymbol{J}_r^{(j)} \coloneqq \boldsymbol{V}^*\circ D\Forward(m^{(j)}) \circ \Decoder \in \mathbb{R}^{\ObsDim \times \RedCoordDim}\quad m^{(j)}\sim\mu. &&& (\text{whitened latent Jacobian})
\end{align}
For notational convenience, we use $\boldsymbol{J}_r^{(j)}$ to denote the i.i.d. samples of the \emph{whitened latent Jacobian sample} of the PtO map.

\subsection{Surrogate-driven lazy map variational inference in the latent space}\label{surrogate_lazy_inference}
We use a trained ridge function surrogate $\FullRidgeForward_{\boldsymbol{w}}\circ\Projector$ to replace the PtO map $\Forward$ and its Jacobian $D\Forward$ evaluations during stochastic optimization of the latent space transport $\Map_{\MapParams}$. Specifically, a single-sample estimate of the surrogate rKL $\widetilde{{\mathcal{L}}}^{\Obs}_1:\ParamSpace\times\R^{d_{\MapParams}}\to\R$, replacing ~\cref{eq:mtvi_loss}, and its gradient, replacing ~\cref{eq:mtvi_gradient}, can be equivalently represented in the latent space as $\widetilde{\mathcal{L}}^{\Obs}_{1,r}:\R^{\RedCoordDim}\times\R^{d_{\MapParams}}\to\R$:
\begin{subequations}
\begin{align} \label{eq:surrogate_transport_loss}
     \widetilde{\mathcal{L}}^{\Obs}_1(m, \MapParams;\boldsymbol{w})&\equiv\widetilde{\mathcal{L}}^{\Obs}_{1,r}(\Encoder m, \MapParams;\boldsymbol{w})\\
     \widetilde{\mathcal{L}}^{\Obs}_{1,r}(\refCoord, \MapParams;\boldsymbol{w})&= \frac{1}{2}\norm{(\RidgeForward_{\boldsymbol{w}}\circ\Map_{\MapParams})(\refCoord) - \boldsymbol{V}^*\Obs}^2 + \frac{1}{2}\norm{ \Map_{\MapParams}(\refCoord)}^2 -\log\det\grad_{\refCoord} \Map_{\MapParams}(\refCoord^{(j)}),\\
    \nabla_{\MapParams}\widetilde{\mathcal{L}}^{\Obs}_{1,r}(\refCoord, \MapParams;\boldsymbol{w}) &= \nabla_{\MapParams} \Map_{\MapParams}(\refCoord)^{\transpose}\left(\nabla_{\refCoord}\RidgeForward_{\boldsymbol{w}}\circ\Map_{\MapParams})(\refCoord)^{\transpose}\left((\RidgeForward_{\boldsymbol{w}}\circ\Map_{\MapParams})(\refCoord) - \boldsymbol{V}^*\Obs\right)\right) \nonumber \\
    &\qquad\qquad\qquad\qquad+ \nabla_{\MapParams}\Map_{\MapParams}(\refCoord^{(j)})^{\transpose}\Map_{\MapParams}(\refCoord^{(j)}) -
    \nabla_{\MapParams}(\log\det\nabla_{\refCoord} \Map_{\theta})(\refCoord^{(j)})
\end{align}
\end{subequations}
Consequently, the surrogate-driven training of lazy maps proceeds entirely in the parameter latent space. After the latent space transport map $\Map_{\MapParams}$ is optimized, we use the map to produce latent space posterior samples $\boldsymbol{x}^{(j)}\sim \Map_{\MapParams}(\refCoord^{(j)})$, and they can be lifted to the full parameter space via sampling the prior as in \cref{eq:sampling_fullpost}.

%% file: lazydino.tex
\section{The \texttt{LazyDINO} method} \label{section:lazy_dino_amortize}

In this section, we present a high-level overview of the steps involved in \texttt{LazyDINO} using schematics and brief descriptions. We refer the reader to \ref{appendix:detailed_lazydino_algorithm} for a more detailed exposition of these steps. 

\subsection{Offline phase: \texttt{RB-DINO} surrogate construction}

\begin{figure}[h]
\center
\includegraphics[width = \textwidth]{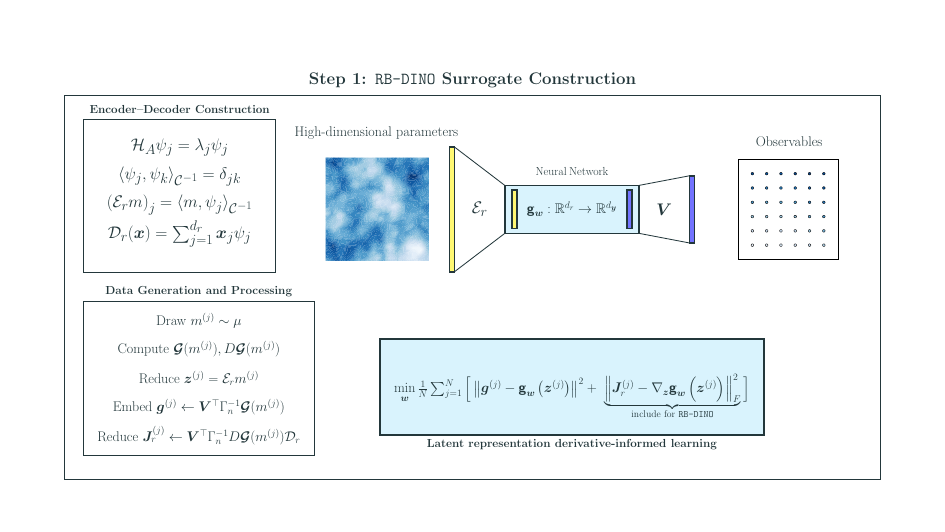}
\caption{Overview of the \texttt{RB-DINO} construction.}\label{fig:dino_construction}
\end{figure}

In Figure \ref{fig:dino_construction}, we provide a schematic for the offline surrogate construction. In this phase, one first defines the prior $\mu$ and PtO map $\Forward$ that determines the class of BIPs to be solved by \texttt{LazyDINO}. Subsequently, encoders and decoders for the parameter are constructed as delineated in \cref{sec:AS_dim_reduction}. The training samples are then generated and reduced to their latent representations as in \cref{eq:latent_parameter,eq:reduced_jacobian}. In particular, the full PtO map Jacobian samples are never formed. Instead, they are compressed matrix-free using the parameter decoder $\Decoder$. The training sample generation is often computationally costly, as the PtO map evaluations often require model solutions, and the PtO map Jacobian actions require computing the forward or adjoint model sensitivity. Once the training samples are collected, a given neural latent representation $\RidgeForward_{\boldsymbol{w}}$ is trained using the derivative-informed learning method in the latent space \cref{eq:h1_learning_latent}. We refer to \cite{o2024derivative, cao2024efficient} for more implementation details and theory on \texttt{RB-DINO} surrogate construction.

\subsection{Online phase: Rapid \texttt{LazyDINO} transport map construction}

\begin{figure}[H]
\includegraphics[width = \textwidth]{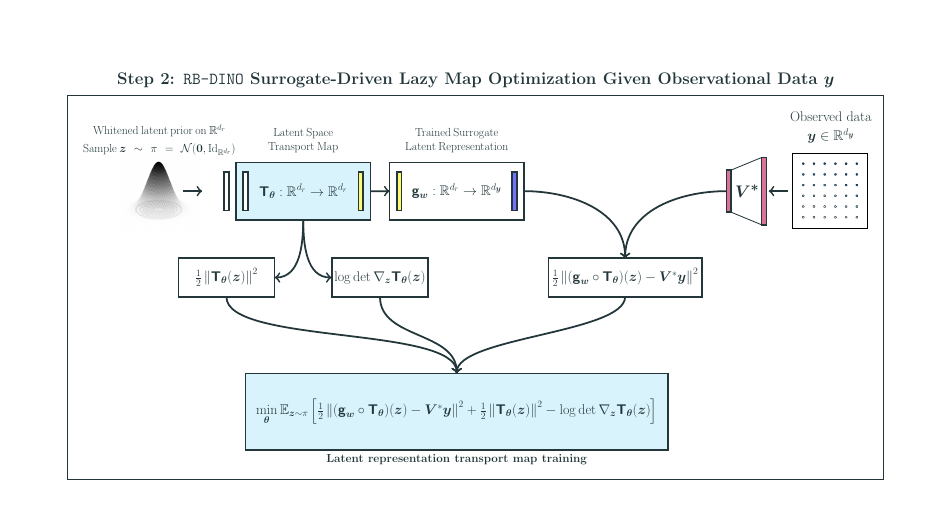}
\caption{Overview of the latent representation lazy map construction.}\label{fig:lazy_dino_construction}
\end{figure}

Once a \texttt{RB-DINO} surrogate PtO map is constructed, it is used in place of the PtO map in the lazy map training, which removes the computational bottleneck of model solutions and makes rapid online inference possible. The process for lazy map construction for a single instance of observational data $\Obs$, is shown in Figure \ref{fig:lazy_dino_construction}. This process involves defining the transport map architecture $\Map_\MapParams$, and the associated latent space rKL objective \cref{eq:surrogate_transport_loss}. The transport map is then optimized with respect to its map parameters $\MapParams$. The trained latent space transport map pushes the whitened latent prior $\pi$ in \cref{eq:whitened_subsapce_prior} to an approximation of the latent posterior induced by the optimal ridge function ${\Map_{\MapParams}}_{\sharp}\pi \approx \widetilde{\pi}^{\Obs}_{\text{opt}}$ as in \cref{theorem:lazy_map_kl}. These samples can then be decoded to generate samples in the full space using the prior $\mu$.

A major point of emphasis for \texttt{LazyDINO} is that the computationally expensive aspects of the method are limited to the offline phase. The \texttt{RB-DINO} surrogate replaces the expensive-to-evaluate and often implicitly defined PtO map with a fast-to-evaluate explicit function. In practice, this leads to potentially enormous speedups for all computations associated with the likelihood and its gradient evaluations. 
 
Likewise, the transport map approximation (Figure \ref{fig:lazy_dino_construction}) occurs in a relatively low-dimensional parameter latent space. It extensively uses highly optimized modern computing kernels such as batch computations, fast sampling of white noise, automatic differentiation, and compile-optimized explicit calculations where all operations are known a priori. 

\subsection{\texttt{LazyDINO} as an amortized inference method}
The latent space transport maps can be rapidly constructed during the online phase, making \texttt{LazyDINO} a compelling method for real-time inference and a competitive alternative to simulation-based amortized inference (SBAI)  methods~\cite{Ganguly_2023} for solving BIPs with the same PtO map but difference instances of observational data. 

\begin{figure}[H]
\center
\includegraphics[width = \textwidth]{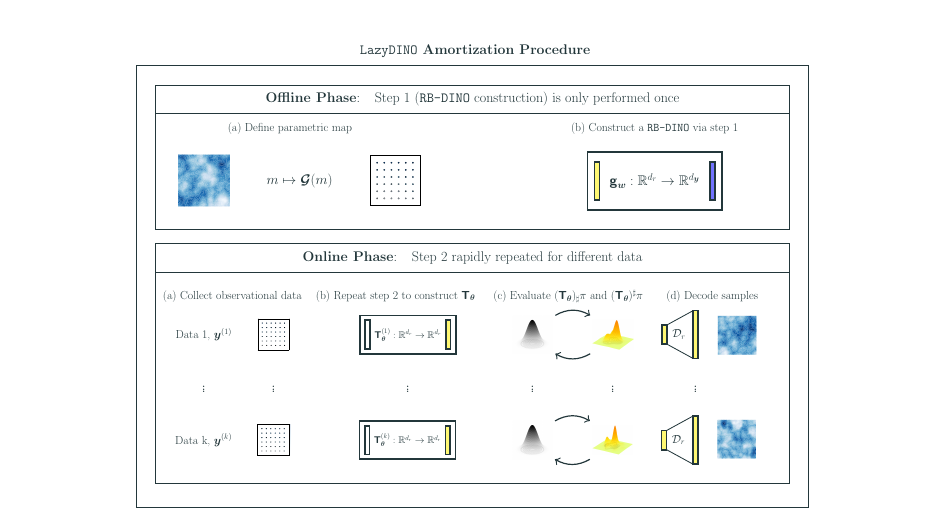}
\caption{Overview of \texttt{LazyDINO} amortization procedure.}\label{fig:lazy_dino_amortization}
\end{figure}

Since we create a surrogate for the PtO map, \emph{and not} the likelihood, the cost of \texttt{RB-DINO} surrogate construction can be amortized by instancing a new likelihood for any new observational data. In this amortization process (see Figure \ref{fig:lazy_dino_amortization}), the construction of the \texttt{RB-DINO} surrogate is amortized over many different BIPs defined by the same PtO map and prior. With this in mind, \texttt{LazyDINO} is an ideal method for settings where many BIPs are solved for the same system. This class of problems can be found in predictive digital twins, state estimation, and Bayesian optimal experimental design.

\paragraph{Comparison of \texttt{LazyDINO} to other SBAI methods} Given joint samples of the latent prior and simulated observational data, the SBVI methods optimize for a \emph{conditional} transport map that matches the pullback distributions, $\Obs^{(j)} \mapsto \Map_{\MapParams}(\Obs^{(j)}, \cdot)^{\sharp}\RedCoordPrior$, to posteriors at the simulated data samples using an fKL objective. The sample generation and the transport map construction are both performed offline. When the observational data $\Obs^\dagger$ is available, the approximate posterior sampling is performed using inversion at latent prior samples $\Map_{\MapParams}(\Obs^\dagger, \cdot)^{-1}(\refCoord^{(j)})$, $\refCoord^{(j)}\iid\RedCoordPrior$, without the computational bottleneck of model simulations, making it an amortized inference method. Simulating observational data for transport map training in SBAI incurs a similar cost compared to \texttt{RB-NO} surrogate construction, requiring PtO map samples with noise perturbation, i.e., $\Obs^{(j)}:= \Forward(\FullParam^{(j)}) + \bn^{(j)}$ with $m^{(j)}\iid\mu$ and $\boldsymbol{n}\iid\mathcal{N}(0, \NoiseCov)$. The training also only requires easy-to-evaluate quantities related to the latent prior; see~\cite{Ganguly_2023}.

Despite ostensibly relaxed training requirements, SBAI is much lower in sample efficiency than \texttt{LazyDINO}. Given limited PtO map samples, SBAI attempts to directly approximate all posteriors (i.e., posteriors for all instances of observational data), whereas our approach invests these samples in surrogate construction to gain almost unlimited access to all surrogate-approximated posteriors. As a result, our approach leads to a much smaller transport map approximation error than SBAI. Furthermore, due to our efficient \texttt{RB-DINO} surrogate construction, the transport map approximation error of SBAI is much higher than the surrogate posterior approximation error in \texttt{LazyDINO}, leading to more than two orders of magnitude lower sample efficiency of SBAI observed in our numerical results in \cref{sec:numerical_results}.

All optimization problems of SBAI are solved offline, which is often regarded as an advantage. In contrast, \texttt{LazyDINO} requires solving an optimization problem for online posterior sampling at each instance of observational data. However, for popular transport map parametrizations such as conditional normalizing flows, sampling requires solving a root-finding problem in $\R^{d_r}$ for map inversion, which incurs a non-negligible cost in practice. The inversion-to-sample approach of SBAI can be more costly than the optimize-to-sample approach of \texttt{LazyDINO} in some situations, e.g., when a large number of approximate posterior samples is needed for each instance of observational data. We provide concrete numerical evidence to support these claims in~\cref{SBAIvLazyDINOTiming}.

\begin{table}[!htbp]
    \centering
    \begin{tabular}{||l|c||c|c|c|c|c||}\hline

          \diagbox{\makecell{\begin{tabular}[t]{@{}l@{}} \textbf{Desired} \\\textbf{Algorithm Characteristics} \end{tabular}}}{\makecell{\begin{tabular}[t]{@{}r@{}} \textbf{Posterior}\\ \textbf{estimate}\\ \end{tabular}}}
         & \rotatebox[origin=c]{90}{\makecell{Ground truth\\ (MCMC)}} & \rotatebox[origin=c]{90}{LA-baseline} & \rotatebox[origin=c]{90}{\texttt{LazyMap}} & \rotatebox[origin=c]{90}{SBAI} & \rotatebox[origin=c]{90}{\texttt{LazyNO}} & \rotatebox[origin=c]{90}{\texttt{LazyDINO}} \\ \hline \hline
        Parallel sampling from posterior estimate & {\color{BrickRed} \xmark} & {\color{ForestGreen} \cmark} & {\color{ForestGreen} \cmark} & {\color{ForestGreen} \cmark} & {\color{ForestGreen} \cmark} & {\color{ForestGreen} \cmark} \\ \hline
 Direct sampling, no inversion $\Map_{\MapParams}^{-1}$ required &-- & -- & {\color{ForestGreen} \cmark} & {\color{BrickRed} \xmark} & {\color{ForestGreen} \cmark} & {\color{ForestGreen} \cmark} \\ \hline
        Uses a neural surrogate PtO map $\RidgeForward_w$ & {\color{BrickRed} \xmark} & {\color{BrickRed} \xmark} & {\color{BrickRed} \xmark} & {\color{BrickRed} \xmark} & {\color{ForestGreen} \cmark} & {\color{ForestGreen} \cmark} \\ \hline
        Parallelly-sampled training data outputs & -- & -- & -- & $\Obs^{(j)}$ & $\RidgeForwardOptApprox^{(j)}$ & \makecell{$\RidgeForwardOptApprox^{(j)}$, \\$\bJ_r^{(j)}$} \\ \hline

        Amortizes inversion of posteriors $\FullPost$ & {\color{BrickRed} \xmark} & {\color{BrickRed} \xmark} & {\color{BrickRed} \xmark} & {\color{ForestGreen} \cmark} & {\color{BrickRed} \xmark} & {\color{BrickRed} \xmark} \\ \hline
        
  Uses $D\Forward$ evaluations & {\color{ForestGreen} \cmark} & {\color{ForestGreen} \cmark} & {\color{ForestGreen} \cmark} & {\color{BrickRed} \xmark} & {\color{BrickRed} \xmark} & {\color{ForestGreen} \cmark} \\ \hline

        Amortizes evaluation of $\Forward, D\Forward$ & {\color{BrickRed} \xmark} & {\color{BrickRed} \xmark} & {\color{BrickRed} \xmark} & {\color{ForestGreen} \cmark} & {\color{ForestGreen} \cmark} & {\color{ForestGreen} \cmark} \\ \hline
        Embarrassingly parallel $\Forward$ evaluation & {\color{BrickRed} \xmark} & {\color{BrickRed} \xmark} & {\color{orange} \cmark} & {\color{ForestGreen} \cmark} & {\color{ForestGreen} \cmark} & {\color{ForestGreen} \cmark} \\ \hline

    \end{tabular}
    \caption{\textbf{Posterior estimate comparison.} ({\color{ForestGreen} \cmark}/{\color{BrickRed} \xmark}) refer to (true/false).  The first two posterior estimates (left two columns) are the ground truth and Laplace approximation baseline (LA-Baseline).
    The MCMC method we use exhibits none of our desired algorithm characteristics, but its samples serve as a trustworthy ground truth. All competing posterior estimation methods (right four columns) use the same parameter dimension reduction and offer parallel i.i.d.\ approximate posterior sampling once trained. We use an orange checkmark {\color{orange} \cmark} to highlight that \texttt{LazyMap} can only partially exploit embarrassingly parallel evaluation of the PtO map and its Jacobians, in particular only within each stochastic optimization iteration. We note that SBVI fully amortizes Bayesian inversion, while \texttt{LazyNO} and \texttt{LazyDINO} only offer amortization of $\Forward, D\Forward$ across the estimation of all BIPs.
    }
    \label{tab:transposed_baselines}
\end{table}

%% file: results.tex
\section{Setup of the numerical studies}\label{section:numerical_study_setup}

In this section, we describe the setup of the numerical examples, including the description of the two BIP examples, neural network architectures, training procedures, and posterior error measures. First, we define the reduced basis dimension, the architecture and training of the \texttt{RB-DINO} surrogate, and the architecture and training of the transport map in~\cref{red_basis_dim}, and \ref{transport_map_specs}. Then, we introduce the error measures for operator learning and posterior approximations in~\cref{quantifying_error_specs}. We proceed by first defining the two PDE problems and associated inverse problems. In both cases, we consider nonlinear elliptic PDEs defined on 2D rectangular domains $\Omega \subset \mathbb{R}^2$. 

For both problems, we define our Gaussian priors using Mat\'ern covariance operators given by an elliptic operator on $\ParamSpace\coloneqq L^2(\Omega)$:
\begin{equation} \label{eq:covariance_definition}
    \mathcal{C} = (-\gamma \nabla \cdot(\bA \nabla ) + \delta \identity_{\ParamSpace})^{-2}.
\end{equation}
Here $\gamma,\delta>0$ are scalar parameters that control the marginal variance and correlation lengths of the random field samples, and $\bA \in \mathbb{R}^{2\times2}$ is a symmetric positive definite matrix that induces anisotropy in the random field samples. 
In both cases, we employ Robin boundary conditions to control boundary artifacts in the samples; see \cite[Equation 37]{VillaPetraGhattas21} and \cite{villa2024note}.

In both cases, $\ParamSpace$ is approximated using linear triangular finite elements, while the state spaces $\StateSpace \subset H^1(\Omega;\mathbb{R}), H^1(\Omega;\mathbb{R}^2)$, respectively, are both approximated using quadratic triangular finite elements. In both cases, the reference solution maps utilize Newton--Raphson methods using sparse direct solvers for each Newton iteration. Notably, sparse direct solvers lead to efficient computation of Jacobian training data; see \cite{o2024derivative,cao2024efficient} for more information on Jacobian computation.

The PtO map  $\Forward$ is defined by composing the PDE solution operator $\SolOp:\ParamSpace\to\StateSpace$ and an observation operator $\ObsOp:\StateSpace\to\ObsSpace$. The inverse problems arise from generating four synthetic observations, each obtained by sampling the prior, evaluating the PtO map, and applying additive white noise.

\subsection{Example I: Inference of the diffusivity field in a nonlinear reaction--diffusion PDE}\label{example1}

For our first example, we consider the following nonlinear reaction--diffusion PDE for $u:\Omega = [0,1]^2\to\R$:
\begin{subequations}
\begin{align}
    - \nabla_{\bs} \cdot \exp(m(\bs)) \nabla u(\bs) + u(\bs)^3 &= 0, \quad \bs \in (0,1)^2, \\
    \exp(m(\bs)) \nabla u(\bs) \cdot \bn &= 0, \quad \bs \in \Gamma_{\mathrm{left}} \cup \Gamma_{\mathrm{right}}, \\
    u(\bs) &= 1, \quad \bs \in \Gamma_{\mathrm{top}}, \\
    u(\bs) &= 0, \quad \bs \in \Gamma_{\mathrm{bottom}},
\end{align}
\end{subequations}
where $\Gamma_{\mathrm{left}}$, $\Gamma_{\mathrm{right}}$, $\Gamma_{\mathrm{top}}$, and $\Gamma_{\mathrm{bottom}}$ denote the left, right, top and bottom boundaries of the unit square, and $\bn$ is the outward unit normal vector. The inverse problem is to find the log-diffusivity field $m:(0,1)^2\to\R$ that best matches noisy observations of $u$ at a set of spatial positions. 

\begin{figure}[H]
    \centering
    \renewcommand{\arraystretch}{0.05}
        \addtolength{\tabcolsep}{-5pt}
    \begin{tabular}{D E E E}
       &\makecell{Synthetic parameter\\$\FullParam\sim\FullPrior$} & \makecell{Noisy observation\\$\Obs\sim \mathcal{N}(\Forward(m), \NoiseCov)$} & \makecell{MAP estimate\\ $m_{\text{MAP}}^{\Obs}$}  \\
       \makecell{BIP\\ \#1}&\includegraphics[width=0.16\textwidth]{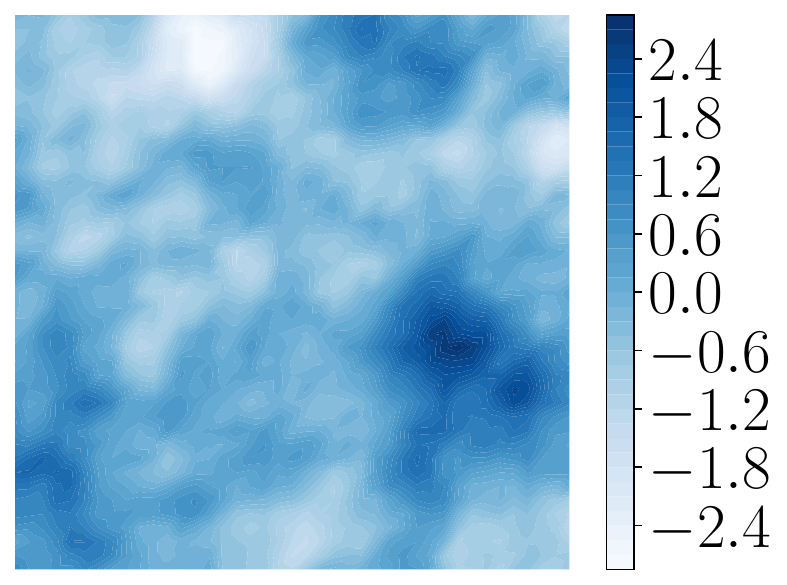}  &  \includegraphics[width=0.16\textwidth]{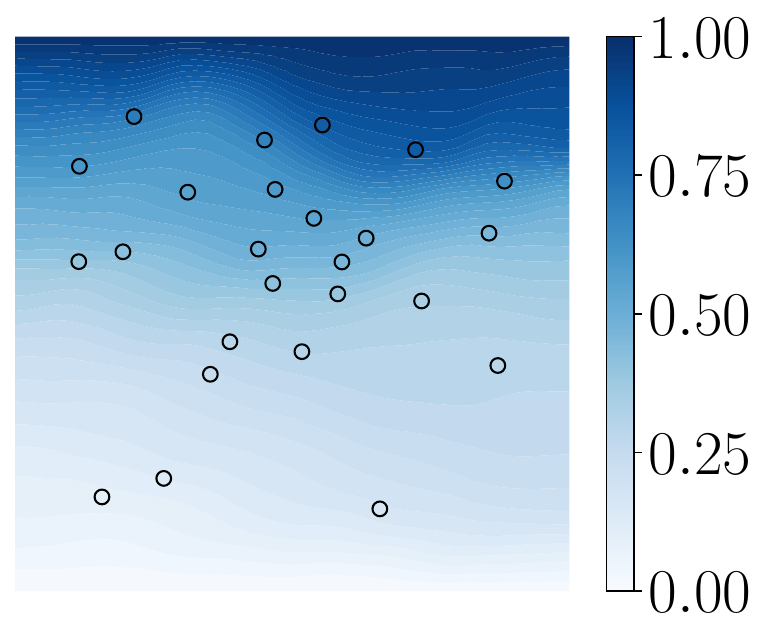} & \includegraphics[width=0.16\textwidth]{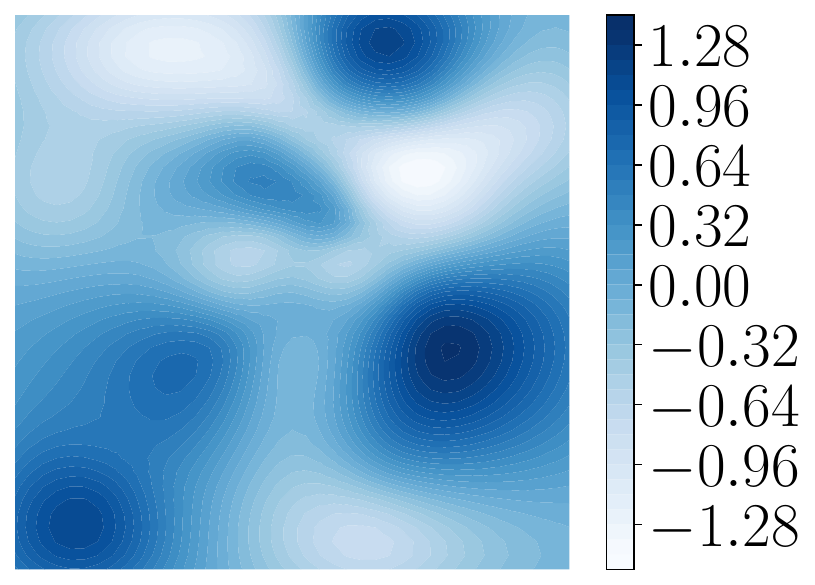} \\
       \makecell{BIP\\ \#2}&\includegraphics[width=0.16\textwidth]{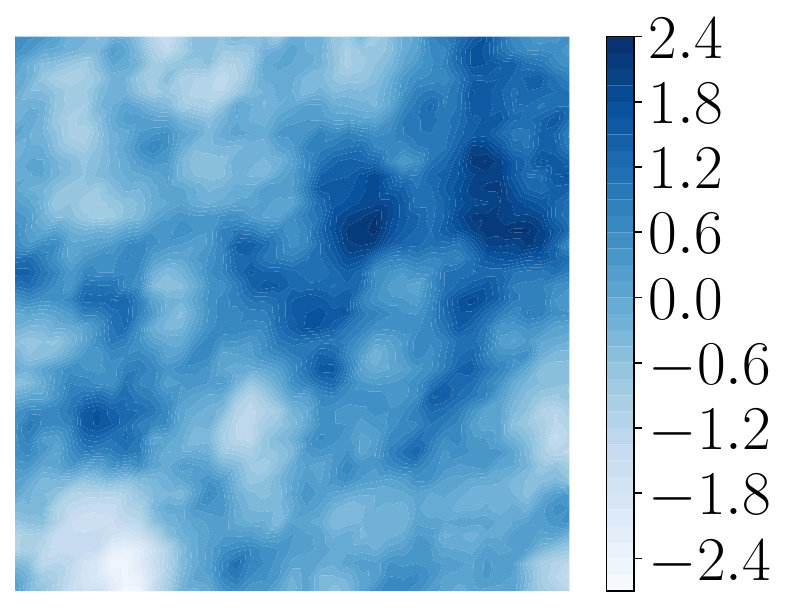}  & \includegraphics[width=0.16\textwidth]{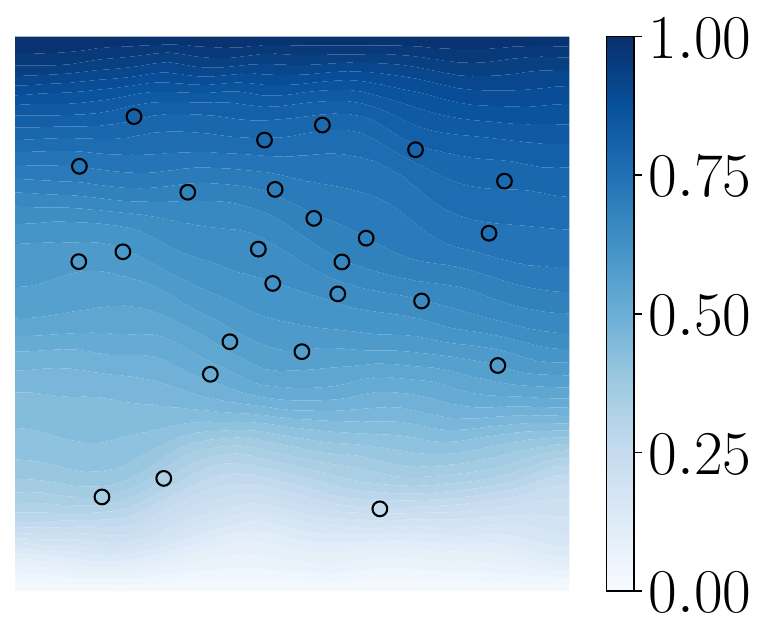} & \includegraphics[width=0.16\textwidth]{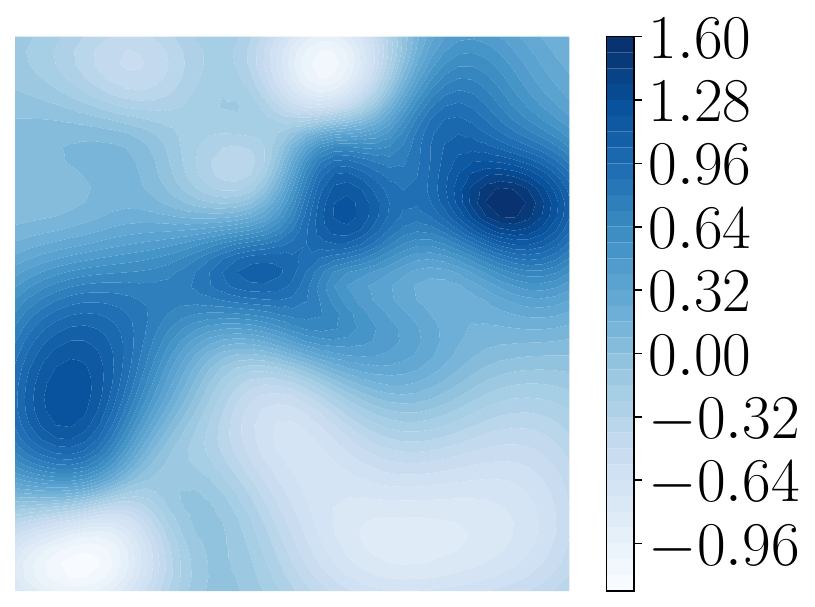} \\
       \makecell{BIP\\ \#3}&\includegraphics[width=0.16\textwidth]{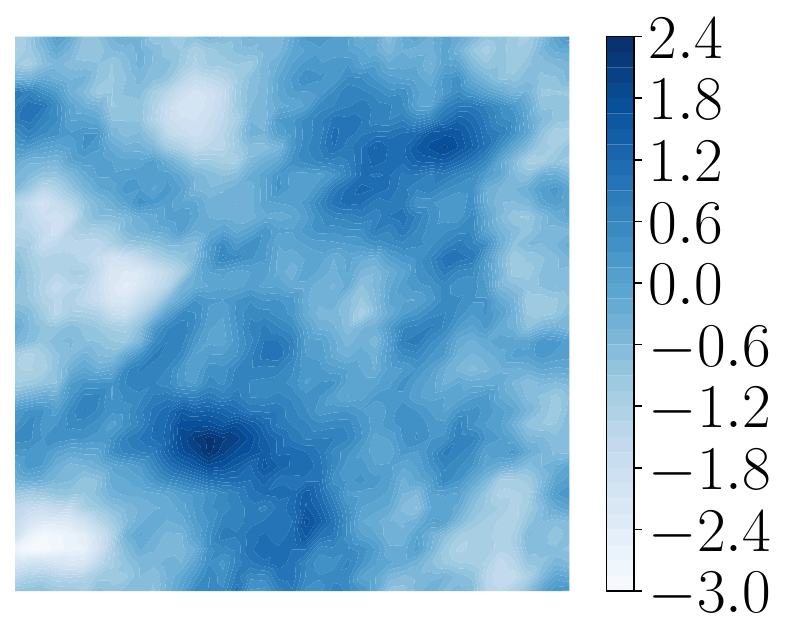}  & \includegraphics[width=0.16\textwidth]{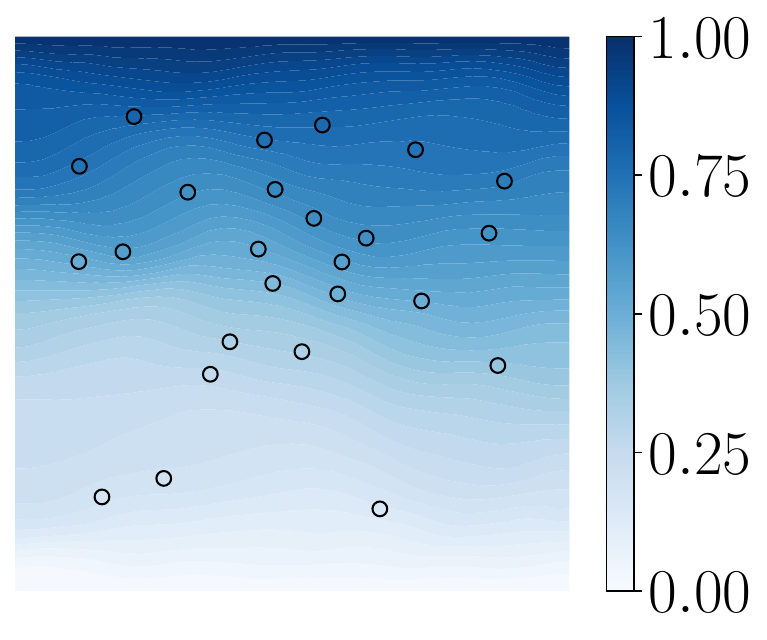} & \includegraphics[width=0.16\textwidth]{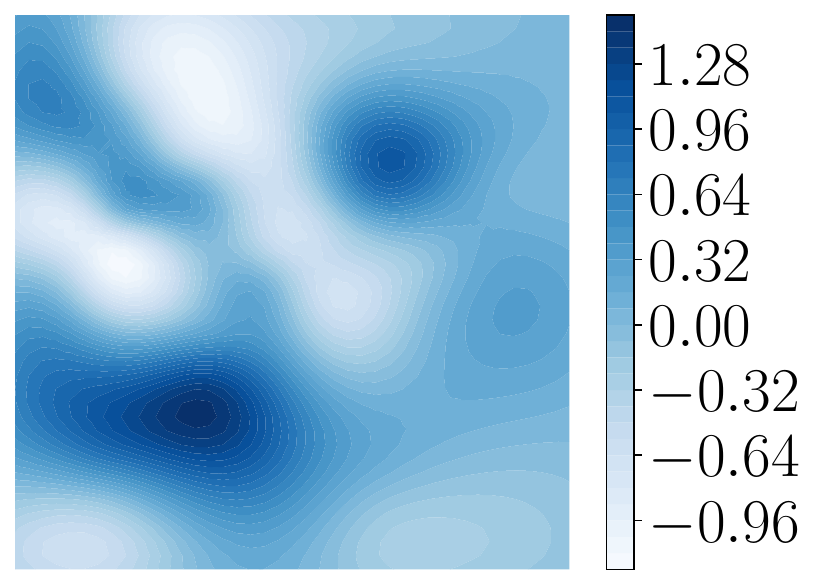} \\
       \makecell{BIP\\ \#4}&\includegraphics[width=0.16\textwidth]{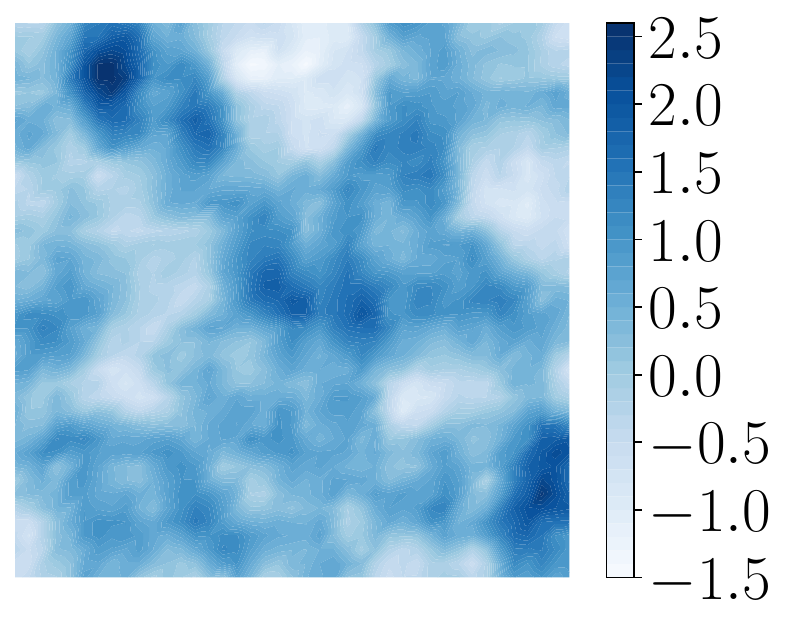}  &  \includegraphics[width=0.16\textwidth]{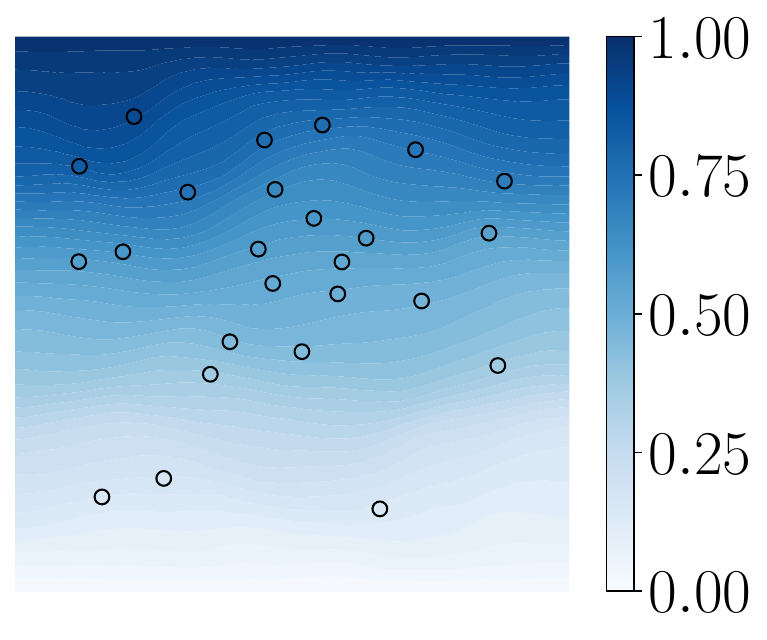} & \includegraphics[width=0.16\textwidth]{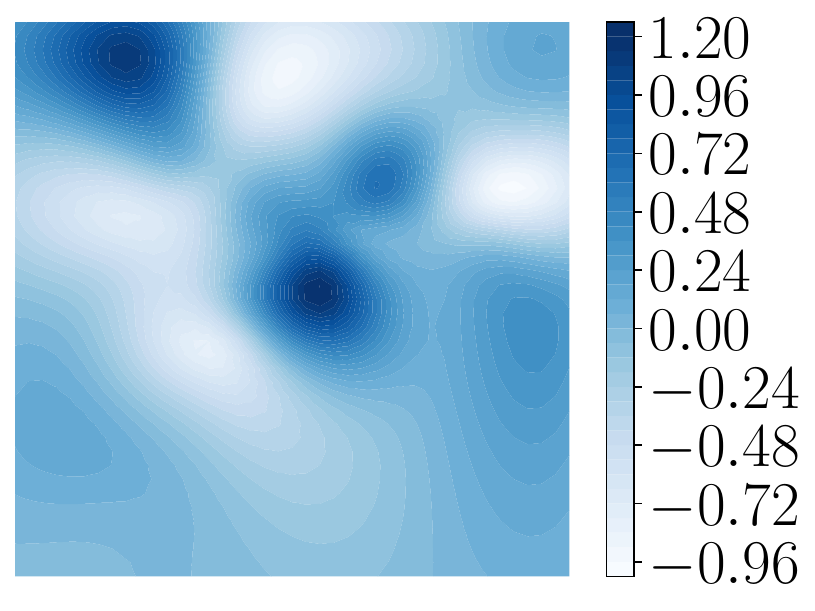}
    \end{tabular}
        \addtolength{\tabcolsep}{5pt}
        \renewcommand{\arraystretch}{1.0}
    \caption{\textbf{Example I}. Setup for inferring the diffusivity field in a nonlinear reaction--diffusion PDE detailed in \cref{example1}. For each BIP (\#1--4), we show the data-generating synthetic parameter $\FullParam$ drawn from the prior, the synthetic data $\Obs$ placed on top of the PDE solution at $\FullParam$, and the MAP estimate $m_{\text{MAP}}^{\Obs}$, i.e., the solution of the deterministic inverse problem.}
    \label{fig:nrd_bip_setup}
\end{figure}

We use a regular grid with $40\times40$ cells. The choice of linear triangular finite elements for $\ParamSpace$ leads to $1,681$ degrees of freedom (DoFs). The choice of quadratic triangular finite elements for $\StateSpace$ leads to $3,362$ DoFs. For the prior covariance~\eqref{eq:covariance_definition}, we choose $\gamma = 0.03$, $\delta = 3.33$, which leads to a point-wise marginal variance around $9$ and a spatial correlation length of around $0.1$. We take $\bA = \identity_{\R^2}$.
We define the observation operator using $\ObsDim=25$ randomly sampled interior points $\{\bs_{\text{obs}}^{(j)}\}_{j=1}^{\ObsDim}$, i.e., $ (\bdmc{O}(u))_{j} = \int_{\bB_\epsilon(\bs^{(j)}_{\text{obs}})} u(\bs)\dd\bs$, where $\bB_{\epsilon}(\bs)\subset(0,1)^2$ is a ball around $\bs$ with a small radius $\epsilon>0$. The noise distribution has covariance $\NoiseCov = 1.94\times 10^{-3} \identity_{\R^{\ObsDim}}$, which corresponds to a signal-to-noise ratio of around $500$.
We visualize the synthetic data set in \cref{fig:nrd_bip_setup}.

\subsection{Example II: Inference of a heterogeneous hyperelastic material property}\label{example2}
For our second example, we consider the uniaxial tensile test of a hyperelastic thin film. The inverse problem aims to recover Young's modulus field, which characterizes spatially varying material strength, from measurements of the material deformation. This problem is of interest to the characterization of heterogeneous material properties from deformation data, see for example \cite{kirchhoff2024inference}. Similar Bayesian inverse problems have been considered in \cite{cao2023residual,cao2024efficient}.

Let $\Omega = (0, 2)\times(0,1)$ be a normalized material domain. The material coordinates $\bs\in\Omega$ of the reference configuration are mapped to the spatial coordinates $\bs + \bu(\bs)$ of the deformed configuration, where $\bu:\Omega\to\R^2$ is the material displacement. The strain energy of the hyperelastic material $\calW_e$ depends on the deformation gradient, i.e.,  $\calW_e=\calW_e(\bF)$ where $\bF = \identity_{\R^{2\times 2}} + \nabla \bu$. We consider the neo-Hookean model for the strain energy density:
\begin{equation}
	\calW_e(\bF) = \frac{\mu_e}{2} (\mathrm{tr}(\bF^{\transpose}\bF) - 3) + \frac{\lambda_e}{2} \left(\text{ln\,det}(\bF)\right)^2 - \mu_e \text{ln\,det}(\bF).
\end{equation}
Here, $\lambda_e$ and $\mu_e$ are the Lam\'e parameters, and they are related to Young's modulus $E_{Y}$ and Poisson ratio $\nu_P$ under the plain strain assumption:
\begin{equation}
    \lambda_e = \frac{E_Y \nu_P}{(1+\nu_P)(1-2\nu_P)}\,, \qquad \mu_e = \frac{E_Y}{2(1+\nu_P)}.
\end{equation}
We assume $\nu_P = 0.4$, and a spatially-varying normalized Young's modulus, $E_Y:\Omega\to({{E}_{Y}}_{\text{min}}, {{E}_{Y}}_{\text{max}})$ with $0<{{E}_{Y}}_{\text{min}}<{{E}_{Y}}_{\text{max}}$. We represent $E_Y$ through a parameter field $m:\Omega\to\R$ as follows
\begin{equation*}
    E_Y(m(\bs)) = \frac{1}{2}\left({{E}_{Y}}_{\text{max}} - {{E}_{Y}}_{\text{min}}\right)\left(\text{erf}(m(\bs)) + 1\right) + {{E}_{Y}}_{\text{min}},
\end{equation*}
where $\text{erf}:\R\to(-1, 1)$ is the error function. We use ${{E}_{Y}}_{\text{min}} = 1$ and ${{E}_{Y}}_{\text{max}} = 7$. The first Piola--Kirchhoff stress tensor is given by $\bP_e(m, \bF) = 2\partial\calW_e(m, \bF)/\partial \bF$. Assuming a quasi-static model with negligible body forces, the balance of linear momentum leads to the following nonlinear PDE:
\begin{subequations}
\begin{align}\label{eq:linmombal}
	\nabla_{\bs} \cdot \bP_e(m, \bF)(\bs) &= \bzero, && \bs\in \Omega;\\
	\bu(\bs) &= \bzero\,, && \bs\in \Gamma_{\text{left}};\\
    \bu(\bs) & = 3/2\,, && \bs\in \Gamma_{\text{right}};\\
	\bP_e(m, \bF)(\bs) \cdot \bn  &= \bzero, && \bs\in \Gamma_{\text{top}}\cup\Gamma_{\text{bottom}};
\end{align}
\end{subequations}
where $\Gamma_{\text{top}}$, $\Gamma_{\text{right}}$, $\Gamma_{\text{bottom}}$, and $\Gamma_{\text{left}}$ denote the material domain's top, right, bottom, and left boundary.

\begin{figure}[H]
    \centering
    \addtolength{\tabcolsep}{-5pt}
    \renewcommand{\arraystretch}{0.2}
    \resizebox{1.0\textwidth}{!}{ 
    \begin{tabular}{D E F E E}
       &\makecell{Synthetic parameter\\$\FullParam\sim\FullPrior$}  & \makecell{Deformed configuration \\$\boldsymbol{s} + \boldsymbol{u}(\boldsymbol{s})$} & \makecell{Noisy observation\\$\Obs\sim \mathcal{N}(\Forward(m), \NoiseCov)$} & \makecell{MAP estimate\\ $m_{\text{MAP}}^{\Obs}$}  \\
       \makecell{BIP\\ \#1}&\includegraphics[width=0.18\textwidth]{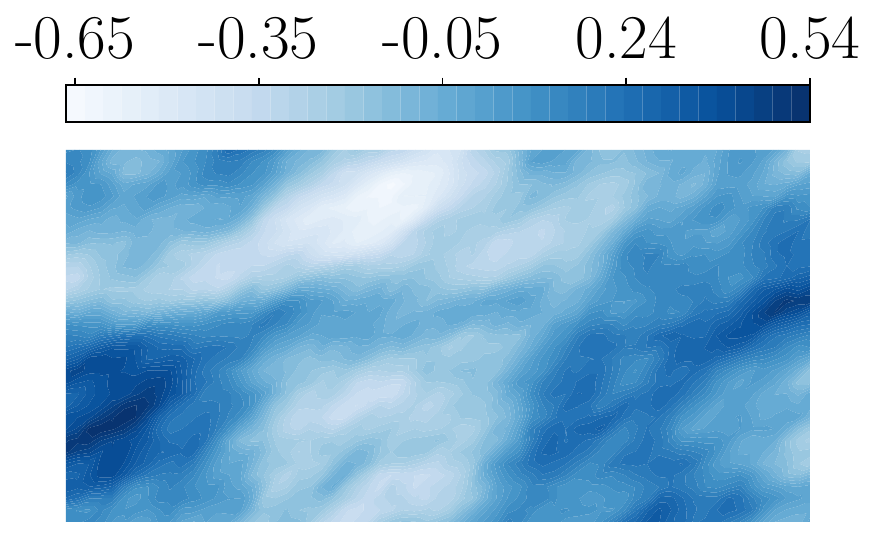}  & \includegraphics[width=0.21\textwidth]{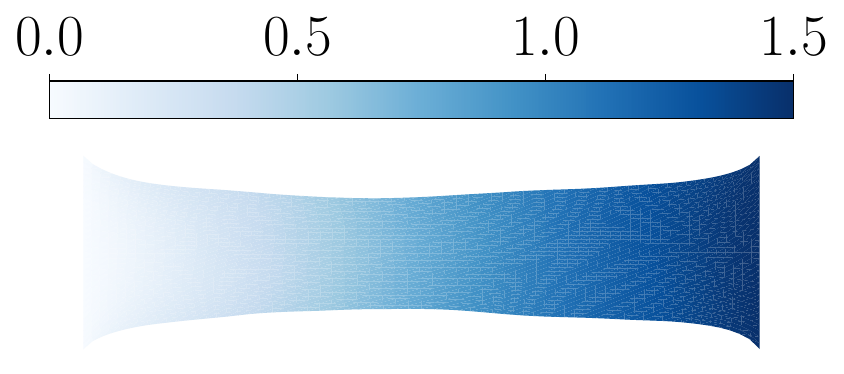} & \includegraphics[width=0.18\textwidth]{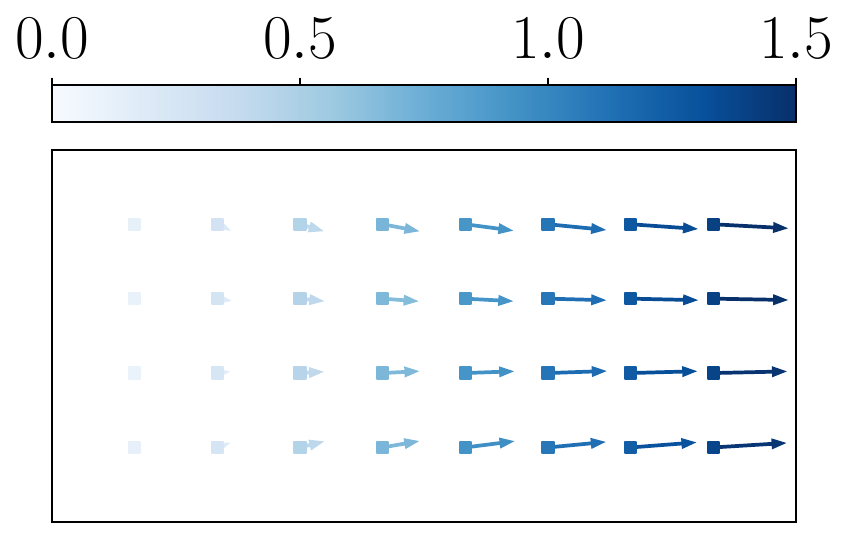} & \includegraphics[width=0.18\textwidth]{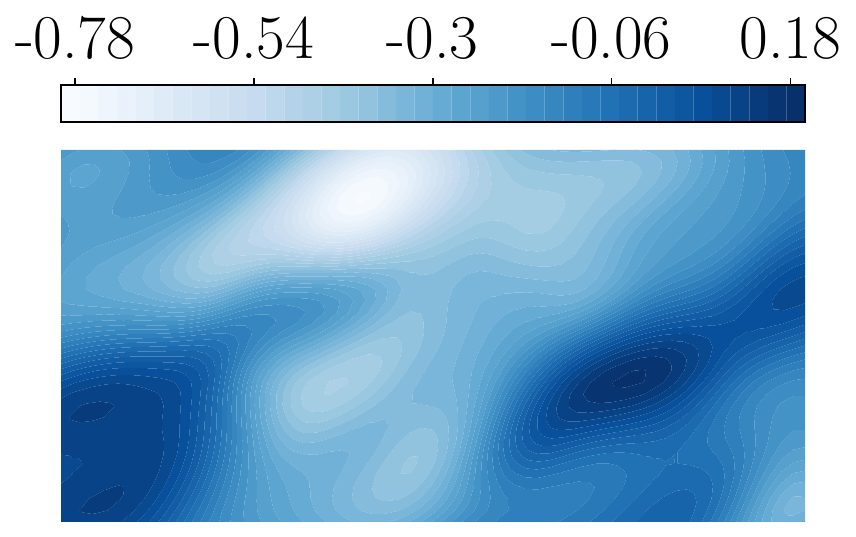} \\
       \makecell{BIP\\ \#2}&\includegraphics[width=0.18\textwidth]{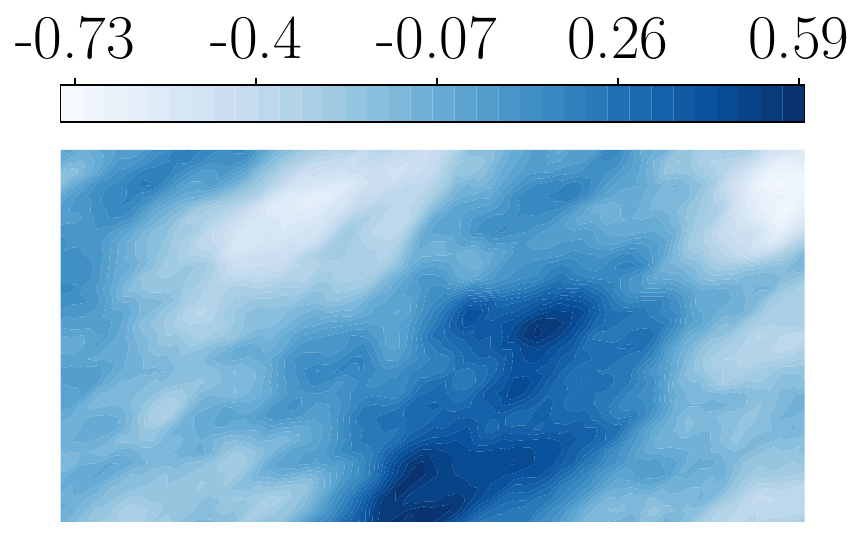}  & \includegraphics[width=0.21\textwidth]{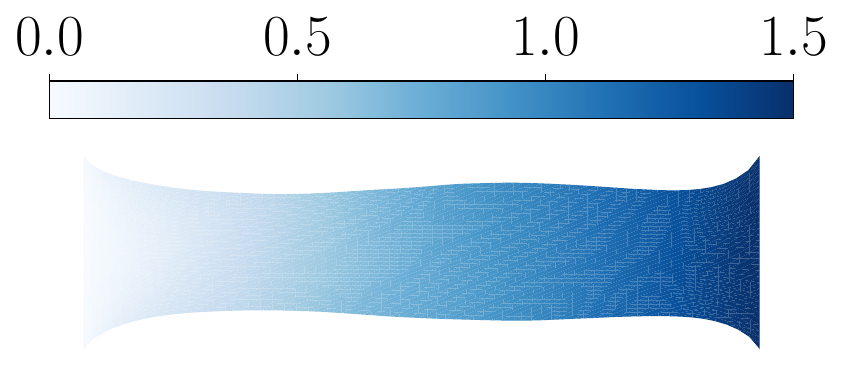} & \includegraphics[width=0.18\textwidth]{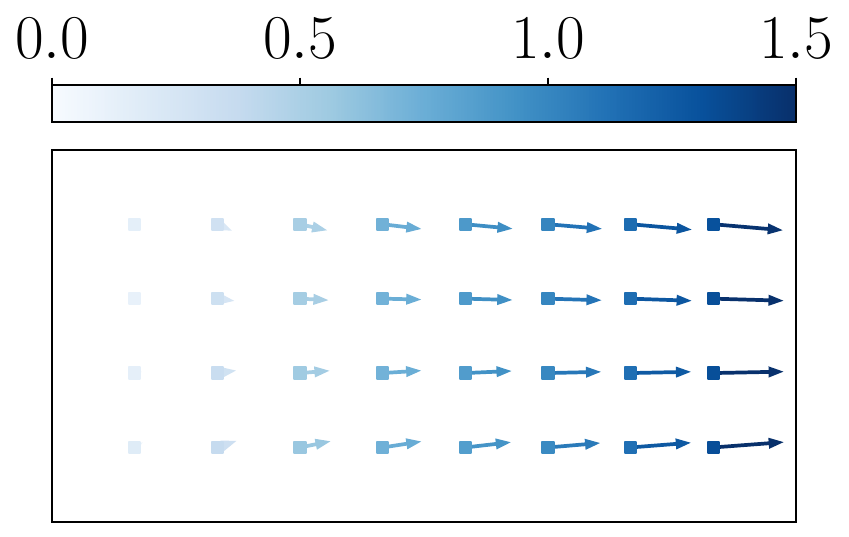} & \includegraphics[width=0.18\textwidth]{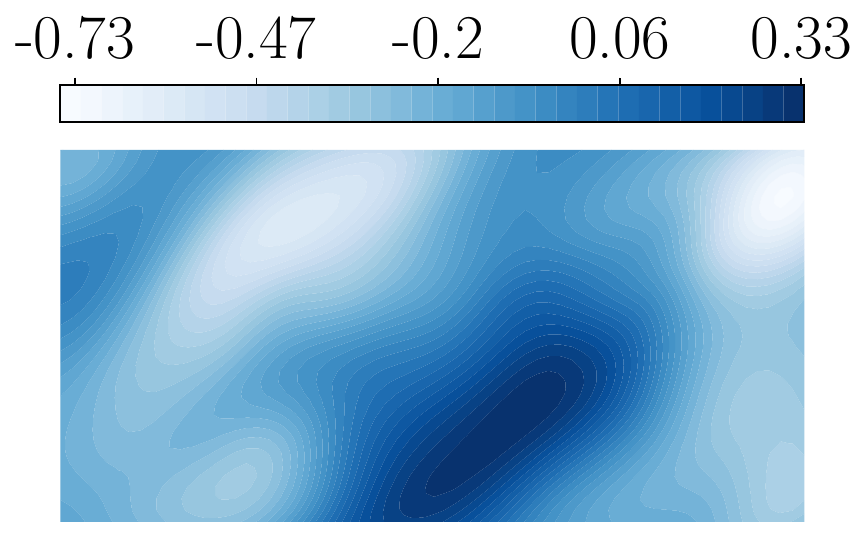} \\
       \makecell{BIP\\ \#3}&\includegraphics[width=0.18\textwidth]{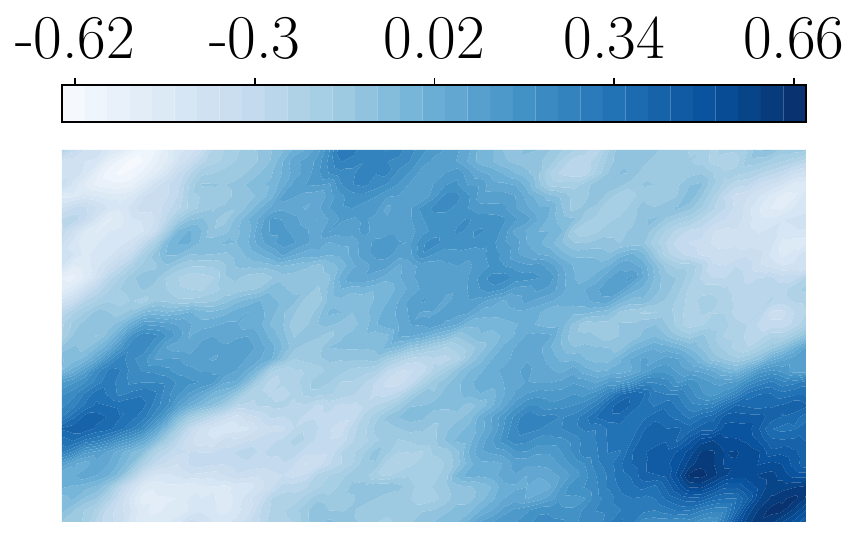}  & \includegraphics[width=0.21\textwidth]{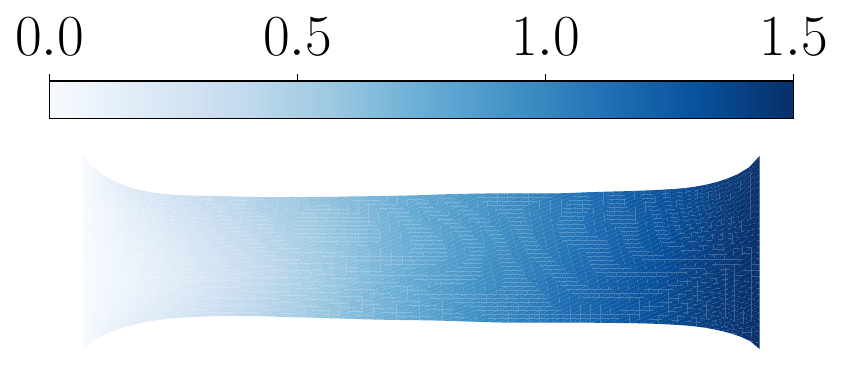} & \includegraphics[width=0.18\textwidth]{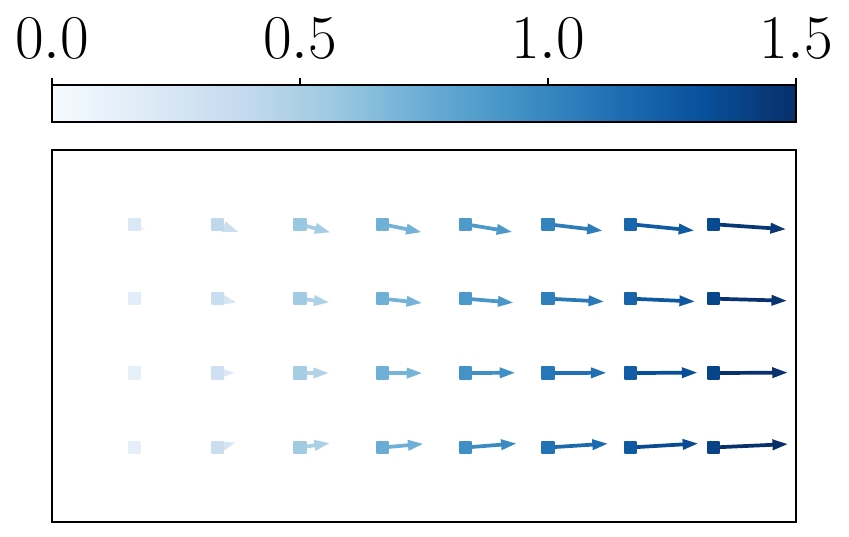} & \includegraphics[width=0.18\textwidth]{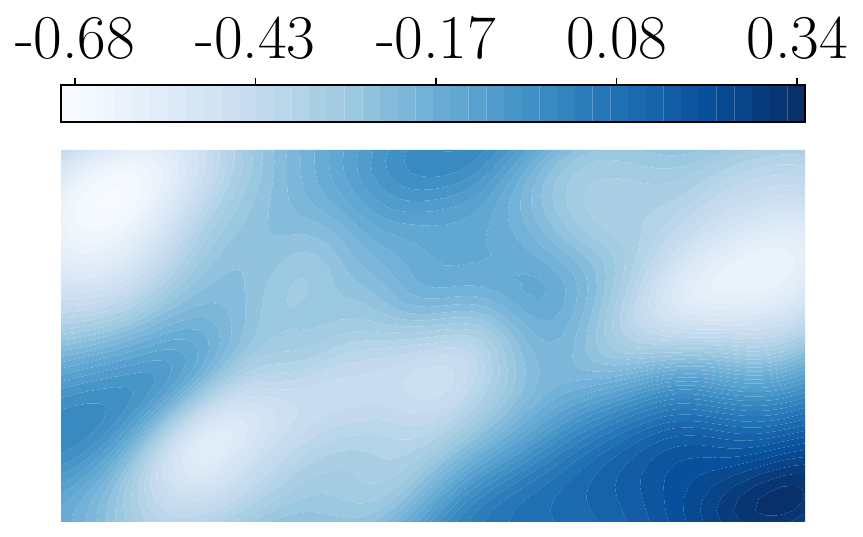} \\
       \makecell{BIP\\ \#4}&\includegraphics[width=0.18\textwidth]{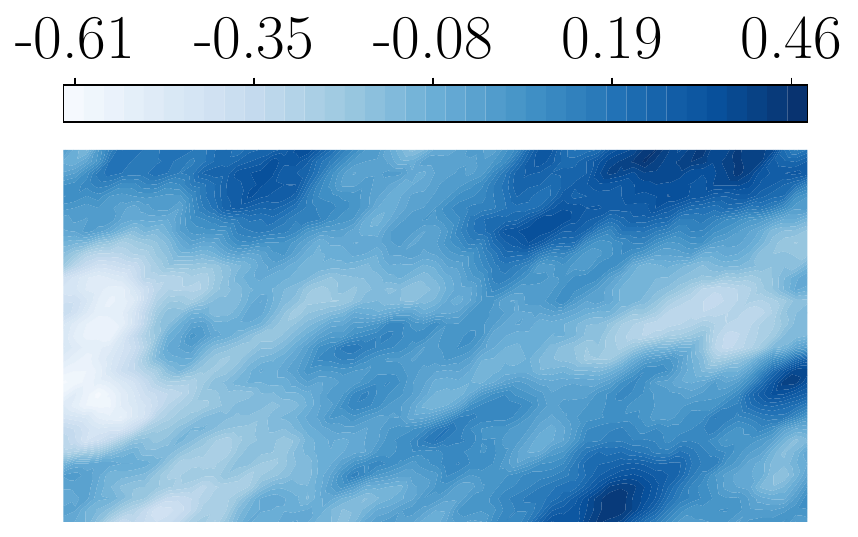}  & \includegraphics[width=0.21\textwidth]{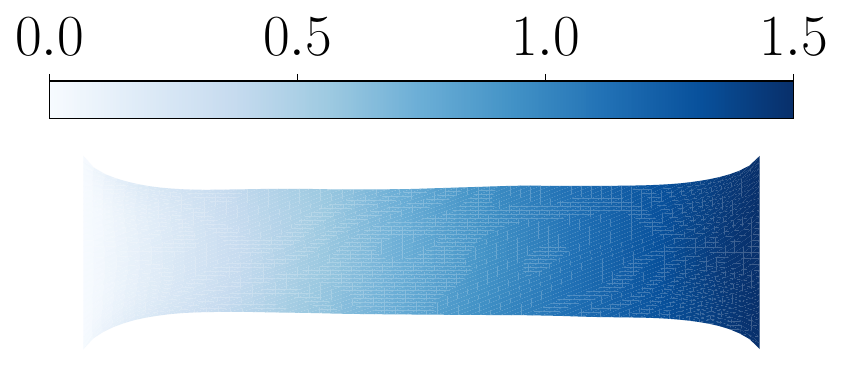} & \includegraphics[width=0.18\textwidth]{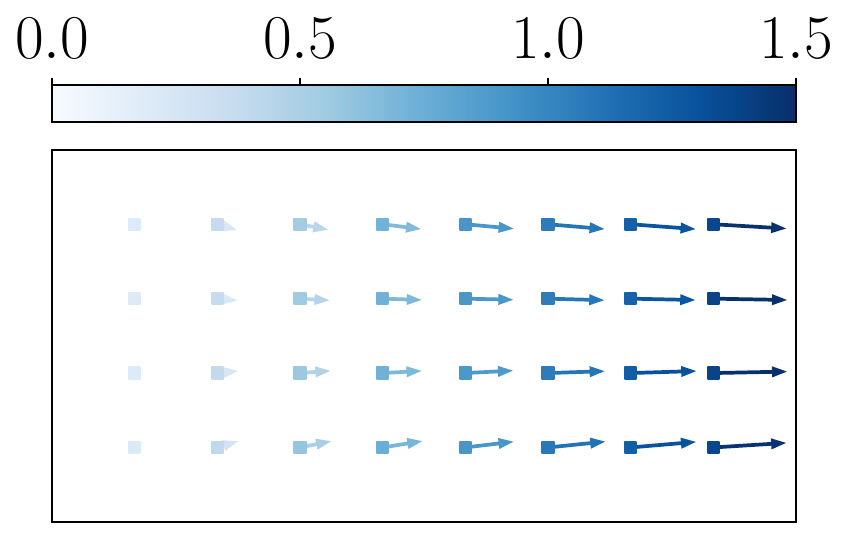} & \includegraphics[width=0.18\textwidth]{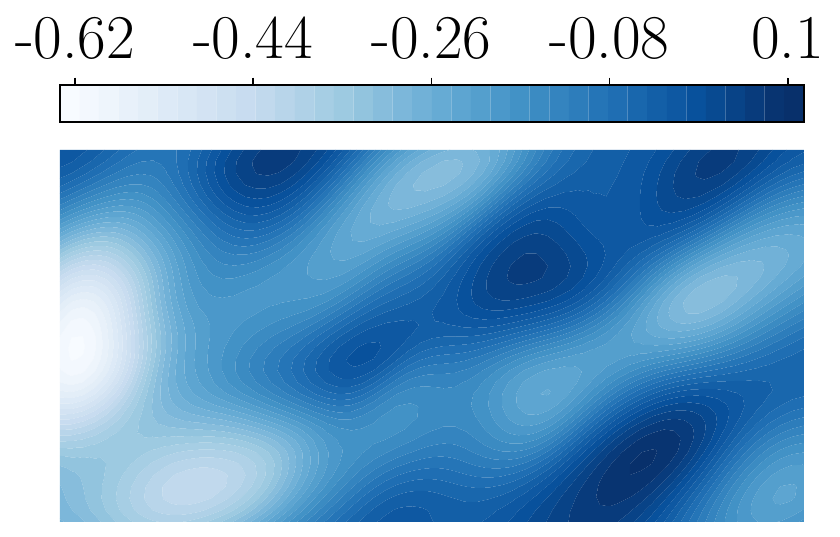} \\
    \end{tabular}
        \renewcommand{\arraystretch}{1.0}
        }
        \addtolength{\tabcolsep}{5pt}
    \caption{\textbf{Example II.} Setup for inferring a heterogeneous hyperelastic material property detailed in \cref{example2}. For each BIP (\#1--4), we visualize the synthetic parameter $\FullParam$ drawn from the prior, the corresponding deformed configuration, the synthetic displacement data $\Obs$, and the MAP estimate $m_{\text{MAP}}^{\Obs}$, i.e., the solution of the deterministic inverse problem.}
    \label{fig:hyper_bip_setup}
\end{figure}

We use a regular grid with $64\times32$ cells. The choice of linear triangular finite elements for $\ParamSpace$ leads to $2,145$ DoFs. The choice of quadratic triangular vector finite elements for $\StateSpace$ leads to $16,770$ DoFs. The Newton--Raphson method is initialized with the homogenous deformation field. For the prior covariance \eqref{eq:covariance_definition}, we induce spatial anisotropy via the following matrix
\begin{gather*}
    \bA = \begin{bmatrix}
        \theta_1\sin(\alpha)^2 + \theta_2\cos(\alpha)^2 & (\theta_1-\theta_2)\sin(\alpha)\cos(\alpha)\\
        (\theta_1-\theta_2)\sin(\alpha)\cos(\alpha) & \theta_1\cos(\alpha)^2 + \theta_2\sin(\alpha)^2
    \end{bmatrix},
\end{gather*}
where $\theta_1 = 2$ and $\theta_2=1/2$, $\alpha = \arctan(2)$. Additionally, we choose $\gamma = 0.3$ and $\delta = 3.3$ which leads to a point-wise marginal variance around $1$, and a spatial correlation of around $2$ and $0.5$ respectively, perpendicular to and along the left bottom to top right diagonal of the material domain. We define the observation operator using $32$ equally spaced interior points $\{\bs_{\text{obs}}^{(j)}\}_{j=1}^{32}$, similar to Example I. This leads to $\ObsDim=64$. The noise distribution has covariance $\NoiseCov = 2.86\times 10^{-3} \identity_{\ObsDim}$, which corresponds to a signal-to-noise ratio of around $500$. We visualize the synthetic data set in \cref{fig:hyper_bip_setup}.

\subsection{Surrogate architecture and training }\label{red_basis_dim}
\paragraph{Reduced basis}
For both problems, we chose the dimension of the parameter latent space to be $d_r=200$ and used 1000 MC samples to compute the reduced basis. We fixed the reduced basis dimension for all studied variational inference methods. Our numerical examples are focused on comparing our proposed \texttt{LazyDINO} with other variational inference methods. Therefore, the effects of varying $d_r$ and MC sample sizes for reduced basis construction are not studied in this work. The eigenvalue decay and basis functions for both examples are visualized in \cref{app:additional_results}. 
\paragraph{Neural network architecture and training} For both examples, we choose 
a dense multi-layer perceptron (MLP) as $\RidgeForward_w$ with 7 hidden layers, each with a width of 400 and a Gaussian error linear unit (GELU) activation function. We train $\RidgeForward_w$ as described in \cref{alg:dino_training} using Adam with $1,500$ epochs and a batch size 25. For \texttt{RB-DINO} training, we used a learning rate of $1\times10^{-3}$ and decreased the learning rate to $3\times10^{-4}$ for the final $375$ epochs. We found that many other training tricks, such as batch normalization or learning rate decay scheduling, were not necessary to produce good generalization for \texttt{RB-DINO}. 

\begin{remark}
    To maintain training stability and prevent overfitting for conventional $L^2_{\mu}$ training, we needed to decrease the learning rate and modify the number of epochs, depending on the training data size. Details reported in~\cref{convention_no_training_details}.
\end{remark}

Implementations of the neural networks and training procedures can be found in \texttt{dinox}, a \texttt{JAX} derivative-informed neural operators library. 

\subsection{Transport map architecture and training} 
\label{transport_map_specs}
\paragraph{Architecture} We chose inverse autoregressive flow (IAF)~\cite{kingma2017improvingvariationalinferenceinverse} as our the transport map $\Map_{\MapParams}:\R^{d_r}\to\R^{d_r}$, with $N_{\Map} = 30$ transport map layers $\tau_i:\R^{d_r}\to\R^{d_r}$ and random input permutation, i.e. $\Map_{\MapParams} = \Map_{\MapParams_{N_{\Map}}}\circ\dots \circ\Map_{\MapParams_1}$ with $\Map_{\MapParams_i} = {\boldsymbol{\mathsf{\tau}}}_{\MapParams_i}\circ \bP_i$ where $\bP_i: \RedCoordParamSpace \to \RedCoordParamSpace$ is the random permutation. Each transport map ${\boldsymbol{\mathsf{\tau}}}_{\MapParams_i}$ is an autoregressive MLP, where inputs are Boolean masked to ensure triangular dependence, i.e., ${\boldsymbol{\mathsf{\tau}}}_{\MapParams_i}(\refCoord_j) = {\boldsymbol{\mathsf{\tau}}}_{\MapParams_i}(\refCoord_1,\ldots, \refCoord_{j-1})$. 
Each autoregressive MLP has 4 hidden layers, each with a width of 400 and GELU activation. In total, the trained parameters are $\MapParams =(\MapParams_1, \ldots, \MapParams_{N_{\Map}})$. For SBAI, the conditional normalizing flow is constructed with a transport map $\Map_{\MapParams}:\R^{d_r}\times\R^{\ObsSpace}$ with a masked autoregressive flow (MAF) with the same architecture, but with a larger input dimension to account for conditioning on observations. We also use the tanh activation function since it empirically produced more stable results. All TMVI methods are implemented and trained via \texttt{lazydinox}, a \texttt{JAX} library for \texttt{LazyDINO} algorithms. 

\paragraph{Training}
The training procedure for all compared methods are taken to be as similar as possible. For \texttt{LazyDINO} and \texttt{LazyNO}, we train $\Map_{\MapParams}$ using Adamax~\cite{kingma2017adammethodstochasticoptimization} with 5 batch sizes, as defined in~\cref{alg:LazyDINO_training}. For minimization $j$, we use $I_j$ iterations using a $B_{j}$--sample MC gradient estimator and learning rate $\alpha_j$, labeled here as $(I_j,B_k,\alpha_j)$: $\Big\{ (5k,200,5\times 10^{-3}), (1k,500,5\times 10^{-3}), (1k,2,000,5\times 10^{-3}), (1k,5,000,5\times 10^{-3}), (1k,7,500,5\times 10^{-4})\Big\}$,
where we decrease the learning rate slightly when we reach the final stochastic approximation batch sample size $7,500$. 
In contrast, for \texttt{LazyMap}, since each Adamax iteration involves PtO map evaluations (referred to as training samples in our results), we use only one batch size, $(I_0,B_0,\alpha_0)=(200,640, 5\times 10^{-3})$ for Example I and $(I_0,B_0,\alpha_0)=(200,100, 5\times 10^{-3})$ for Example II, for a total of 128,000 and 20,000 PtO map evaluations, respectively. For comparison in the proceeding section, error measures are recorded after steps \#($5, 10,\ldots, 320, 640)$, which equal training sample sets of size $(1,000, 2,000, \ldots, 64,000, 128,000)$, respectively. For SBAI, we train with batches of size 100 sampled-without-replacement over epochs with a fixed learning rate of $5\times 10^{-4}$. We terminated optimization when the validation error had not decreased in 10 epochs.

\subsection{Error measures} \label{quantifying_error_specs}

\paragraph{Surrogate approximation error}
The approximation error of the neural ridge function surrogate $\FullRidgeForward_{\boldsymbol{w}}\circ\Projector$ for the PtO map approximation, $\boldsymbol{E}_{\RidgeForward}$, and the PtO map latent Jacobian approximation, $\boldsymbol{E}_{\nabla \RidgeForward}$, are defined as follows.

\begin{align*} 
\boldsymbol{E}_{\RidgeForward} &= \sqrt{\frac{1}{N_{\text{MC}}}\sum_{j=1}^{N_{\text{MC}}}\left[\frac{\|\RidgeForward_w(\refCoord^{(j)})- \boldsymbol{g}^{(j)}\|^2}{\|\boldsymbol{g}^{(j)}\|^2}\right]} &&(\text{Relative PtO map error})\\
 \boldsymbol{E}_{\nabla\RidgeForward} &= \sqrt{\frac{1}{N_{\text{MC}}}\sum_{j=1}^{N_{\text{MC}}}\left[\frac{\|\boldsymbol{J}_r^{(j)}-\nabla\RidgeForward_w(\refCoord^{(j)})\|_{F}^2}{\|\boldsymbol{J}_r^{(j)}\|_{F}^2}\right]}
&&(\text{Relative latent Jacobian error})
\end{align*}
We compute the errors by using $N_{\text{MC}}=5,000$ joint samples of the prior, whitened PtO evaluations and its latent Jacobian evaluations. 

\paragraph{Posterior approximation error}

Posterior approximation accuracy can be assessed in many ways. It is important to ensure posterior accuracy where probability mass is more present, i.e., measures of central concentration or tendency, such as central moments or modes. Accuracy can also be assessed via probability divergences, which measure the overall deviation from the posterior. Accounting for these various forms of accuracy measurement, we consider the quality of posterior approximation under a nonlinear transport map $\FullMap$ via \emph{moment discrepancies} and \emph{density-based diagnostics}, which are described as follows.
\begin{enumerate}
    \item \textbf{Moment discrepancies.} Let $\overline{\FullParam}^{\Obs}\in\ParamSpace$, $\mathcal{C}^{\Obs}\in \HS(\ParamSpace)$ denotes the mean and covariance of $\FullPost$ and ${\mathcal{S}_{25}^{\Obs}}\in\R^{25\times 25\times 25}$ denotes the skewness of $\FullPost$ in the leading 25 latent space coordinates. Let $\overline{\FullParam}^{\FullMap}$, $\mathcal{C}^{\FullMap}$, and $\mathcal{S}_{25}^{\FullMap}$ denote the same quantities for $\FullMap_{\sharp}\mu$. We consider the following relative error in the moments
    \begin{align*}
        \boldsymbol{E}_{\text{mean}} &= \norm{\overline{\FullParam}^{\Obs} - \overline{\FullParam}^{\mathcal{T}}}_{\ParamSpace}/\norm{\overline{\FullParam}^{\Obs}}_{\ParamSpace} && (\text{Relative mean error})\\
        \boldsymbol{E}_{\text{cov}} &= \norm{\mathcal{C}^{\Obs} - \mathcal{C}^{\FullMap}}_{\HS(\ParamSpace)}/\norm{\mathcal{C}^{\Obs}}_{\HS(\ParamSpace)} && (\text{Relative covariance error})\\
        \boldsymbol{E}_{\text{skew}} &= \norm{{\mathcal{S}_{25}^{\Obs}}-{\mathcal{S}_{25}^{\mathcal{T}}}}_{F}/\norm{{\mathcal{S}_{25}^{\Obs}}}_{F} && (\text{Relative skewness error})
    \end{align*}
    The central moments of both $\FullPost$ and $\FullMap_{\sharp}\FullPrior$ are estimated using samples, where samples from $\FullPost$ are obtained using up to $5\times10^6$ MCMC samples using a simplified manifold MCMC method \cite{beskos2008mcmc, cao2024efficient}. The discrepancies are reported in terms of percentages.
    
    Since all central moments must converge as a posterior estimator converges to the posterior, analyzing moment discrepancies of varying orders together is more helpful than analyzing them independently. Estimating higher-order statistics becomes progressively more challenging, so we consider only the first three moments and compute the skewness in the leading dimensions of the latent space.
    
    \item \textbf{Density-based diagnostics.}\label{density_diag} Let $\Phi_{\FullMap}(m) \coloneqq \log\left(\frac{\meas\FullPrior}{\meas(\FullMap_{\sharp}\FullPrior)}(m)\right)$ denote the log density of the prior with respect to the pushforward distribution; see \cref{derivation_of_KL_divergences} for explicit forms. Let $\widetilde{w}(m) = \exp(-2\Misfit(m) + 2\Phi_{\FullMap}(m))$ denote the unnormalized density of the posterior with respect to the pushforward density and $w(m)$ denote its normalization by $\mathbb{E}_{m\sim\FullMap_{\sharp}\FullPrior}\left[\widetilde{w}(m)\right]$. We consider the following quantities related to the quality of each posterior approximation:
    \begin{align*}
        \boldsymbol{E}_{\text{fKL}} &= \mathbb{E}_{\FullParam\sim\FullMap_{\sharp}\FullPrior}\left[ \Misfit(m) - \Phi_{\FullMap}(m)\right] + C_1  && (\text{Shifted rKL divergence})\\
        \boldsymbol{E}_{\text{rKL}} &= \mathbb{E}_{\FullParam\sim\FullMap_{\sharp}\FullPrior}\left[ w(m) (-\Misfit(m) + \Phi_{\FullMap}(m))\right] + C_2  && (\text{Shifted ANIS fKL divergence})\\
     \text{ESS}_{N}\% &= \frac{\left(\frac{1}{N}\sum_{j=1}^{N}\widetilde{w}(m^{(j)})\right)^2}{\frac{1}{N}\sum_{j=1}^{N}\widetilde{w}(m^{(j)})^2}\times \frac{100\%}{N},\quad m^{(j)}\iid \FullMap_{\sharp}\FullPrior  &&(\text{ANIS effective sample size percentage})\\
        \boldsymbol{E}_{\text{MAP}} &= \norm{\FullParam_{\textrm{MAP}}^{\Obs} - \FullParam^{\mathcal{T}_{\textrm{MAP}}}}_{\ParamSpace}/\norm{\FullParam_{\textrm{MAP}}^{\Obs}}_{\ParamSpace} && (\text{Relative MAP point error})
    \end{align*}

   We use both rKL and fKL to measure the posterior approximation error. The former measures the optimality gap due to surrogate error in \texttt{LazyNO} and \texttt{LazyDINO}. 
   On the other hand, rKL can be small when the pushforward distribution is overly concentrated. Therefore, considering both fKL and rKL together provides a fuller picture of posterior approximation accuracy. We use the auto-normalized importance sampling weights~\cite{agapiou2017importancesamplingintrinsicdimension} to compute a biased but consistent estimator for the fKL. We note that rKL is shifted by the normalization constant, and both rKL and fKL are additionally shifted by a constant for visualizations in the log scale. We use $10^5$ i.i.d.\ samples from the pushforward to estimate the expectations in the rKL and fKL divergences.

   We also consider the effective sample size percentage estimated using $N$ i.i.d.\ samples from the pushforward distribution, denoted $\text{ESS}_N\%$. We take $N=10^5$ in our numerical studies. Effective sample size (ESS) is a commonly used diagnostic to assess the quality of approximate posterior sampling, and it has been observed that producing large effective sample percentages can be difficult for many numerical methods~\cite{ImportanceSamplingHighDim,MaiaPolo2020EffectiveSS}. It is related to the forward $\chi^2$ divergence $\chi^2(\FullPost||\FullMap_\sharp\FullPrior) \approx 1-\text{ESS}_N\%/100$~\cite{Sanz_Alonso_2020}.

   Lastly, we study convergence in the MAP estimate, where the ground truth MAP point is obtained by LA.

\end{enumerate}

\section{Numerical results}\label{sec:numerical_results}
In this section, we present numerical results for the two problems described in the previous section. First we show the generalization errors accomplished when training \texttt{RB-DINO} and \texttt{RB-NO} ridge function surrogates with different training sample sizes in~\cref{neural_operator_generalization}. These errors are tied to the accuracy of the surrogate-based posterior approximation via \cref{theorem:kl_bound}, and to the optimality gap in surrogate-driven LMVI via \cref{corollary:optimality_gap}.

In the subsequent results, we compare the posterior errors for \texttt{LazyDINO} and \texttt{LazyNO}. As points of comparison, we additionally consider the Laplace approximation as a baseline, SBAI, as well as \texttt{LazyMap}, which utilizes the true PtO map in training instead of the surrogate. 

\begin{remark}
    In the neural operator training results, we compare the computational costs measured in terms of nonlinear PDE solves. The additional costs associated with the \texttt{RB-DINO} training are measured in a cost basis that is relative to the nonlinear PDE solves. Since these additional computational costs associated with the latent Jacobian computations are negligible due to the use of sparse direct solvers \cite[Section 4.3]{cao2024efficient}, this point of comparison is not considered for the posterior approximation error comparisons. 
\end{remark}

\subsection{Neural operator ridge function generalization}\label{neural_operator_generalization}

\begin{figure}[htbp]
    \centering
    \begin{tabular}{c c}
    \multicolumn{2}{c}{\includegraphics[width=0.05\textwidth]{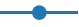}\hspace{0.01\textwidth} \texttt{RB-DINO} (derivative-informed learning)\hspace{0.03\textwidth}\includegraphics[width=0.05\textwidth]{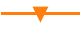} \hspace{0.01\textwidth}  \texttt{RB-NO} (conventional supervised learning)}\\[0.1 in]
    \hspace{0.06 \textwidth}\makecell{PtO map generalization\\ relative error, $\boldsymbol{E}_{\RidgeForward}$}& \hspace{0.06 \textwidth}\makecell{Latent Jacobian generalization\\ relative error, $\boldsymbol{E}_{\nabla\RidgeForward}$} \\
        \includegraphics[width = 0.45\textwidth]{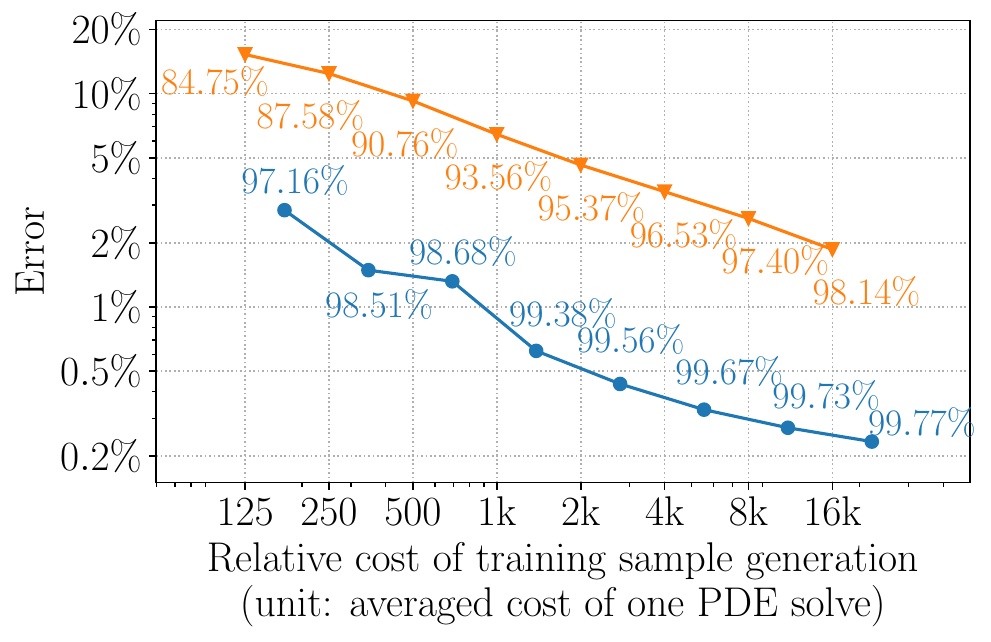} & 
         \includegraphics[width = 0.45\textwidth]{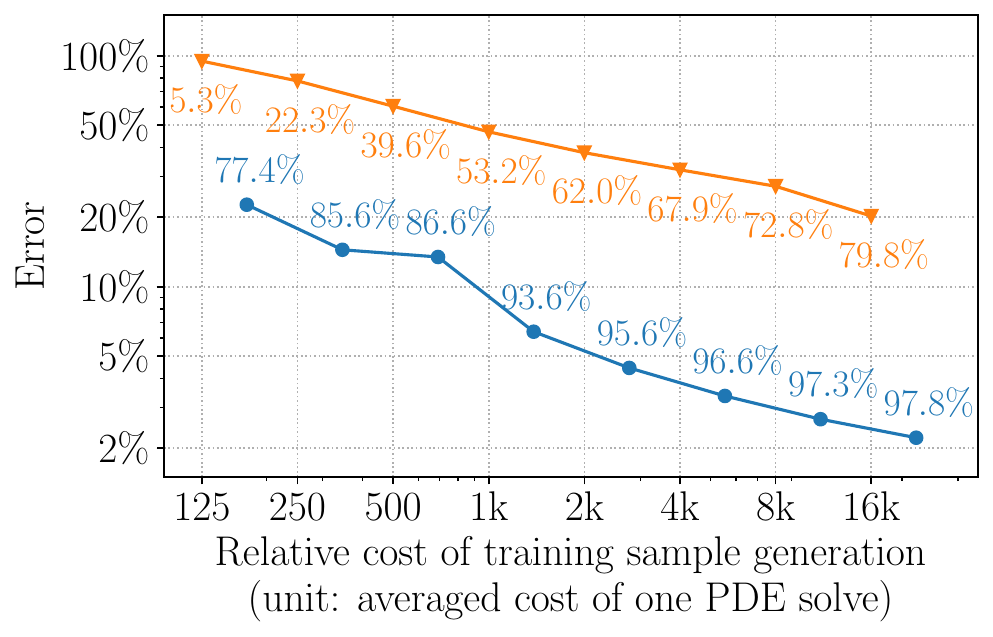}
    \end{tabular}
    \caption{ \textbf{Example I neural ridge function surrogate testing.} Percentage accuracy, $100\% \times \left(1- \text{error}\right)$, is also overlaid. 
    Overall, these results demonstrate a significant cost reduction for achieving any given generalization accuracy in both the PtO map and the latent Jacobian via the derivative-informed learning method. While convergence rates are similar, \texttt{RB-DINO} has a greater than $64\times$ higher cost efficiency measure in relative errors. Due to the bounds in~\cref{H1mu}, we expect this to reflect in increased sample efficiency of \texttt{LazyDINO} over \texttt{LazyNO} in posterior error measures. We note a statistical anomaly in \texttt{RB-DINO} training encountered at 500 training samples, which poses a downstream impact on the posterior error measures in subsequent figures.}
  \label{fig:gen_error_1}
\end{figure}

\begin{figure}[htbp]
    \centering
    \begin{tabular}{c c}

    \multicolumn{2}{c}{\includegraphics[width=0.05\textwidth]{figures/legend_dino_full.pdf}\hspace{0.1 in} \texttt{RB-DINO} (derivative-informed learning)\hspace{0.2 in}\includegraphics[width=0.05\textwidth]{figures/legend_no_full.pdf} \texttt{RB-NO} (conventional supervised learning))}\\[0.1 in]
    \hspace{0.06\textwidth}\makecell{Observable prediction generalization\\relative error, $\boldsymbol{E}_{\RidgeForward}$}& \hspace{0.07\textwidth}\makecell{Latent Jacobian generalization \\relative error, $\boldsymbol{E}_{\nabla\RidgeForward}$} \\
        \includegraphics[width = 0.45\textwidth]{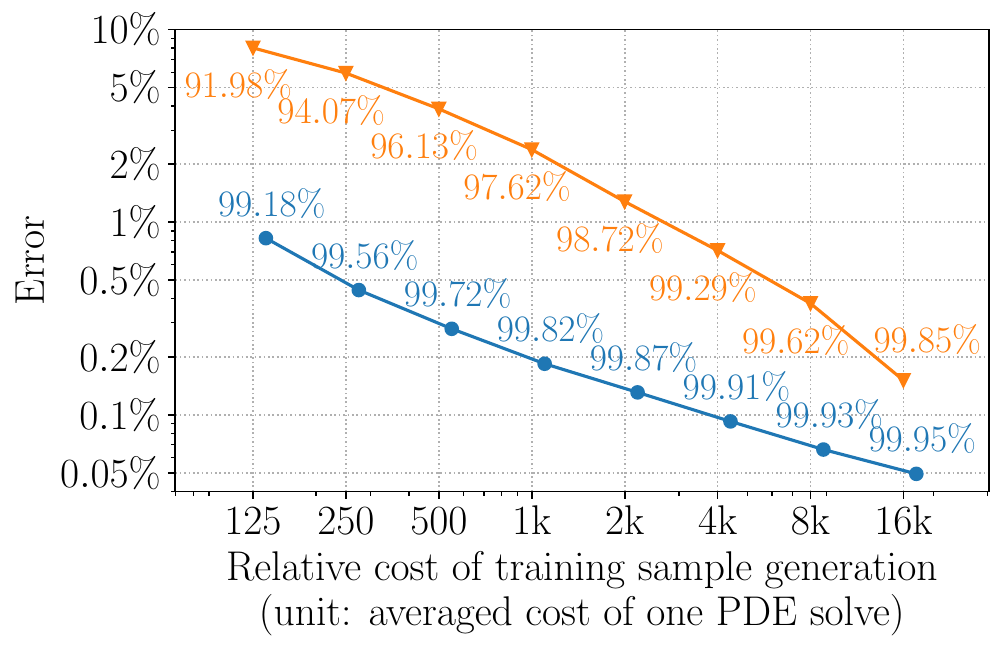} & 
         \includegraphics[width = 0.45\textwidth]{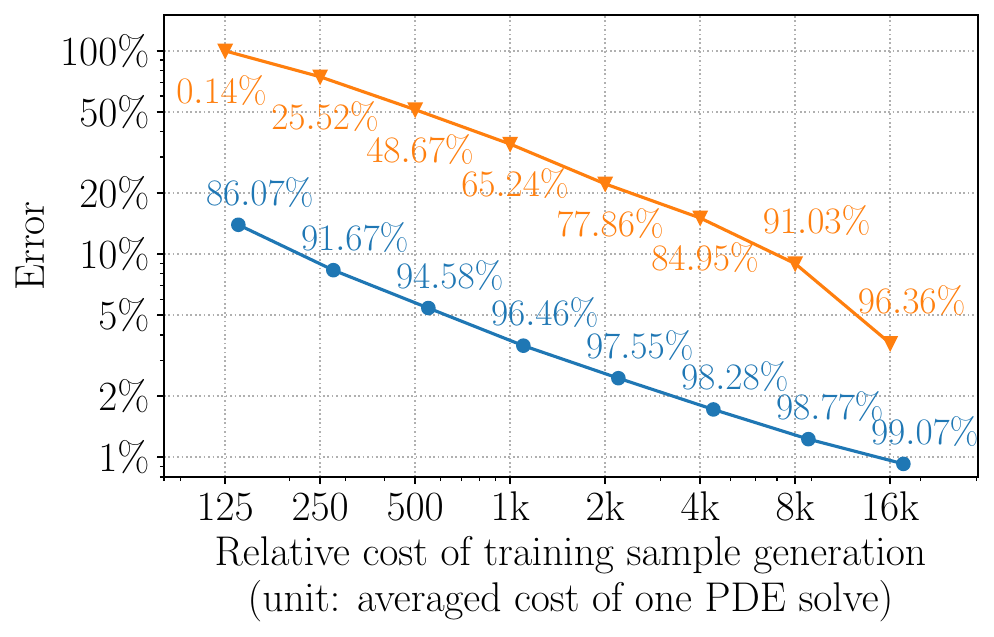}
    \end{tabular}
    \caption{\textbf{Example II neural operator ridge function generalization.} In a similar pattern as~\cref{fig:gen_error_1}, \texttt{RB-DINO} enjoys $8$--$32\times$ lower data generation cost over \texttt{RB-NO} for achieving a given accuracy. We expect this to reflect in increased sample efficiency of \texttt{LazyDINO} over \texttt{LazyNO} in posterior error measures.}
    \label{fig:gen_error2}
\end{figure}

We begin by assessing the results of \texttt{RB-DINO}, and \texttt{RB-NO} surrogate construction. The results for Example I (a nonlinear reaction--diffusion PDE) can be found in ~\cref{fig:gen_error_1}, while the results for Example II (deformation of a hyperelastic thin film) can be found in ~\cref{fig:gen_error2}. Overall, the trend demonstrates that the derivative-informed learning method leads to a significant cost reduction in learning the latent representation of the PtO map and its Jacobian compared to conventional supervised learning. These results are consistent with those in \cite{o2024derivative,cao2024efficient}. We expect that the improvement of the PtO map approximation leads to better fidelity in the surrogate approximated posterior through \cref{theorem:kl_bound}, and the combination of improved PtO map approximation and Jacobian approximation leads to more accurate LMVI optimization through \cref{corollary:optimality_gap}.

\newpage

\subsection{Posterior approximation error}\label{posterior_error}


In this section we investigate the posterior error metrics described in~\cref{quantifying_error_specs}. As a reminder, we compare \texttt{LazyDINO} against \texttt{LazyNO}, the PDE-driven \texttt{LazyMap} as well as two additional baselines: the Laplace approximation and SBAI via conditional transport.

For all plots in this section, the markers for BIP \#1--4 are labeled as \wye[rotate=180], \wye, \wye[rotate=90], and \wye[rotate=270]. The average error is plotted in a darker color, and a line is drawn between each average to visualize the trend. We cut off vertical-axis errors at 200--300\%, depending on the plot, for readability since the scale of the errors across methods varies widely. 

For the LA-baseline, horizontal lines are included to facilitate visual comparison, even though it is only computed once for each BIP--the horizontal axis, i.e., the number of training samples is not meaningful in this context. A conservative estimate of training cost equivalent to 100 training samples, used for MAP estimation and Hessian-inverse covariance estimation, is marked in the plots. To potentially outperform the LA-baseline in amortized Bayesian inversion, a method should achieve lower posterior error as few amortized PtO and Jacobian evaluations as possible. 

\begin{figure}[htbp]

    \centering
    \addtolength{\tabcolsep}{-6pt}
    \begin{tabular}{C C}
    \hspace{0.07\textwidth}Posterior mean relative error, $\boldsymbol{E}_{\text{mean}}$& \hspace{0.07\textwidth}Posterior covariance relative error, $\boldsymbol{E}_{\text{cov}}$ \\  
    \includegraphics[width=0.5\textwidth]{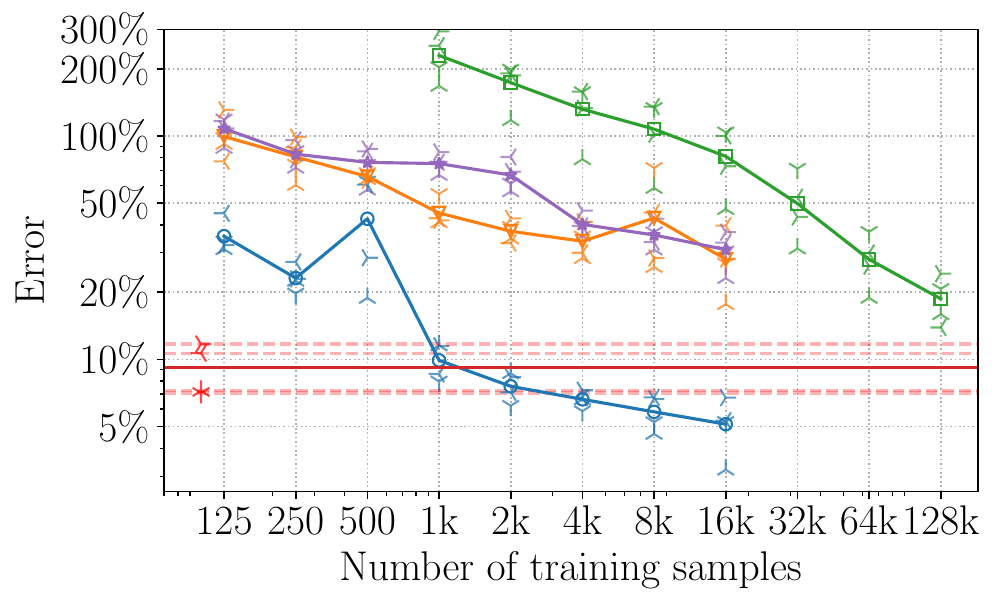} & \includegraphics[width=0.5\textwidth]{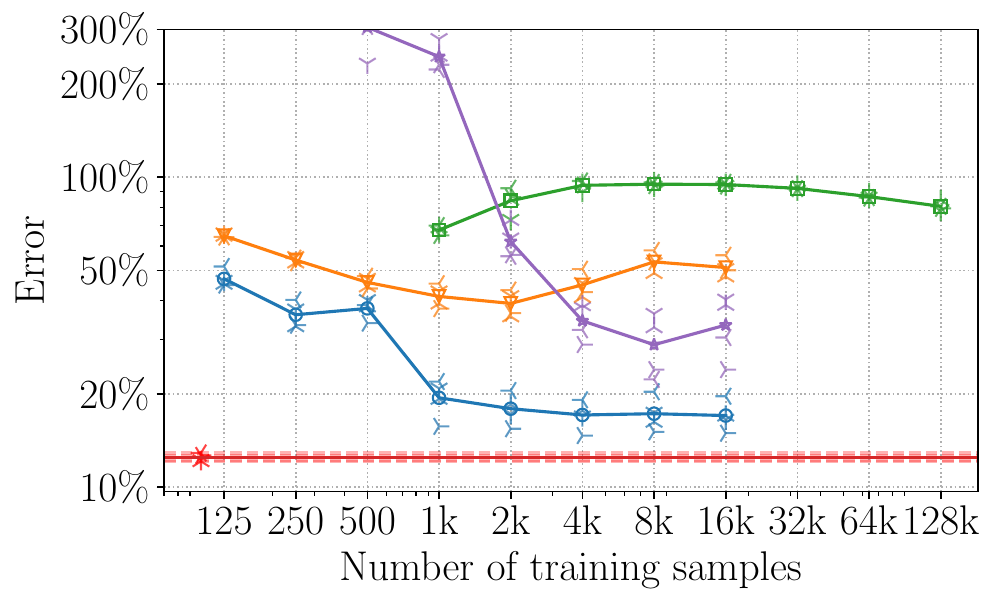}\\
    {\hspace{0.07\textwidth}\makecell{Posterior skewness relative error\\
    in $25$ leading latent space coordinates, $\boldsymbol{E}_{\text{skew}}$}} &
    \\
    \includegraphics[width=0.5\textwidth]{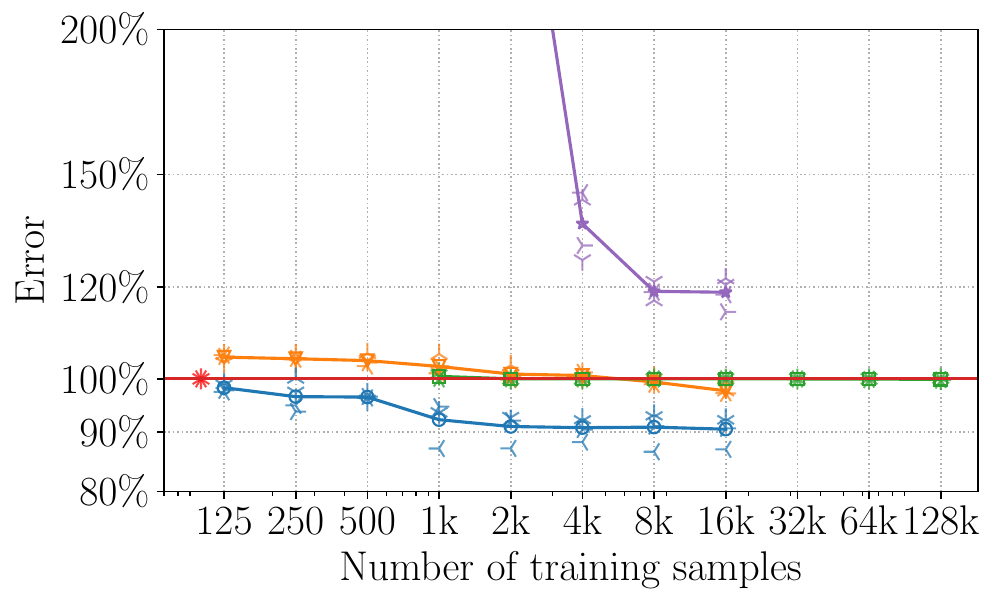} &
    \vspace{-0.02125\textwidth}
             \addtolength{\tabcolsep}{6pt}
             \hspace{0.9cm}{\fbox{
\begin{tabular}{l l}
    \multicolumn{2}{l}{\makecell{Training samples \textit{amortized} over 4 BIPs}}\\
    \includegraphics[width=0.04\textwidth]{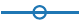} & \texttt{LazyDINO} \\
    \includegraphics[width=0.04\textwidth]{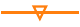} & \texttt{LazyNO} \\
    \includegraphics[width=0.04\textwidth]{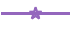} & SBAI\\[0.0386875\textwidth]
    \multicolumn{2}{l}{\makecell{Training samples \textit{repeated} for each BIP}}\\
    \includegraphics[width=0.04\textwidth]{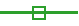} & \texttt{LazyMap} \\
    \includegraphics[width=0.04\textwidth]{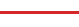} & LA-baseline
\end{tabular}
}}
    
        \addtolength{\tabcolsep}{-6pt}
    \end{tabular}
        \addtolength{\tabcolsep}{6pt}
    \caption{\textbf{Example I moment discrepancies.}  In all cases, a lower error indicates a better posterior approximation. (\texttt{LazyDINO} vs. \texttt{LazyNO}): 
    Apart from the statistical anomaly at 500 samples attributed to the stochasticity in surrogate training seen in~\cref{fig:gen_error_1}, \texttt{LazyDINO} is $64\times$ more sample-efficient measured in mean relative error, over $4\times$ in covariance relative error, and over $64\times$ in skewness relative error. The discrepancy in efficiency is even more pronounced in the higher sample regime ($>500$ samples). (\texttt{LazyDINO} vs. SBAI): \texttt{LazyDINO} is $64\times$ more sample-efficient measured in mean relative error, and SBAI is uncompetitive in covariance and skewness error in all sample regimes. (\texttt{LazyDINO} vs. \texttt{LazyMap}): In all error measures, \texttt{LazyDINO} achieves orders of magnitude higher sample-efficiency compared to \texttt{LazyMap}. We also note that since \texttt{LazyMap} repeats computations of the PtO for each BIP, the number of training samples in total is $4\times$ the number for the other approaches. (\texttt{LazyDINO} vs. LA-baseline): \texttt{LazyDINO} achieves lower error in the mean and the skewness, particularly noticeable in the higher sample size regime. While no approach studied in this work can achieve lower relative error in the covariance compared to LA-baseline, we note that the density-based diagnostics in~\cref{fig:NRD_Density_diagonistics} lend further support to the proposition that \texttt{LazyDINO} improves upon the baseline, even for small sample sizes.}
    \label{fig:NRD_moment_discrepancies}
\end{figure}

\begin{figure}[htbp]

    \centering
    \addtolength{\tabcolsep}{-6pt}
    \begin{tabular}{C C}
    \hspace{0.07\textwidth}Posterior mean relative error, $\boldsymbol{E}_{\text{mean}}$& \hspace{0.07\textwidth}Posterior covariance relative error, $\boldsymbol{E}_{\text{cov}}$ \\  
    \includegraphics[width=0.5\textwidth]{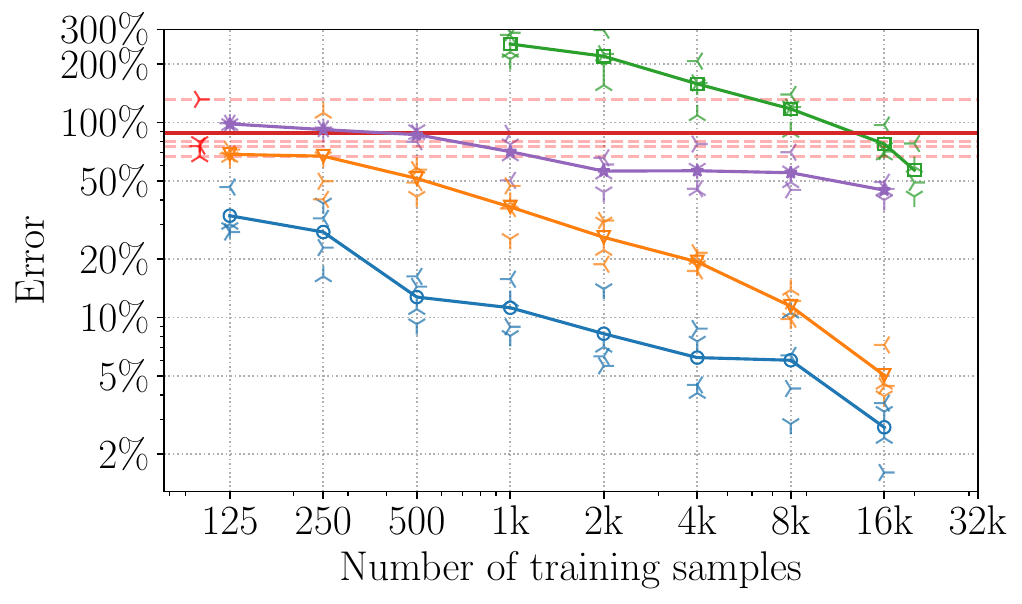} & \includegraphics[width=0.5\textwidth]{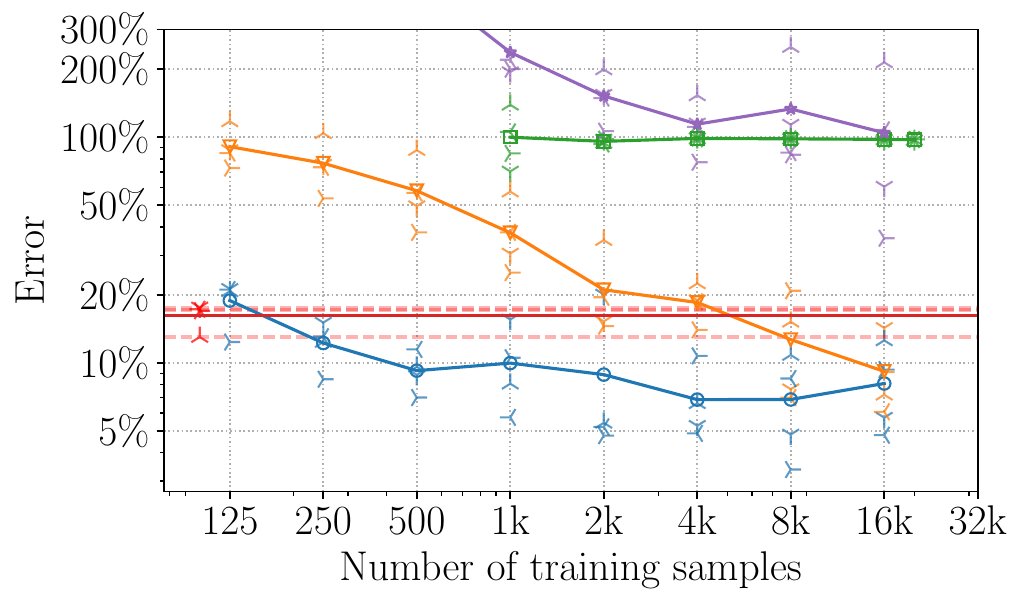}\\
    {\hspace{0.07\textwidth}\makecell{Posterior skewness relative error\\
    in $25$ leading latent space coordinates, $\boldsymbol{E}_{\text{skew}}$}} &
    \\
    \includegraphics[width=0.5\textwidth]{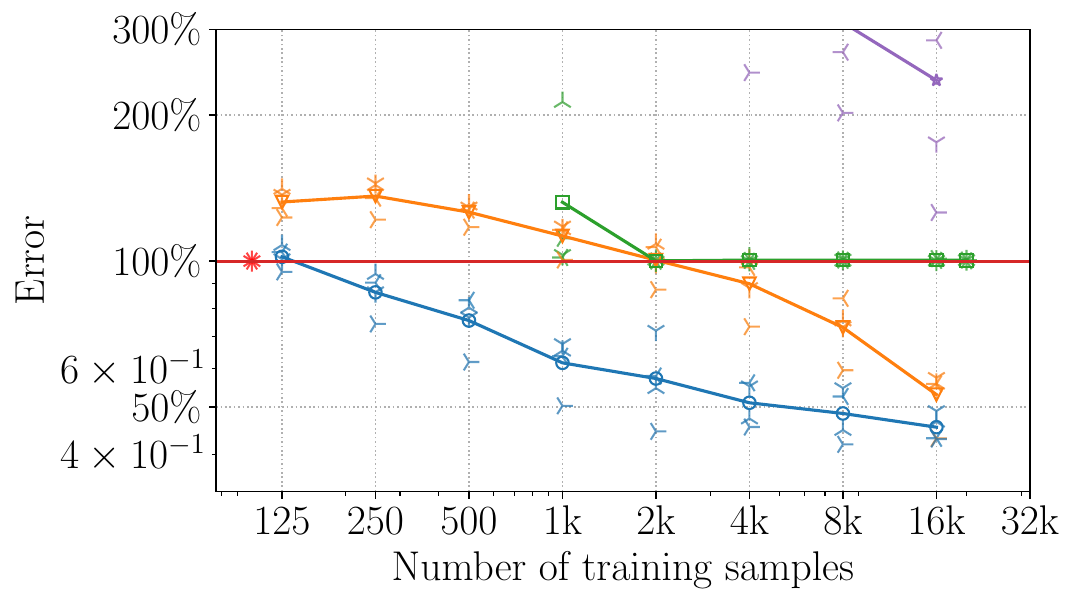} &
    \vspace{-0.018\textwidth}
             \addtolength{\tabcolsep}{6pt}
              \hspace{0.9cm}{\fbox{
\begin{tabular}{l l}
    \multicolumn{2}{l}{\makecell{Training samples \textit{amortized} over 4 BIPs}}\\
    \includegraphics[width=0.04\textwidth]{figures/legend_dino_empty.pdf} & \texttt{LazyDINO} \\
    \includegraphics[width=0.04\textwidth]{figures/legend_no_empty.pdf} & \texttt{LazyNO} \\
    \includegraphics[width=0.04\textwidth]{figures/legend_cm_empty.pdf} & SBAI\\[0.0215\textwidth]
    \multicolumn{2}{l}{\makecell{Training samples \textit{repeated} for each BIP}}\\
    \includegraphics[width=0.04\textwidth]{figures/legend_pde_empty.pdf} & \texttt{LazyMap} \\
    \includegraphics[width=0.04\textwidth]{figures/legend_la.pdf} & LA-baseline
\end{tabular}
}}
    
        \addtolength{\tabcolsep}{-6pt}
    \end{tabular}
        \addtolength{\tabcolsep}{6pt}
    \caption{\textbf{Example II moment discrepancies.}     
    We observe similar trends as in Example I (\cref{fig:NRD_moment_discrepancies}). (\texttt{LazyDINO} vs. \texttt{LazyNO}) The trend of consistent outperformance of \texttt{LazyNO} is clear in this example; the derivative-informed learning of \texttt{RB-DINO} yields $2-16\times$ higher sample efficiency.
    (\texttt{LazyDINO} vs. SBAI) The best-performing SBAI at the high sample regime is still less accurate in posterior approximation compared to the worst \texttt{LazyDINO} at the low sample regime.
    (\texttt{LazyDINO} vs. \texttt{LazyMap}) \texttt{LazyMap} exhausts the training sample budget before performing comparably to \texttt{LazyDINO}. 
    (\texttt{LazyDINO} vs. \texttt{LA-baseline}) At 250 training samples, \texttt{LazyDINO} produces lower error than the LA-baseline in all moment discrepancies.
    }
    \label{fig:hyper_moment_discrepancies}
    \bigskip
\end{figure}

We begin by comparing moment discrepancies for Example I and II in \cref{fig:NRD_moment_discrepancies} and \cref{fig:hyper_moment_discrepancies}, respectively. In general the LA-baseline provided a reasonable baseline point of comparison; and in the case of covariance approximations for Example I, it consistently outperformed each method. In all other cases, however, \texttt{LazyDINO} eventually produced substantially better predictions of moments, particularly given a lot of data (e.g., consider the mean error for \cref{fig:hyper_moment_discrepancies}). Of the remaining methods, \texttt{LazyDINO} produced the best matching moments for each problem in each training sample size. \texttt{LazyNO} was typically the next best performing method, although in some cases SBAI or \texttt{LazyMap} performed comparably or slightly better. A notable phenomenon was the relatively poor performance of SBAI, which was typically more than an order of magnitude worse than \texttt{LazyDINO}. The poor performance relatively of \texttt{LazyMap} is easily explained by the sample intensity required to minimize the latent space transport training objective reliably. In particular, we artificially truncated the sampling budget at $128,000$, while the \texttt{LazyDINO/NO} required 16 million total samples over all iterations. This point of comparison demonstrates the essential benefit of the \texttt{LazyDINO} approach: by first building a reliable PtO surrogate over the prior using a fixed number of samples, we can later enable an optimization algorithm requiring orders of magnitude more samples. 

\begin{figure}[hbtp]
    \centering
    \addtolength{\tabcolsep}{-6pt} 
    \begin{tabular}{C C}
    \hspace{0.06\textwidth}Shifted rKL divergence error, $\boldsymbol{E}_{\text{rKL}}$&  \hspace{0.06\textwidth}\makecell{ANIS Effective Sample Percentage, $\text{ESS}_{100\text{k}}\%$} \\
    \includegraphics[width=0.5\textwidth]{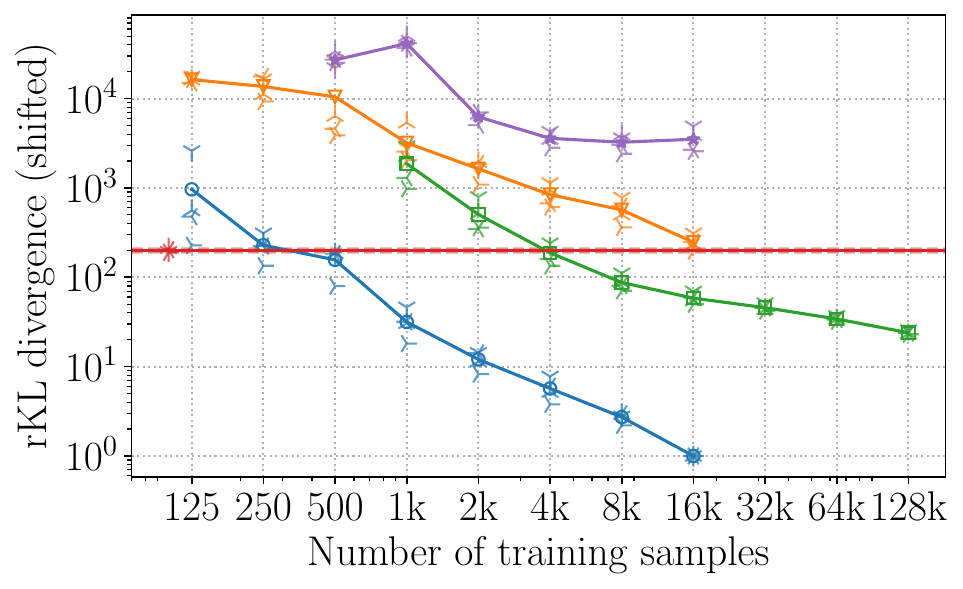} & \includegraphics[width=0.5\textwidth]{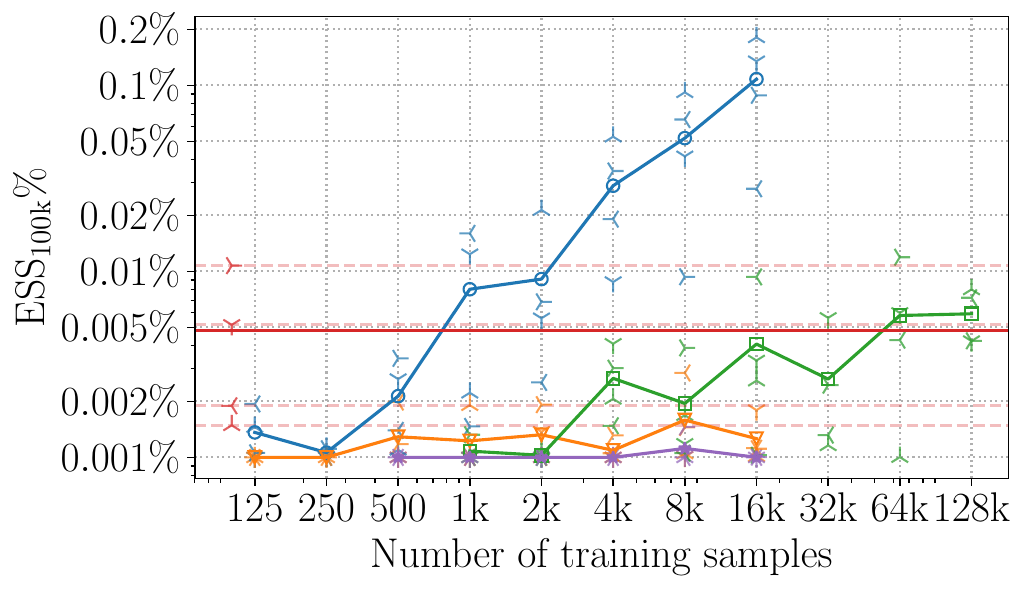}\\
    \hspace{0.06\textwidth} Shifted ANIS fKL divergence error, $\boldsymbol{E}_{\text{fKL}}$ & \hspace{0.06\textwidth}  MAP estimate relative error, $\boldsymbol{E}_{\text{MAP}}$  \\
    \includegraphics[width=0.5\textwidth]{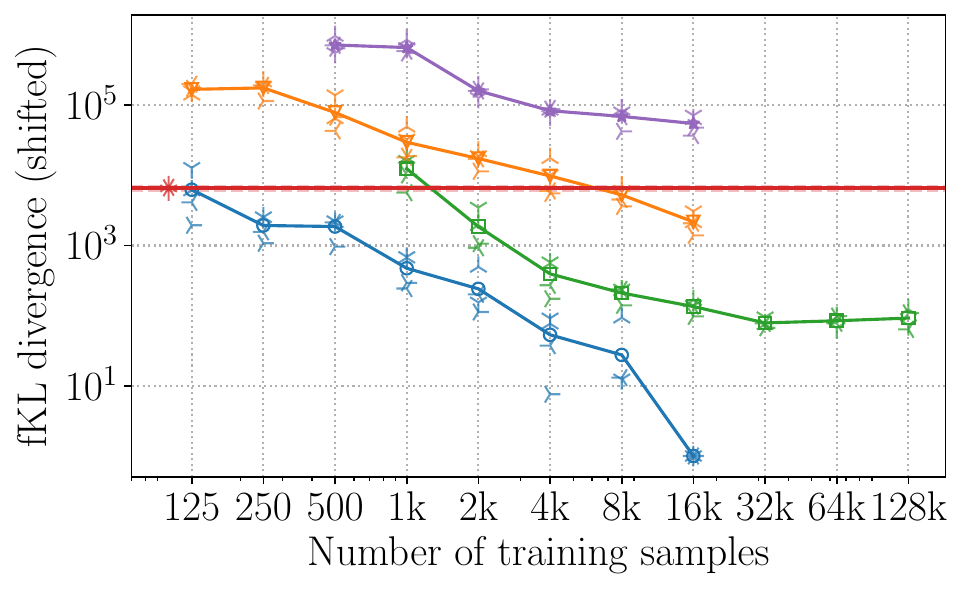}&\includegraphics[width=0.5\textwidth]{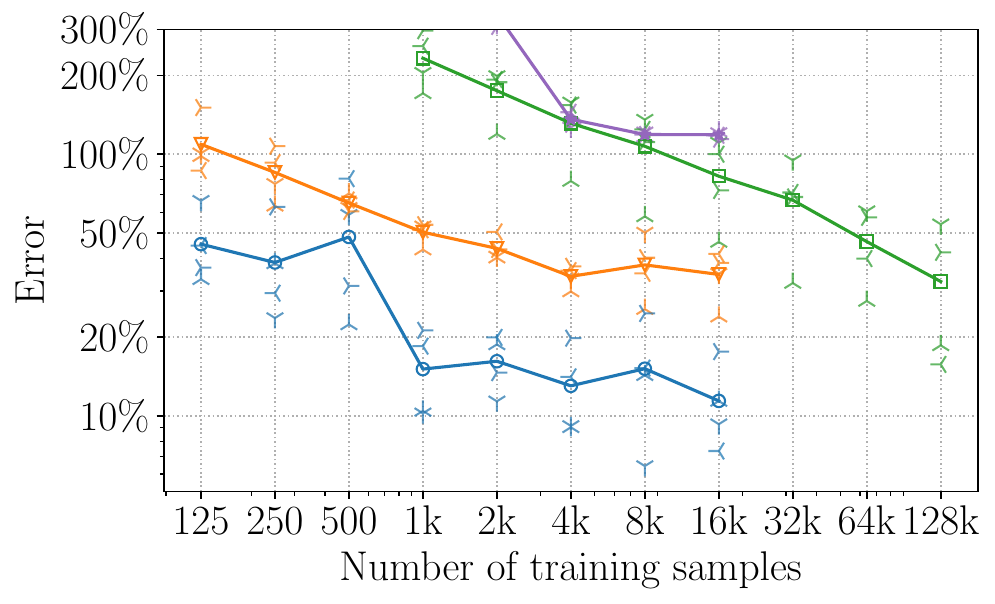}
    \end{tabular}
    \addtolength{\tabcolsep}{6pt}
    
 \begin{adjustwidth}{1.120cm}{}
\begin{tabular}{|l l l l c  | l l|}\hline
    \multicolumn{4}{|l}{\makecell{Training samples \textit{amortized} over 4 BIPs}} & &   \multicolumn{2}{l|}{\makecell{Training samples \textit{repeated} for each BIP}}\\
    \includegraphics[width=0.04\textwidth]{figures/legend_dino_empty.pdf} & \texttt{LazyDINO} & \includegraphics[width=0.04\textwidth]{figures/legend_cm_empty.pdf} & SBAI & &\includegraphics[width=0.04\textwidth]{figures/legend_pde_empty.pdf} & \texttt{LazyMap} \\
    \includegraphics[width=0.04\textwidth]{figures/legend_no_empty.pdf} & \texttt{LazyNO} & & &&\includegraphics[width=0.04\textwidth]{figures/legend_la.pdf} & LA-baseline \\\hline
\end{tabular}
\end{adjustwidth}

    \caption{\textbf{Example I density-based diagnostics.} Higher values of $\text{ESS}_{100\text{K}}\%$ and lower values of all other diagnostics imply better posterior approximation. We observe similar trends to those observed in the moment discrepancy comparisons in~\cref{fig:NRD_moment_discrepancies}. Notably, the \texttt{LazyDINO} eventually yields the best error in each case. \texttt{LazyDINO} enjoys over $8$--$128\times$ higher sample efficiency compared to the other methods. Though $\text{ESS}_{100\text{K}}\%$ is low across the board, \texttt{LazyDINO} produces an impressive $\approx100$ effective samples while other methods only achieve $\approx1$--$6$ effective samples.}
    \label{fig:NRD_Density_diagonistics}
\end{figure}

\begin{figure}[hbtp]
    \centering

    \addtolength{\tabcolsep}{-6pt} 
    \begin{tabular}{C C}
    \hspace{0.06\textwidth}Shifted rKL divergence error, $\boldsymbol{E}_{\text{rKL}}$&  \hspace{0.06\textwidth}\makecell{ANIS Effective Sample Percentage, $\text{ESS}_{100\text{k}}\%$} \\
    \includegraphics[width=0.5\textwidth]{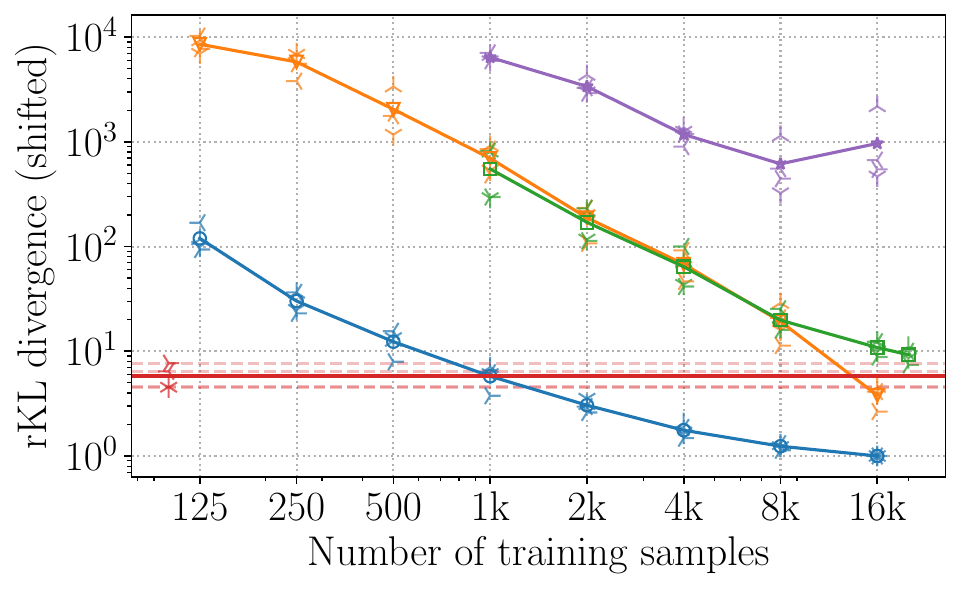} & \includegraphics[width=0.5\textwidth]{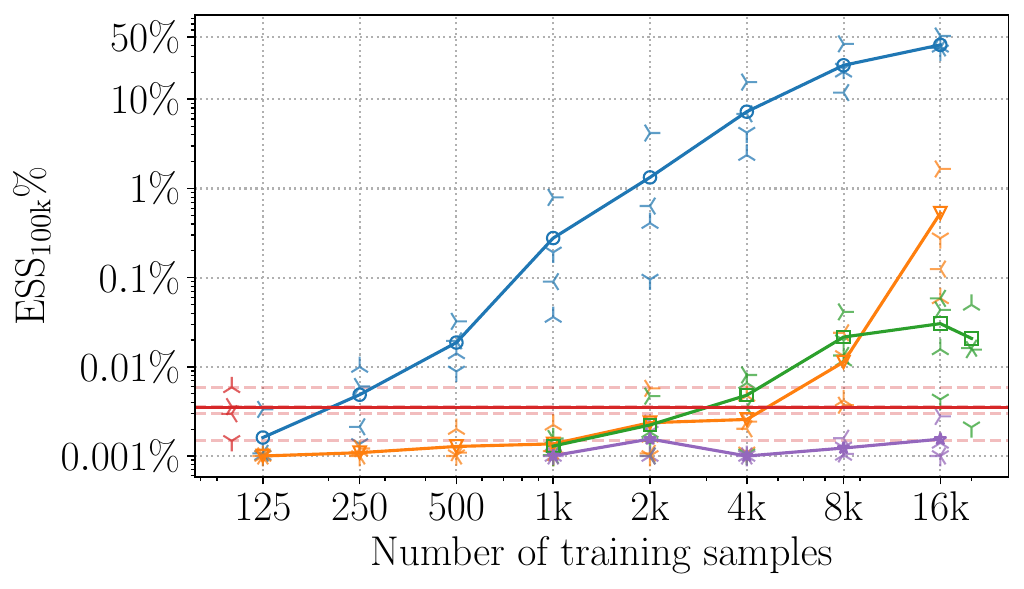}\\
    \hspace{0.06\textwidth} Shifted ANIS fKL divergence error, $\boldsymbol{E}_{\text{fKL}}$ & \hspace{0.06\textwidth}  MAP estimate relative error, $\boldsymbol{E}_{\text{MAP}}$  \\
    \includegraphics[width=0.5\textwidth]{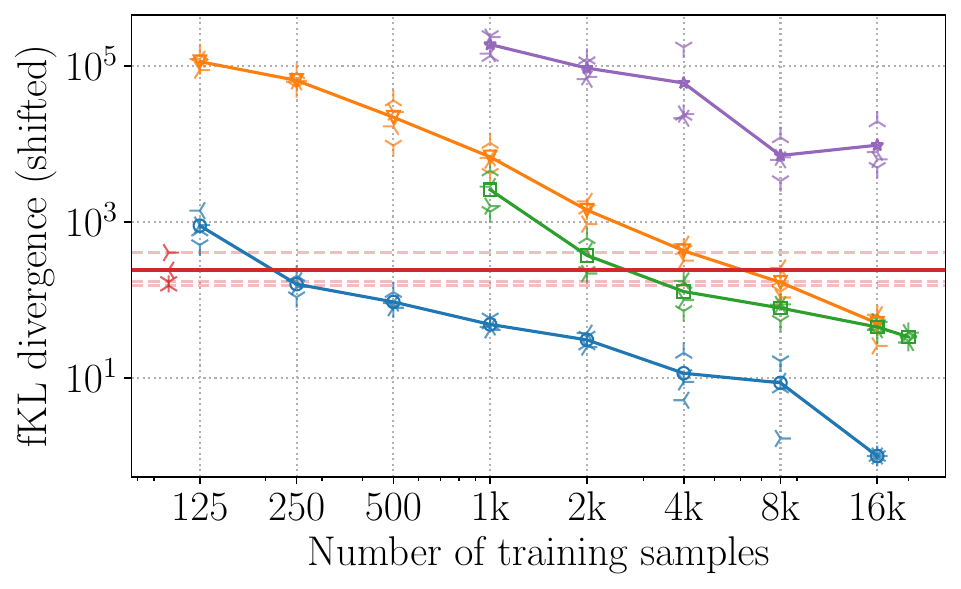}&\includegraphics[width=0.5\textwidth]{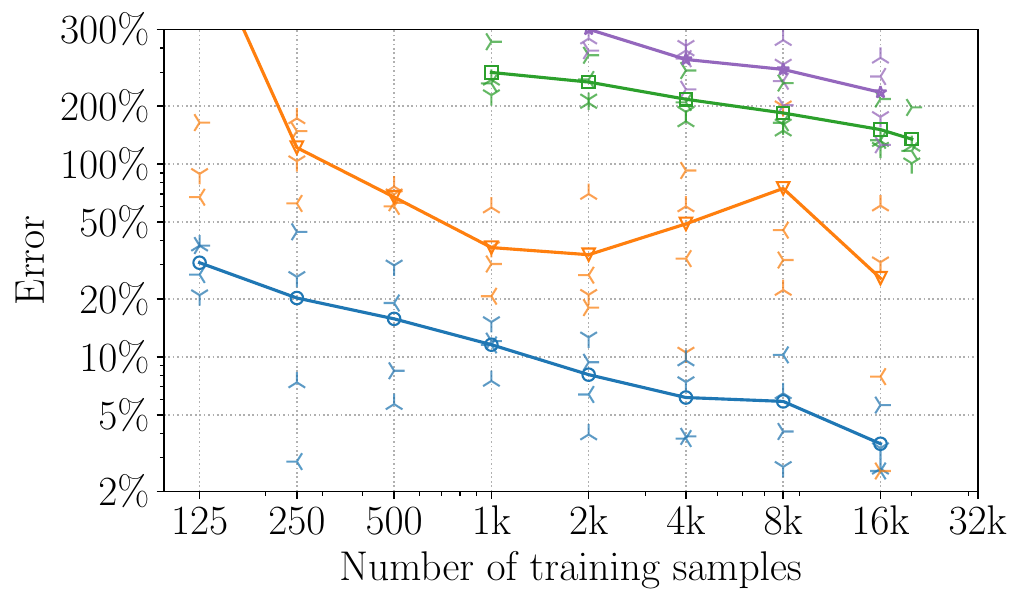}
    \end{tabular}
    \addtolength{\tabcolsep}{6pt}

 \begin{adjustwidth}{1.120cm}{}
\begin{tabular}{|l l l l c  | l l|}\hline
    \multicolumn{4}{|l}{\makecell{Training samples \textit{amortized} over 4 BIPs}} & &   \multicolumn{2}{l|}{\makecell{Training samples \textit{repeated} for each BIP}}\\
    \includegraphics[width=0.04\textwidth]{figures/legend_dino_empty.pdf} & \texttt{LazyDINO} & \includegraphics[width=0.04\textwidth]{figures/legend_cm_empty.pdf} & SBAI & &\includegraphics[width=0.04\textwidth]{figures/legend_pde_empty.pdf} & \texttt{LazyMap} \\
    \includegraphics[width=0.04\textwidth]{figures/legend_no_empty.pdf} & \texttt{LazyNO} & & &&\includegraphics[width=0.04\textwidth]{figures/legend_la.pdf} & LA-baseline \\\hline
\end{tabular}
\end{adjustwidth}

    \caption{\textbf{Example II density-based diagnostics.} We observe similar trends as in Example I. Notably, \texttt{LazyDINO} achieves nearly 50\% ANIS effective sample percentage as a 100,000-sample independent sampler with 16,000 training samples, 50 to 50,000 times the percentage achieved for competing methods. }
    \label{fig:hyper_Density_diagonistics}
\end{figure}

In the next set of results, we consider various density-based diagnostic criteria, which are defined in \cref{quantifying_error_specs}. In \cref{fig:NRD_Density_diagonistics} and \cref{fig:hyper_Density_diagonistics}, we compare the performance of the different methods through their shifted rKL and fKL, ANIS effective sample percentage, and MAP point estimates. As with the moment discrepancy results, we see again the consistent superior performance of \texttt{LazyDINO} compared to the other TMVI methods and the LA-baseline for $>500/1000$ samples for Example I and II,
respectively. Notably, the \texttt{LazyDINO} effective sample size is orders of magnitude higher than the other methods for both examples in the largest sample case.

\subsection{Timing comparisons}\label{amortized_actual_times}

In the previous section, we compared various methods on a sample-cost basis. In this section, we consider additional computational speedups, such as parallelism, to make a comparison on a time-cost basis. 

We begin by making empirical comparisons of \texttt{LazyDINO} and the original \texttt{LazyMap} algorithm, taking into account that \texttt{LazyMap} can be made more efficient with parallelism. We continue by comparing the online evaluation costs for SBAI and \texttt{LazyDINO/NO}, demonstrating that in addition to being much more accurate than SBAI, \texttt{LazyDINO/NO} have a smaller overall cost in an amortized setting. 

\begin{remark}
Compute times will vary based on computing environments, but we note that all GPU computations were performed on Nvidia A100 GPUs with 40 and 80GB of RAM, and all CPU computations were performed on an Intel Xeon Gold 6248R 3.00GHz CPU with 1.2 Terabytes of RAM (CPU computations were compute-bound, not memory bound).
\end{remark}
\paragraph{\texttt{LazyMap} vs. \texttt{LazyDINO}}

The efficiency of \texttt{LazyDINO} is impacted by the compute time for the offline phase. Since \texttt{LazyMap} can exploit parallelism within each iteration, depending on the parallel computing resources available, \texttt{LazyMap} might still be competitive in particular sample size regimes since it does not require first training a surrogate. For this comparison, we investigate the relative performance of \texttt{LazyMap} and \texttt{LazyDINO} for a similar end-to-end computational budget while allowing for parallelism in the \texttt{LazyMap} calculations. 

In ~\cref{tab:alg_timing_1} and~\cref{tab:alg_timing_2}, we provide a comparison of solution time and posterior mean accuracy for the two methods, given access to 20 concurrent CPU evaluations of the PtO map and Jacobian action. All parallel times provided are theoretical (without accounting for communication), computed to two decimals of accuracy, since we lacked a parallel implementation of the PtO map.  We provide actual compute times for training sample data generation for \texttt{LazyDINO} and the training times for \texttt{LazyDINO} and \texttt{LazyMap}. The \texttt{LazyDINO} offline phase times are reported amortized across four instances of observation data. All computations repeated for each instance of observation data are reported as an average over the four instances and rounded up to the nearest 10 seconds. 

\begin{table}[htbp]
    \centering
    {\renewcommand{\arraystretch}{1.5}
    \begin{tabular}{||c|r|c|c|c|c|| }\hline
        & \multirow{3}{*}{\textbf{Category (unit)}} & \multicolumn{4}{c||}{\textbf{LMVI Method}} \\ \cline{3-6}
        & & \makecell{\texttt{LazyMap} \\ (16k)} & \makecell{\texttt{LazyMap} \\ (128k)} & \makecell{\texttt{LazyDINO} \\ (1k)} &  \makecell{\texttt{LazyDINO} \\ (16k)} \\ \hline
        
        \multirow{7}{*}{\rotatebox[origin=c]{90}{Algorithm steps}} 
        & Amortized PtO evaluations (sec) &  -- & -- & 130/4 & 1,950/4 \\ \cline{2-6}
        & Parallel amortized PtO evaluations (sec) & -- & -- & 6.5 /4 & 97.5/4 \\ \cline{2-6}
        & Amortized Jacobian (sec) & -- & -- & 50/4 & 750/4 \\ \cline{2-6}
        & Parallel amortized Jacobian (sec) & -- & -- & 2.5/4 & 37.5/4 \\ \cline{2-6}
        & Amortized DINO training (sec) & -- & -- & 80/4 & 1,220 /4 \\ \cline{2-6}
        & TMVI training$^*$ (sec) & 2,710 & 21,560 & 460 & 460 \\ \cline{2-6}
        & Parallel TMVI training (sec) & 135.5 & 1078 & -- & -- \\ \hline\hline

        \multirow{2}{*}{\rotatebox[origin=c]{90}{Total}}
        & Time per BIP (sec) &  2,710 & 21,560 & 525 & 1,440 \\ \cline{2-6}
        & Parallel time per BIP (sec) &   135.5  & 1078 & 482.45  & 798.75 \\ \hline\hline
        
        \multicolumn{2}{|r|}{Relative mean error achieved (\%)} & 80 & 20 & 10 & 5 \\ \hline 
        
    \end{tabular}
    }
    
    \caption{\textbf{Example I: \texttt{LazyDINO}/\texttt{LazyMap} timing comparison.}  We include two sample sizes (16k and 128k) for \texttt{LazyMap} for comparison with \texttt{LazyDINO} at 1k and 16k training data. The total sequential execution times are similar for \texttt{LazyMap} (16k) and \texttt{LazyDINO} (16k), and the total 20-way parallel execution times are similar for \texttt{LazyMap} (128k) and \texttt{LazyDINO} (16k). In both cases, the relative mean error achieved is much lower for \texttt{LazyDINO}.  Moreover, \texttt{LazyDINO} (1k) achieves smaller relative mean error than \texttt{LazyDINO} (128k) in less time.
    $^*$ denotes the fact that \texttt{LazyDINO} already performs batch-vectorized computation, so that parallel computation of the surrogate PtO map is not applicable.
        }
    \label{tab:alg_timing_1}
\end{table}

\begin{table}[htbp]
    \centering
    {\renewcommand{\arraystretch}{1.5}
    \begin{tabular}{||c|r|c|c|c|c|| }\hline
        & \multirow{3}{*}{\textbf{Category (unit)}} & \multicolumn{4}{c||}{\textbf{Method}} \\ \cline{3-6}
        & & \makecell{\texttt{LazyMap} \\ (1k)} & \makecell{\texttt{LazyMap} \\ (16k)} & \makecell{\texttt{LazyDINO} \\ (1k)} & \makecell{\texttt{LazyDINO} \\ (16k)} \\ \hline

        \multirow{7}{*}{\rotatebox[origin=c]{90}{Algorithm steps}} 
        & Amortized PtO evaluations (sec) & -- & -- & 2,150/4 & 34,200/4 \\ \cline{2-6}
        & Parallel amortized PtO evaluations (sec) & -- & -- & 107.5/4 & 1,710/4 \\ \cline{2-6}
        & Amortized Jacobian (sec) & -- & -- & 220/4 & 3,440/4 \\ \cline{2-6}
        & Parallel amortized Jacobian (sec) & -- & -- & 11/4 & 172/4 \\ \cline{2-6}
        & Amortized DINO training (sec) & -- & -- & 120/4 & 1,840/4 \\ \cline{2-6}
        & TMVI training$^*$ (sec) & 2,390 & 38,500 & 750 & 750 \\ \cline{2-6}
        & Parallel TMVI training (sec) & 119.5 & 1,925 & -- & -- \\ \hline\hline

        \multirow{2}{*}{\rotatebox[origin=c]{90}{Total}} 
        & Time per BIP (sec) & 2,390 & 38,500 & 1,372.5 & 10,620 \\ \cline{2-6}
        & Parallel time per BIP (sec) & 119.5 & 1,925 & 809.625 & 1,680.5 \\ \hline\hline
        \multicolumn{2}{|r|}{Relative mean error achieved (\%)} & 230 & 90 & 12 & 3.5 \\ \hline  
    \end{tabular}
    }
    \caption{\textbf{Example II: \texttt{LazyDINO}/\texttt{LazyMap} timing comparison.} We include two sample sizes (1k and 16k) for comparison. The total sequential execution times for the same sample sizes are less for \texttt{LazyDINO} due to the amortization of the offline phase across four BIPs corresponding to four instances of observational data. The total 20-way parallel execution times are similar for \texttt{LazyMap}(16k) and \texttt{LazyDINO} (16k). For both sample sizes, the relative mean error achieved is much lower for \texttt{LazyDINO}. Moreover, \texttt{LazyDINO} (1k) achieves much smaller relative mean error than \texttt{LazyDINO} (16k) in less time.
    $^*$ denotes the fact that \texttt{LazyDINO} already performs batch-vectorized computation, so that parallel computation of the surrogate PtO map is not applicable.}
    \label{tab:alg_timing_2}
\end{table}

Overall, these results demonstrate that for similar end-to-end computational costs, \texttt{LazyDINO} still performs substantially better than \texttt{LazyMap}, even allowing 20-way parallelism. We additionally note that, while we considered only relative mean error in the tables, \texttt{LazyDINO}'s efficiency gains in other error measures are even higher, as visible in the figures in~\cref{posterior_error}. 

Each iteration of LazyMap is computed with a 200-sample MC gradient estimator, \texttt{LazyMap} (16k) refers to  80 stochastic iterations, and \texttt{LazyMap} (128k) refers to 640 stochastic iterations.  In contrast, since each iteration of \texttt{LazyDINO} LMVI is cheap, the training time reported is for the 16 million surrogate PtO and Jacobian actions resulting from the increasing sample size strategy provided in~\cref{transport_map_specs}. This demonstrates a key takeaway: for extremely query-intensive algorithms such as the expected risk minimization problem arising in transport map training, surrogates are necessary, and since the associated training costs with the high-fidelity model would be so expensive, one can invest significant offline computations and still save orders of magnitude in computational costs.

\paragraph{SBAI vs. \texttt{LazyDINO}}\label{SBAIvLazyDINOTiming}
In typical formulations of SBAI, sampling requires inverting the transport map.\footnote{While this need not necessarily be the case, the alternative requires inverting the transport map \emph{during training}, which is often considered too expensive.} For particular transport map parametrizations, such as the inverse autoregressive flows (IAFs) we used for numerical results, the scalability of this inversion can be preserved. In the case of IAFs, the cost is dominated by $d_r$ 1-dimensional root-finding problems that can be solved relatively quickly via the bisection method. However, in contrast, sampling with \texttt{LazyDINO} involves only explicit evaluations of neural networks and produces samples in significantly less time.

In ~\cref{table:sampling_SBAI}, we report average times to train and compute 1 million i.i.d. approximate posterior samples for four instances of observational data studied in Example I and II. We note that these results depend on transport map architecture and the quality of implementations; however, the overall point concerning the additional expense to invert maps in SBAI is broadly applicable. 
Due to these higher sampling costs, the transport map training time required by \texttt{LazyDINO} \emph{for each } instance of observational data is amortized across the \emph{sampling of the posterior}, such that for large enough sample sizes, it can be less time consuming to first train a \texttt{LazyDINO} and then to subsequently sample, rather than to sample using the SBAI method.

\begin{table}[htbp]

    \begin{center}
    {\renewcommand{\arraystretch}{1.75}
    \begin{tabular}{ ||c | c | c||}\hline
     Example I & \multicolumn{2}{c||}{\textbf{Method}}  \\\hline
    \textbf{Time (sec)} & SBAI (16k) & \makecell{\texttt{LazyDINO}} (16k)\\\hline
    1 million samples & 1130  &  60 \\\hline
   Non-amortized training  & --- & 460\\\hline
   Amortized training & 930 / 4 & 1220/4 \\\hline\hline
  Total: sampling per BIP  & 1362.5 & 825 \\\hline

    \end{tabular}
    }
    \end{center}
    \begin{center}
    {\renewcommand{\arraystretch}{1.75}
    \begin{tabular}{ ||c | c | c||}\hline
     Example II & \multicolumn{2}{c||}{\textbf{Method}} \\\hline
    \textbf{Time (sec)} & SBAI (16k) & \makecell{\texttt{LazyDINO}} (16k)\\\hline
     1 million samples &  1150  &  60 \\\hline
   Non-amortized training  & --- & 750\\\hline
   Amortized training & 960 / 4 & 1840/4\\\hline\hline
   Total: sampling per BIP  & 1390 & 1270 \\\hline

    \end{tabular}
    }
    \end{center}
    \caption{\textbf{SBAI vs \texttt{LazyDINO} sampling times.} For large enough sample sizes, the inversion-to-sample approach of SBAI is more costly than the optimize-to-sample appraoch of \texttt{LazyDINO}. Since previous comparisons demonstrated that \texttt{LazyDINO} is consistently much more accurate on a sample cost basis, this adds to the argument for utilizing \texttt{LazyDINO} over SBAI in the setting where amortized offline computations are desirable to facilitate real-time solutions.
    To facilitate direct comparison, we employ identical IAF architectures for SBAI and \texttt{LazyDINO}, except that the input dimension for the conditional transport map for SBAI is $\mathbb{R}^{d_r} \times \mathbb{R}^{d_{\boldsymbol{y}}}$ rather than $\mathbb{R}^{d_r}$ in $\texttt{LazyDINO}$. The architecture used is the same as the one described with all LMVI methods in~\cref{transport_map_specs}. Times are averaged across the four BIPs corresponding to four instances of observational data and rounded up to the nearest 10 seconds. We report the times for 16k training samples. However, the times are similar across all sample sizes. Since the transport map architecture was chosen to be identical, sampling times for the two examples were essentially the same. The additional time to sample with SBAI stems from the need to apply the inverse of the transport map when sampling. Amortized training refers to conditional transport map training for SBAI and \texttt{RB-DINO} training for \texttt{LazyDINO}.}
    \label{table:sampling_SBAI}
\end{table}
 
\subsection{Visualization of discrepancy in marginals}
\begin{figure}[htbp]
    \centering
    {\begin{tabular}{l l l l l l} \includegraphics[width=0.04\textwidth]{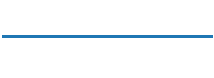}& \texttt{LazyDINO} & \includegraphics[width=0.04\textwidth]{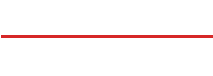} & LA-baseline&\includegraphics[width=0.04\textwidth]{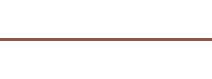} & True posterior via MCMC
    \end{tabular}}  
        \addtolength{\tabcolsep}{-6pt}
    \begin{tabular}{c c c}
    \hspace{0.04\textwidth}At $250$ training samples & \hspace{0.04\textwidth}At $2$k training samples & \hspace{0.04\textwidth}At $16$k training samples\\
    \includegraphics[width=0.33\textwidth]{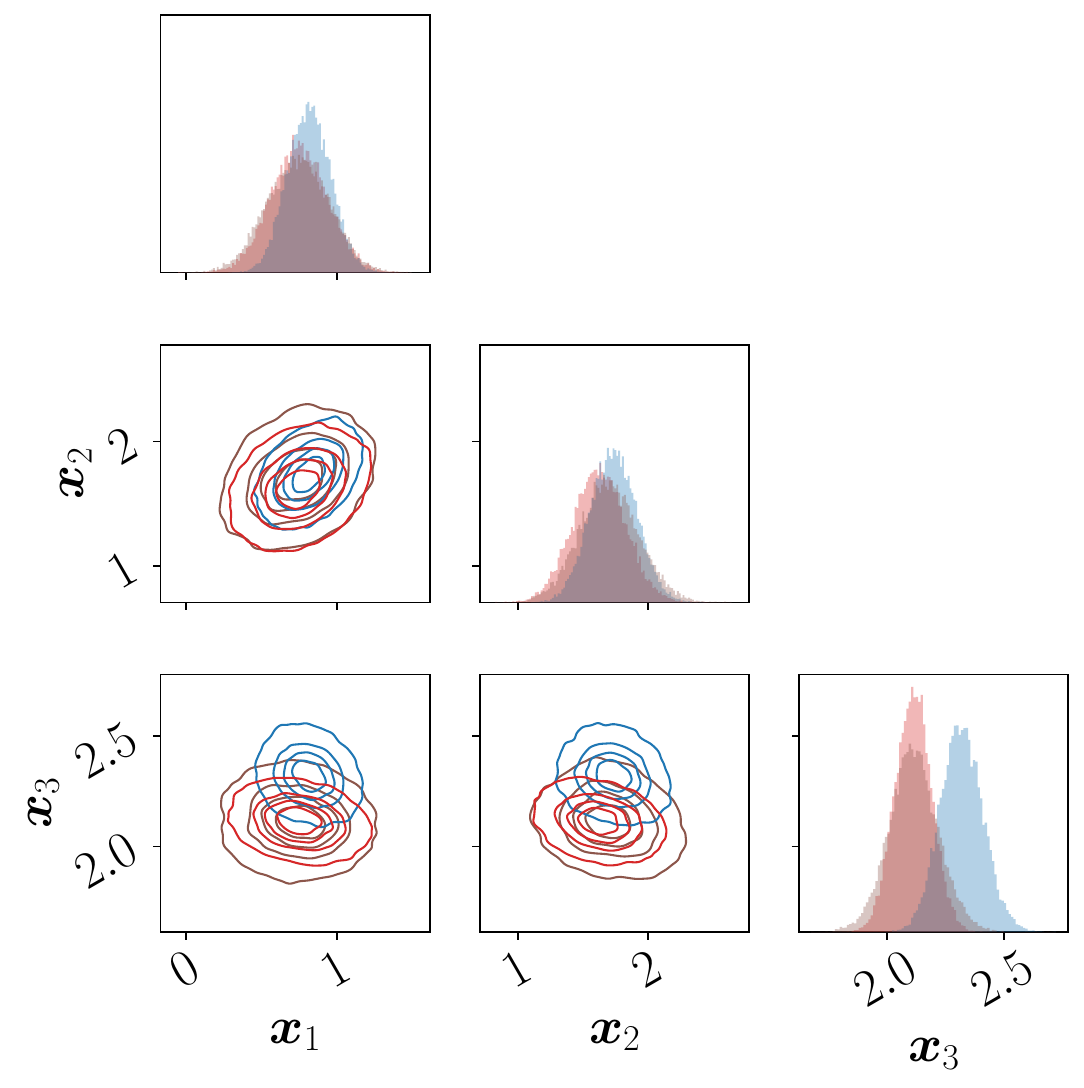}     & \includegraphics[width=0.33\textwidth]{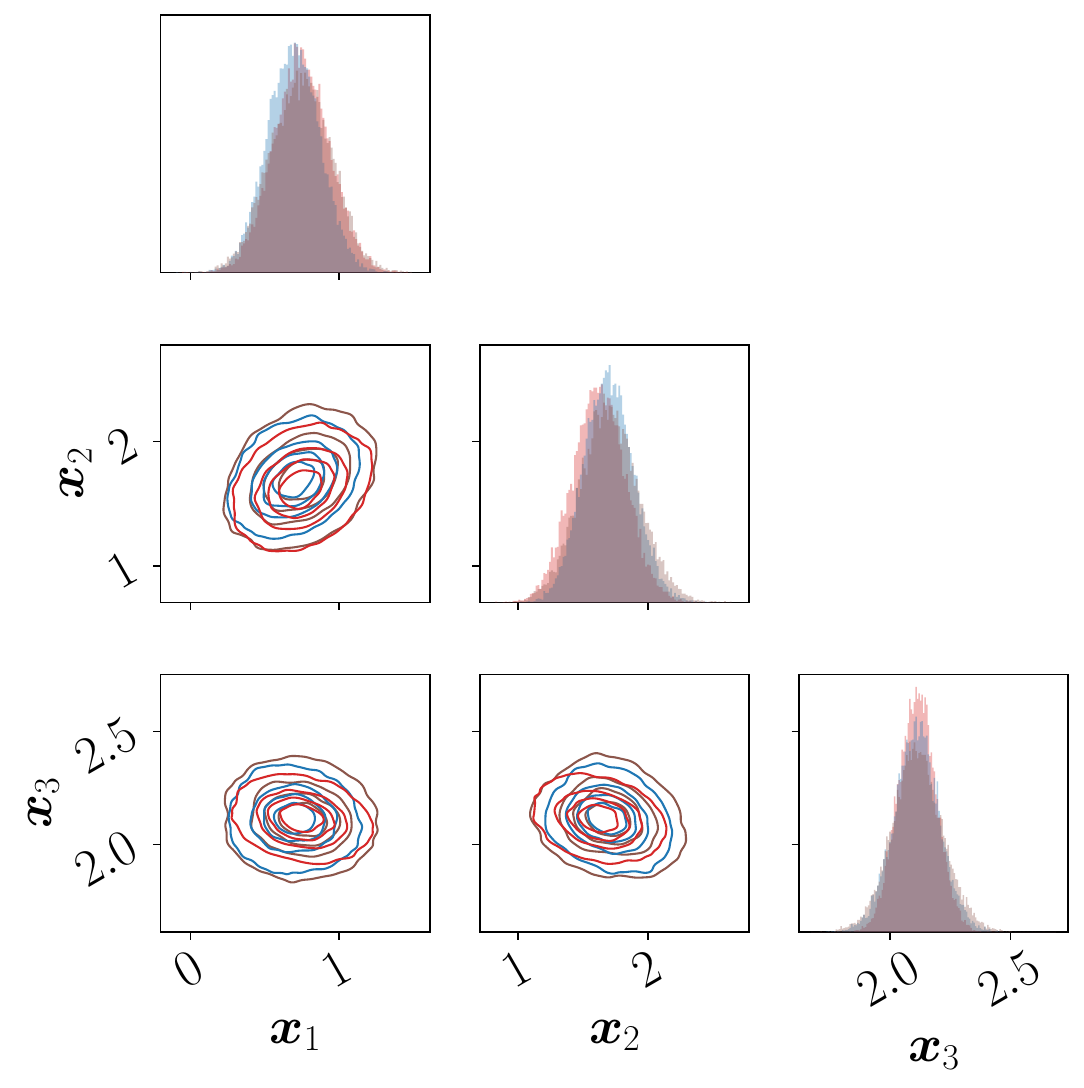} & \includegraphics[width=0.33\textwidth]{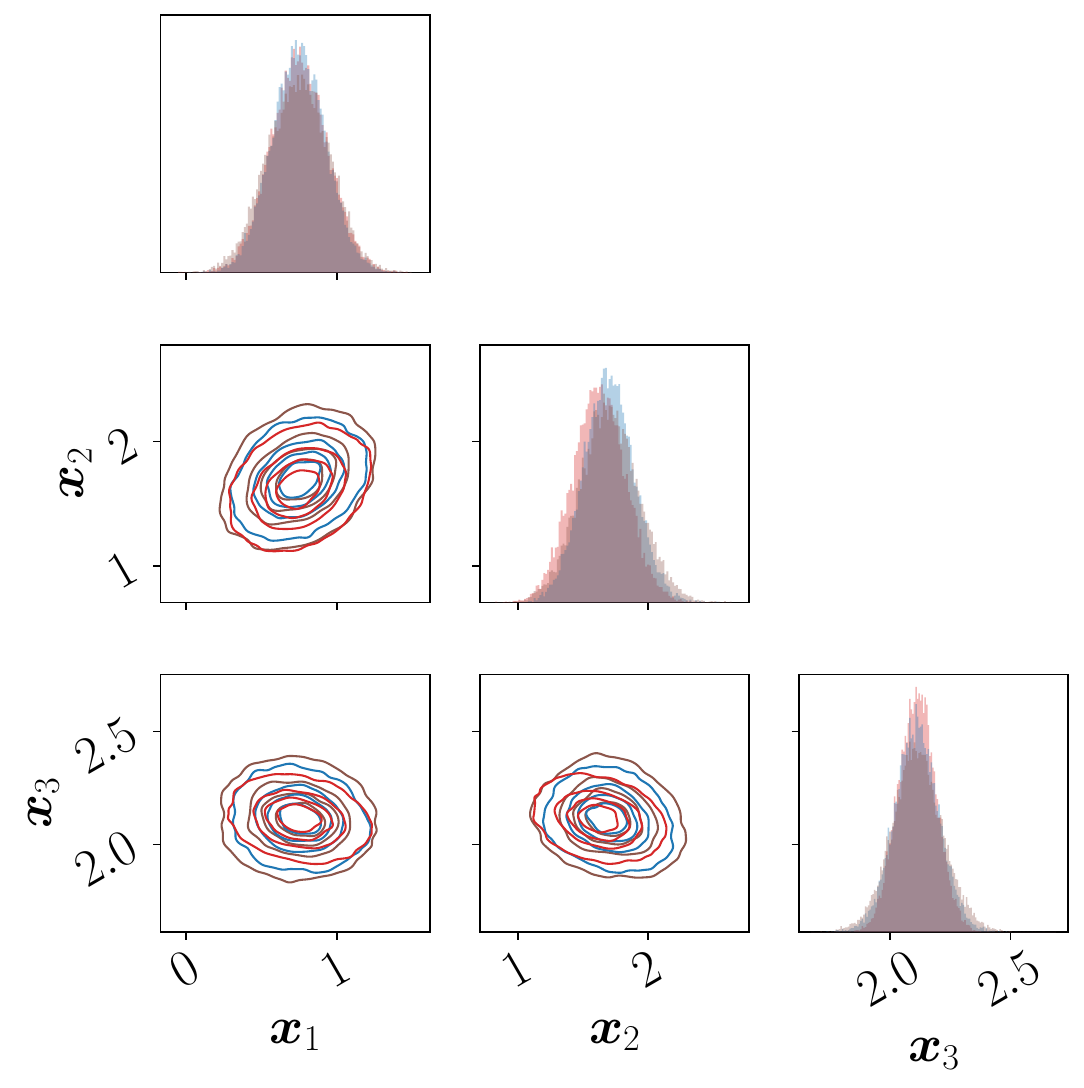}  \\
    \end{tabular}
    \addtolength{\tabcolsep}{6pt}
    \caption{\textbf{Example I: \texttt{LazyDINO} vs Laplace approximation marginals. } At 250 training samples, we see clear deviation in the marginals for both approaches, though LA-baseline has contours that match the posterior marginals more closely. By 2k training samples, $\texttt{LazyDINO}$ seems to match the marginals better, especially in the `tail' contours of the true posterior.}
    \label{fig:ndr_marg1}

    {\begin{tabular}{l l l l l l} \includegraphics[width=0.04\textwidth]{figures/legend_dino_line.pdf}& \texttt{LazyDINO} & \includegraphics[width=0.04\textwidth]{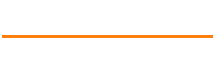} & \texttt{LazyNO} &\includegraphics[width=0.04\textwidth]{figures/legend_mcmc_line.pdf} & True posterior via MCMC  
    \end{tabular}}  
        \addtolength{\tabcolsep}{-6pt}
    \begin{tabular}{c c c}
    \hspace{0.04\textwidth}At $250$ training samples & \hspace{0.04\textwidth}At $2$k training samples & \hspace{0.04\textwidth}At $16$k training samples\\
    \includegraphics[width=0.33\textwidth]{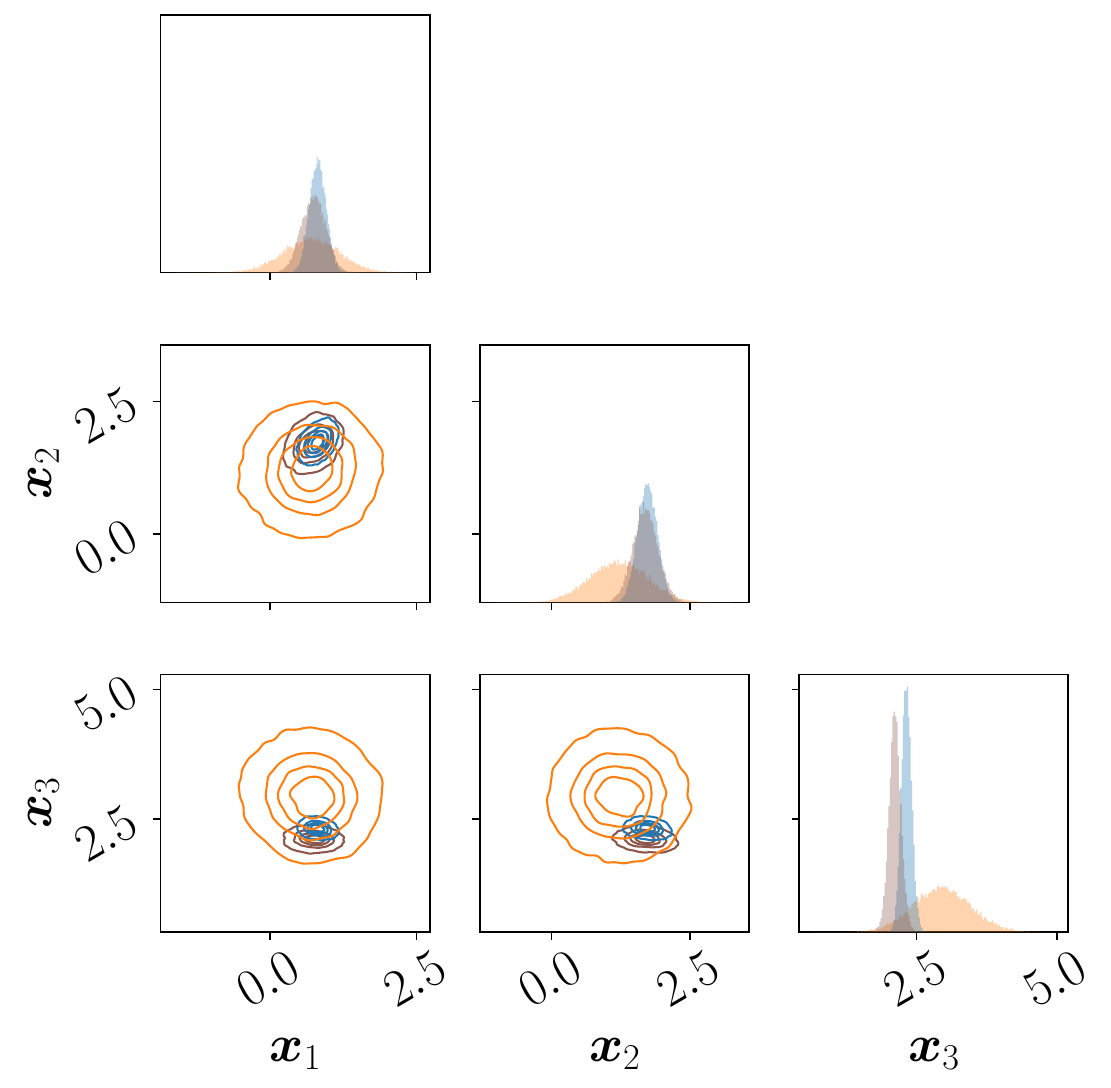}     & \includegraphics[width=0.33\textwidth]{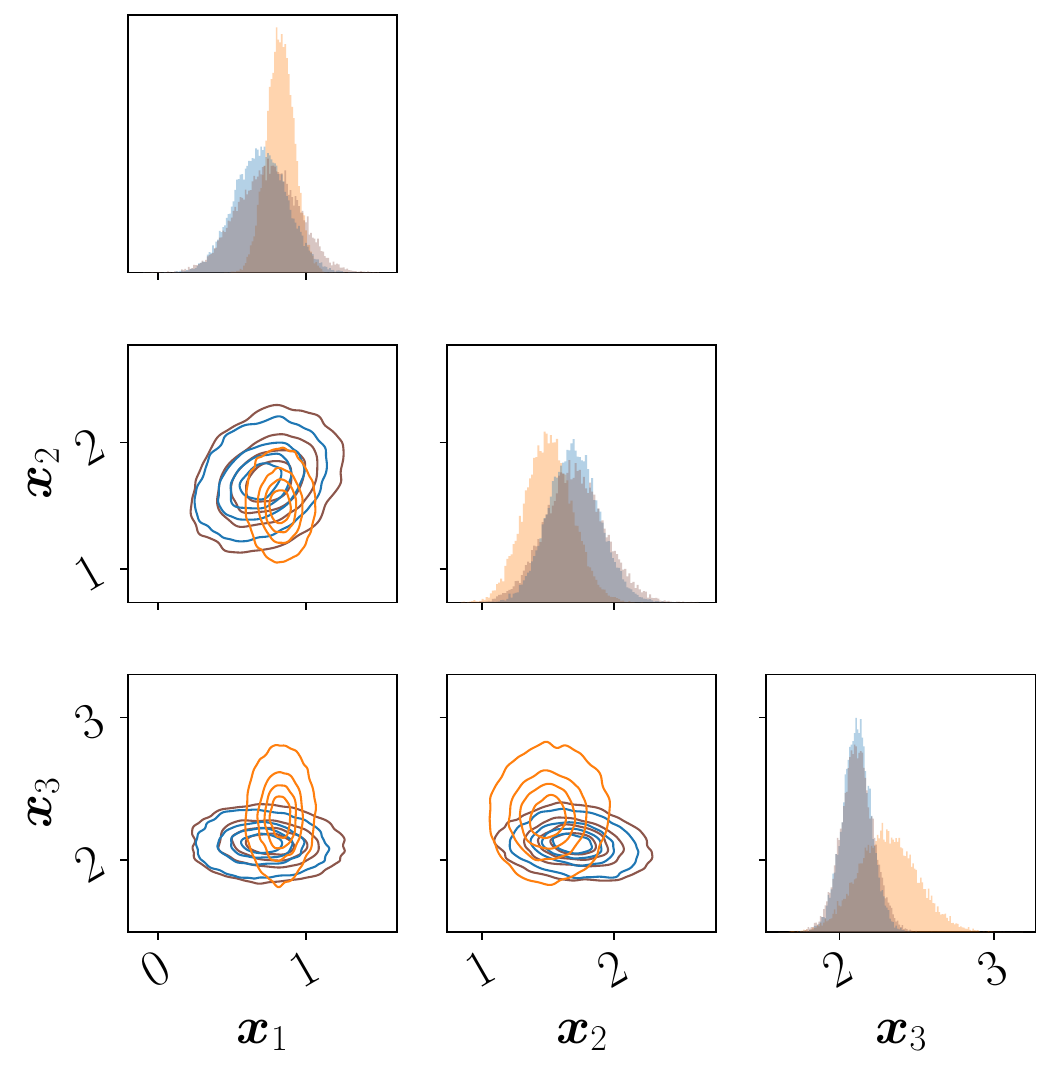} & \includegraphics[width=0.33\textwidth]{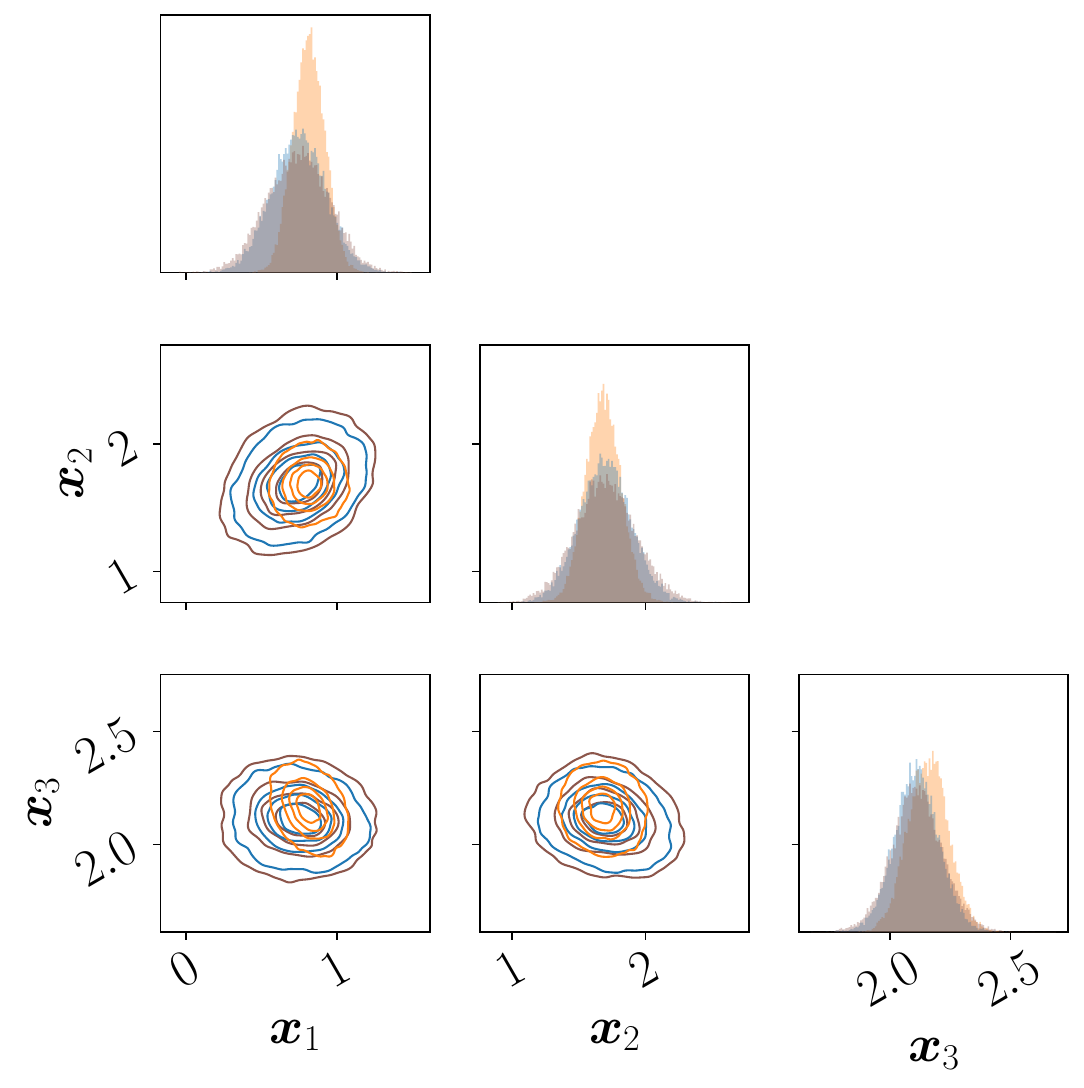}  \\
    \end{tabular}
    \addtolength{\tabcolsep}{6pt}
    \caption{\textbf{Example I: \texttt{LazyDINO} vs \texttt{LazyNO} marginals.} Consistent with~\cref{fig:gen_error_1} and~\cref{fig:NRD_Density_diagonistics}, we see the \texttt{LazyNO} fail to capture the posterior marginals well in all sample sizes shown. In comparison, \texttt{LazyDINO} closely matches the posterior marginals at 250 training samples (\emph{left})}
    \label{fig:ndr_marg2}
\end{figure}

\begin{figure}[htbp]
    \centering

    {\begin{tabular}{l l l l l l} \includegraphics[width=0.04\textwidth]{figures/legend_dino_line.pdf}& \texttt{LazyDINO} & \includegraphics[width=0.04\textwidth]{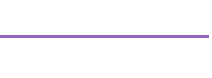} & SBAI &\includegraphics[width=0.04\textwidth]{figures/legend_mcmc_line.pdf} & True posterior via MCMC  
    \end{tabular}}  
        \addtolength{\tabcolsep}{-6pt}
    \begin{tabular}{c c c}
    \hspace{0.04\textwidth}At $250$ training samples & \hspace{0.04\textwidth}At $2$k training samples & \hspace{0.04\textwidth}At $16$k training samples\\
    \includegraphics[width=0.33\textwidth]{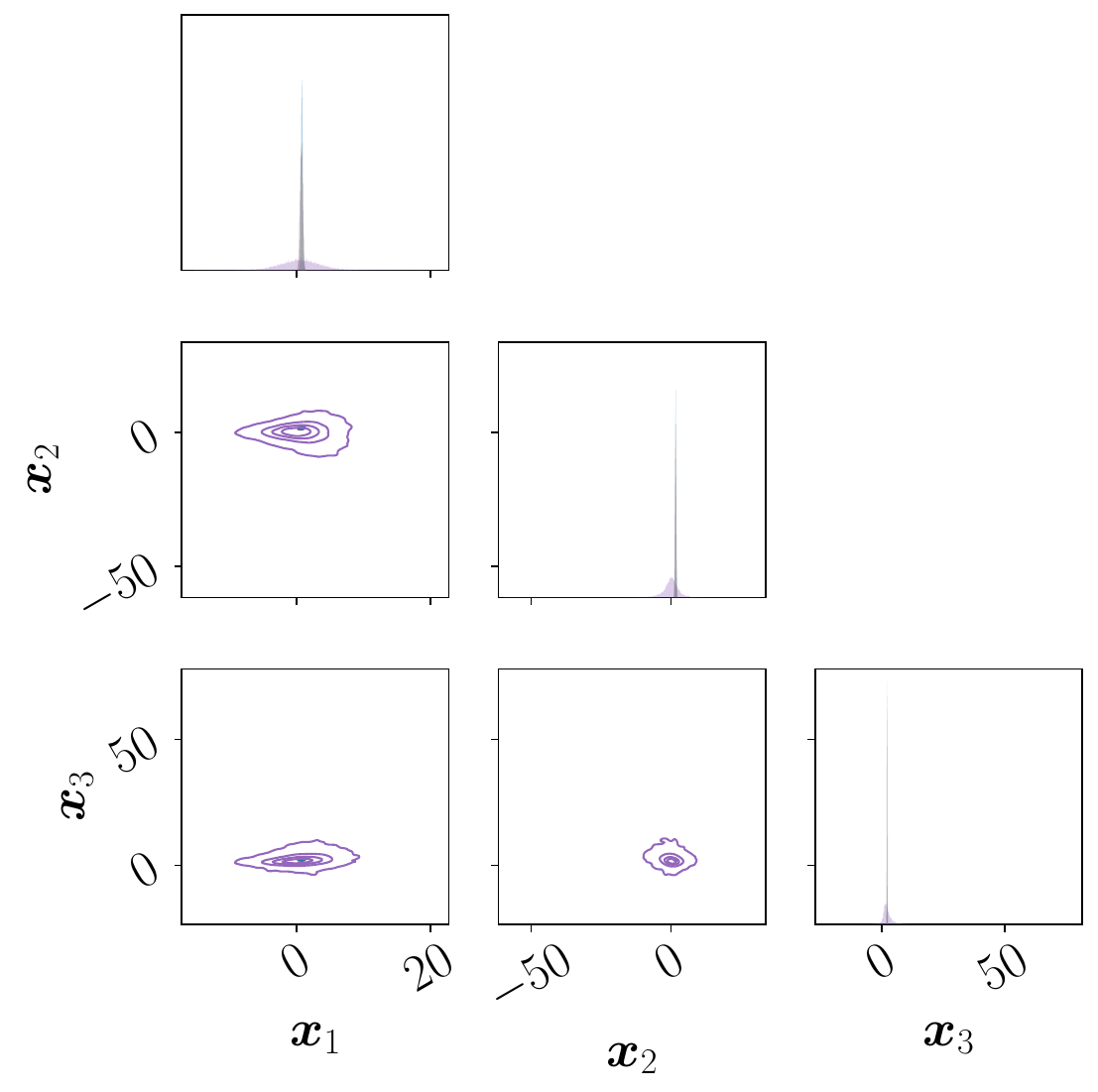}     & \includegraphics[width=0.33\textwidth]{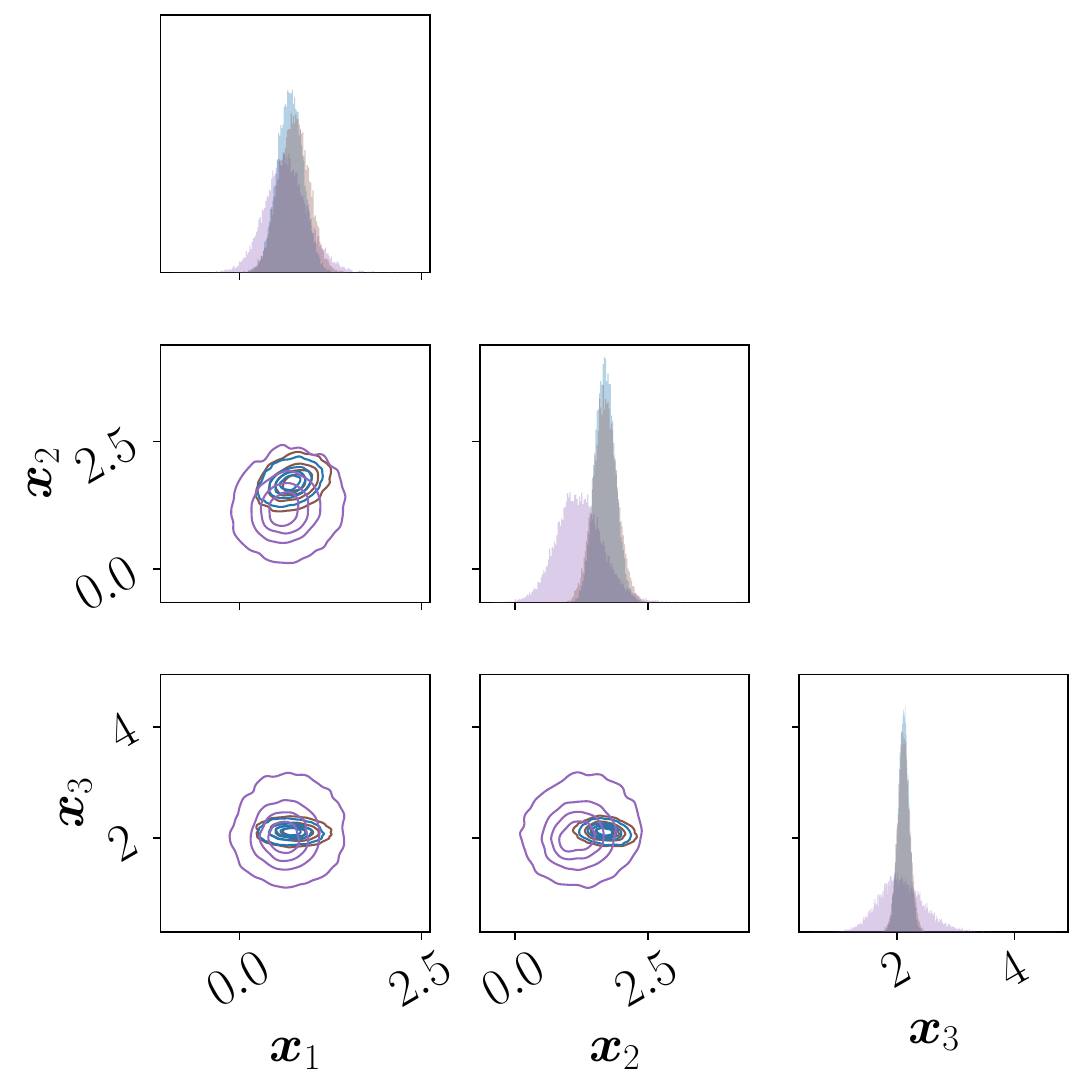} & \includegraphics[width=0.33\textwidth]{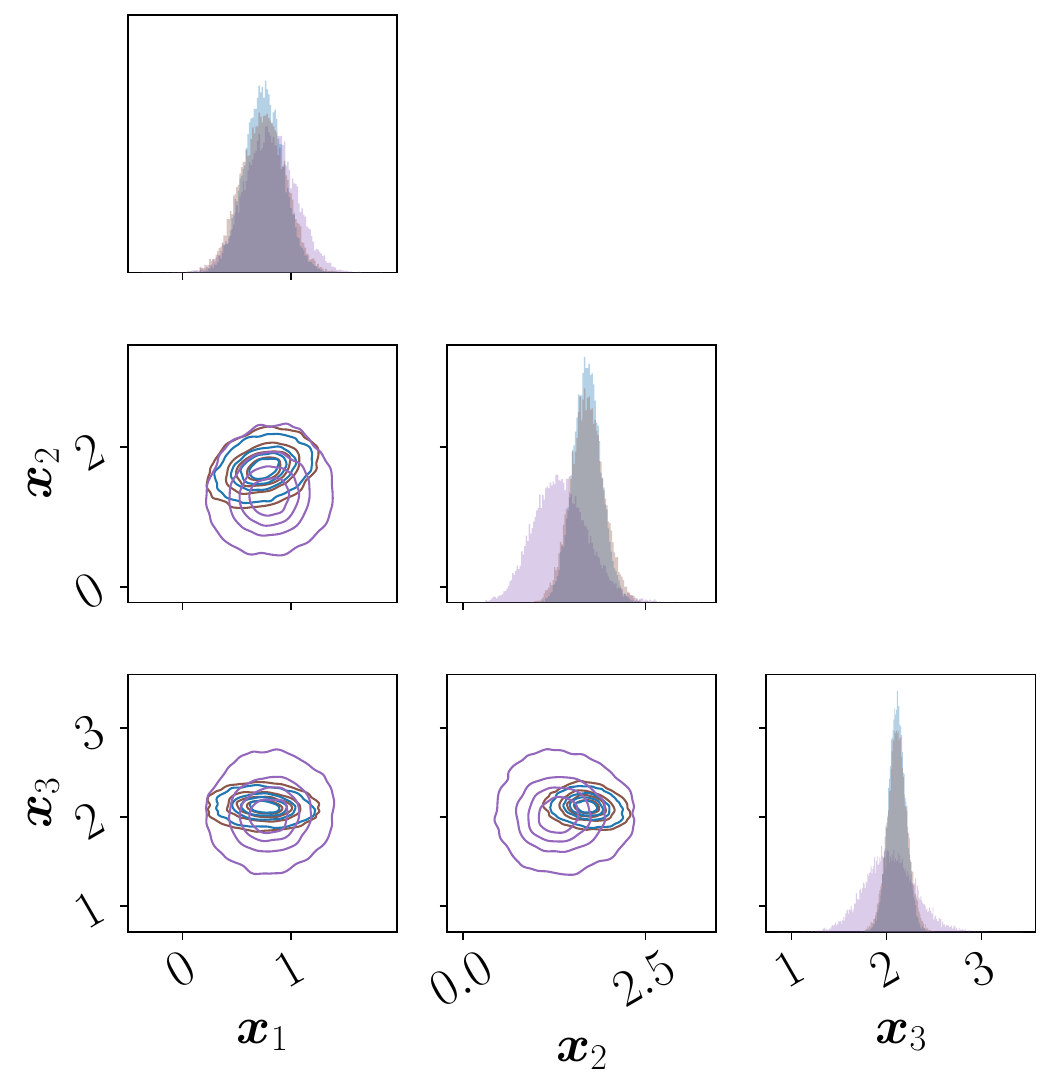}  \\
    \end{tabular}
    \addtolength{\tabcolsep}{6pt}
    \caption{\textbf{Example I: \texttt{LazyDINO} vs SBAI marginals.} The marginals produced via SBAI are much further from the true posterior marginals than the other approaches.
    Notably, SBAI consistently overestimates the uncertainty in posterior reconstruction and still yields a poor reconstruction of posterior marginals for $16,000$ samples. }
    \label{fig:ndr_marg3}
\end{figure}

\begin{figure}[htbp]
    \centering

    {\begin{tabular}{l l l l l l} \includegraphics[width=0.04\textwidth]{figures/legend_dino_line.pdf}& \texttt{LazyDINO} & \includegraphics[width=0.04\textwidth]{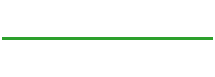} & \texttt{LazyMap} &\includegraphics[width=0.04\textwidth]{figures/legend_mcmc_line.pdf} & True posterior via MCMC  
    \end{tabular}}  
        \addtolength{\tabcolsep}{-6pt}
    \begin{tabular}{c c c}
    \hspace{0.04\textwidth}At $1$k training samples & \hspace{0.04\textwidth}At $4$k training samples & \hspace{0.04\textwidth}At $16$k training samples\\
    \includegraphics[width=0.33\textwidth]{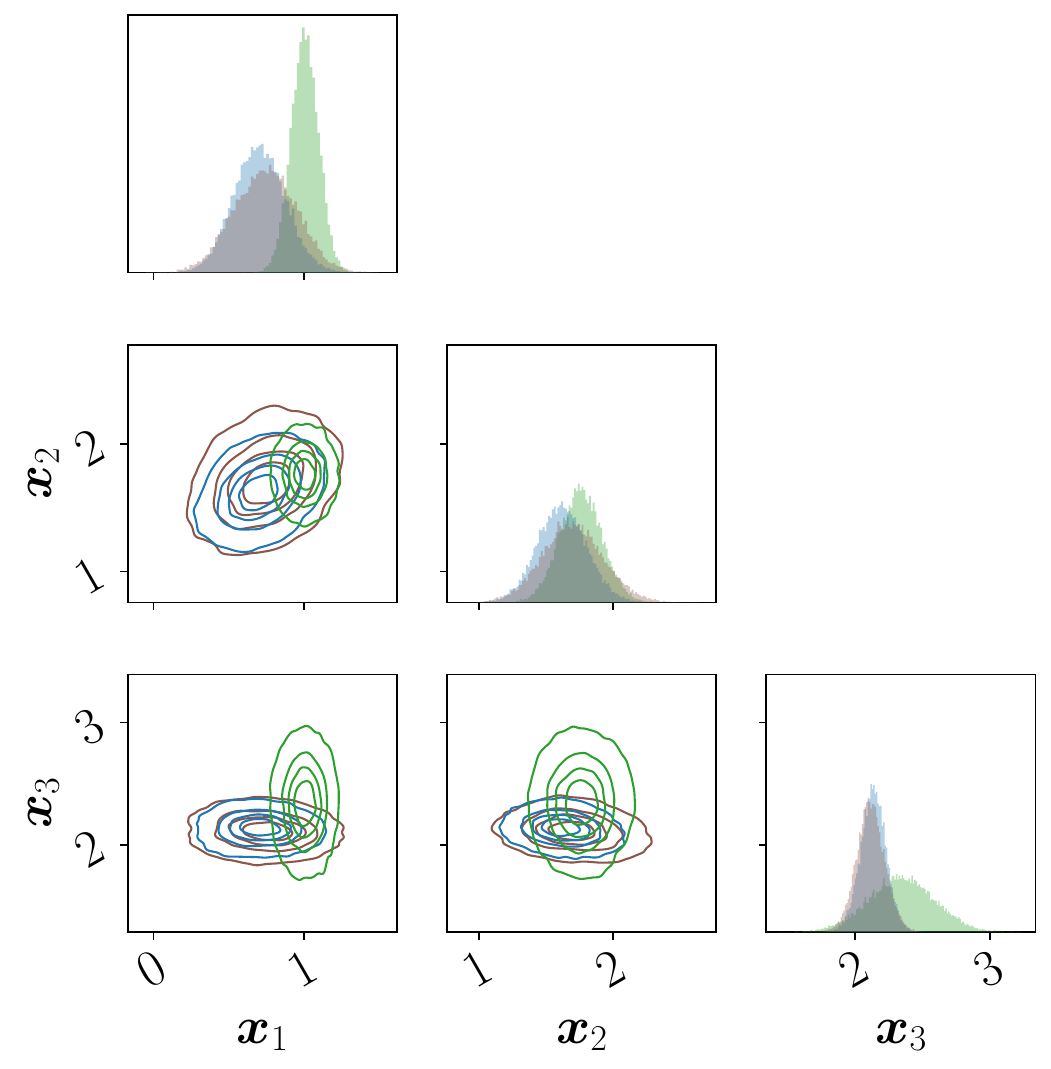}     & \includegraphics[width=0.33\textwidth]{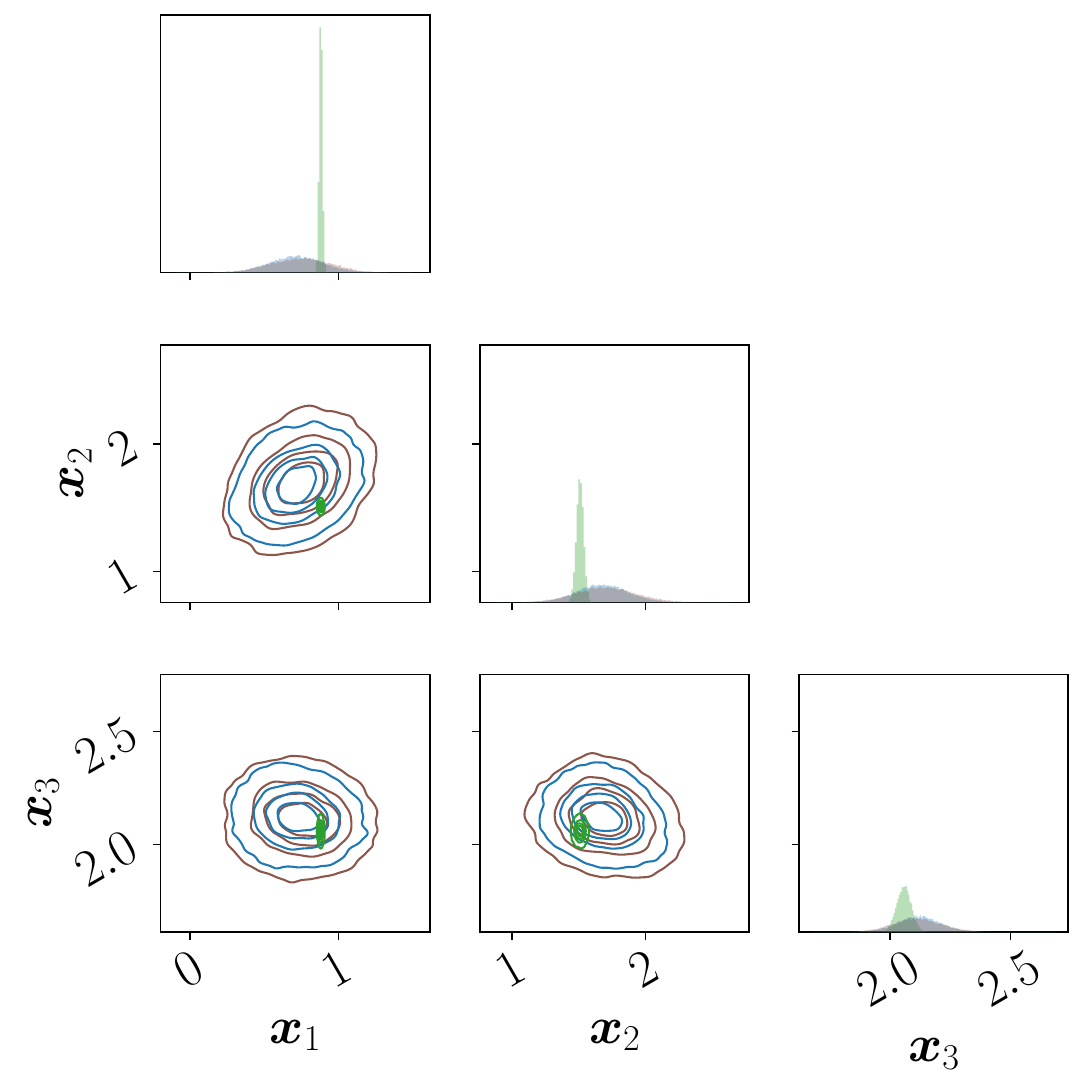} & \includegraphics[width=0.33\textwidth]{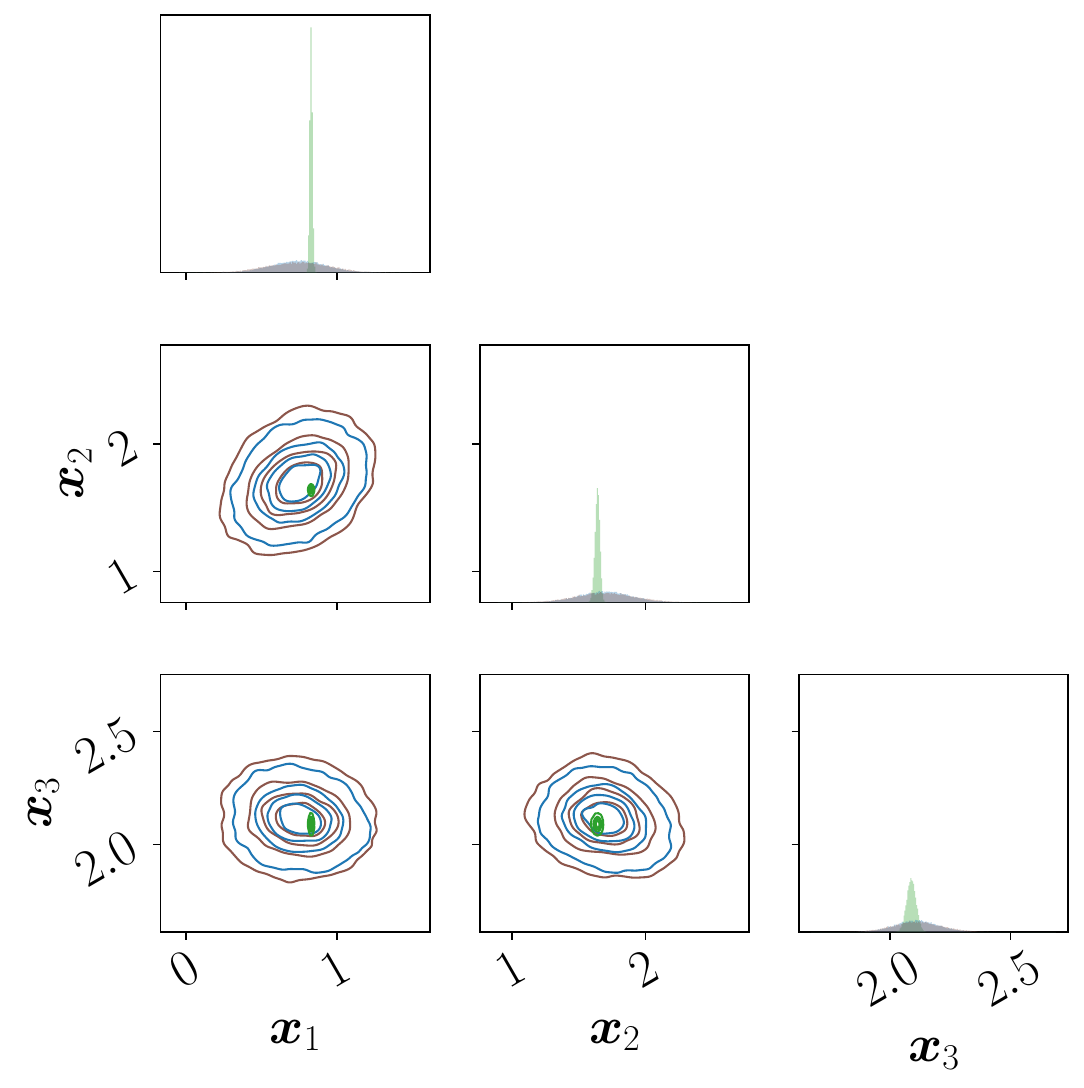}  \\
    \end{tabular}
    \addtolength{\tabcolsep}{6pt}
    \caption{\textbf{Example I: \texttt{LazyDINO} vs \texttt{LazyMap} marginals.} The \text{LazyMap} marginals are quite poor for 1k training samples and only get more concentrated as the number of training samples increases. The equivalent sample-cost \texttt{LazyDINO} posterior marginals match the ground truth substantially better. 
    Notably, \texttt{LazyMap} is highly concentrated, leading to na\"ive underestimation of the uncertainty in the posterior approximation. In comparison, \texttt{LazyDINO} yields a faithful approximation for only $250$ training data. }
    \label{fig:ndr_marg4}
\end{figure}

In this section, we provide visual evidence of the relative performances of different methods by investigating pairwise 2D marginal kernel density estimates and 1D marginal histograms of samples produced via the posterior sampling. We use a geometric MCMC method \cite[$\infty$-mMALA in section 3.2]{beskos2017geometric} to produce samples from the ground truth posterior. Marginal discrepancies provide a diagnostic measure for the quality of each posterior approximation---matching marginals is only a necessary, not sufficient, condition for matching the joint distribution. We plot a progression of marginals increasing sample size from left to right. We investigate the different posterior marginal comparisons for Example I in \cref{fig:ndr_marg1,fig:ndr_marg2,fig:ndr_marg3,fig:ndr_marg4}. Consistent with the previous results, \texttt{LazyDINO} produces consistently better approximations than the other TMVI methods, and overtakes the LA-baseline for a relatively small amount of training samples.

\begin{figure}[htbp]
    \centering
    {\begin{tabular}{l l l l l l} \includegraphics[width=0.04\textwidth]{figures/legend_dino_line.pdf}& \texttt{LazyDINO} (Ours) & \includegraphics[width=0.04\textwidth]{figures/legend_la_line.pdf} & LA-baseline&\includegraphics[width=0.04\textwidth]{figures/legend_mcmc_line.pdf} & True posterior via MCMC  
    \end{tabular}}  
    \addtolength{\tabcolsep}{-6pt}
    \begin{tabular}{c c c}
    \hspace{0.06\textwidth}At $250$ training samples & \hspace{0.06\textwidth}At $2$k training samples & \hspace{0.06\textwidth}At $16$k training samples\\
    \includegraphics[width=0.33\textwidth]{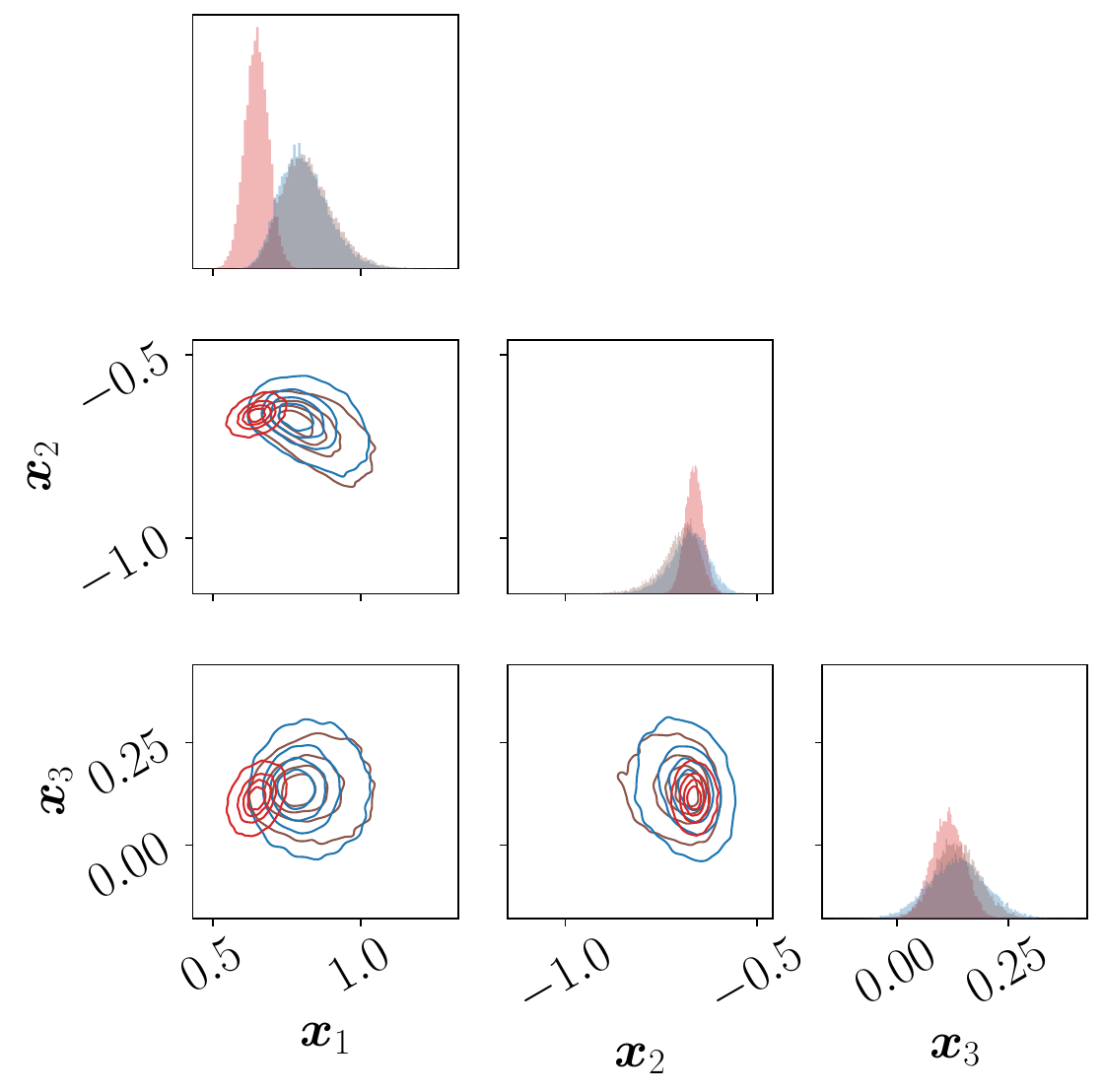}     & \includegraphics[width=0.33\textwidth]{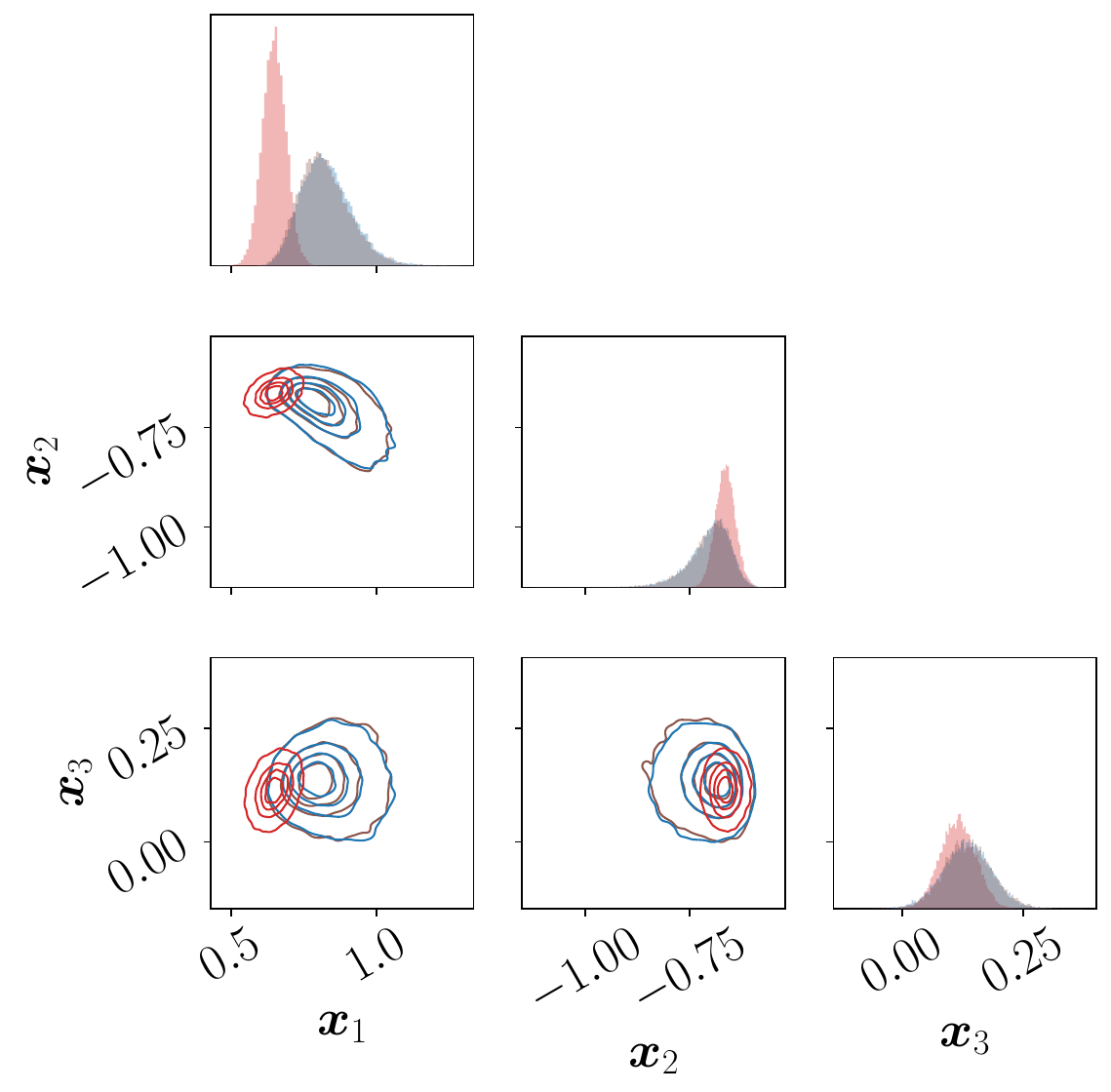} & \includegraphics[width=0.33\textwidth]{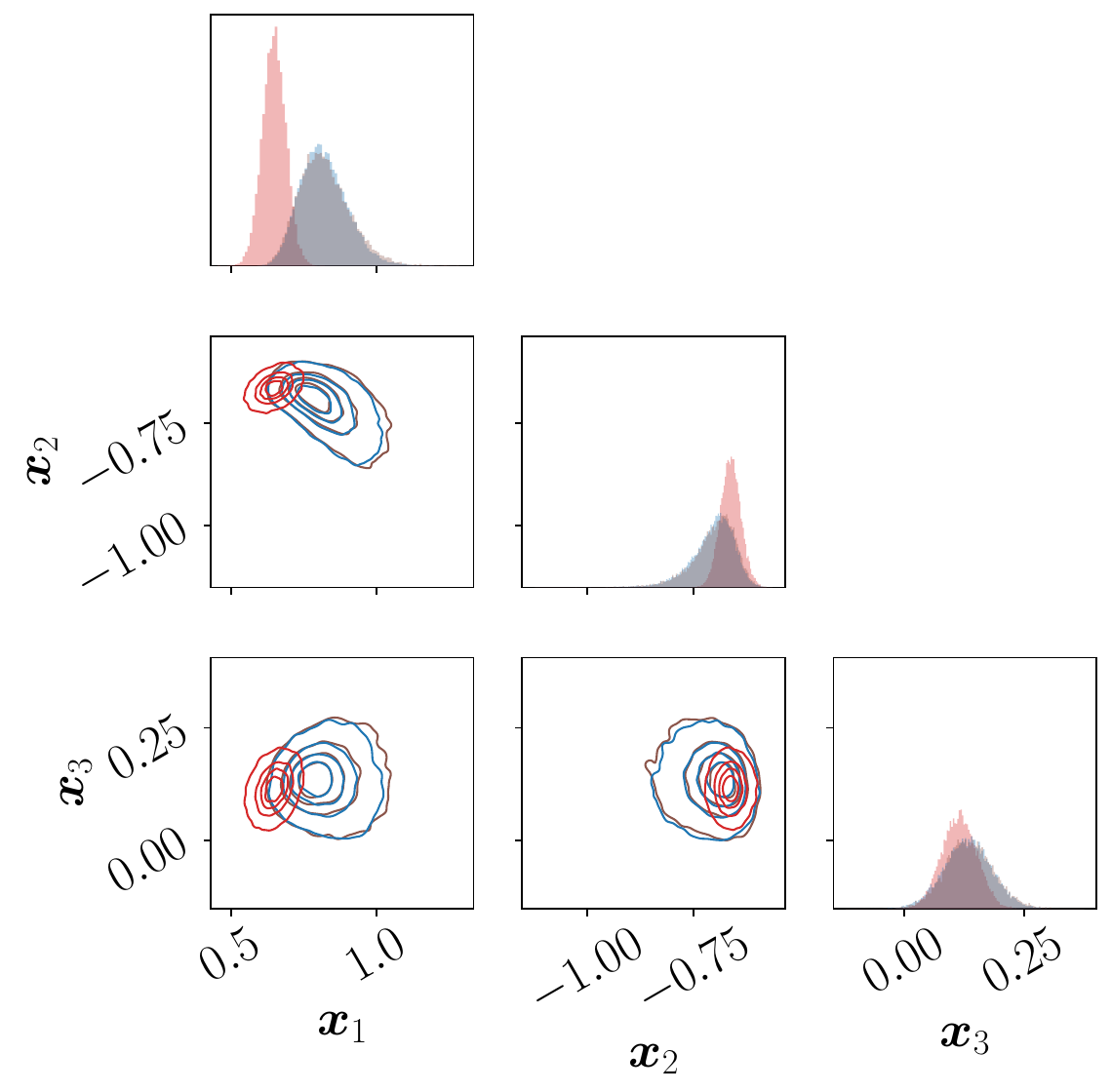}  \\
    \end{tabular}
    \addtolength{\tabcolsep}{6pt}
    \caption{
    {\textbf{Example II: \texttt{LazyDINO} vs Laplace approximation marginals.} Due to the non-Gaussianity of the posterior, Laplace approximation struggles to capture the overall behavior of the marginals depicted, whereas the \texttt{LazyDINO} marginals are close to the posterior marginals at 2k training samples. We note that a marginal of the MAP estimate is not the same as the MAP of a marginal in general, which may explain the apparent inconsistency that this marginal of the Laplace approximation is not visually centered on the region of highest probability of these particular posterior marginals.}}
    \label{fig:hyper_marg1}
\end{figure}

\begin{figure}[htbp]
    \centering
    {\begin{tabular}{l l l l l l} \includegraphics[width=0.04\textwidth]{figures/legend_dino_line.pdf}& \texttt{LazyDINO} (Ours) & \includegraphics[width=0.04\textwidth]{figures/legend_no_line.pdf} & \texttt{LazyNO} &\includegraphics[width=0.04\textwidth]{figures/legend_mcmc_line.pdf} & True posterior via MCMC  
    \end{tabular}}  
    \addtolength{\tabcolsep}{-6pt}
    \begin{tabular}{c c c}
    \hspace{0.06\textwidth}At $250$ training samples & \hspace{0.06\textwidth}At $2$k training samples & \hspace{0.06\textwidth}At $16$k training samples\\
    \includegraphics[width=0.33\textwidth]{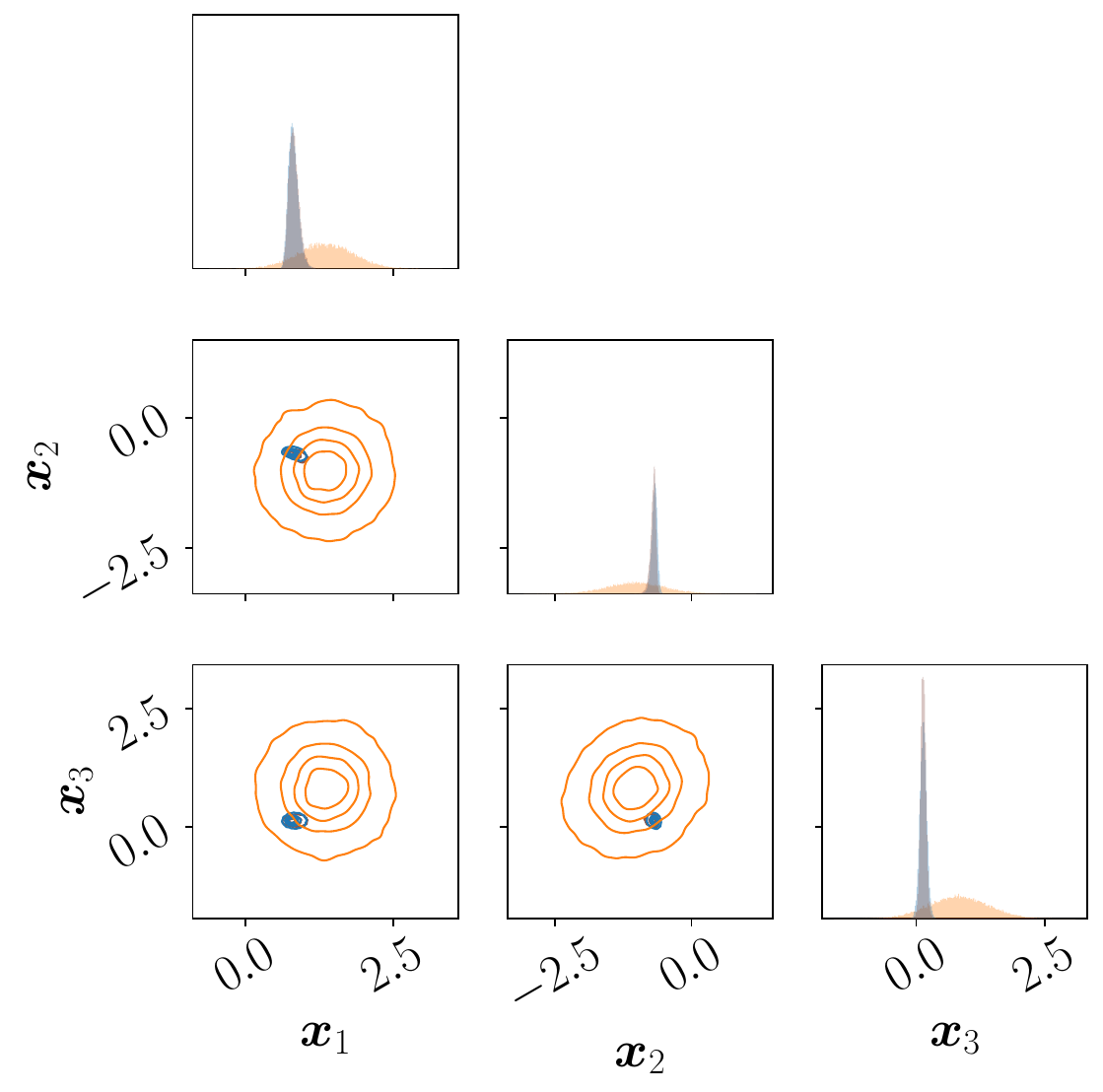}     & \includegraphics[width=0.33\textwidth]{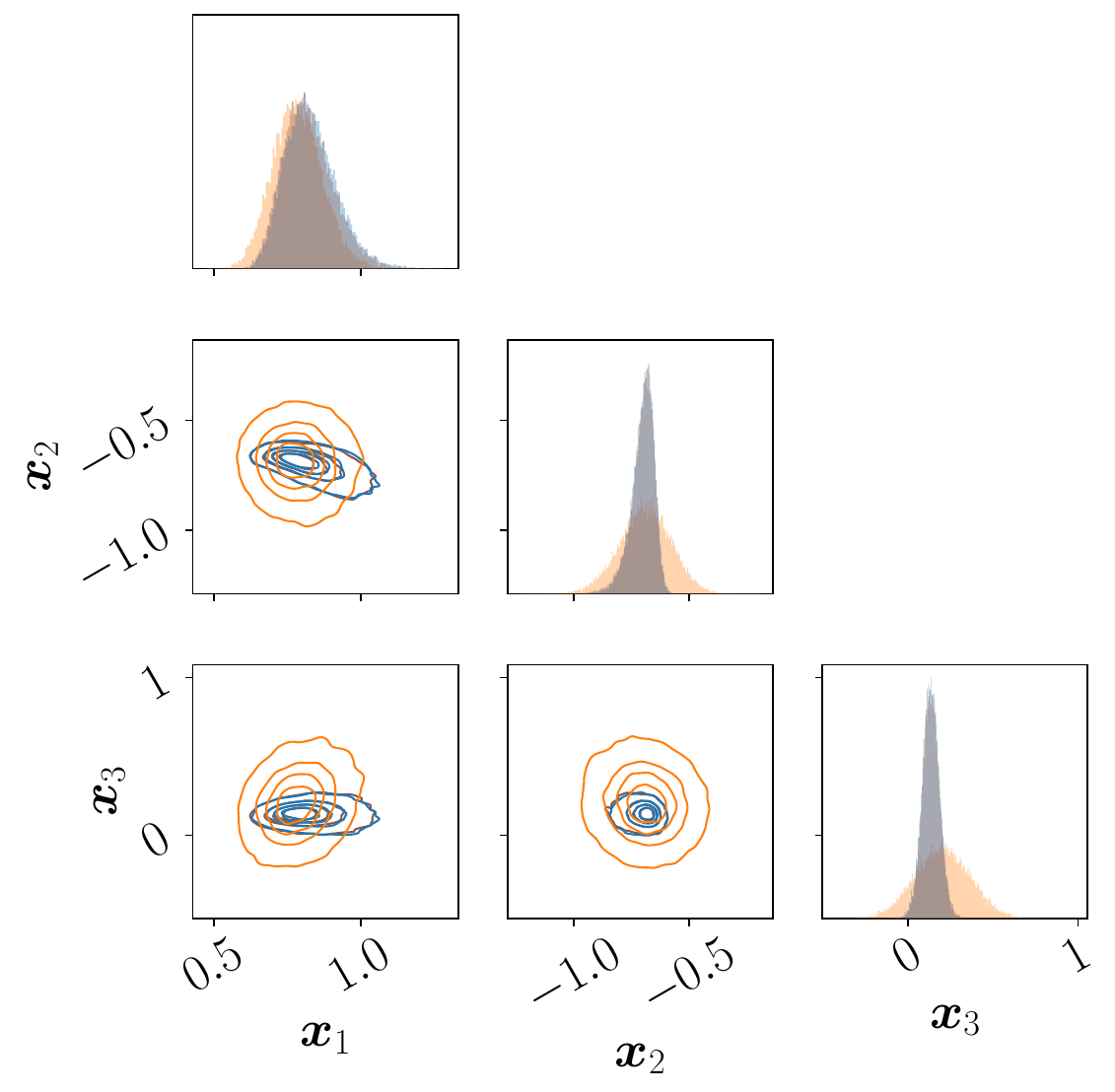} & \includegraphics[width=0.33\textwidth]{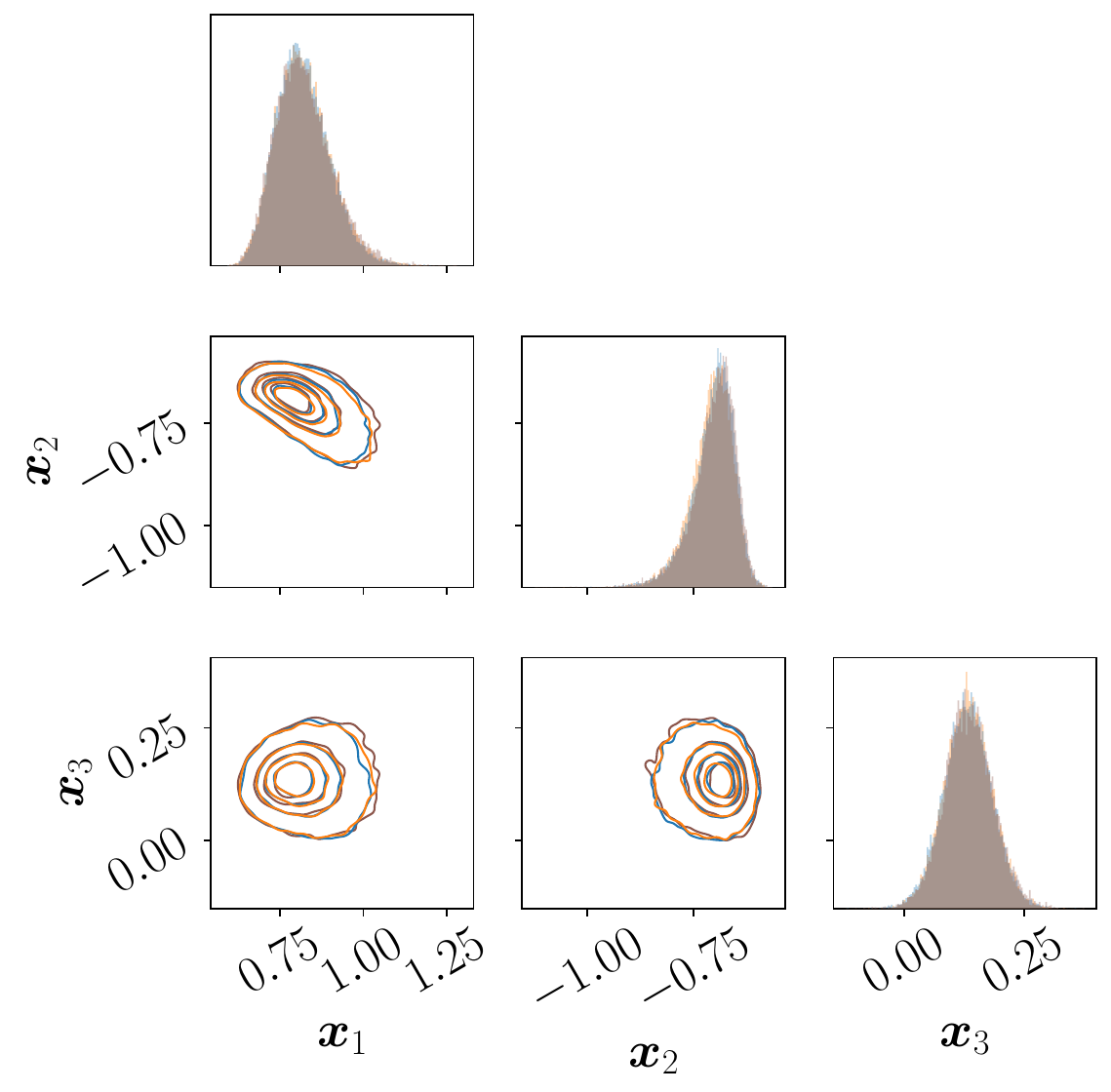}  \\
    \end{tabular}
    \addtolength{\tabcolsep}{6pt}
    \caption{\textbf{Example II: \texttt{LazyDINO} vs \texttt{LazyNO}  marginals.} With 250 training samples, \texttt{LazyNO} produces far underconcentrated samples in these marginals. With 16k training samples, both methods produce similar marginal distributions. The notable takeaway is that \texttt{LazyNO} requires multiple orders of magnitude more samples to have comparable performance to \texttt{LazyDINO} with $250$ samples. }
    \label{fig:hyper_marg2}
\end{figure}

\begin{figure}[htbp]
    \centering
    {\begin{tabular}{l l l l l l} \includegraphics[width=0.04\textwidth]{figures/legend_dino_line.pdf}& \texttt{LazyDINO} (Ours)& \includegraphics[width=0.04\textwidth]{figures/legend_cm_line.pdf} & SBAI &\includegraphics[width=0.04\textwidth]{figures/legend_mcmc_line.pdf} & True posterior via MCMC  
    \end{tabular}}  
    \addtolength{\tabcolsep}{-6pt}
    \begin{tabular}{c c c}
    \hspace{0.06\textwidth}At $250$ training samples & \hspace{0.06\textwidth}At $2$k training samples & \hspace{0.06\textwidth}At $16$k training samples\\
    \includegraphics[width=0.33\textwidth]{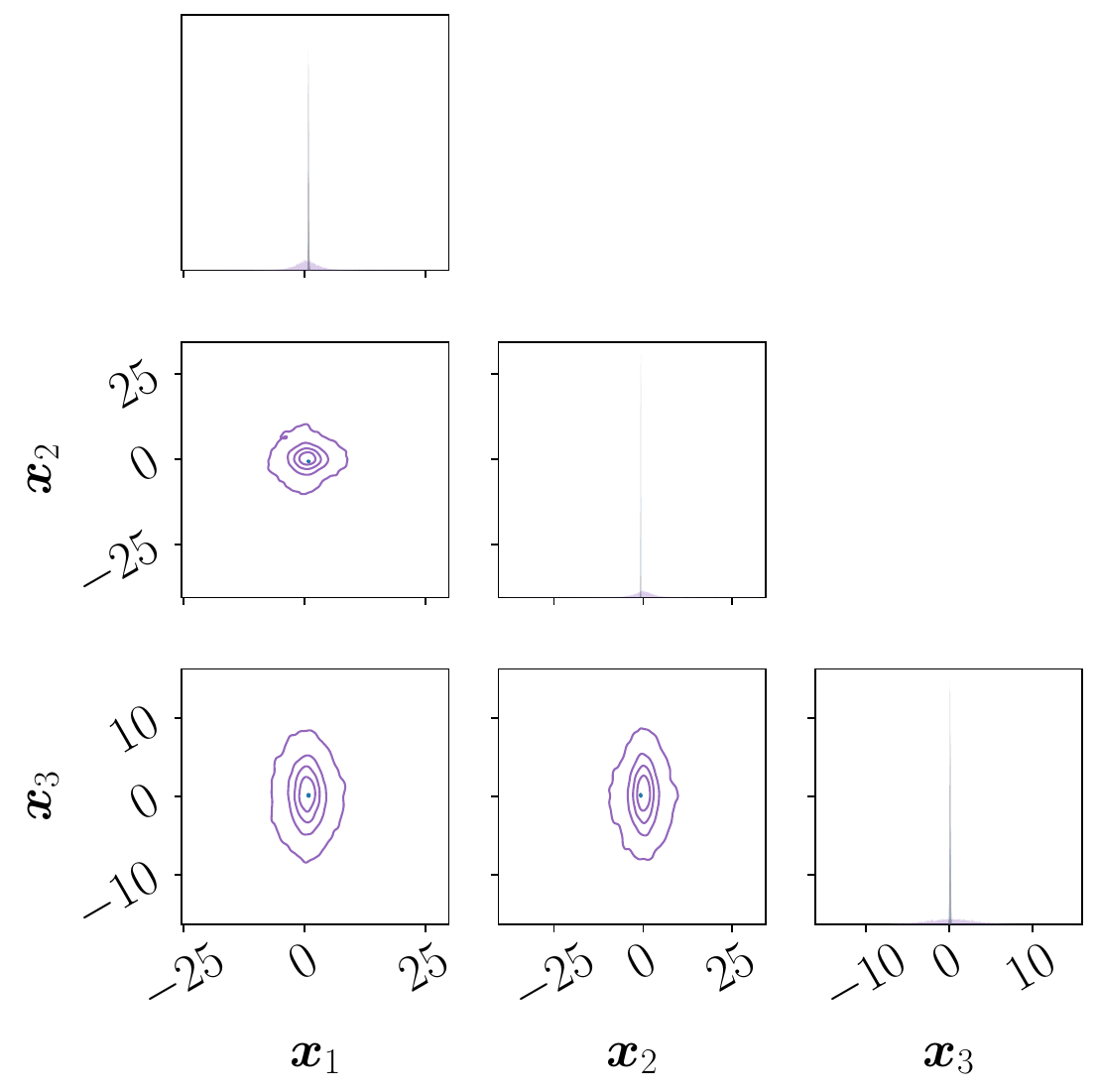}     & \includegraphics[width=0.33\textwidth]{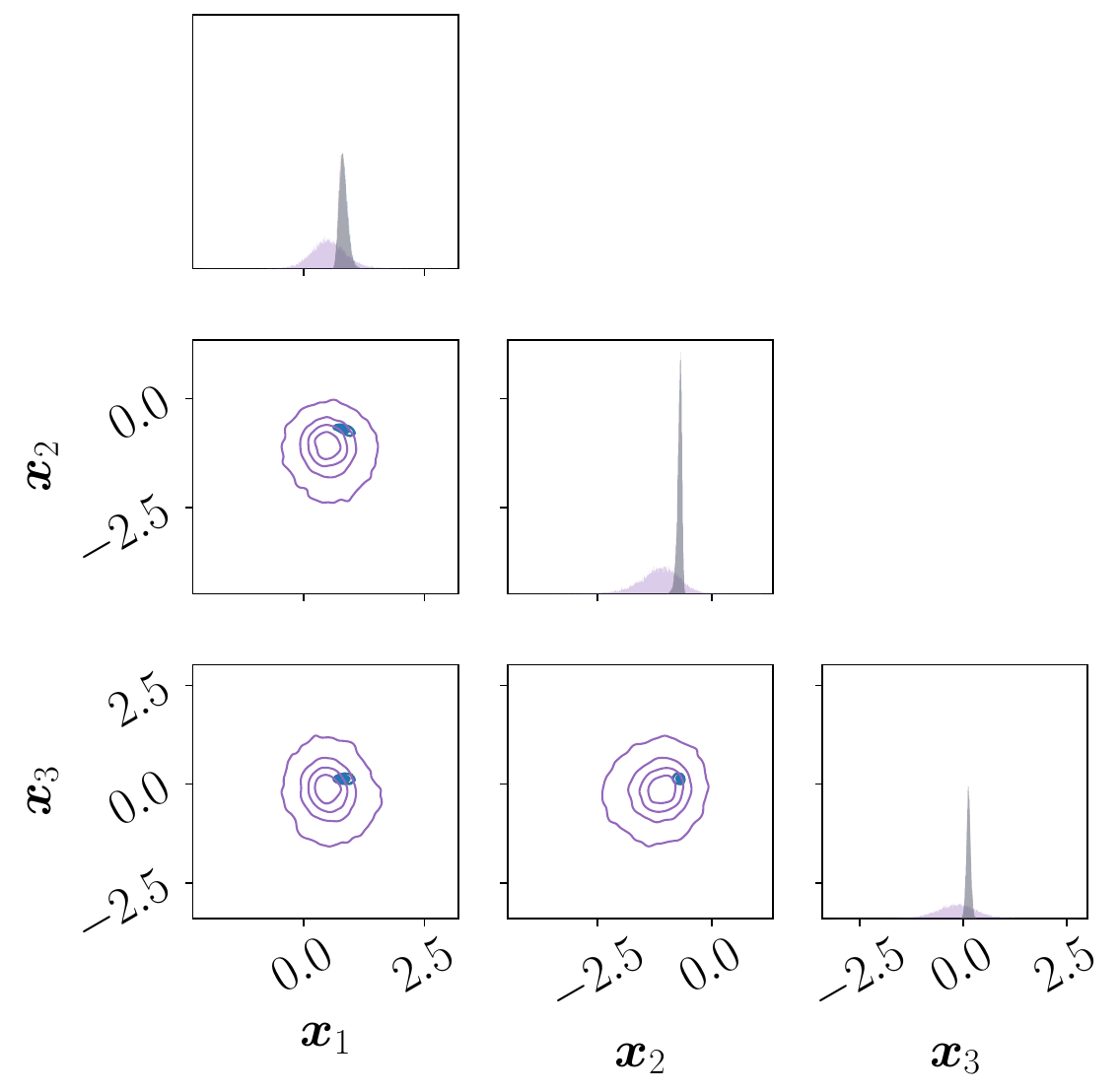} & \includegraphics[width=0.33\textwidth]{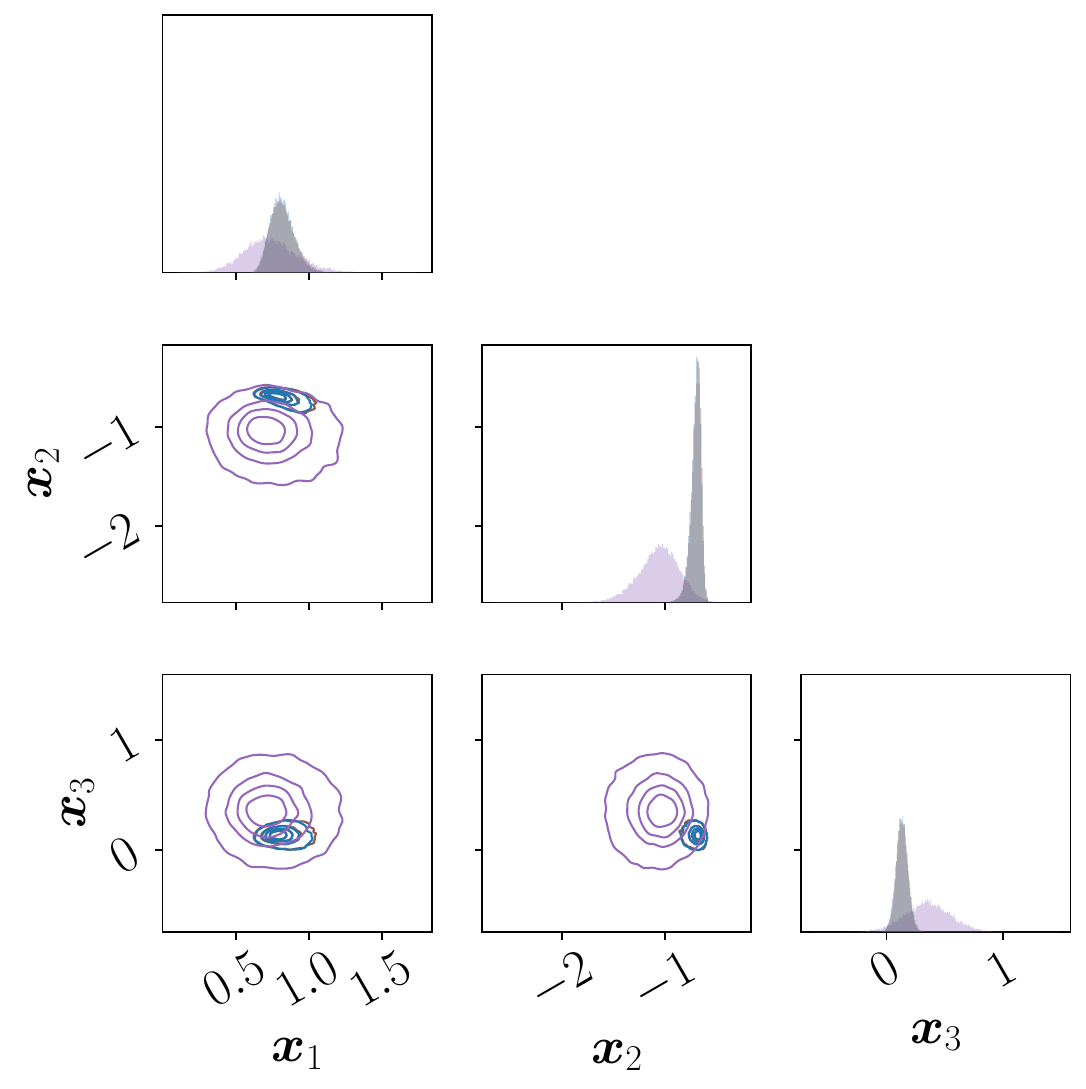}  \\
    \end{tabular}
    \addtolength{\tabcolsep}{6pt}
    \caption{\textbf{Example II: \texttt{LazyDINO} vs SBAI marginals.} SBAI produces samples that are highly under-concentrated. SBAI is simultaneously off in capturing the peak locations of the marginals while overestimating the uncertainty in the parameter.}
    \label{fig:hyper_marg3}
\end{figure}

\begin{figure}[htbp]
    \centering
    {\begin{tabular}{l l l l l l} \includegraphics[width=0.04\textwidth]{figures/legend_dino_line.pdf}& \texttt{LazyDINO} (Ours)& \includegraphics[width=0.04\textwidth]{figures/legend_pde_line.pdf} & \texttt{LazyMap} &\includegraphics[width=0.04\textwidth]{figures/legend_mcmc_line.pdf} & True posterior via MCMC  
    \end{tabular}}  
    \addtolength{\tabcolsep}{-6pt}
    \begin{tabular}{c c c}
    \hspace{0.06\textwidth}At $1k$ training samples & \hspace{0.06\textwidth}At $4$k training samples & \hspace{0.06\textwidth}At $16$k training samples\\
    \includegraphics[width=0.33\textwidth]{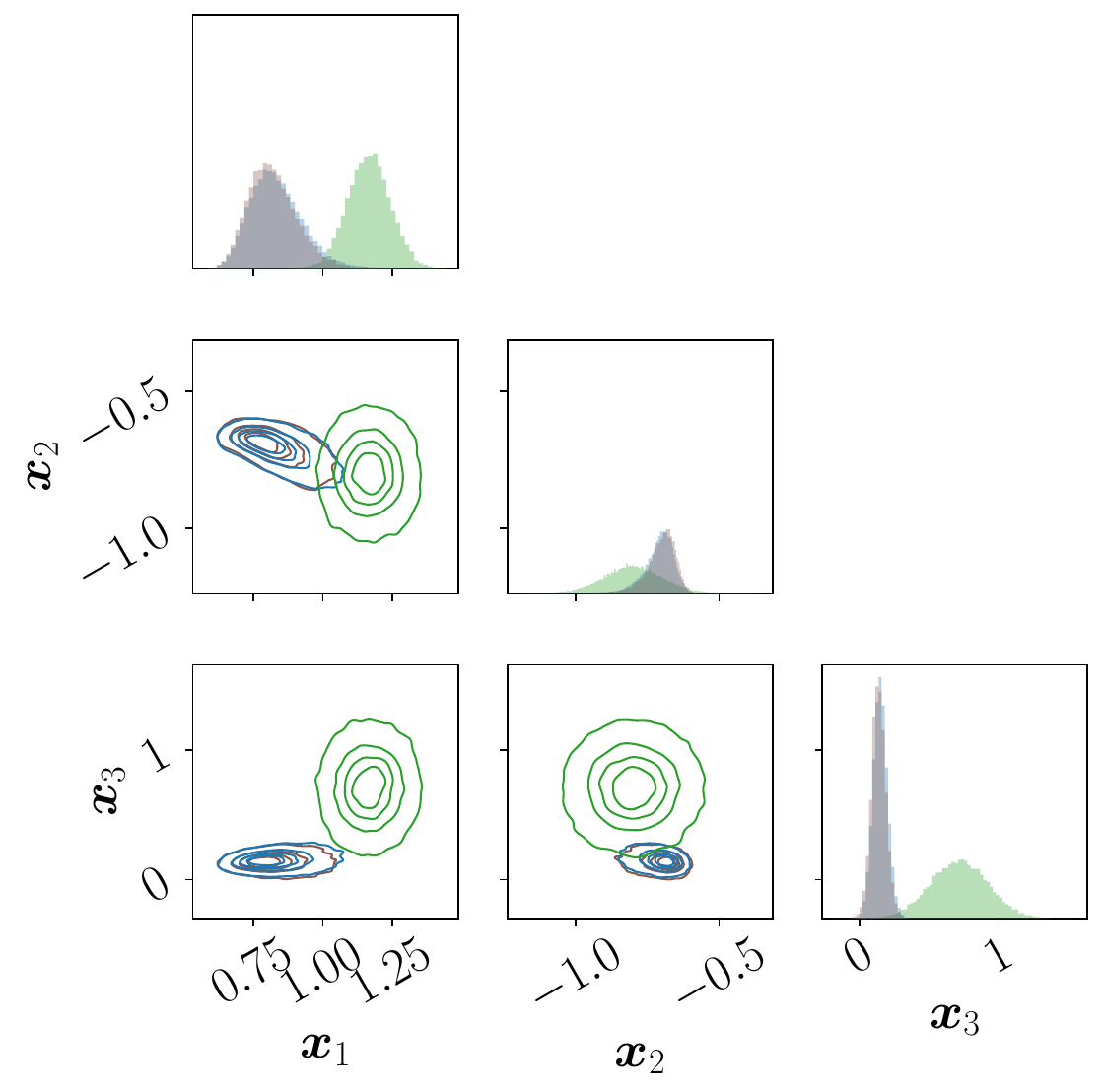}     & \includegraphics[width=0.33\textwidth]{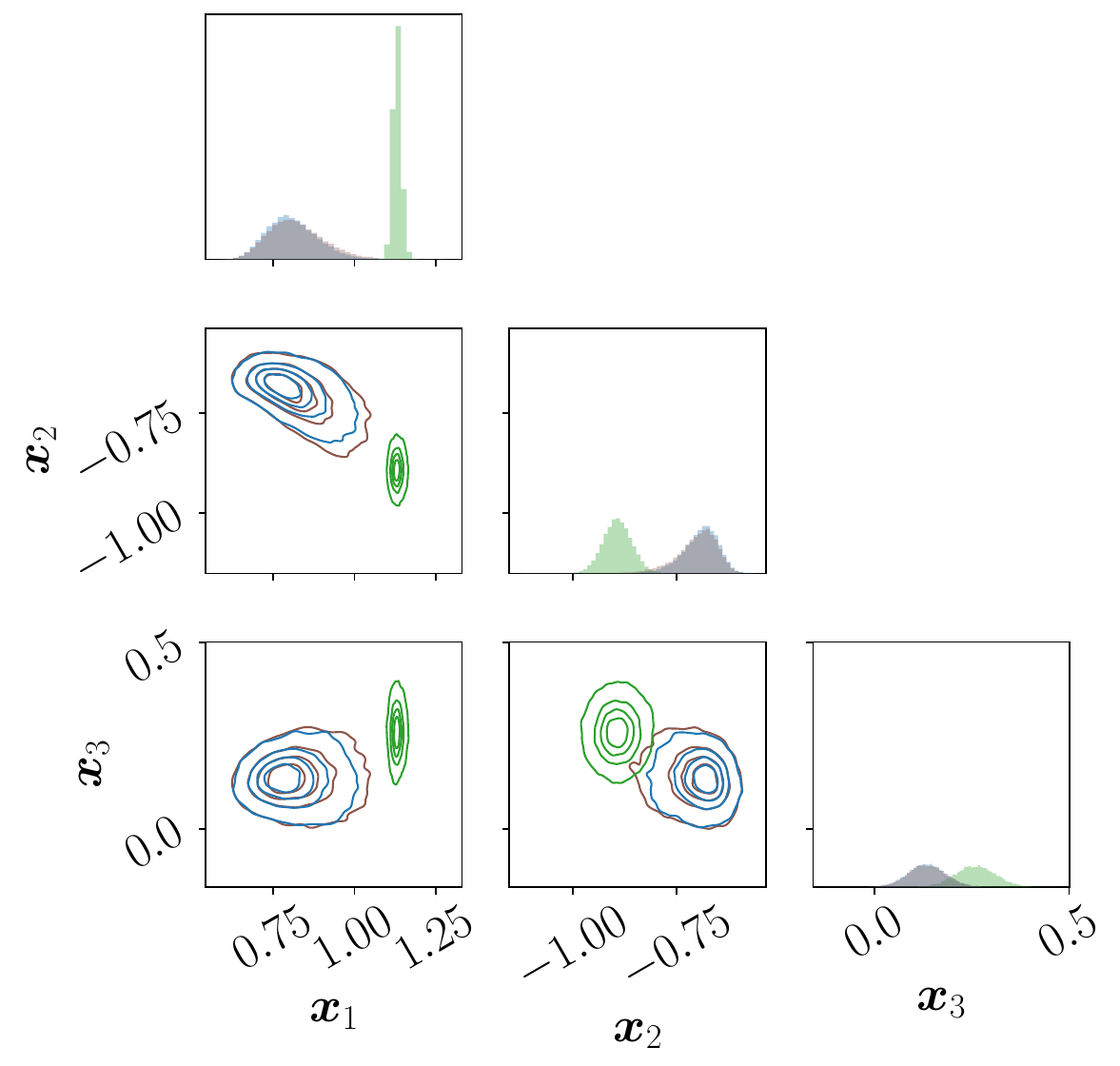} & \includegraphics[width=0.33\textwidth]{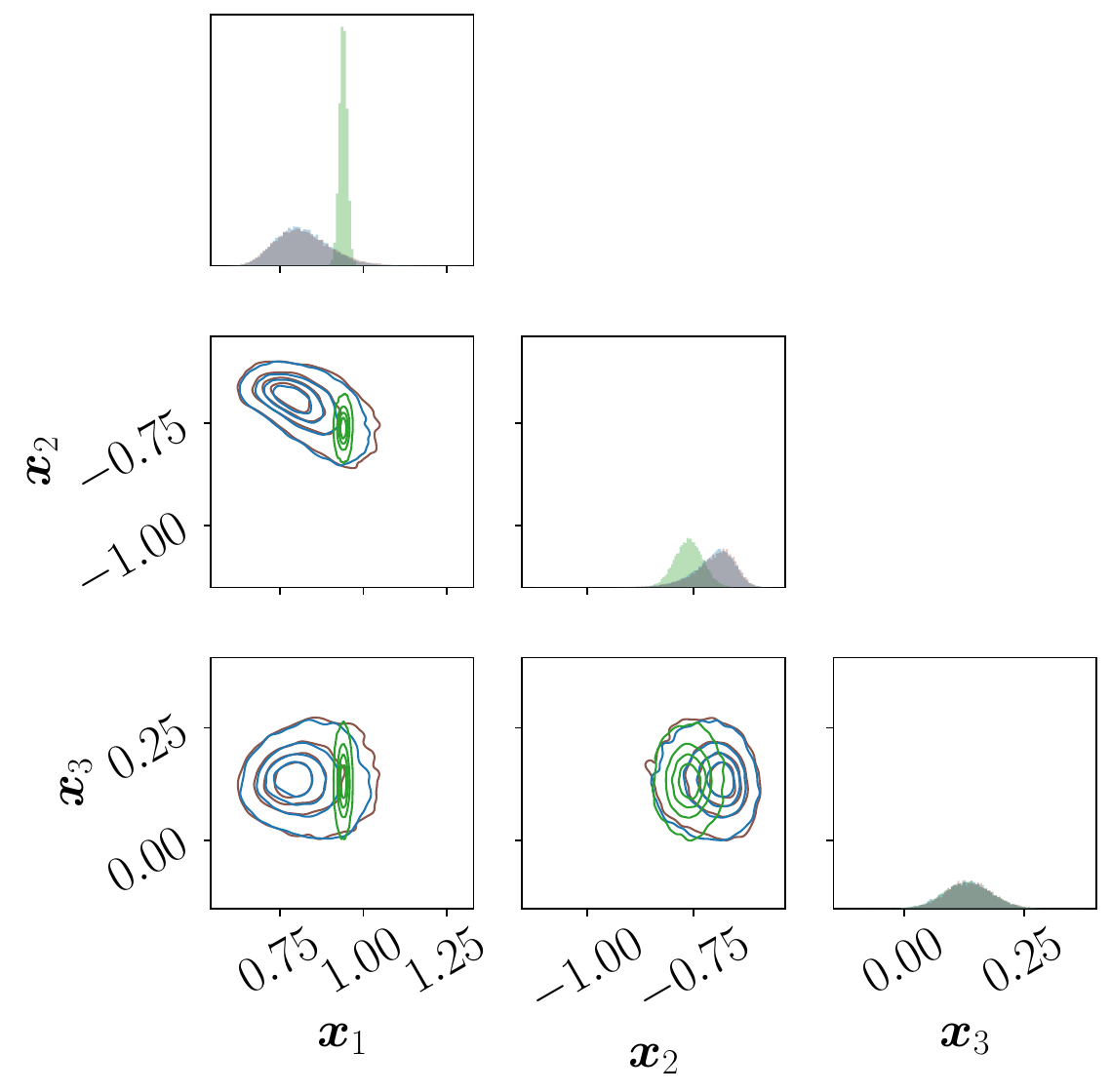}  \\
    \end{tabular}
    \addtolength{\tabcolsep}{6pt}
    \caption{\textbf{Example II: \texttt{LazyDINO} vs. \texttt{LazyMap} marginals.} \texttt{LazyMap} fails to capture the contours of the posterior marginals faithfully. \texttt{LazyMap} is off by a constant error in capturing the peak of the marginals. As more data is available, it tends to over-concentrate, leading to na\"ive underestimation of risk.}
    \label{fig:hyper_marg4}
\end{figure}

We now investigate posterior marginals for Example II in \cref{fig:hyper_marg1,fig:hyper_marg2,fig:hyper_marg3,fig:hyper_marg4}. For this set of results, the target marginals exhibit more non-Gaussianity than the previous example. In this set of results, \texttt{LazyDINO} consistently outperforms all other methods. Notably, due to the non-Gaussianity of the problem, the LA-baseline does not produce accurate marginal approximations, allowing the nonlinearly parametrized \texttt{LazyDINO} to overtake it for very limited sample data. 


\subsection{Visualization of discrepancy in mean, MAP point, and point-wise variance }

\begin{figure}[htbp]
    \centering
    \addtolength{\tabcolsep}{-6pt}
    \renewcommand{\arraystretch}{0.35}
    
    Ground truth mean\\
    \includegraphics[width=0.25\textwidth]{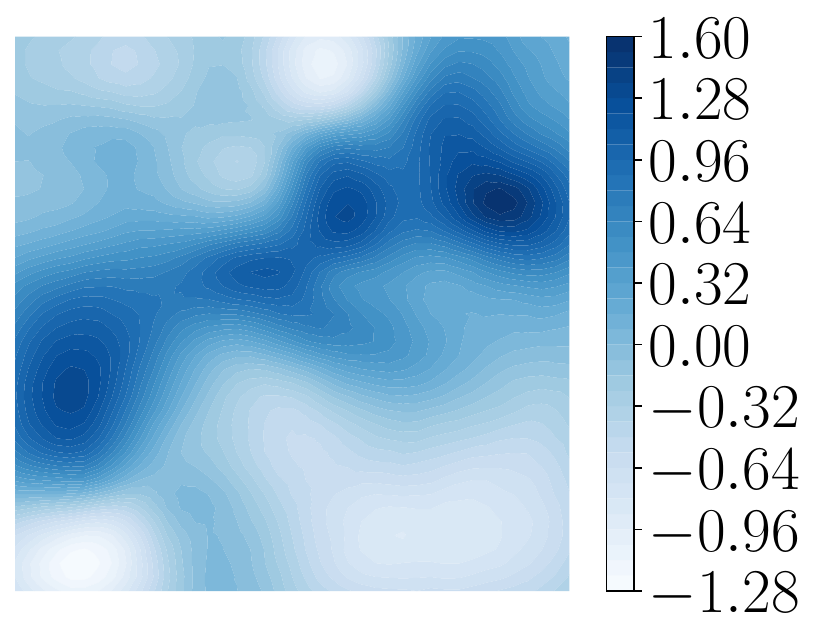}
    
    \begin{tabular}{|>{\centering\arraybackslash}p{0.75cm}|c|c|c|c|}
    \hline
    \multicolumn{5}{|c|}{\bf Posterior mean estimators (BIP \#2)} \\
    \hline \hline
    
    & \multicolumn{3}{c|}{Mean estimates} & {\makecell[c]{point-wise \\ absolute\\ error\\\includegraphics[width=0.1\textwidth]{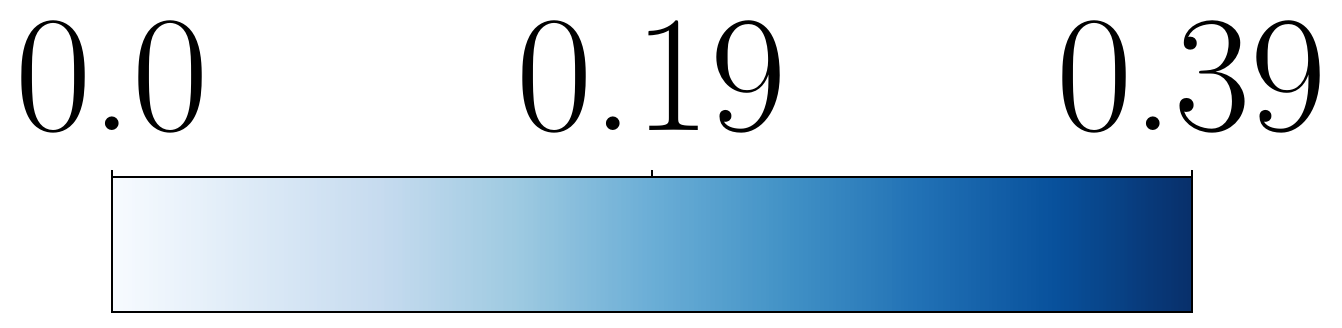}}}
    
    \\ \hline
       {$N$}& 250 & 2k & 16k &  \makecell[c]{16k}  \\ \hline 
     
     \makebox[0pt][c]{\rotatebox[origin=l]{90}{\makebox[0.1\textwidth][c]{\texttt{LazyDINO}}}} 
     & 
     \includegraphics[width=0.1\textwidth]{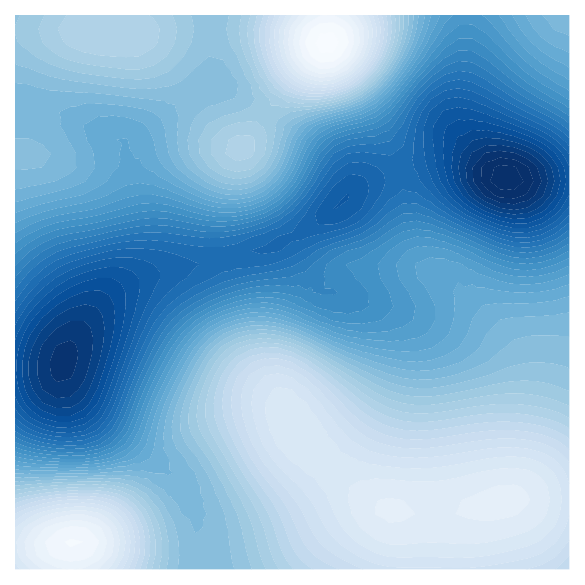} &
     \includegraphics[width=0.1\textwidth]{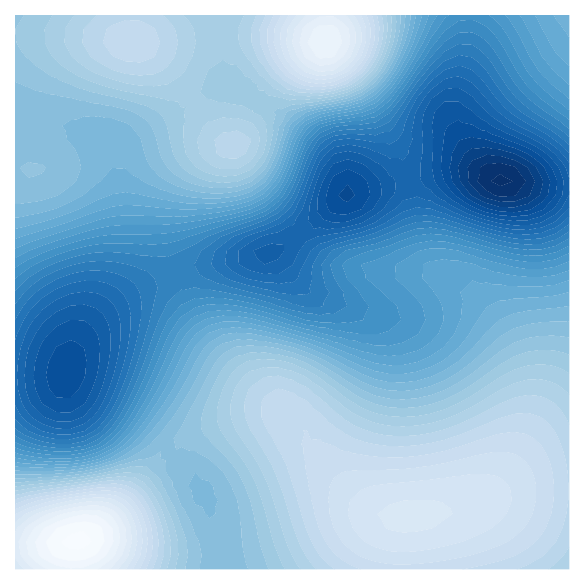} &
     \includegraphics[width=0.1\textwidth]{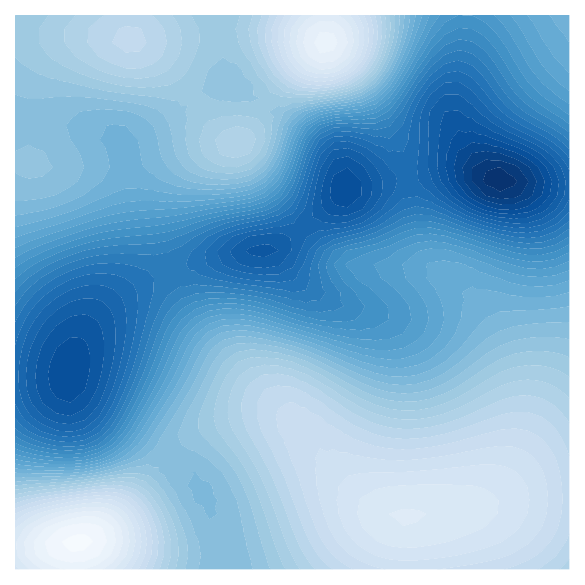} &
     \includegraphics[width=0.1\textwidth]{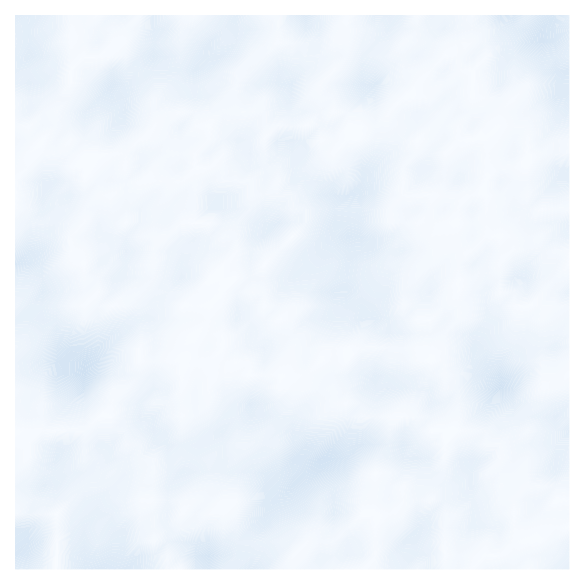} \\ \hline
     
     \makebox[0pt][c]{\rotatebox[origin=l]{90}{\makebox[0.1\textwidth][c]{\texttt{LazyNO}}}} 
     &
     \includegraphics[width=0.1\textwidth]{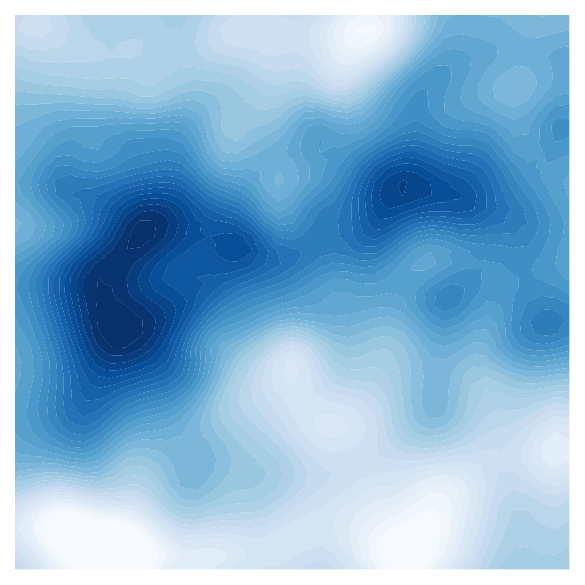} &
     \includegraphics[width=0.1\textwidth]{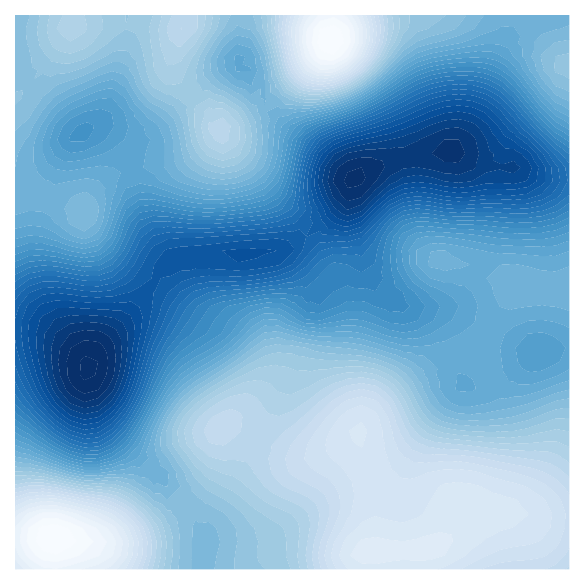} &
     \includegraphics[width=0.1\textwidth]{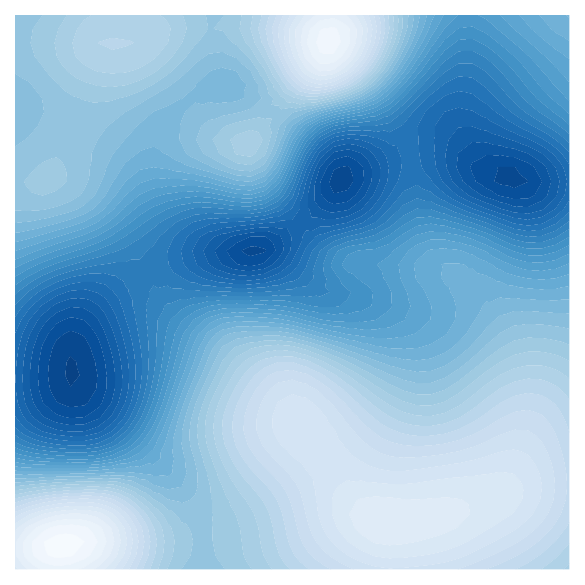} &
     \includegraphics[width=0.1\textwidth]{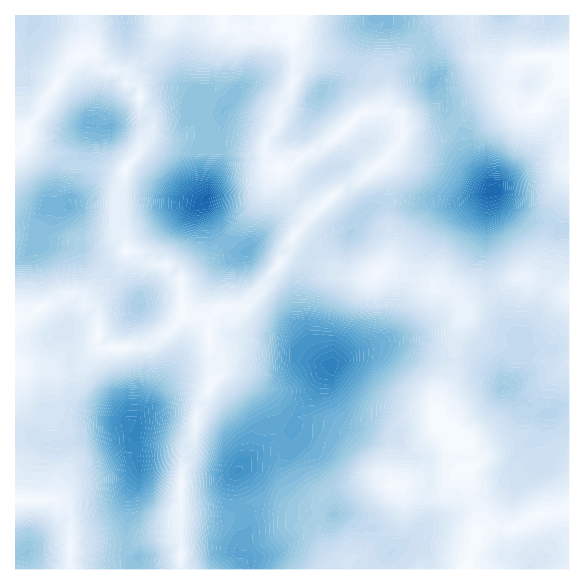} \\ \hline
     
    \makebox[0pt][c]{\rotatebox[origin=l]{90}{\makebox[0.1\textwidth][c]{SBAI}}} 
     &
     \includegraphics[width=0.1\textwidth]{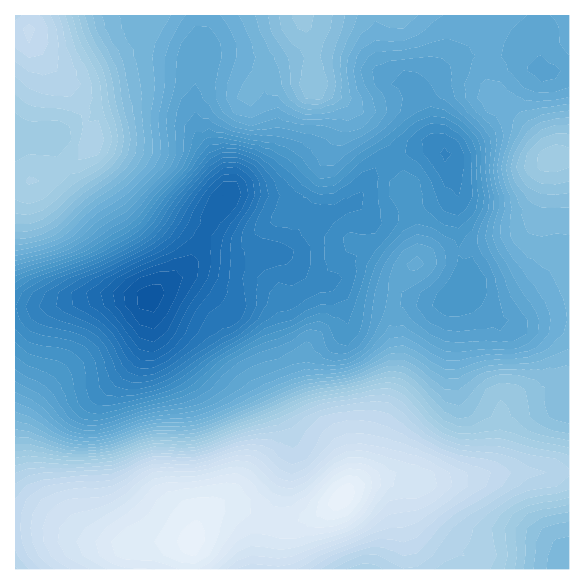} &
     \includegraphics[width=0.1\textwidth]{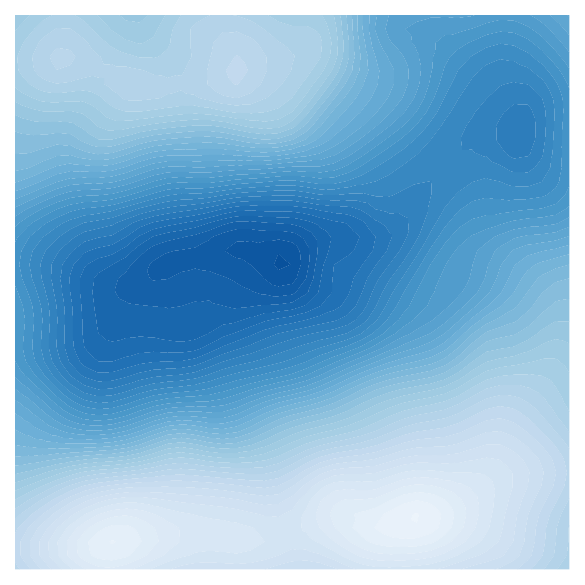} &
     \includegraphics[width=0.1\textwidth]{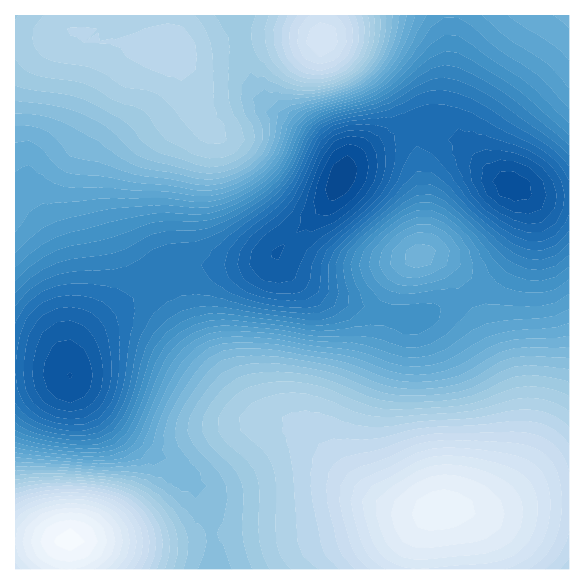} &
     \includegraphics[width=0.1\textwidth]{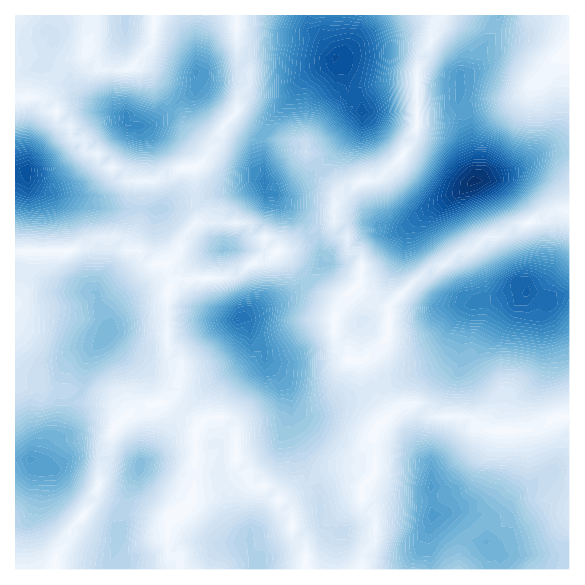} \\ 
     \hline
       {$N$}& 2k & 16k & 20k & \makecell[c]{128k}\\
       \hline
    \makebox[0pt][c]{\rotatebox[origin=l]{90}{\makebox[0.1\textwidth][c]{\texttt{LazyMap}}}} 
     &
     \includegraphics[width=0.1\textwidth]{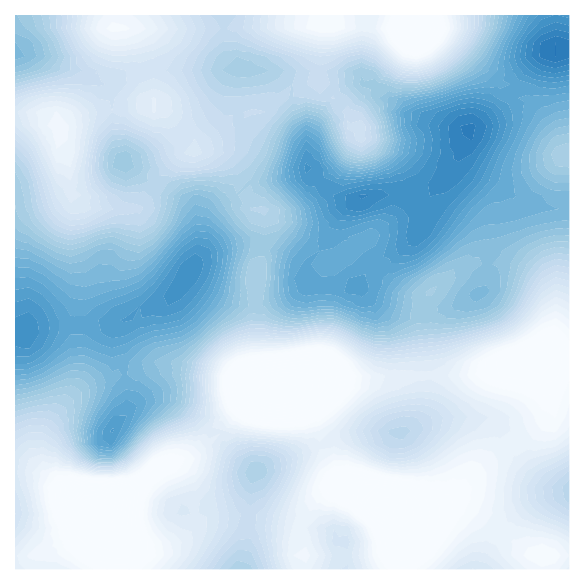} &
     \includegraphics[width=0.1\textwidth]{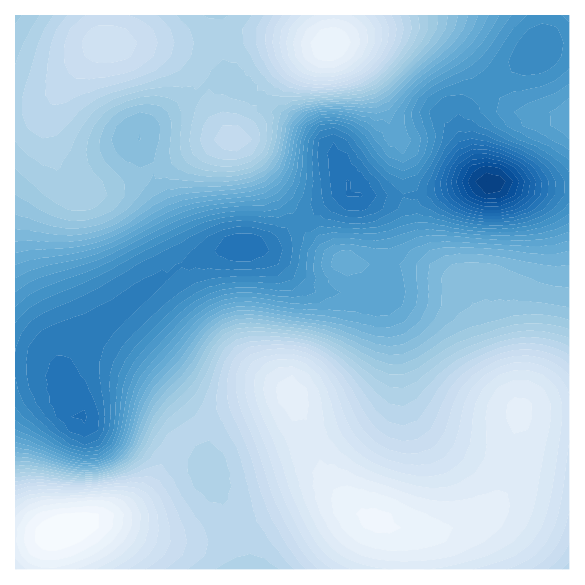} &
     \includegraphics[width=0.1\textwidth]{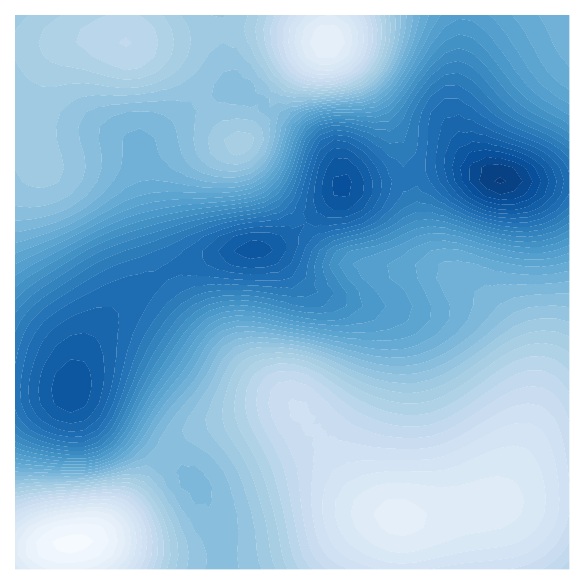} &
     \includegraphics[width=0.1\textwidth]{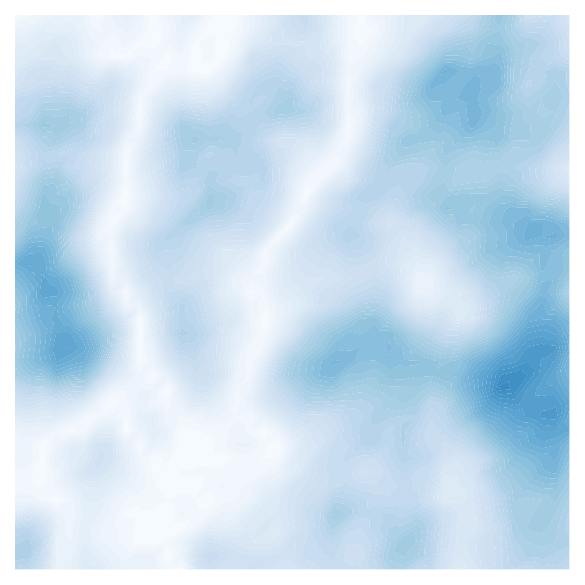} \\ \hline
\end{tabular}

     \caption{\textbf{Example I: progression of mean estimators.} \texttt{LazyDINO} already visually captures the mean well with $250$ training samples and leads to significantly smaller point-wise errors at $16,000$ samples. The next best performing method is \texttt{LazyNO}, which struggles to resolve the essential features of the mean until it has $16,000$ training samples and still results in substantially higher point-wise absolute errors than \texttt{LazyDINO} at this amount of training data. The SBAI and \texttt{LazyMap} mean estimates are poor and dominated by artifacts. The mean constructions are only reasonable at the largest training sample size and still yield substantial absolute point-wise errors. 
     }
     \label{fig:ex1_progression_mean_estimator}
\end{figure}

\begin{figure}[htbp]
    \centering
    \addtolength{\tabcolsep}{-6pt}
    \renewcommand{\arraystretch}{0.35}
    
    Ground truth MAP\\
    \includegraphics[width=0.25\textwidth]{figures/ndr_map_bip1.pdf}
    
    \begin{tabular}{|>{\centering\arraybackslash}p{0.75cm}|c|c|c|c|}
        \hline
    \multicolumn{5}{|c|}{\bf Posterior MAP estimators (BIP \#2)} \\
    \hline \hline
    
    & \multicolumn{3}{c|}{MAP estimates} & {\makecell[c]{point-wise \\ absolute\\ error\\\includegraphics[width=0.1\textwidth]{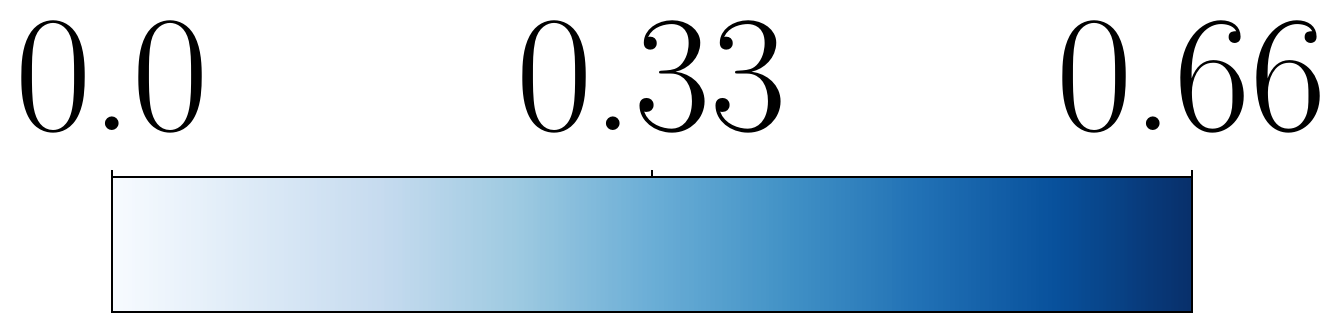}}}
    
    \\ \hline
       {$N$}& 250 & 2k & 16k &  \makecell[c]{16k}  \\ \hline 
     \makebox[0pt][c]{\rotatebox[origin=l]{90}{\makebox[0.1\textwidth][c]{\texttt{LazyDINO}}}} 
     & 
     \includegraphics[width=0.1\textwidth]{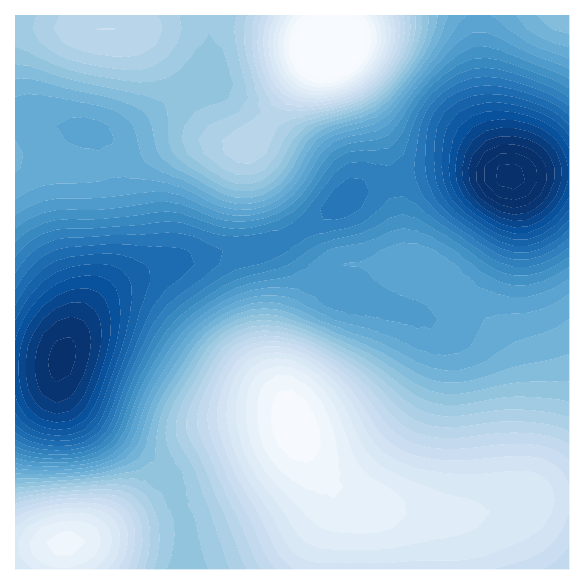} &
     \includegraphics[width=0.1\textwidth]{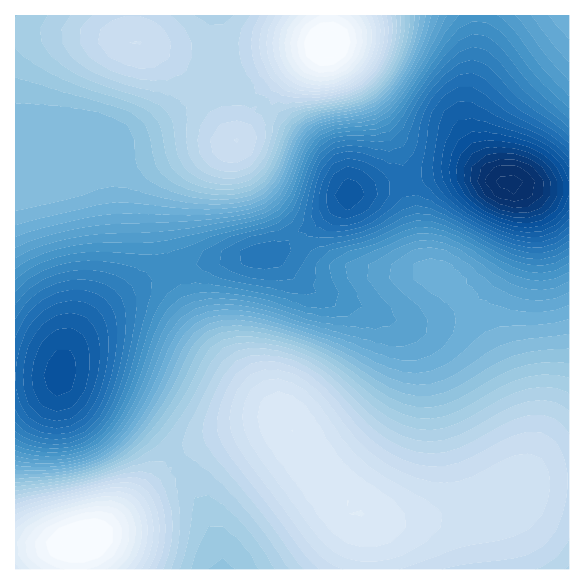} &
     \includegraphics[width=0.1\textwidth]{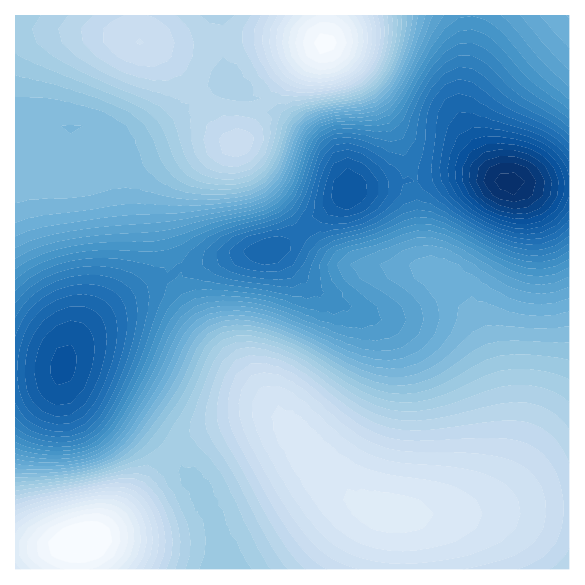} &
     \includegraphics[width=0.1\textwidth]{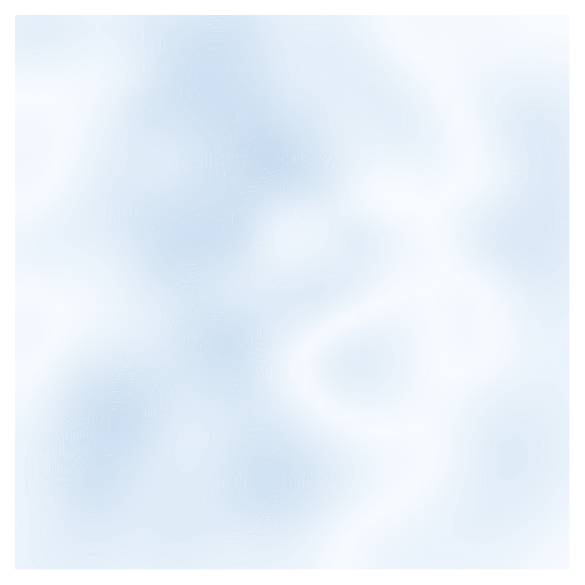} \\ \hline
     
     \makebox[0pt][c]{\rotatebox[origin=l]{90}{\makebox[0.1\textwidth][c]{\texttt{LazyNO}}}} 
     &
     \includegraphics[width=0.1\textwidth]{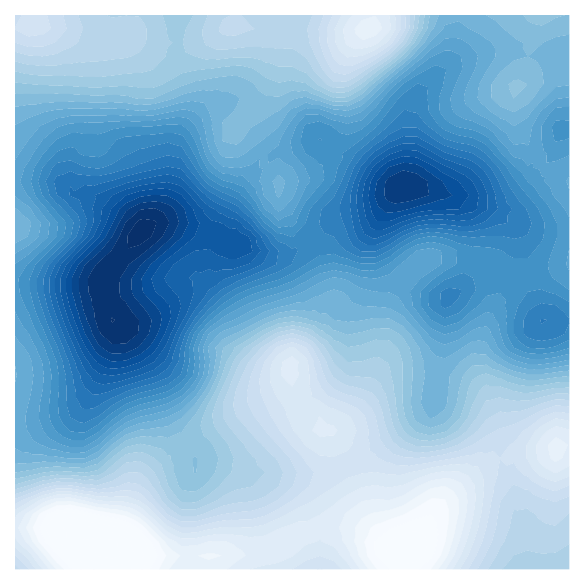} &
     \includegraphics[width=0.1\textwidth]{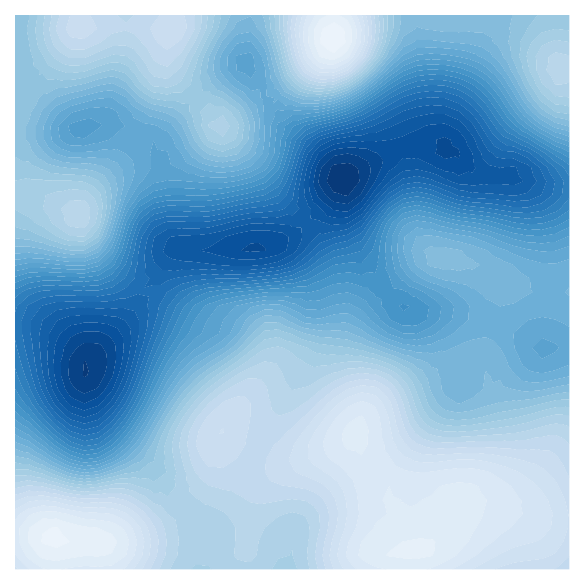} &
     \includegraphics[width=0.1\textwidth]{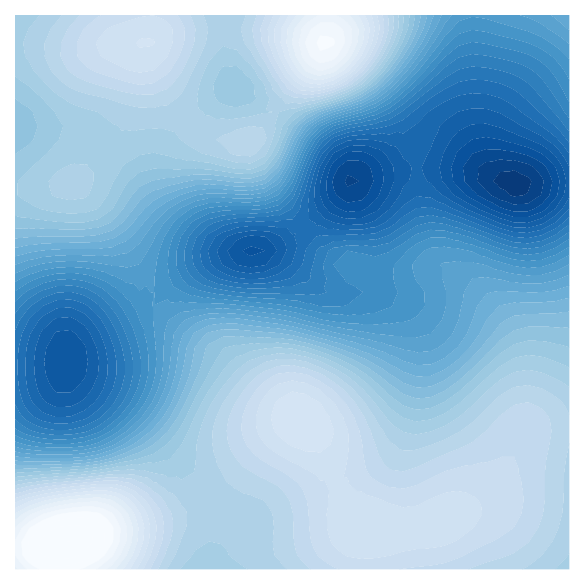} &
     \includegraphics[width=0.1\textwidth]{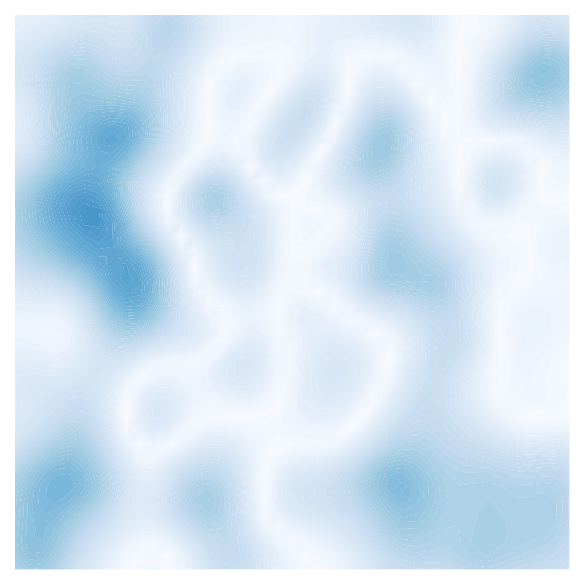} \\ \hline
     
    \makebox[0pt][c]{\rotatebox[origin=l]{90}{\makebox[0.1\textwidth][c]{SBAI}}} 
     &
     \includegraphics[width=0.1\textwidth]{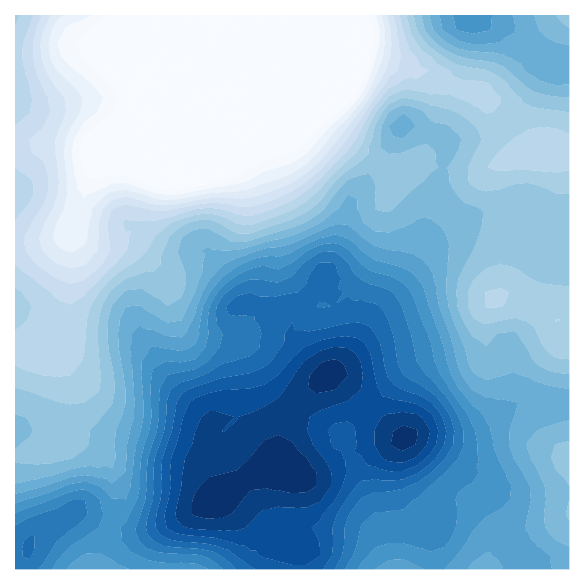} &
     \includegraphics[width=0.1\textwidth]{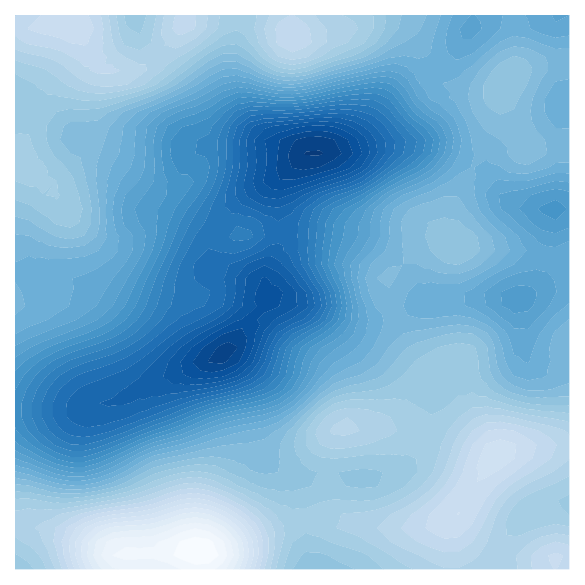} &
     \includegraphics[width=0.1\textwidth]{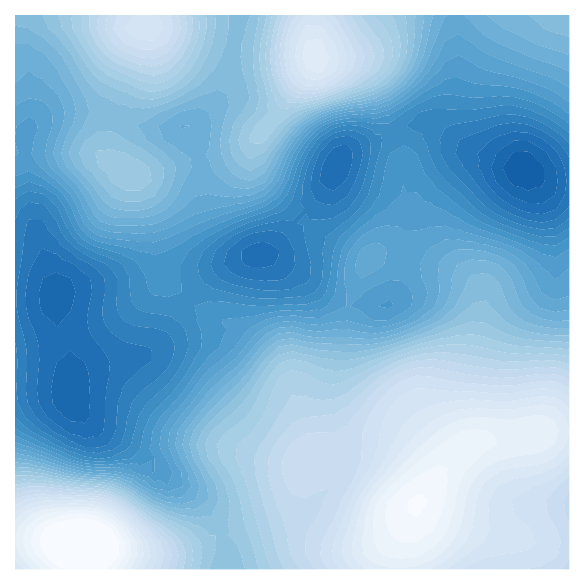} &
     \includegraphics[width=0.1\textwidth]{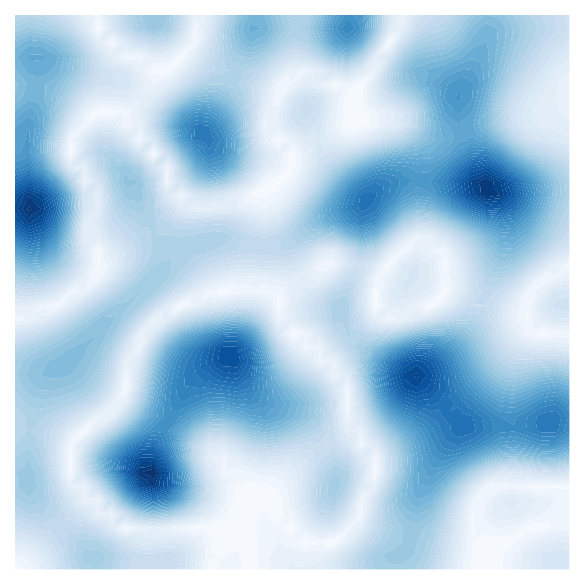} \\ 
     \hline
       {$N$}& 2k & 16k & 128k & \makecell[c]{128k}\\
       \hline
    \makebox[0pt][c]{\rotatebox[origin=l]{90}{\makebox[0.1\textwidth][c]{\texttt{LazyMap}}}} 
     &
     \includegraphics[width=0.1\textwidth]{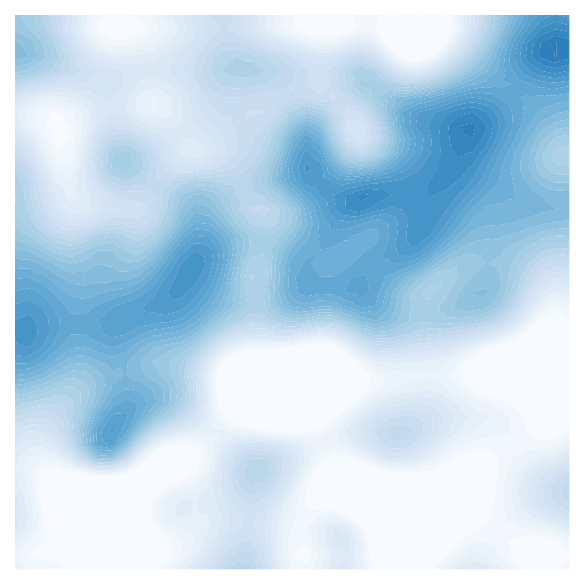} &
     \includegraphics[width=0.1\textwidth]{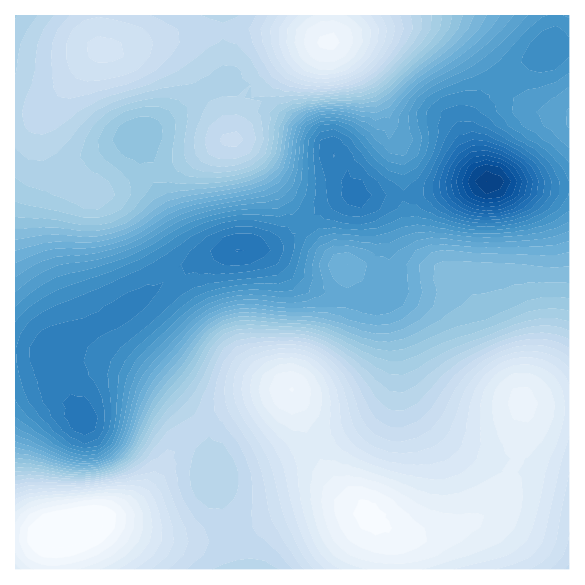} &
     \includegraphics[width=0.1\textwidth]{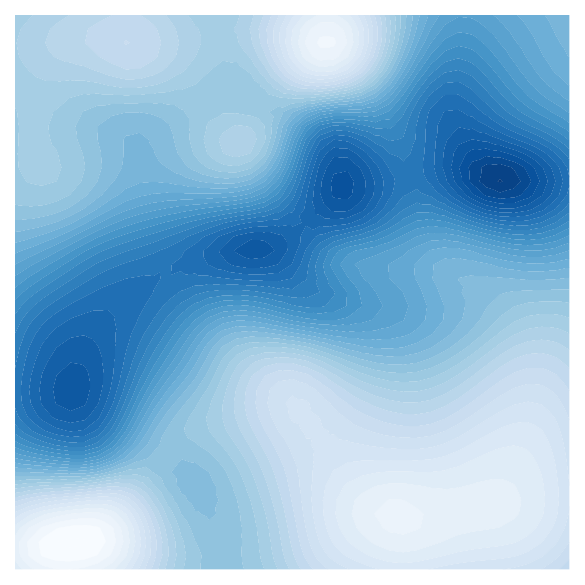} &
     \includegraphics[width=0.1\textwidth]{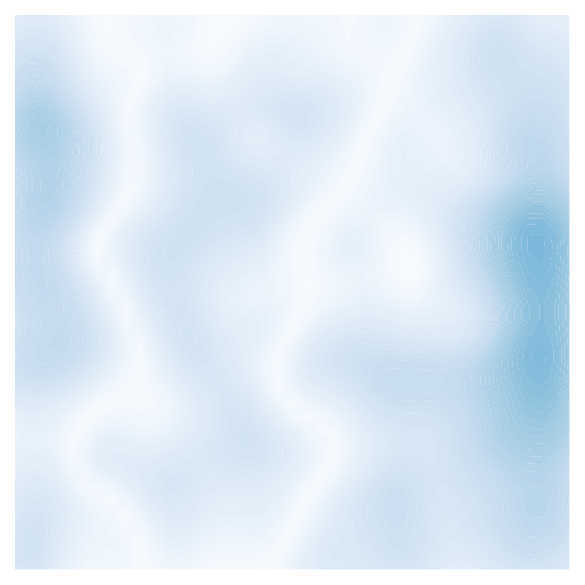} \\ \hline
     
    \end{tabular}
     \caption{\textbf{Example I: progression of MAP estimators.} In a similar story to \cref{fig:ex1_progression_mean_estimator} the \texttt{LazyDINO} MAP reconstruction is already quite accurate for $250$ training data, and consistently outperforms the other methods. \texttt{LazyNO} and \texttt{LazyMap} yield reasonable MAP reconstructions given $16,000$ and $128,000$ samples respectively, albeit with higher absolute point-wise errors than \texttt{LazyDINO}. The SBAI MAP point estimate is poor even for $16,000$ training data. This is evident from the absolute point-wise errors. 
     }\label{fig:ex1_progression_map}

\end{figure}

\begin{figure}[htbp]
    \centering
    \addtolength{\tabcolsep}{-6pt}
    \renewcommand{\arraystretch}{0.35}
    
        Ground truth point-wise \\marginal variance\\
        \includegraphics[width=0.25\textwidth]{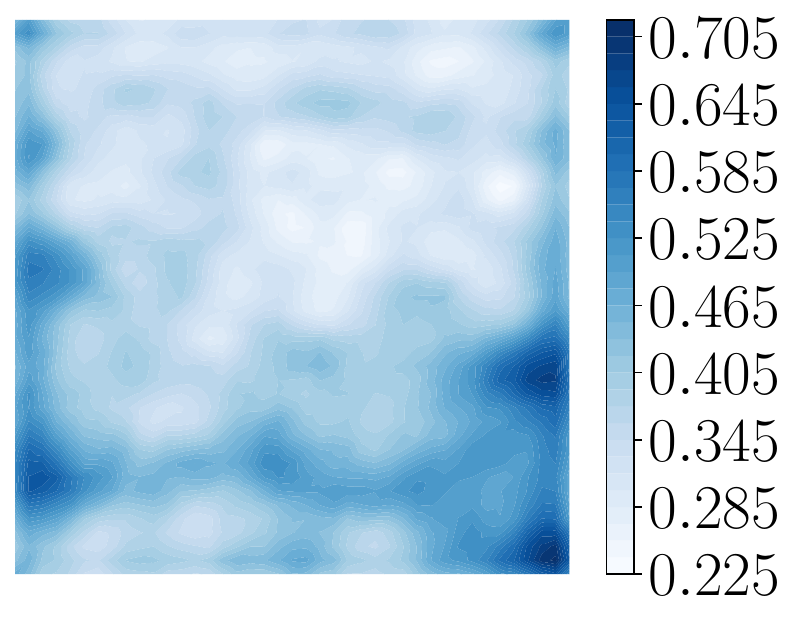}
        
    \begin{tabular}{|>{\centering\arraybackslash}p{0.75cm}|c|c|c|c|}
        \hline
    \multicolumn{5}{|c|}{\bf \makecell{Posterior point-wise marginal \\variance estimators (BIP \#2)}} \\
    \hline \hline
    
    & \multicolumn{3}{c|}{\makecell[c]{Point-wise marginal \\variance estimates}} & {\makecell[c]{point-wise \\ absolute\\ error\\\includegraphics[width=0.1\textwidth]{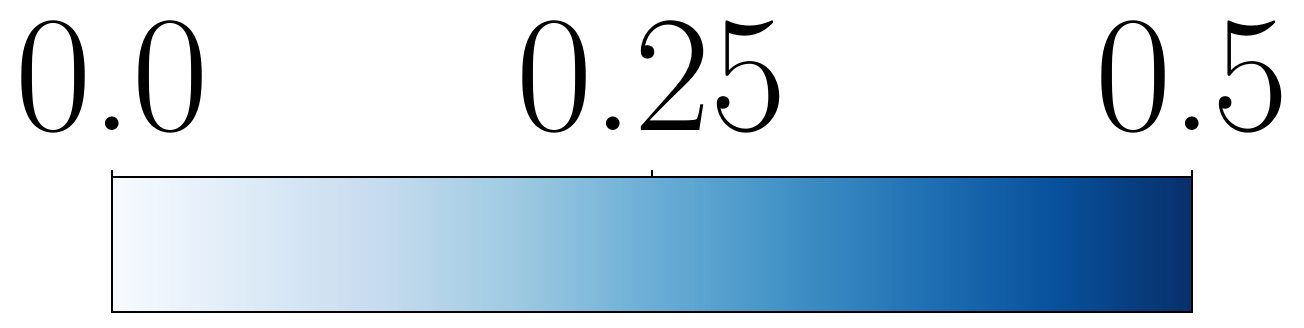}}}
    
    \\ \hline
       {$N$}& 250 & 2k & 16k &  \makecell[c]{16k}  \\ \hline 
     \makebox[0pt][c]{\rotatebox[origin=l]{90}{\makebox[0.1\textwidth][c]{\texttt{LazyDINO}}}} 
     & 
     \includegraphics[width=0.1\textwidth]{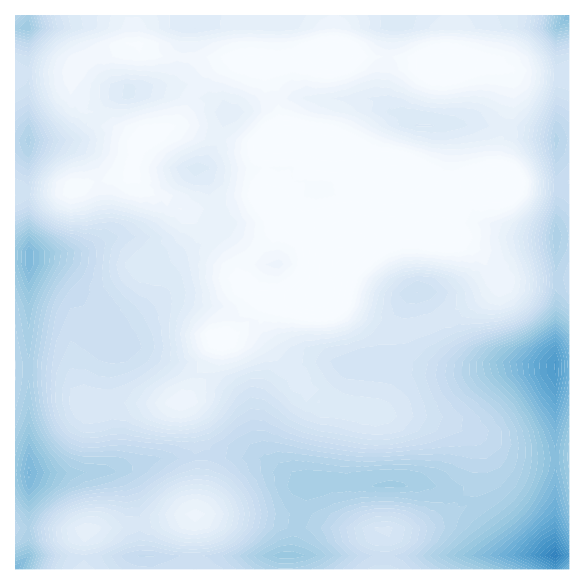} &
     \includegraphics[width=0.1\textwidth]{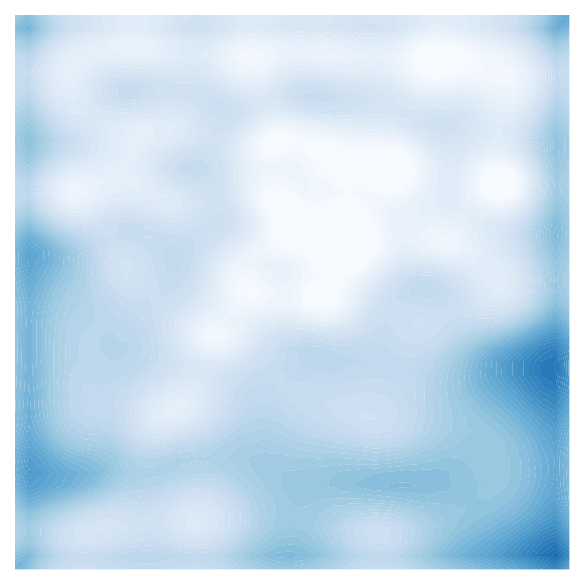} &
     \includegraphics[width=0.1\textwidth]{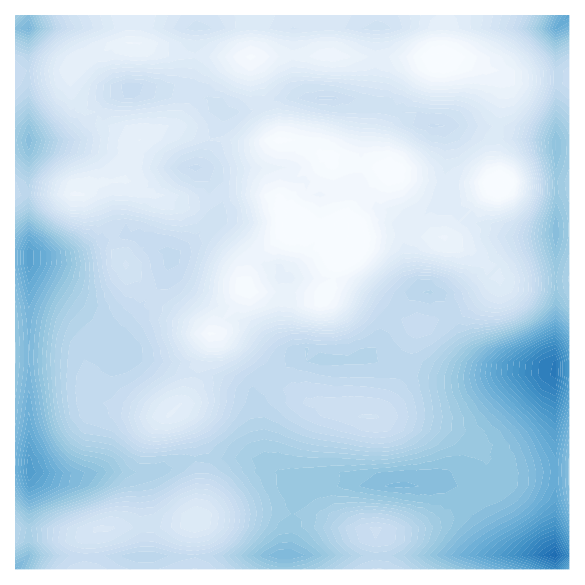} &
     \includegraphics[width=0.1\textwidth]{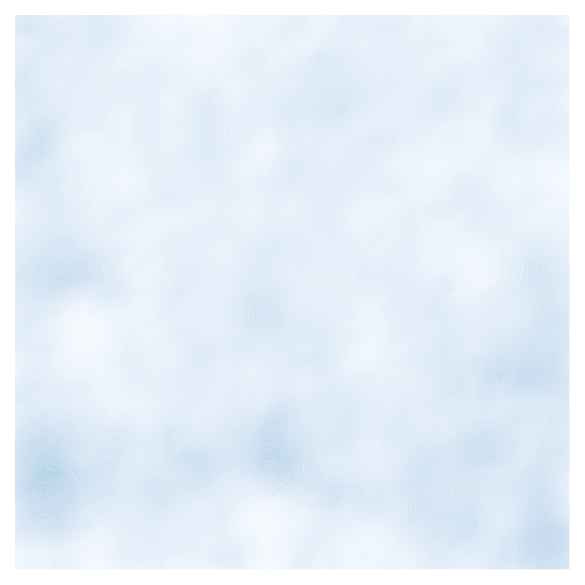} \\ \hline
     
     \makebox[0pt][c]{\rotatebox[origin=l]{90}{\makebox[0.1\textwidth][c]{\texttt{LazyNO}}}} 
     &
     \includegraphics[width=0.1\textwidth]{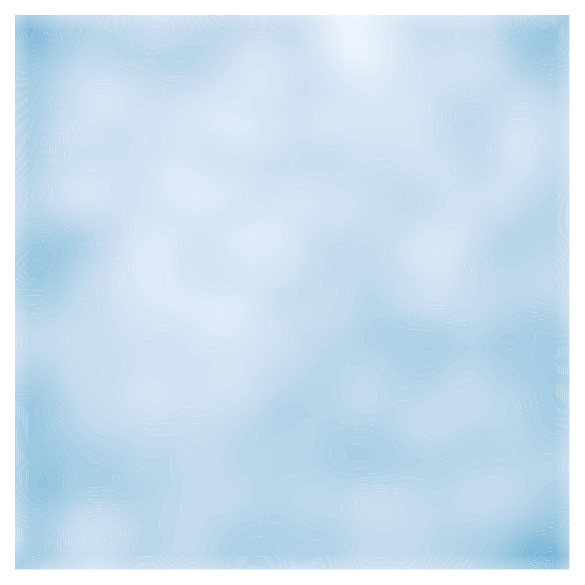} &
     \includegraphics[width=0.1\textwidth]{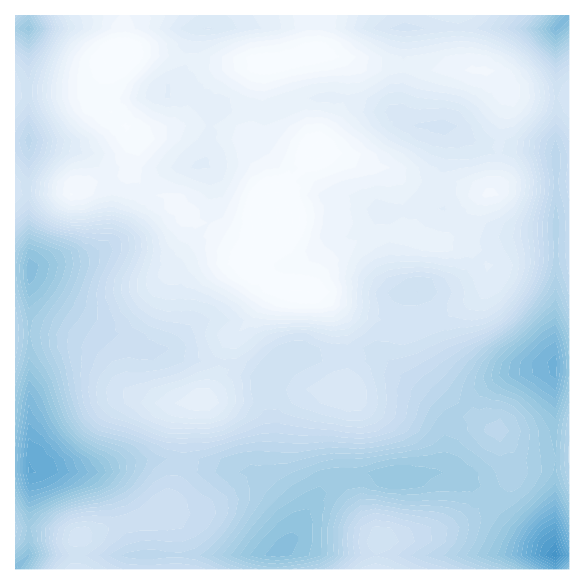} &
     \includegraphics[width=0.1\textwidth]{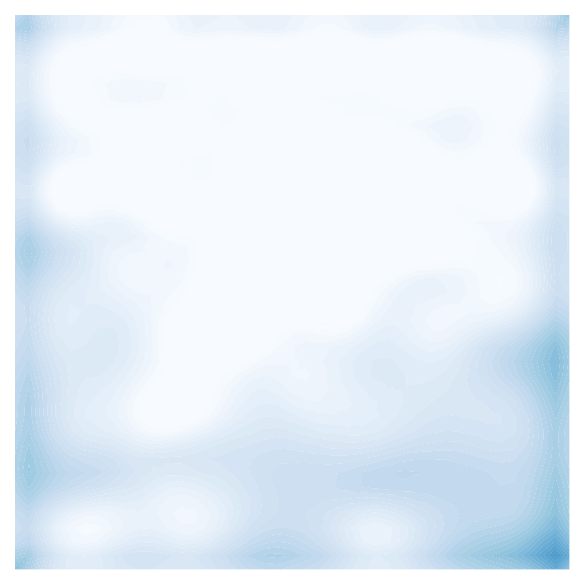} &
     \includegraphics[width=0.1\textwidth]{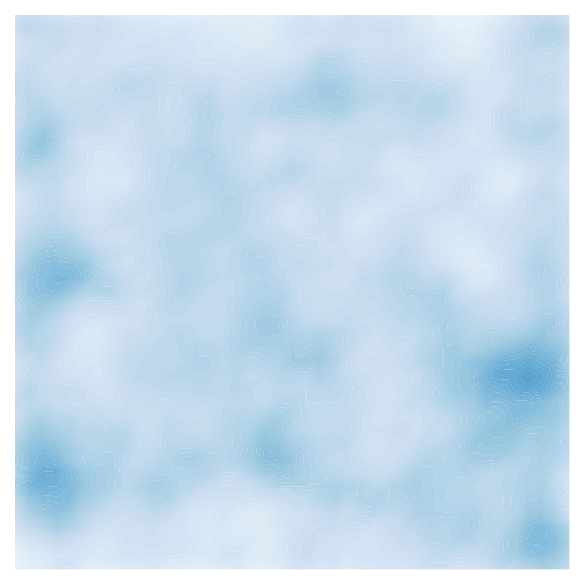} \\ \hline
     
    \makebox[0pt][c]{\rotatebox[origin=l]{90}{\makebox[0.1\textwidth][c]{SBAI}}} 
     &
     \includegraphics[width=0.1\textwidth]{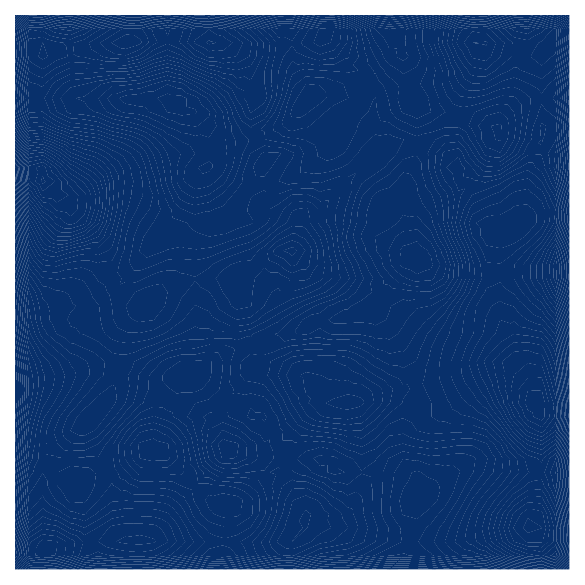} &
     \includegraphics[width=0.1\textwidth]{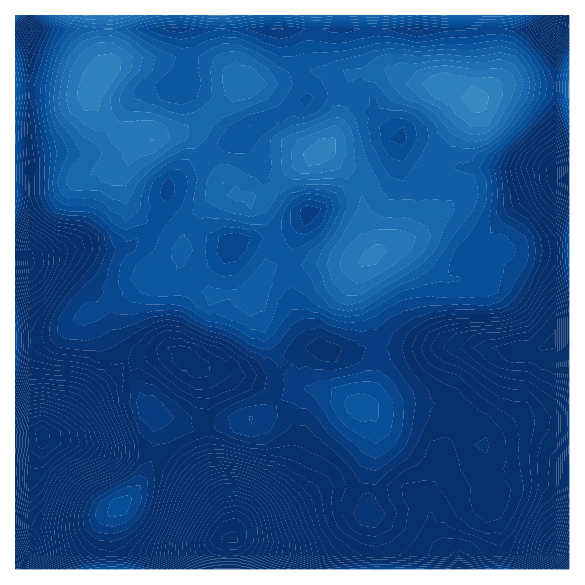} &
     \includegraphics[width=0.1\textwidth]{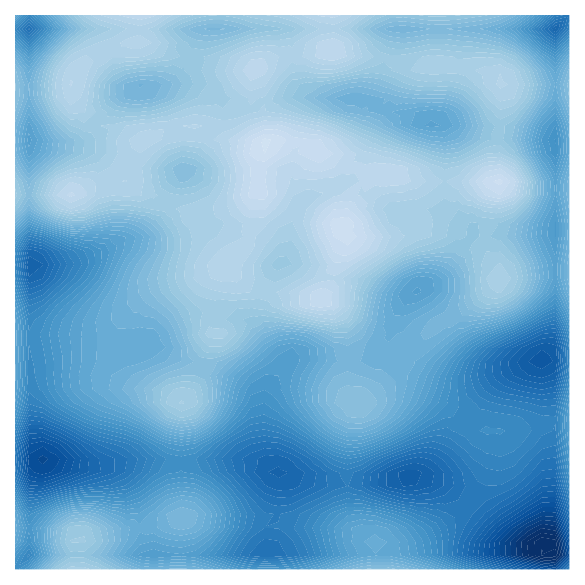} &
     \includegraphics[width=0.1\textwidth]{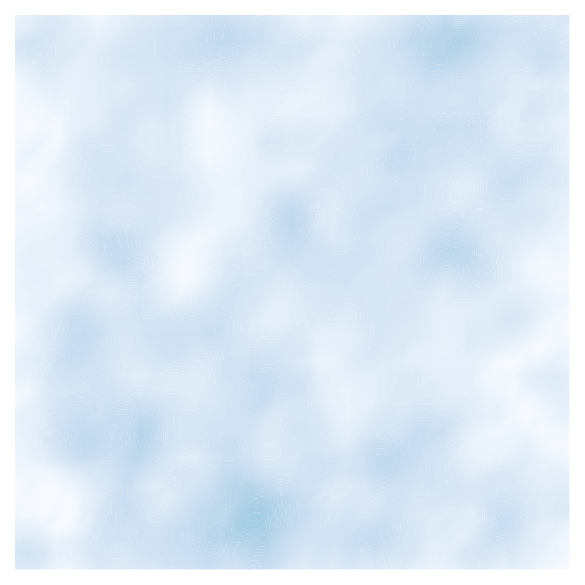} \\ 
     \hline
       {$N$}& 2k & 16k & 20k & \makecell[c]{128k}\\
       \hline
    \makebox[0pt][c]{\rotatebox[origin=l]{90}{\makebox[0.1\textwidth][c]{\texttt{LazyMap}}}} 
     &
     \includegraphics[width=0.1\textwidth]{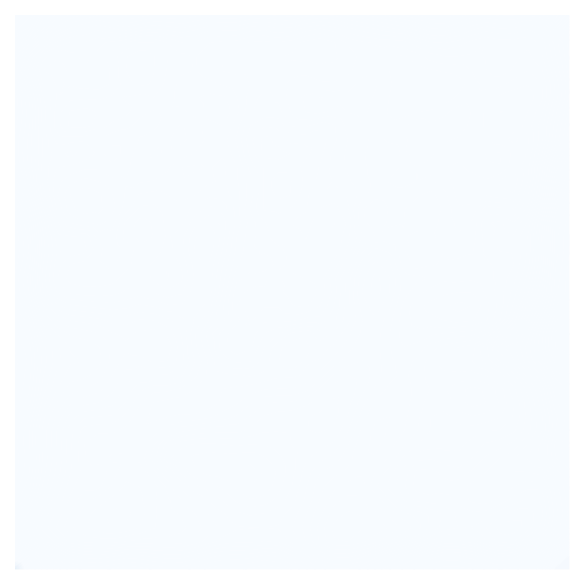} &
     \includegraphics[width=0.1\textwidth]{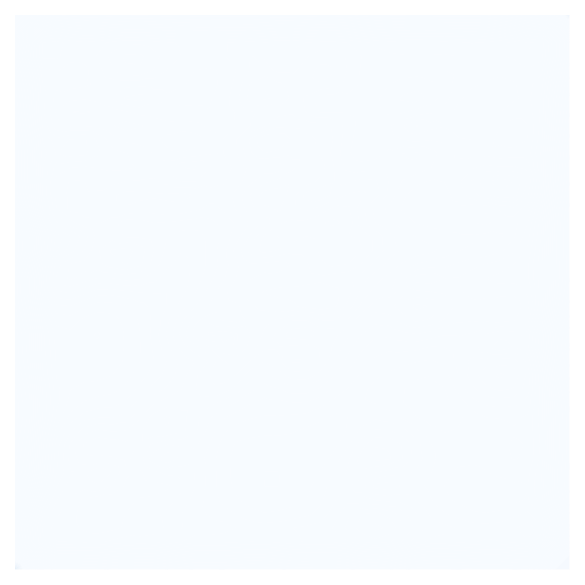} &
     \includegraphics[width=0.1\textwidth]{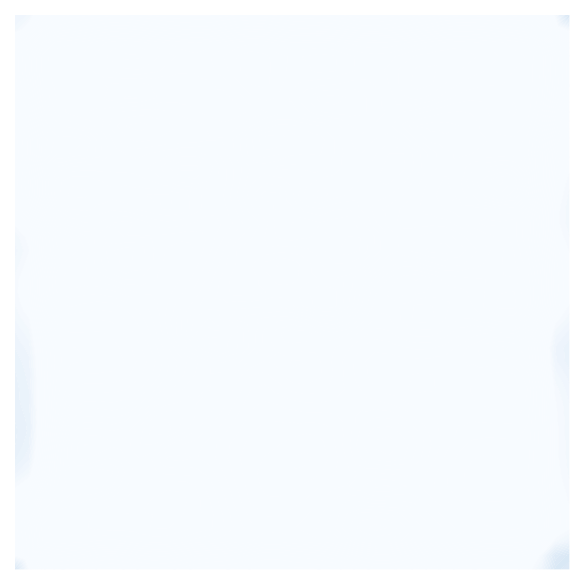} &
     \includegraphics[width=0.1\textwidth]{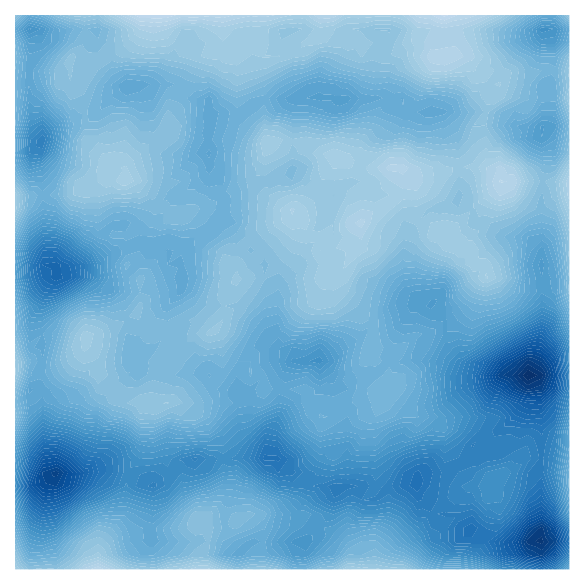} \\ \hline
     
    \end{tabular}
     \caption{\textbf{Example I: progression of point-wise marginal variance estimators.} 
     Again, the \texttt{LazyDINO} estimate of point-wise marginal variance is superior to the other methods as in the previous studies of mean and MAP estimation. Notably, SBAI overestimates the point-wise marginal variance, while \texttt{LazyMap} significantly underestimates the marginal variance. This is consistent with the evidence provided in \cref{fig:ndr_marg3,fig:ndr_marg4}, which demonstrates a similar phenomenon in the marginals: SBAI is spread out while \texttt{LazyMap} is highly concentrated. 
     }\label{fig:ex_1_progression_variance}
\end{figure}

In this section, we proceed with additional visualizations of the discrepancies in the mean, MAP point estimates and point-wise variances, and associated point-wise absolute errors arising from the different methods for varying amounts of training data used. These comparisons allow the methods to be visually differentiated in their ability to resolve features in the parameter reconstruction via the BIP. We begin by visualizing the mean, MAP estimates, and point-wise variances for Example I in \cref{fig:ex1_progression_mean_estimator,fig:ex1_progression_map,fig:ex_1_progression_variance}, respectively. The general trend is consistent with the previous numerical studies: \texttt{LazyDINO} leads to superior approximation than the other methods. Notably, \texttt{LazyDINO} can capture the mean, MAP, and point-wise marginal variance somewhat faithfully for $250$ samples, while other methods struggle with orders of magnitude more samples. 


\clearpage

\begin{figure}[htbp]
    \centering
        \addtolength{\tabcolsep}{-6pt}
        \renewcommand{\arraystretch}{0.35}

   Ground truth mean\\
    \includegraphics[width=0.35\textwidth]{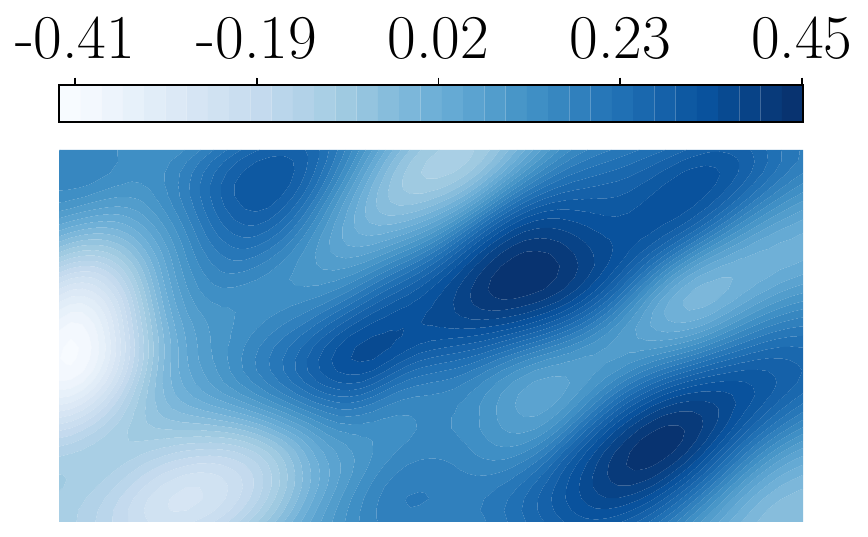}
    
  \begin{tabular}{|>{\centering\arraybackslash}p{0.75cm}|c|c|c|c|}
    \hline
    \multicolumn{5}{|c|}{\bf \makecell{Progression of posterior mean estimator (BIP \#4)}} \\
    \hline \hline
    
    & \multicolumn{3}{c|}{Mean estimates} & {\makecell[c]{point-wise \\ absolute error\\\includegraphics[width=0.2\textwidth]{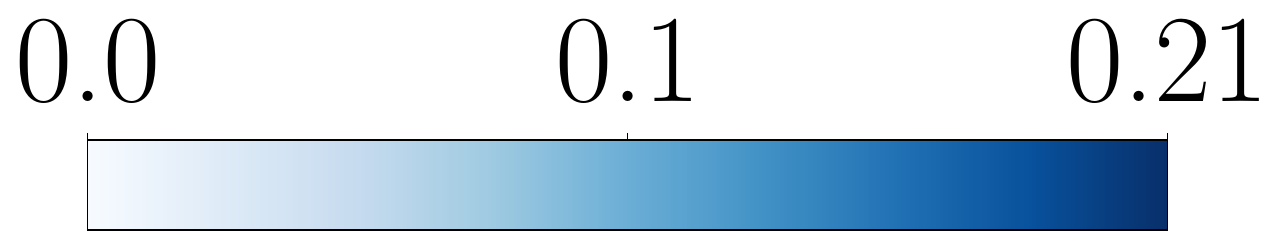}}}
    
    \\ \cline{1-5}
     $N$& 250 & 2k & 16k &  \makecell[c]{16k}  \\ \hline
     
     \makebox[0pt][c]{\rotatebox[origin=l]{90}{\makebox[0.1\textwidth][c]{\texttt{LazyDINO}}}} 
     & 
     \includegraphics[width=0.2\textwidth]{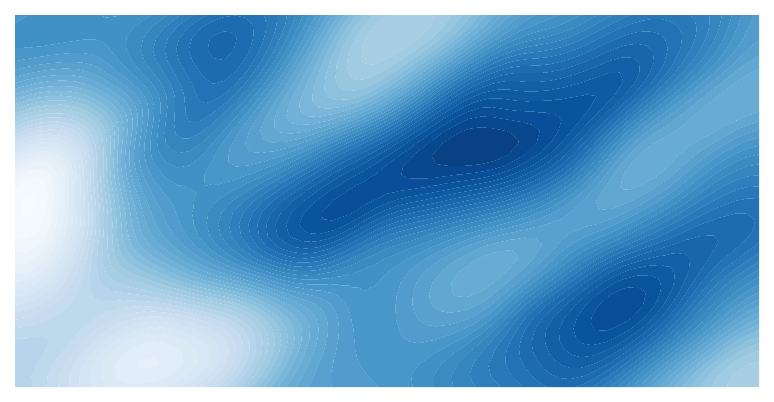} &
     \includegraphics[width=0.2\textwidth]{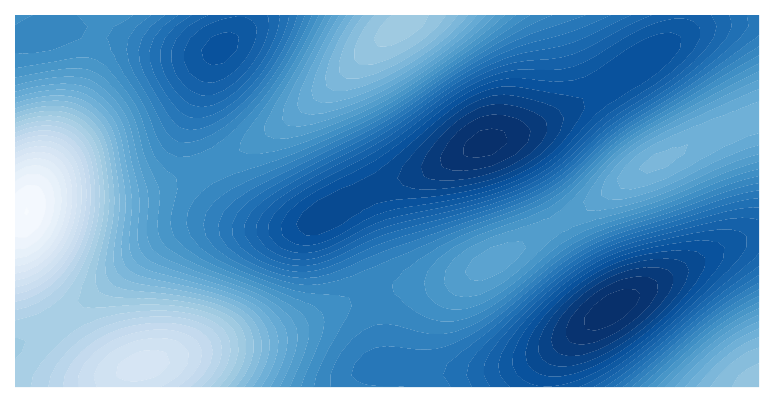} &
     \includegraphics[width=0.2\textwidth]{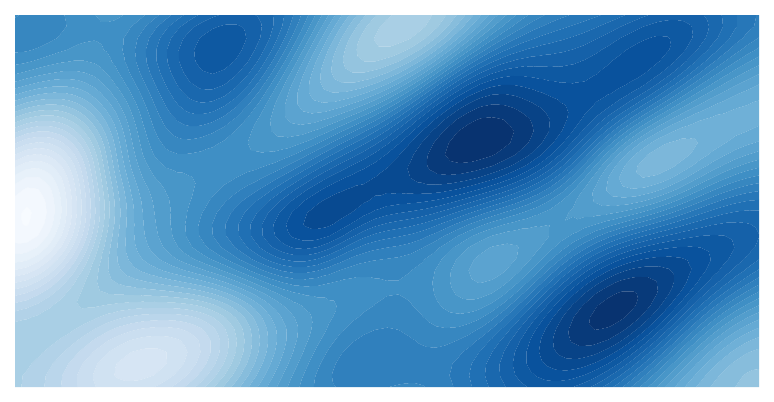} &
     \includegraphics[width=0.2\textwidth]{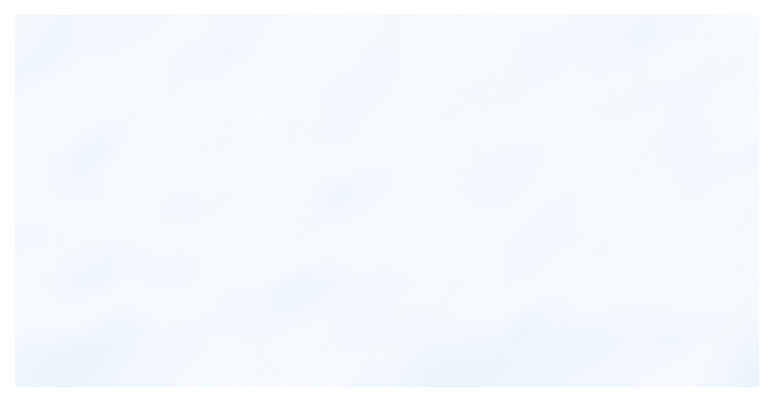} \\ \hline

   \makebox[0pt][c]{\rotatebox[origin=l]{90}{\makebox[0.1\textwidth][c]{\texttt{LazyNO}}}} 
     &
     \includegraphics[width=0.2\textwidth]{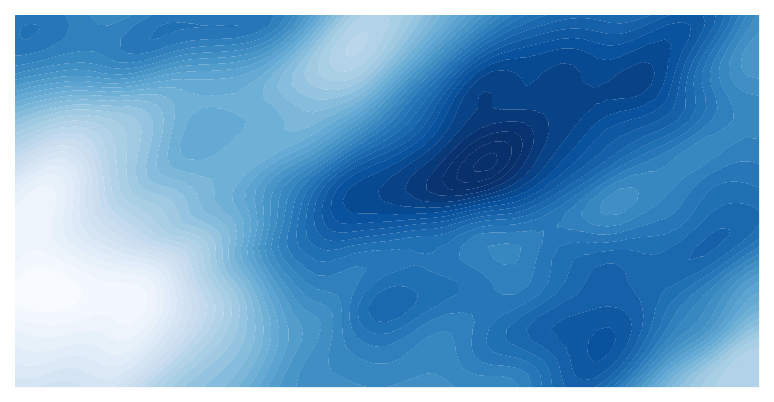} &
     \includegraphics[width=0.2\textwidth]{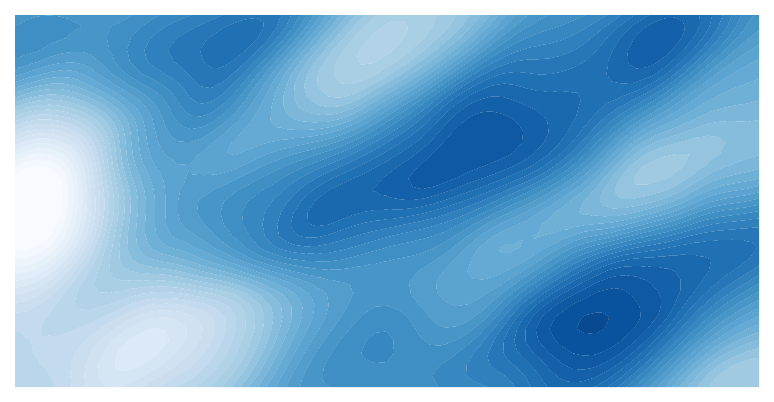} &
     \includegraphics[width=0.2\textwidth]{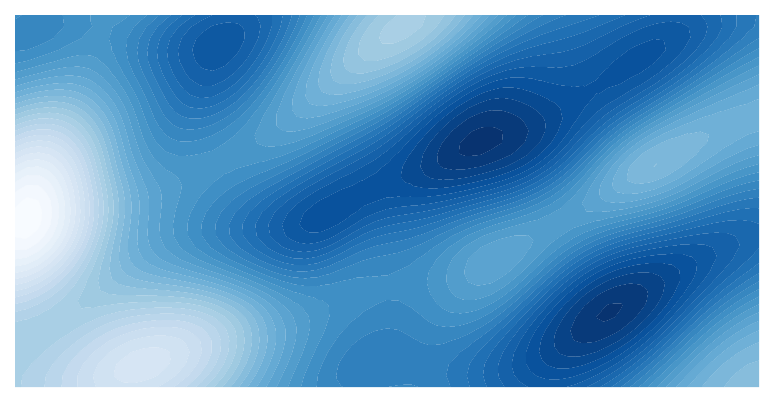} &
     \includegraphics[width=0.2\textwidth]{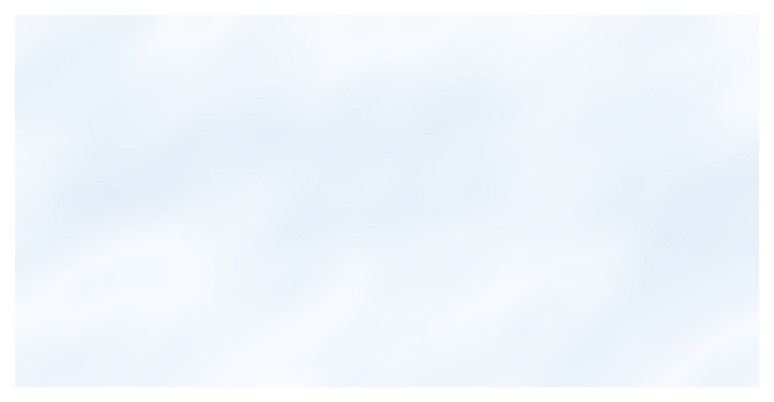} \\ \hline
     
    \makebox[0pt][c]{\rotatebox[origin=l]{90}{\makebox[0.1\textwidth][c]{SBAI}}} 
     &
     \includegraphics[width=0.2\textwidth]{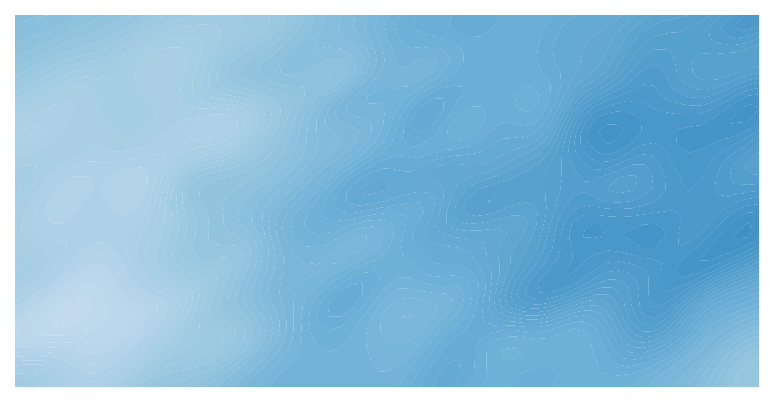} &
     \includegraphics[width=0.2\textwidth]{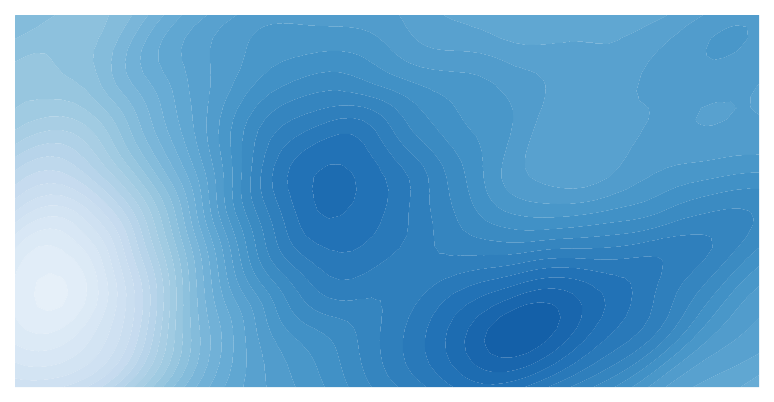} &
     \includegraphics[width=0.2\textwidth]{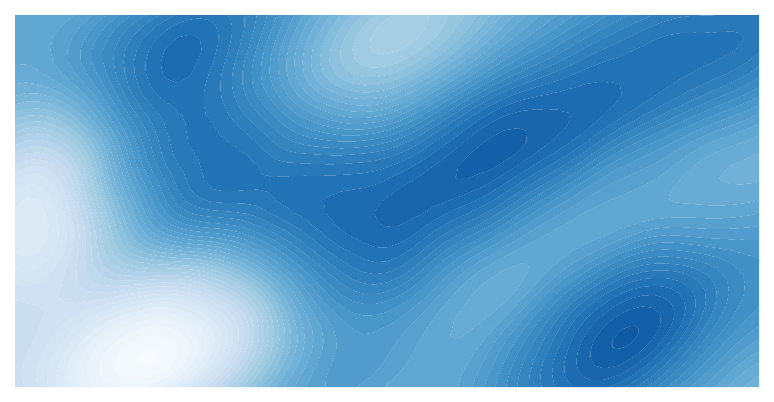} &
     \includegraphics[width=0.2\textwidth]{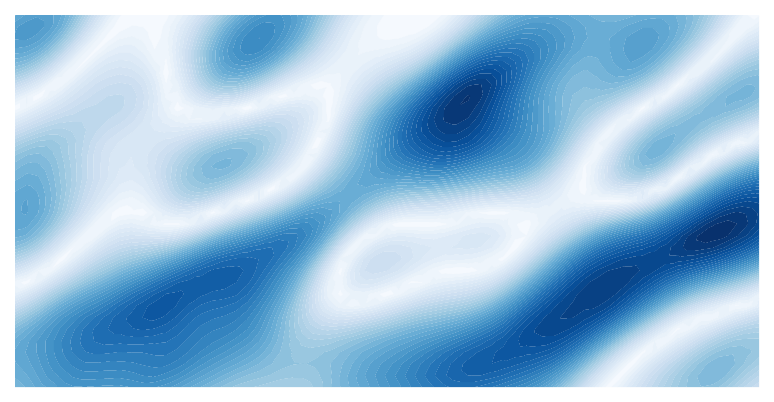} \\ 
     \hline
       {$N$}& 2k & 16k & 20k & \makecell[c]{20k}\\
       \hline
    \makebox[0pt][c]{\rotatebox[origin=l]{90}{\makebox[0.1\textwidth][c]{\texttt{LazyMap}}}} 
     &
     \includegraphics[width=0.2\textwidth]{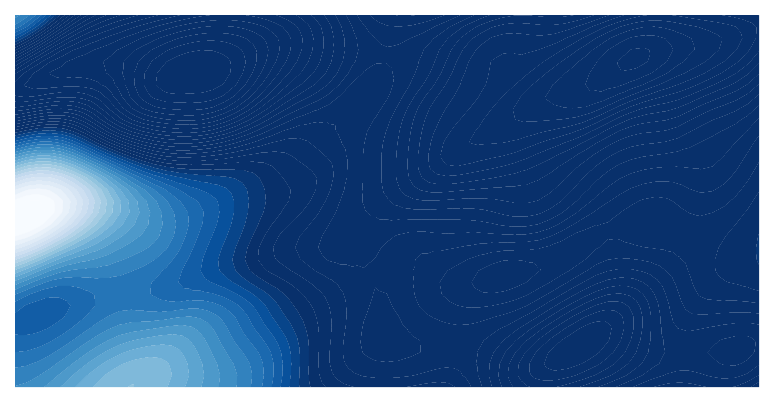} &
     \includegraphics[width=0.2\textwidth]{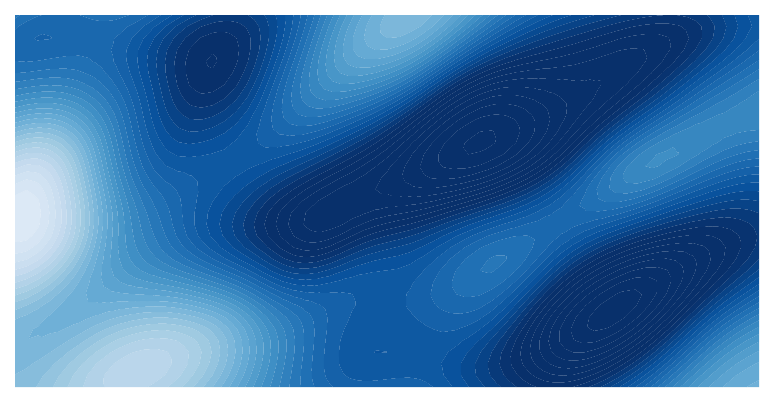} &
     \includegraphics[width=0.2\textwidth]{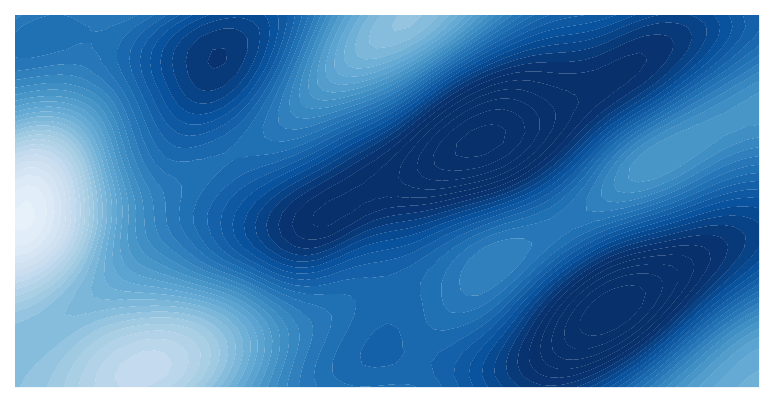} &
     \includegraphics[width=0.2\textwidth]{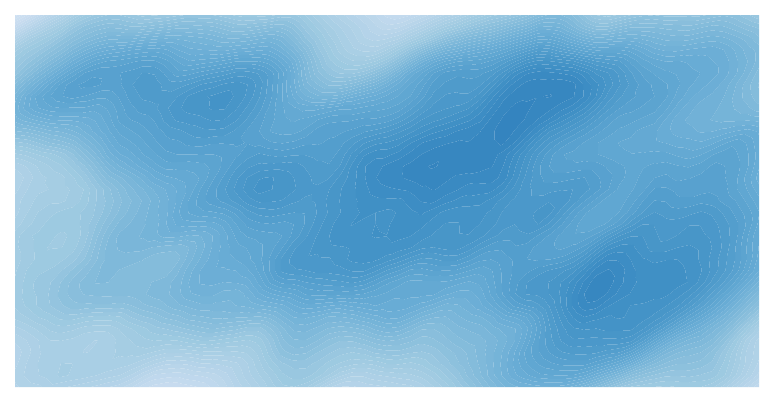} \\ \hline
     
\end{tabular}
            \renewcommand{\arraystretch}{1}
        \addtolength{\tabcolsep}{6pt}
    \caption{\textbf{Example II: progression of mean estimators with $N$.} 
    As with \cref{fig:ex1_progression_mean_estimator} we see a similar trend, where \texttt{LazyDINO} yields superior approximation than the other methods. Notably, \texttt{LazyDINO} achieves a reasonably accurate mean estimate with $250$ samples. \texttt{LazyNO} eventually also yields a faithful estimate but requires an order of magnitude more samples to do so. SBAI and \texttt{LazyMAP} both struggle substantially with this problem and yield large point-wise absolute errors, even for the largest amount of training samples utilized in their construction. 
    }
    \label{fig:ex_2_mean_estimator_progression}
\end{figure}

\begin{figure}
    \centering
    \addtolength{\tabcolsep}{-6pt}
    \renewcommand{\arraystretch}{0.35}

     Ground truth MAP\\
    \includegraphics[width=0.35\textwidth]{figures/hyper_map_bip3.pdf}
    
    \begin{tabular}{|>{\centering\arraybackslash}p{0.75cm}|c|c|c|c|}
    \hline
    \multicolumn{5}{|c|}{\bf \makecell{Progression of posterior MAP estimator (BIP \#4)}} \\
    \hline \hline
    
    & \multicolumn{3}{c|}{MAP estimates} & {\makecell[c]{point-wise \\ absolute error\\\includegraphics[width=0.2\textwidth]{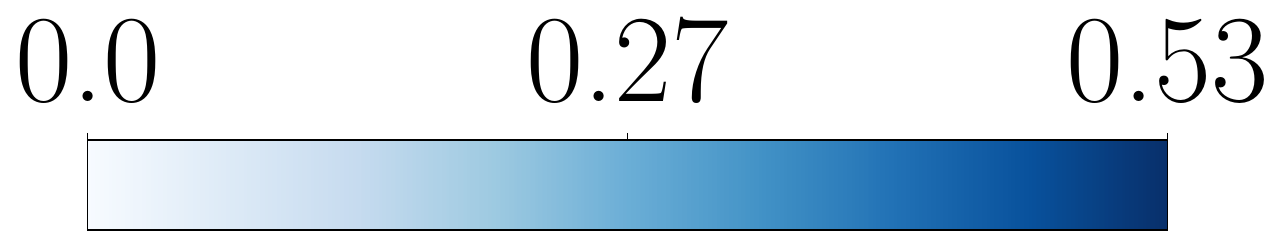}}}
    
    \\ \cline{1-5}
     $N$& 250 & 2k & 16k &  \makecell[c]{16k} \\ \hline
     
     \makebox[0pt][c]{\rotatebox[origin=l]{90}{\makebox[0.1\textwidth][c]{\texttt{LazyDINO}}}} 
     & 
     \includegraphics[width=0.2\textwidth]{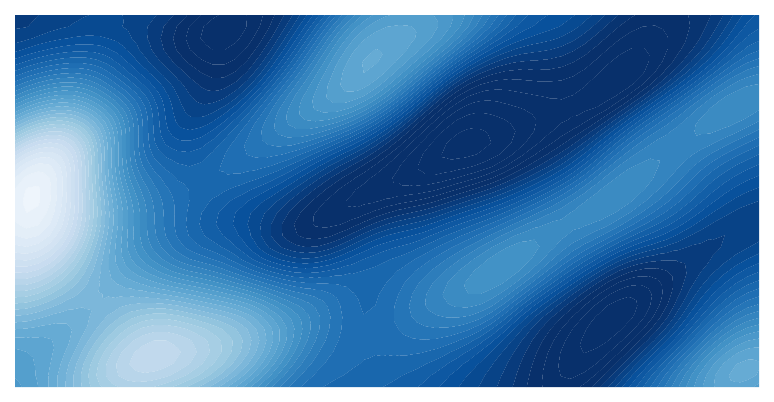} &
     \includegraphics[width=0.2\textwidth]{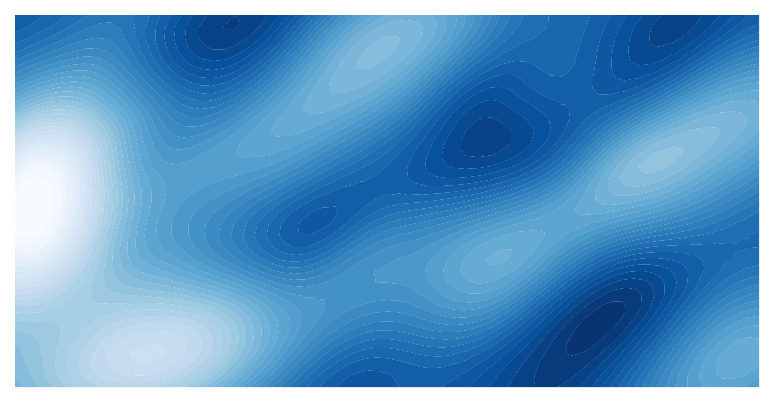} &
     \includegraphics[width=0.2\textwidth]{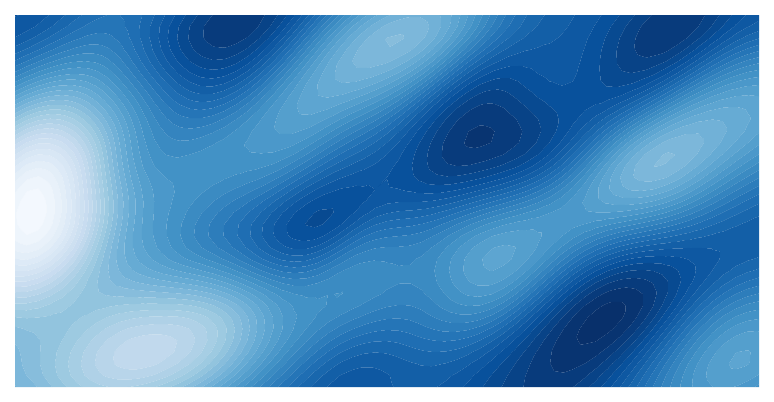} &
     \includegraphics[width=0.2\textwidth]{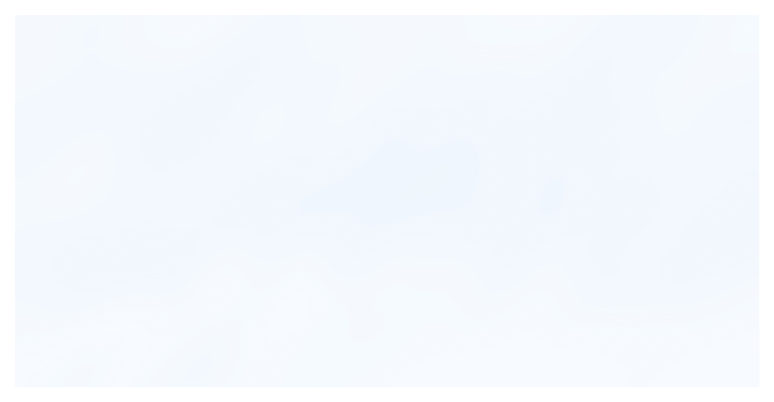} \\ \hline

   \makebox[0pt][c]{\rotatebox[origin=l]{90}{\makebox[0.1\textwidth][c]{\texttt{LazyNO}}}} 
     &
     \includegraphics[width=0.2\textwidth]{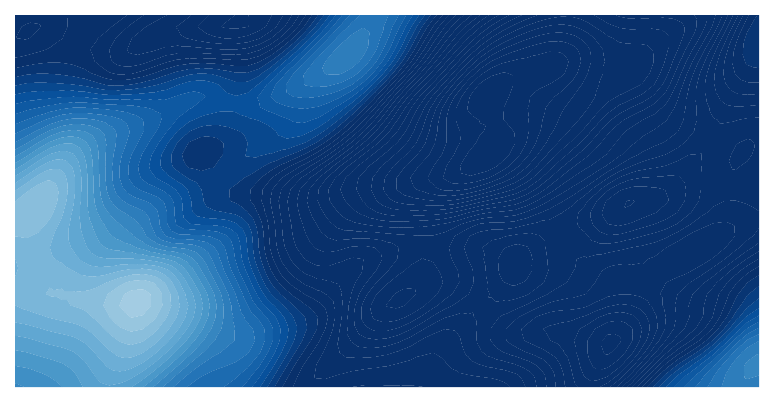} &
     \includegraphics[width=0.2\textwidth]{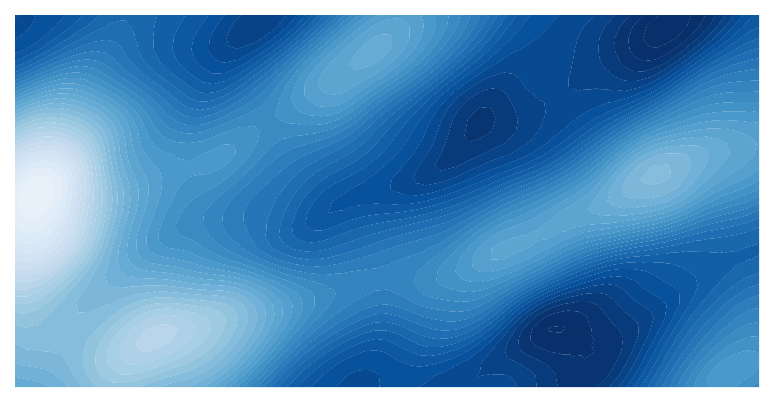} &
     \includegraphics[width=0.2\textwidth]{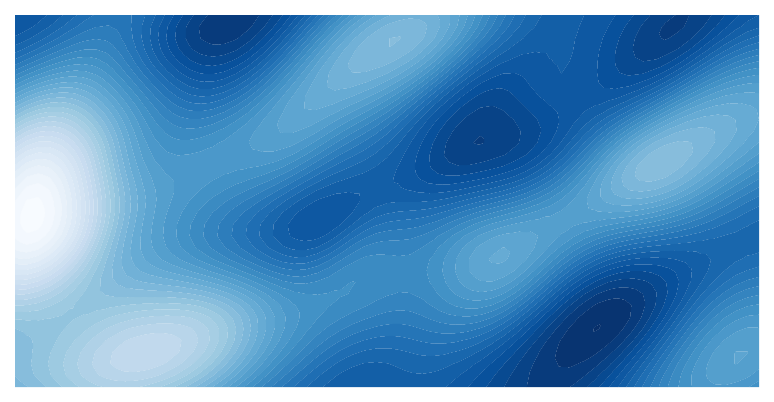} &
     \includegraphics[width=0.2\textwidth]{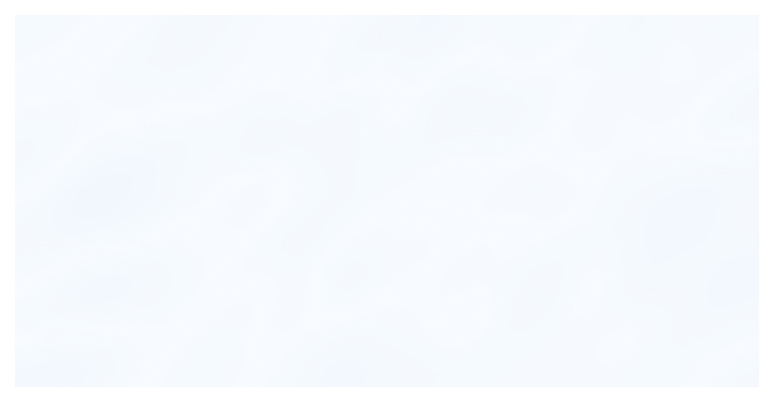} \\ \hline
     
    \makebox[0pt][c]{\rotatebox[origin=l]{90}{\makebox[0.1\textwidth][c]{SBAI}}} 
     &
     \includegraphics[width=0.2\textwidth]{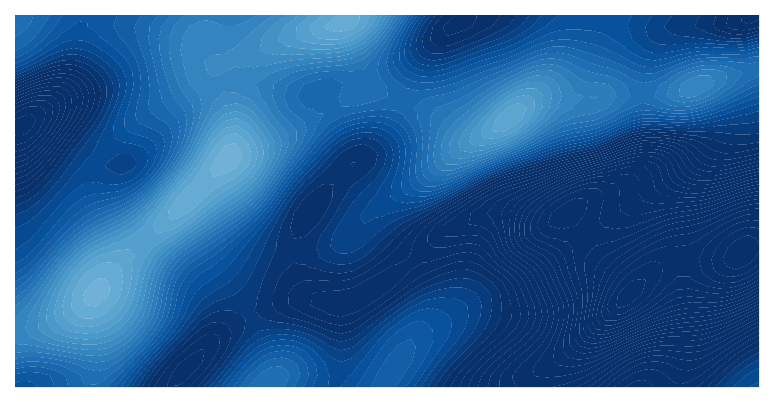} &
     \includegraphics[width=0.2\textwidth]{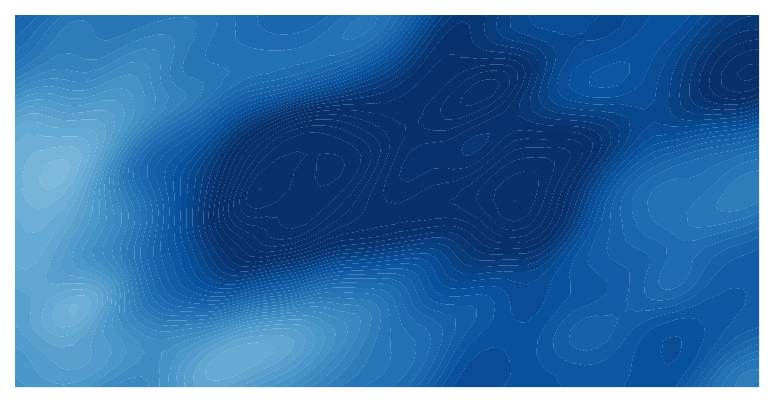} &
     \includegraphics[width=0.2\textwidth]{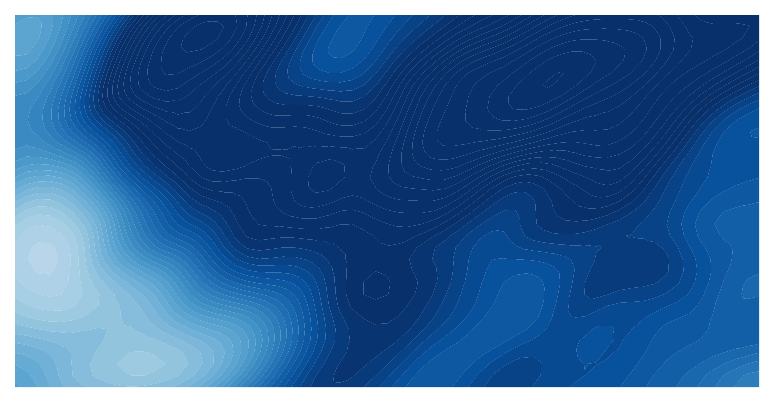} &
     \includegraphics[width=0.2\textwidth]{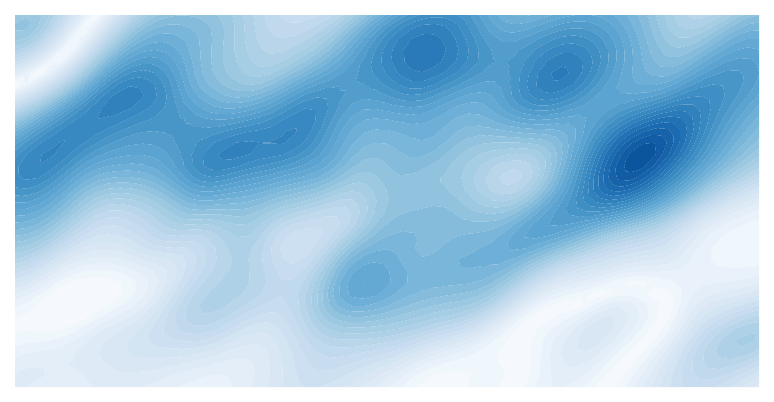} \\ 
     \hline
       {$N$}& 2k & 16k & 20k & \makecell[c]{20k}\\
       \hline
    \makebox[0pt][c]{\rotatebox[origin=l]{90}{\makebox[0.1\textwidth][c]{\texttt{LazyMap}}}} 
     &
     \includegraphics[width=0.2\textwidth]{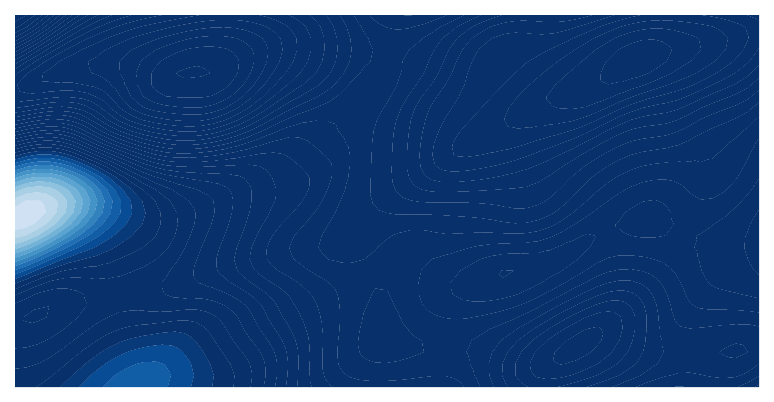} &
     \includegraphics[width=0.2\textwidth]{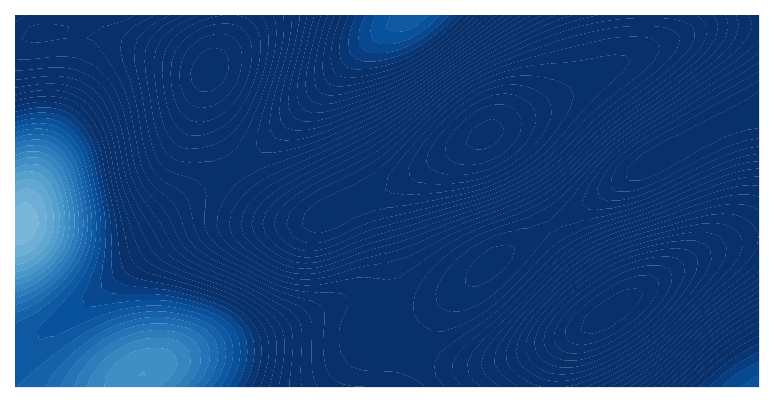} &
     \includegraphics[width=0.2\textwidth]{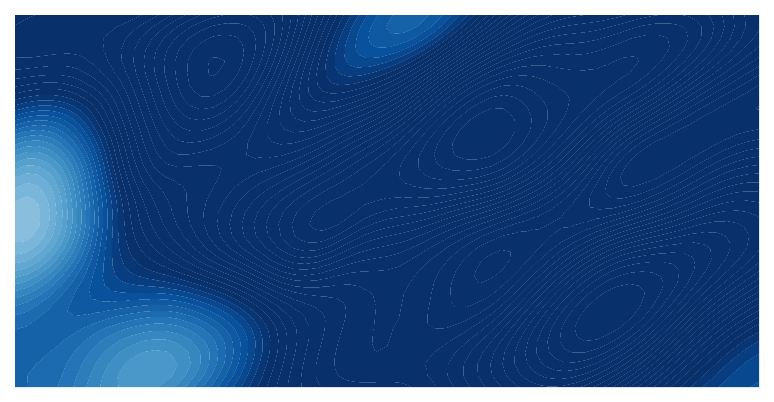} &
     \includegraphics[width=0.2\textwidth]{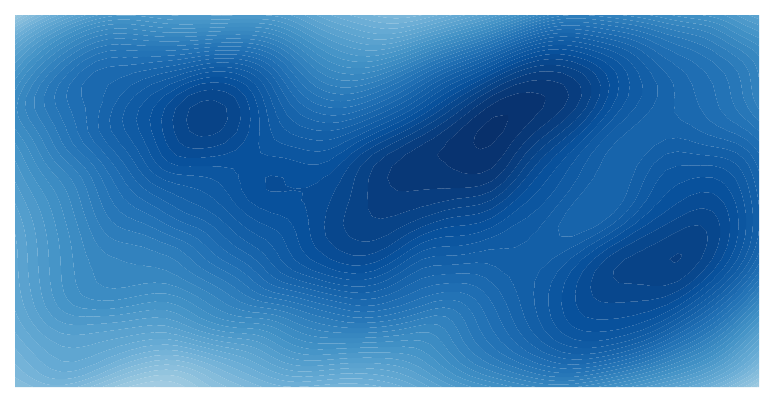} \\ \hline
     
\end{tabular}
\caption{\textbf{Example II: progression of MAP estimators.} As with \cref{fig:ex_2_mean_estimator_progression}, we see a similar story, where \texttt{LazyDINO} yields a superior approximation of the MAP point, particularly given fewer samples for its construction. Eventually, \texttt{LazyNO} yields a comparable approximation but requires an order of magnitude more samples to catch up to \texttt{LazyDINO}. Both SBAI and \texttt{LazyMap} yield poor reconstructions for the range of training samples considered in this study. 
}
    \label{fig:ex_2_map_progression}
\end{figure}

\begin{figure}[htbp]
    \centering
    \addtolength{\tabcolsep}{-6pt}
    \renewcommand{\arraystretch}{0.35}

     Ground truth point-wise\\marginal variance\\
    \includegraphics[width=0.35\textwidth]{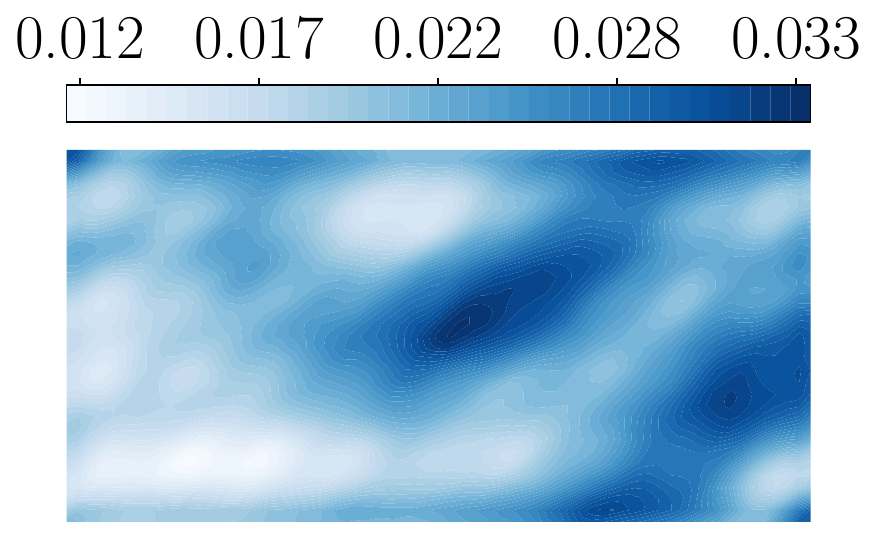}
    
    \begin{tabular}{|>{\centering\arraybackslash}p{0.75cm}|c|c|c|c|}
    \hline
    \multicolumn{5}{|c|}{\bf \makecell{Progression of posterior point-wise \\marginal variance (BIP \#4) }} \\
    \hline \hline
    
    & \multicolumn{3}{c|}{Point-wise marginal variance estimates} & {\makecell[c]{point-wise \\ absolute error\\\includegraphics[width=0.2\textwidth]{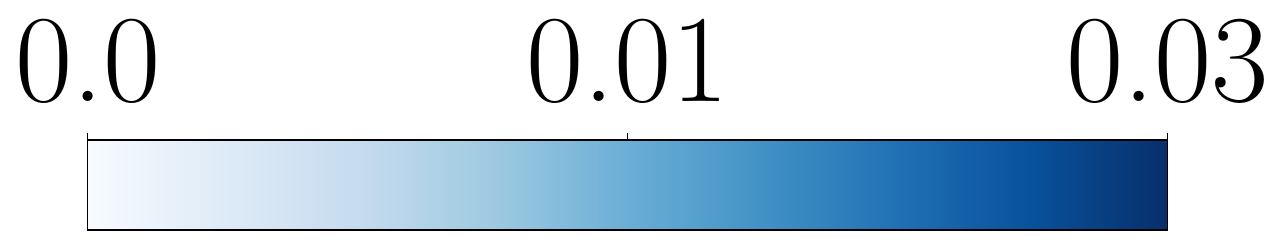}}}
    
    \\ \cline{1-5}
     $N$& 250 & 2k & 16k & \makecell[c]{16k}  \\ \hline
     
     \makebox[0pt][c]{\rotatebox[origin=l]{90}{\makebox[0.1\textwidth][c]{\texttt{LazyDINO}}}} 
     & 
     \includegraphics[width=0.2\textwidth]{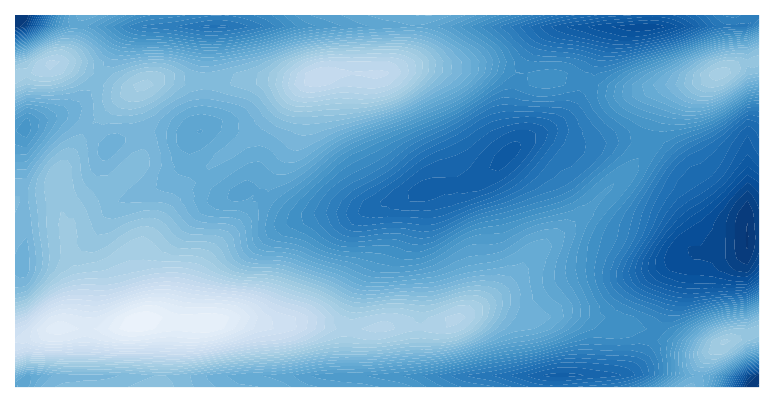} &
     \includegraphics[width=0.2\textwidth]{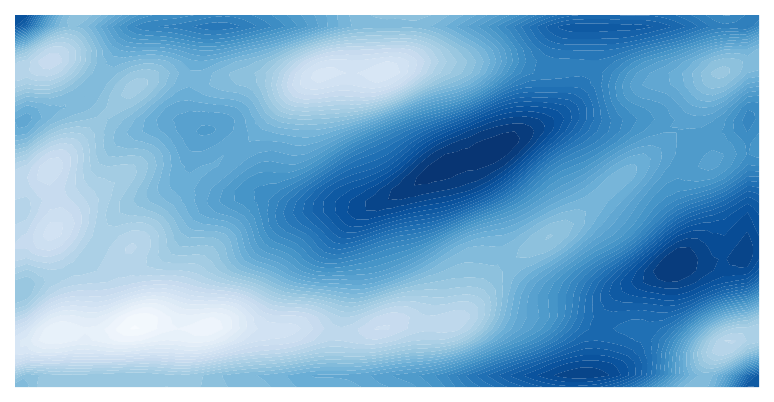} &
     \includegraphics[width=0.2\textwidth]{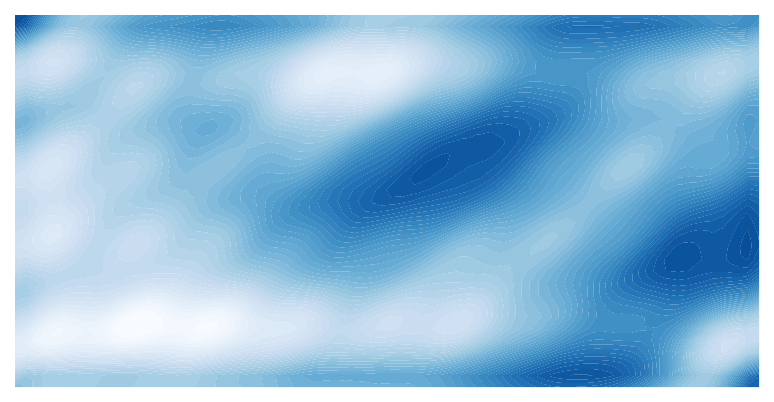} &
     \includegraphics[width=0.2\textwidth]{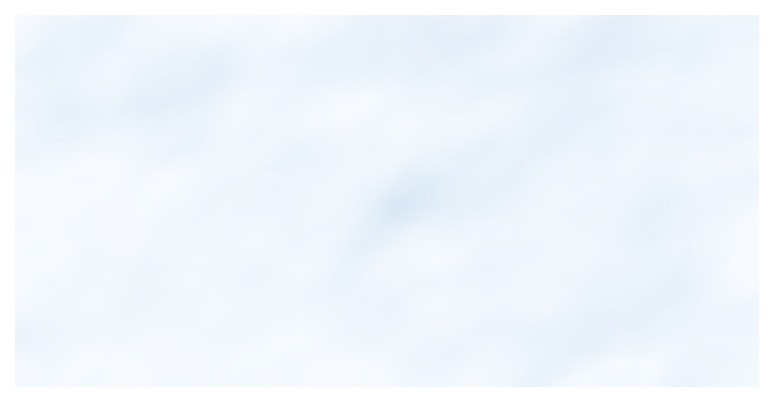} \\ \hline

   \makebox[0pt][c]{\rotatebox[origin=l]{90}{\makebox[0.1\textwidth][c]{\texttt{LazyNO}}}} 
     &
     \includegraphics[width=0.2\textwidth]{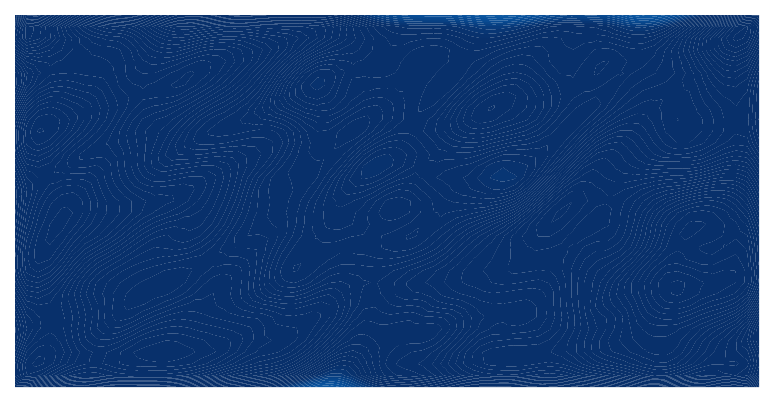} &
     \includegraphics[width=0.2\textwidth]{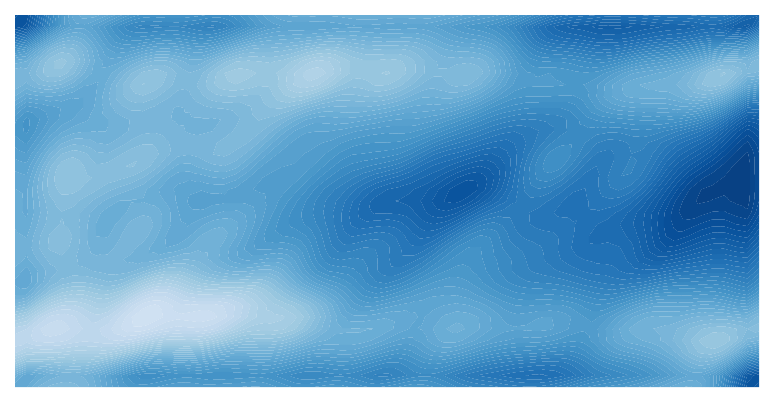} &
     \includegraphics[width=0.2\textwidth]{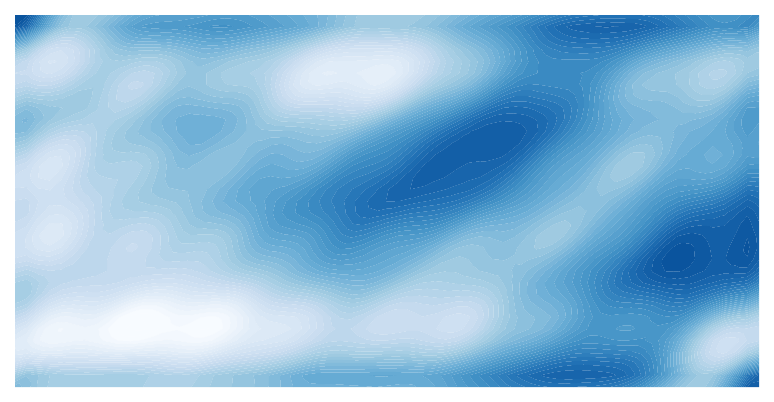} &
     \includegraphics[width=0.2\textwidth]{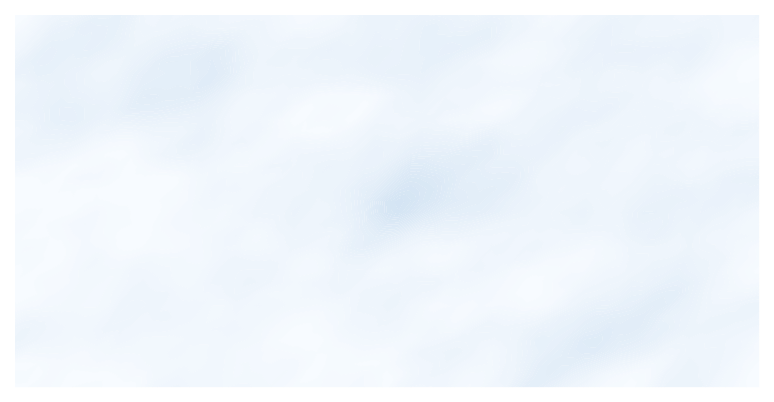} \\ \hline
     
    \makebox[0pt][c]{\rotatebox[origin=l]{90}{\makebox[0.1\textwidth][c]{SBAI}}} 
     &
     \includegraphics[width=0.2\textwidth]{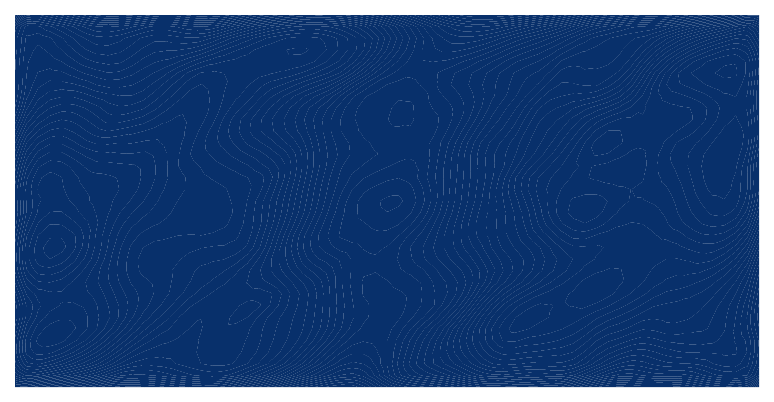} &
     \includegraphics[width=0.2\textwidth]{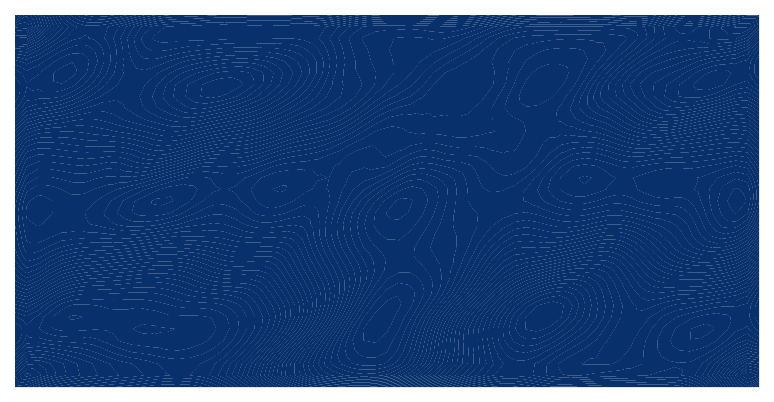} &
     \includegraphics[width=0.2\textwidth]{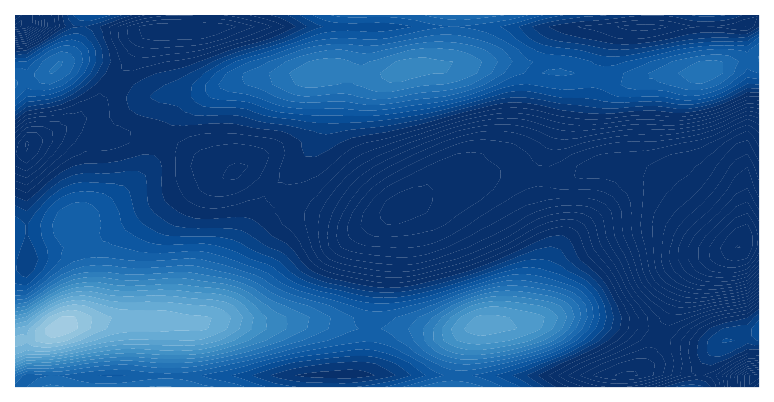} &
     \includegraphics[width=0.2\textwidth]{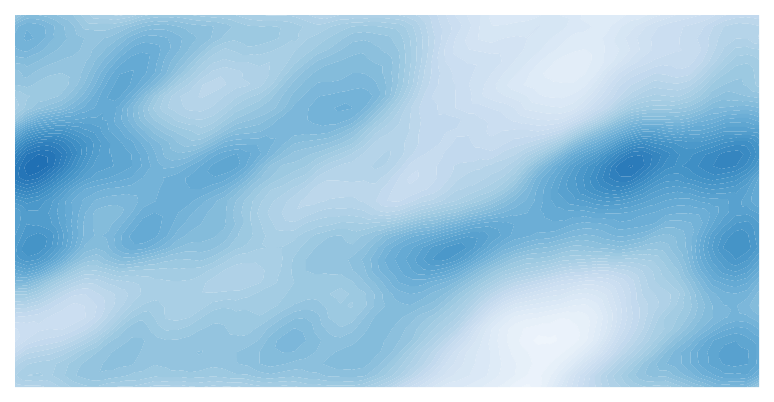} \\ 
     \hline
       {$N$}& 2k & 16k & 20k & \makecell[c]{20k}\\
       \hline
    \makebox[0pt][c]{\rotatebox[origin=l]{90}{\makebox[0.1\textwidth][c]{\texttt{LazyMap}}}} 
     &
     \includegraphics[width=0.2\textwidth]{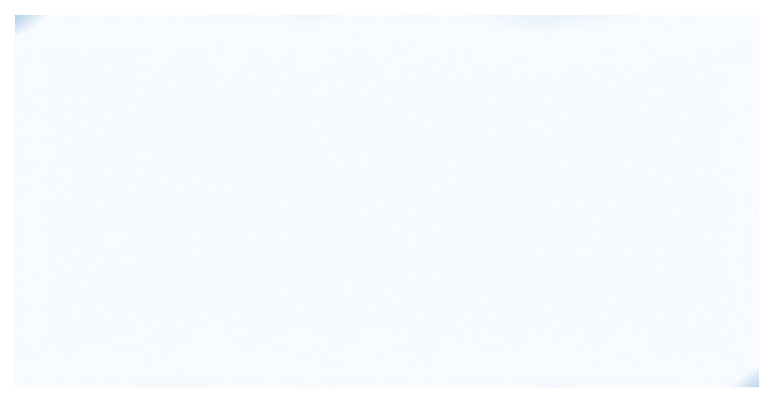} &
     \includegraphics[width=0.2\textwidth]{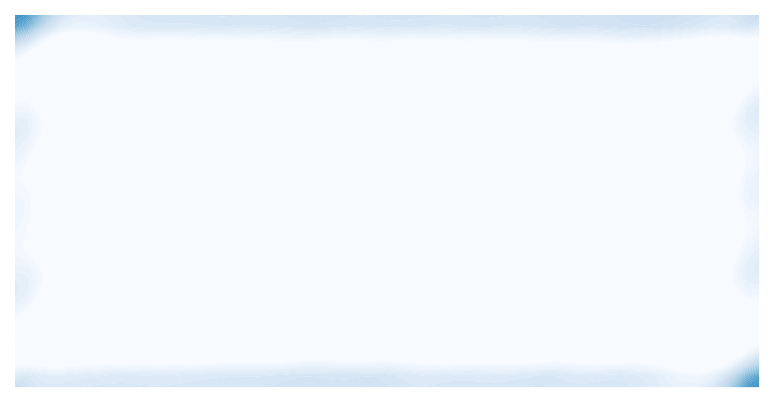} &
     \includegraphics[width=0.2\textwidth]{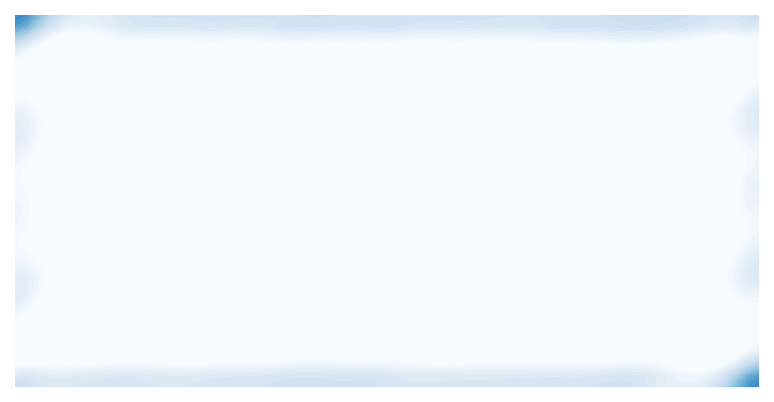} &
     \includegraphics[width=0.2\textwidth]{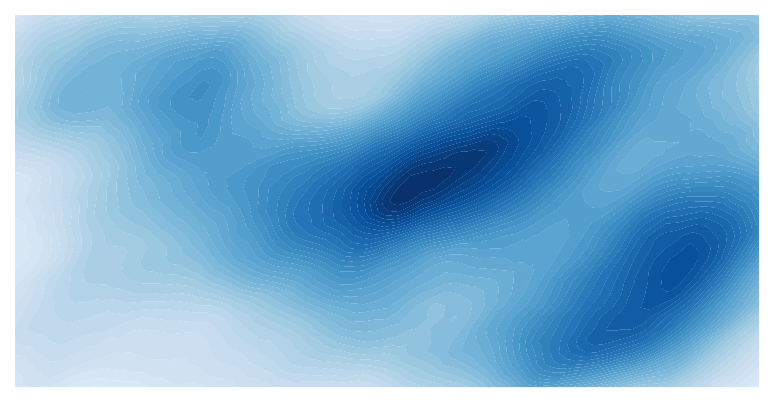} \\ \hline
     
\end{tabular}
\caption{\textbf{Example II: progression of point-wise marginal variance estimators.} 
We see, yet again, the same trend of the previous studies: \texttt{LazyDINO} yields a superior approximation of the point-wise marginal variance and, notably, can deliver a faithful approximation for as little as $250$ training data. \texttt{LazyNO} requires an order of magnitude more training data to catch up to the accuracy of \texttt{LazyDINO}, eventually yielding similar performance at $16,000$ training data, where the approximation capabilities of this fixed architecture may have saturated. SBAI and \texttt{LazyMap}, respectively, over and underestimate point-wise variance. This observation is consistent with the spreads seen in the marginal visualizations for SBAI and \texttt{LazyMap} in \cref{fig:hyper_marg3,fig:hyper_marg4}, respectively.
}
    \label{fig:ex_2_variance_progression}
\end{figure}

We proceed with the study of mean, MAP point, and point-wise marginal variance estimation for Example II in \cref{fig:ex_2_mean_estimator_progression,fig:ex_2_map_progression,fig:ex_2_variance_progression}. The overall story is similar to all preceding numerical studies: \texttt{LazyDINO} provides superior approximation to the other methods, and is notably able to give faithful approximations for limited samples, while the other methods require orders of magnitude more samples to achieve similar accuracy. . 

This concludes our extensive numerical comparison. In almost every point of comparison \texttt{LazyDINO} yielded the most accurate estimation of the posterior distribution as evidenced by moment discrepancies, density-based diagnostics, posterior marginals, mean, MAP and point-wise marginal variance estimations. Notably, \texttt{LazyDINO} gives faithful posterior estimates for orders of magnitude fewer samples than the alternative TVMI methods. While the LA-baseline did perform well in some metrics given limited samples (e.g., MAP estimate and covariance), this approximation assumes posterior Gaussianity. It leads to constant irreducible error, making it unviable for complex nonlinear BIPs. 

%% file: conclusion.tex
\section{Conclusion}\label{sec:conclusion}

In this work, we present \texttt{LazyDINO}, a fast, scalable, and efficiently amortized method for high-dimensional Bayesian inversion with expensive PtO maps. The method is composed of offline and online phases. During the offline phase, we generate join samples of the PtO map and its Jacobian to construct a \texttt{RB-DINO} surrogate of the PtO map via derivative-based dimension reduction and derivative-informed learning methods. During the online phase, when observational data is given, we seek rapid posterior approximation via surrogate-driven optimization of lazy maps, i.e., structure-exploiting transport maps with relatively low-dimensional nonlinearity. The trained lazy map is used for approximate posterior sampling and density estimation. 

We provide theoretical results demonstrating that the \texttt{RB-DINO} surrogate construction is optimized for amortized Bayesian inversion via lazy map variational inference. In \cref{theorem:kl_bound}, we show that the conventional supervised learning of the \texttt{DIPNet} surrogate architecture minimizes the upper bound on the expected error in posterior approximation when the ridge function surrogate replaces the PtO map. This architecture constricts the surrogate approximation to the parameter subspace that captures prior-to-posterior. In \cref{theorem:kl_gradient,corollary:optimality_gap}, We show that the derivative-informed learning of the surrogate minimizes the expected gradient error and optimality gap due to surrogate-driven transport map optimization.  This result reflects that the surrogate Jacobian accuracy affects the quality of the trained lazy map and thus directly influences the posterior approximation accuracy.

The \texttt{LazyDINO} method has several desirable traits. 
\begin{enumerate}
    \item \textit{Scalability}. The surrogate and transport map training in \texttt{LazyDINO} are independent of the parameter dimension as their latent representations reside in the same relatively low-dimensional derivative-informed subspace (\cref{fig:dino_construction,fig:lazy_dino_construction}). 
    \item \textit{Fast online inference}. Using a cheap-to-evaluate surrogate rKL objective for transport map optimization, \texttt{LazyDINO} fully exploits GPU-based
accelerations to rapidly approximate posteriors (\cref{tab:alg_timing_1,tab:alg_timing_2}). While our method requires solving an optimization problem to sample, we demonstrate that it leads to faster online posterior sampling than the typical inversion-to-sample approach of SBAI in the large sample size regime (\cref{table:sampling_SBAI}). 
    \item \textit{High posterior accuracy at low offline cost}. First, the RB-DINO surrogate and the lazy map are co-designed to exploit the structure of the BIP efficiently. Second, the derivative-informed learning method is highly cost-efficient and outperforms conventional supervised learning by one to two orders of magnitude (\cref{fig:gen_error_1,fig:gen_error2}). Consequently, \texttt{LazyDINO} requires a much smaller cost in offline computation to achieve high accuracy in online posterior approximation across multiple instances of observational.
\end{enumerate}

We studied two challenging infinite-dimensional PDE-constrained BIPs, each with four different instances of observational data: (i) inferring diffusivity field in a nonlinear reaction--diffusion PDE, and (ii) inferring Young's modulus field in a hyperelastic material thin film under deformation. In both cases, we observed one to two orders of magnitude of offline cost reduction for achieving similar accuracy in posterior approximation compared to alternative amortized inference methods such as \texttt{LazyNO} and SBAI via conditional transport. Moreover, LazyDINO consistently outperforms Laplace approximation at a small offline training sample regime (250--1,000), except for covariance approximation for Example I. In contrast, LazyNO and SBAI struggle to outperform Laplace approximation and, in some cases, failed at 16,000 offline training samples.

\texttt{LazyDINO} is a powerful method for settings requiring the repeated solution of BIPs defined by the same PtO map and prior. The efficiency gains achieved via \texttt{LazyDINO} motivate further study. First, \texttt{LazyDINO} can be applied to the case of posterior approximation for multiple independent observations simultaneously. Given the difficulty associated with concentration-of-measure for posteriors for many independent observations, the potential efficiency gains could be significant. Lastly, we will explore using \texttt{LazyDINO} in real-time uncertainty quantification for risk-averse decision-making, such as optimal experimental design and optimization under uncertainty for complex physical systems.

%% file: acknowledgment.tex
\section*{Ackowledgement}

This work was partially supported by the U.S. Department of Energy, Office of Science, Office of Advanced Scientific Computing Research under awards DE-SC0021239, DE-SC0023171, and DE-SC0023187. Thomas O'Leary-Roseberry, Joshua Chen, and Omar Ghattas are partially supported by the National Science Foundation under awards OAC-2313033 and DMS-234643, and the Air Force Office of Scientific Research under MURI grant FA9550-21-1-0084. The work of Lianghao Cao was partially supported by a Department of Defense Vannevar Bush Faculty Fellowship held by Andrew M. Stuart, and by the SciAI Center, funded by the Office of Naval Research (ONR), under Grant Number N00014-23-1-2729.

%% file: appendix_glossary.tex
\section{Glossary of terminology}

\begin{multicols}{2}
\begin{itemize}[leftmargin=*, label={}]

    \item \textbf{amortized cost}: Discounted cost, by spreading it out over the solution of additional problems. The more problems solved, the cheaper the amortized cost.

    \item \textbf{ANIS}: Auto-normalized importance sampling, also known as self-normalized importance sampling, ~\cite{agapiou2017importancesamplingintrinsicdimension}.

    \item \textbf{SBAI}: Simulation-based amortized inference~\cite{Ganguly_2023}.

    \item \textbf{BIP}: Bayesian inverse problem, also known as Bayesian inference problem. \emph{Inverse} problems are typically governed by physics.

    \item \textbf{ESS}: Effective sample size, an estimation of the number of independent samples drawn from a proposal distribution to obtain an expectation estimator with equivalence variance as standard Monte Carlo.

    \item \textbf{fKL}: forward Kullback-Leibler. The fKL divergence is the Kullback-Leibler divergence of the approximate distribution $\mu$ from the target distribution $\nu$, $\Dkl(\nu || \mu)$.

    \item \texttt{LazyDINO/NO}: The name of our algorithm, referring to a hybrid of \texttt{LazyMap} and the (DI)NO, (derivative-informed) neural operator trained surrogate PtO map.

    \item \texttt{LazyMap}: A transport map that acts only in a subspace of the input, by way of a linear projection-based dimension reduction~\cite{brennan2020greedyinferencestructureexploitinglazy}. It is also used to refer to the algorithm to train it.

    \item \textbf{LMVI}: Lazy map variational inference. Variational inference using structure-exploiting transport map with relatively low-dimensional non-linearity~\cite{brennan2020greedyinferencestructureexploitinglazy}.

    \item \textbf{MAP}: Maximum A-Posteriori, a/the MAP estimate is a/the point of highest probability concentration of a distribution.

    \item \textbf{MC}: Monte Carlo, which will always refer specifically to an i.i.d. sampling based approach to approximating an expectation.

    \item \textbf{MCMC}: Markov chain Monte Carlo.

    \item \textbf{PtO}: parameter-to-observable. A PtO map is a function taking the parameter we wish to infer to the non-noisy \emph{observable}. In contrast, \emph{observations} are measurements of the observable corrupted via noise.

    \item \textbf{rKL}: reverse Kullback-Leibler. The rKL divergence is the Kullback-Leibler divergence of the target distribution $\nu$ from an approximate distribution $\mu$, $\Dkl(\mu || \nu)$.

    \item \textbf{TMVI}: transport map variational inference, as in~\cite{marzouk2016introduction}. Though most modern variational inference is transport map variational inference, we use this term to \emph{distinguish} from other forms of variational inference, i.e., inference within families of distributions that do not originate via the transport of a reference distribution.
    
\end{itemize}
\end{multicols}

%% file: appendix_inf_dimensions.tex
\section{The single-sample estimate of the gradient of the optimal ridge function}\label{gradient_conditional_expectation_commute}
\begin{proof}Given $\mu = \mu_r \otimes \mu_{\perp}$, where $\mu_r = \Projector_{\sharp}\mu$ and $\mu_{\perp} = (\identity_{\ParamSpace}-\Projector)_{\sharp}\mu$, we reminder the reader that the projector is defined in terms of the $\ParamCM$-orthonormal reduced basis $\Psi_r=\{\psi_j\}_{j=1}^{\RedCoordDim}$, with $\text{span}(\Psi_r)=\text{Im}(\Projector).$ Using the compact notation 
$\mathbb{E}[ \Forward| \sigma( \Projector ) ](m_r) := \mathbb{E}_{m_{\perp}\sim \mu_{\perp}}\left[\Forward(m_r + m_{\perp})\right]$ for  $m_r \in \text{Im}(\Projector)$ for the conditional expectation, and noting that $\Decoder\refCoord \overset{d}{=}m_r \overset{d}{=}\Projector m$, and $\boldsymbol{V}^*\FullRidgeForward_{\text{opt}}(\Decoder \refCoord) = \RidgeForward_{\text{opt}}( \refCoord, \boldsymbol{w})$ we have that

\begin{align}
  &\mathbb{E}_{m_r\sim\mu_r}\Big[\left\|  \mathbb{E}[ \Forward| \sigma( P ) ](m_r) - \FullRidgeForward_{\boldsymbol{w}}(m_r)\right\|^2_{\NoiseCov^{-1}} + \left\|D_H \left(\mathbb{E}[ \Forward| \sigma( P ) ](m_r)\right) - D_{H} \FullRidgeForward_{\boldsymbol{w}}(m_r)\right\|^2_{\HS(\ParamCM,\ObsCM)} \Big]\\
     =&\mathbb{E}_{m\sim\mu}\Big[\left\|\Forward_{\text{opt}}(\Projector m) - \FullRidgeForward_{\boldsymbol{w}}(\Projector m)\right\|^2_{\NoiseCov^{-1}} + \left\|D_H \Forward_{\text{opt}}(\Projector m) - D_{H} (\FullRidgeForward_{\boldsymbol{w}}\circ\Projector)(m)\right\|^2_{\HS(\ParamCM,\ObsCM)} \Big]\\ 
\overset{d}{=}&\mathbb{E}_{\refCoord\sim\pi}\Big[\left\|\RidgeForward_{\text{opt}}( \refCoord)- \RidgeForward_{\boldsymbol{w}}(\refCoord)\right\|^2 + \left\| \mathbb{E}[ D_H \Forward| \sigma( P ) ](\Decoder \refCoord) - D_{H} \FullRidgeForward_{\boldsymbol{w}}(\Decoder \refCoord)\right\|^2_{\HS(\ParamCM,\ObsCM)} \Big]
\end{align}

For the last line, since $\Forward\in H^1_\mu(\ParamSpace; \ObsCM)$, the conditional expectation also belongs to the space, $\mathbb{E}[\Forward | \sigma(\Projector)]\in H^1_\mu(\ParamSpace; \ObsCM)$, by an isometry property (Proposition 1.2.8,~\cite{nualart2006malliavin}), and the usual commuting property between conditional expectation and differentiation apply---the Malliavin derivative of the conditional expectation with respect to the sigma algebra generated by the $\Projector$ is the conditional expectation of the Malliavin derivative. 
Then, $D_H \Forward_{\text{opt}}(\Projector m) := D_H  \mathbb{E}[\Forward | \sigma(\Projector)](\Projector m) = \mathbb{E}[D_H \Forward | \sigma(\Projector)](\Projector m)$ a.e. Lastly, we use a change-of-variables.  

In particular, we have the connection to the $\pi$-weighted Sobolev norm $H^1_{\pi}(\R^{\RedCoordDim}; \ObsSpace)$ objective: 
\begin{align}
&\mathbb{E}_{\refCoord\sim \pi} \left[\left\|\RidgeForwardOpt(\refCoord) - \RidgeForward_{\boldsymbol{w}}(\refCoord) \right\|^2 + \left\|\grad\RidgeForwardOpt(\refCoord) - \nabla_{\refCoord}\RidgeForward_{\boldsymbol{w}}(\refCoord) \right\|_F^2\right]\\
=&\mathbb{E}_{\refCoord\sim\pi}\Big[\left\|\boldsymbol{V}^*\mathbb{E}_{m_{\perp}\sim\mu_{\perp}}\left[\Forward(\Decoder\refCoord + m_{\perp})\right] - \RidgeForward_{\boldsymbol{w}}(\refCoord)\right\|^2
   + \big\| \boldsymbol{V}^* \mathbb{E}_{m_{\perp}\sim\mu_{\perp}}\left[D_H \Forward(\Decoder \refCoord + m_{\perp})\right]\circ \Decoder - \nabla_{\refCoord}\RidgeForward_{\boldsymbol{w}}(\refCoord)\big\|^2_{F} \Big]
\end{align}
where the second term follows from a chain rule, which produces $\Decoder$. Now, approximating this objective with a single sample $m_{\perp} \sim \mu_{\perp}$ for each $z \sim \pi$ produces our desired single-sample estimate result.
\end{proof}

\section{Transport map variational inference on Hilbert spaces \cite[section 6.6]{bogachev1998gaussian}}\label{inf_map}

Let us consider the rKL objective using nonlinear transformation of Gaussian measure $\mu=\mathcal{N}(0,\PriorCov)$ on a separable Hilbert space $\ParamSpace$. Let $\FullMap = \identity_{\ParamSpace} + \calK$ be the transport map where $\calK:\ParamSpace\to\ParamCM$ is nonlinear operator with a stochastic derivative $D_{H}\calK: \ParamSpace\to \text{HS}(\ParamCM)$ such that $D_{H}\FullMap\in \text{HS}(\ParamCM)$ is invertible $\mu$-a.e. Then we have
\begin{align*}
    \mathcal{D}_{\text{KL}}(\FullPrior||\FullMap^{\#}\FullPrior^{\Obs})=\int_{\ParamSpace} \log\left(\frac{\meas\mu}{\meas\left(\FullPrior^{\Obs}\circ\FullMap\right)}(\FullParam)\right)\, \meas\FullPrior(\FullParam)\,.
\end{align*}
Here, we derive the density between the pullback measure and the prior. Let $\scrA$ be any measurable subset of $\scrM$, we have
\begin{align*}
    (\mu^{\by}\circ\calT)(\scrA) &= \int_{\calT(\scrA)}\meas \mu^{\by}(\FullParam)
    =\frac{1}{Z^{\Obs}}\int_{\scrA}\exp\left(-(\Phi^{\by}\circ\calT)(\FullParam)\right) \dd(\mu\circ\calT)(\FullParam)
\end{align*}
The existence and the formula of the Radon--Nikodym derivative between $\FullPrior\circ\FullMap$ and $\mu$ is non-trivial; see \cite[section 6.6]{bogachev1998gaussian} for details. In the case where $D_{H}\calK(\FullParam)$ is a trace class operator on $\ParamCM$ $\mu$-a.e., we have 
\begin{equation*}
    \frac{\meas(\FullPrior\circ\FullMap)}{\meas\FullPrior}(\FullParam) = \text{det}_{\ParamCM}(D_H\FullMap)\exp(-\left\langle\calK(\FullParam), m\right\rangle_{\PriorCov^{-1}} -\frac{1}{2}\norm{\calK( m)}^2_{\PriorCov^{-1}}),
\end{equation*}
where the determinant is taken as the product of eigenvalues with $\ParamCM$-orthonormal eigenbases.

Therefore, the transport objective is given by
\begin{equation*}
    \mathcal{D}_{\text{KL}}(\FullPrior||\FullMap^{\#}\FullPrior^{\Obs})=\mathbb{E}_{{\FullParam}\sim\mu}\left[(\Misfit\circ \FullMap)(\FullParam) - \log\text{det}_{\ParamCM}(D_H\FullMap) + \left\langle\calK m, m\right\rangle_{\PriorCov^{-1}} +\frac{1}{2}\norm{\calK(\FullParam)}^2_{\PriorCov^{-1}} \right] + C_1
\end{equation*}

Let us consider the perturbation of the identity $\calK$ only acting on the latent space $\text{Im}(\calE_r)=\R^{d_r}$ of the projection $\mathcal{P}=\calD_r\circ\calE_r$ with $\calE_r\circ\calD_r = \identity_{\R^r}$. Then we have the following alternative definition of lazy map through $\calK$:
\begin{equation*}
    \calK = \underbrace{\calD_r\circ(\Map - \identity_{\RedCoordParamSpace})\circ\calE_r}_{\substack{\text{Perturbation of the identity} \\ \text{transport in } \text{Im}(\Projector)}},\quad \FullMap = \identity_{\ParamSpace} +\calK,
\end{equation*}
where $\Decoder$ and $\Encoder$ consists of $\ParamCM$-orthonormal reduced bases.
In this case, we have
\begin{align*}
    \log\text{det}(D_H\FullMap(\FullParam)) &\implies \log\text{det}(\nabla \Map(\Encoder m)),\\
    \left\langle\calK(\FullParam), m\right\rangle_{\PriorCov^{-1}}& \implies (\Encoder m)^{\transpose} \Map(\Encoder m) - \norm{\Encoder m}^2,\\
    \frac{1}{2}\norm{\calK(\FullParam)}^2_{\PriorCov^{-1}} &\implies \frac{1}{2}\norm{\Map(\Encoder m)}^2 + \frac{1}{2}\norm{\Encoder m}^2 - (\Encoder m)^{\transpose} \Map(\Encoder m).
\end{align*}
Therefore we have
\begin{equation}
\begin{aligned}
    \mathcal{D}_{\text{KL}}(\FullPrior||\FullMap^{\#}\FullPrior^{\Obs})&= \mathbb{E}_{{\FullParam}\sim\FullPrior}\left[(\Misfit\circ \FullMap)(\FullParam) - \log\text{det}(\grad \Map(\Encoder m)) + \frac{1}{2}\norm{\Map(\Encoder m)}^2\right] + C_2\\
    &= \mathbb{E}_{(\refCoord,m_{\perp})\sim\pi\otimes\mu_{\perp}}\left[(\Phi^{\by}\left((\calD_r\circ \Map\right)(\refCoord) + m_{\perp}) -\log\text{det}(\grad \Map(\refCoord)) + \frac{1}{2}\norm{ \Map(\refCoord)}^2\right] + C_2,
\end{aligned}
\end{equation}
where $\mu_{\perp} = (\calI-\mathcal{P})_{\#}\mu$ is the pushforward measure in the complimentary space. 

%% file: appendix_proof.tex
\section{Proofs of \cref{theorem:kl_bound,theorem:kl_gradient} and \cref{corollary:optimality_gap} }
\subsection{Proof of \cref{theorem:kl_bound}}\label{app:kl_bound_proof}
\begin{proof}
Proposition 4.1 by Cui and Zahm \cite{cui2021data} gives us the following equality:
\begin{equation*}
    \mathbb{E}_{\Obs\sim\gamma} \left[ \KL( \FullPost||\widetilde{\FullPrior}^{\Obs})\right] = \frac{1}{2}\mathbb{E}_{\FullParam\sim\FullPrior}\left[\norm{\Forward(m) - \FullRidgeForward(\Projector m)}^2_{\NoiseCov^{-1}}\right] - \KL(\gamma||\widetilde{\gamma}).
\end{equation*}
Using triangle and Cauchy--Schwarz inequalities, we have 
\begin{equation*}
        \frac{1}{2}\mathbb{E}_{\FullParam\sim\FullPrior}\left[\norm{\Forward(m) - \FullRidgeForward(\Projector m)}^2_{\NoiseCov^{-1}}\right]\leq \mathbb{E}_{\FullParam\sim\FullPrior}\left[\norm{\Forward(m) - \FullRidgeForwardOpt(\Projector m)}^2_{\NoiseCov^{-1}}\right] + \mathbb{E}_{\refCoord\sim\pi}\left[\norm{ \RidgeForwardOpt(\refCoord)-\RidgeForward(\refCoord)}^2\right].
\end{equation*}
Furthermore, Proposition 7 in Cao et al.~\cite{cao2024efficient} gives the following upper bound:
\begin{equation*}
    \mathbb{E}_{\FullParam\sim\FullPrior}\left[\norm{\Forward(m) - \FullRidgeForwardOpt(\Projector m)}^2_{\NoiseCov^{-1}}\right] \leq \Tr_{\ParamCM}\left(\left(\identity_{\ParamCM}-\Projector\right)\mathcal{H}_A\left(\identity_{\ParamCM}-\Projector\right)\right),
\end{equation*}
which completes the proof.
\end{proof}
\subsection{Proof}\label{app:kl_gradient_proof}
\begin{proof}[Proof of \cref{theorem:kl_gradient}]
First, since the density between ${\Map_{\MapParams}}_{\sharp}\pi$ and $\pi$ is essentially bounded we have the following bound for any $f\in L^1(\pi)$ due to a change-of-variables formula and the H\"older inequality: 
\begin{align*}
\mathbb{E}_{\boldsymbol{x}\sim{\Map_{\MapParams}}_{\sharp}\pi}\left[|f(\boldsymbol{x})|\right] &\leq C_1\mathbb{E}_{\refCoord\sim\pi}\left[|f(\refCoord)|\right],
\end{align*}
where $C_1$ is the essential supremum of the density.
Next, we have
    \begin{align*}
        \nabla_{\MapParams}\mathcal{L}_1^{\Obs}(\refCoord, \MapParams) - \nabla_{\MapParams}\widetilde{\mathcal{L}}^{\Obs}_1(\refCoord, \MapParams) &=\nabla_{\MapParams} \Map_{\MapParams}(\refCoord)^{\transpose}(\nabla\RidgeForwardOpt\circ\Map_{\MapParams})(\refCoord)^{\transpose}\left((\RidgeForwardOpt\circ\Map_{\MapParams})(\refCoord)-\boldsymbol{V}^*\Obs\right)\\
        &\qquad- \nabla_{\MapParams} \Map_{\MapParams}(\refCoord)^{\transpose}(\nabla\RidgeForward\circ\Map_{\MapParams})(\refCoord)^{\transpose}\left((\RidgeForward\circ\Map_{\MapParams})(\refCoord)-\boldsymbol{V}^*\Obs\right)\\
        & = \nabla_{\MapParams} \Map_{\MapParams}(\refCoord)^{\transpose}\left((\nabla\RidgeForwardOpt\circ\Map_{\MapParams})(\refCoord) - (\nabla\RidgeForward\circ\Map_{\MapParams})(\refCoord)\right)^{\transpose}\left((\RidgeForwardOpt\circ\Map_{\MapParams})(\refCoord)-\boldsymbol{V}^*\Obs\right)\\
        &\qquad + \nabla_{\MapParams} \Map_{\MapParams}(\refCoord)^{\transpose}(\nabla\RidgeForward\circ\Map_{\MapParams})(\refCoord)^{\transpose}\left((\RidgeForwardOpt\circ\Map_{\MapParams})(\refCoord)-(\RidgeForward\circ\Map_{\MapParams})(\refCoord)\right).
    \end{align*}
By Jensen's, triangle, and H\"older's inequalities, we have
\begin{align*}
     \left(\mathbb{E}_{\Obs\sim\gamma}\left[\norm{\nabla_{\MapParams}\mathcal{L}^{\Obs}(\MapParams) - \nabla_{\MapParams}\widetilde{\mathcal{L}}^{\Obs}(\MapParams)}\right]\right)^2&\leq\left(\mathbb{E}_{(\Obs,\refCoord)\sim\gamma\otimes\pi}\left[ \norm{\nabla_{\MapParams}\mathcal{L}_1^{\Obs}(\refCoord, \MapParams) - \nabla_{\MapParams}\widetilde{\mathcal{L}}^{\Obs}_1(\refCoord,\MapParams)}\right]\right)^2\\
    &\leq \max\{C_2, C_3\}\mathbb{E}_{\boldsymbol{x}\sim{\Map_{\MapParams}}_{\sharp}\pi}\left[\norm{\RidgeForwardOpt(\boldsymbol{x}) - \RidgeForward(\boldsymbol{x})}^2 + \norm{\nabla\RidgeForwardOpt(\boldsymbol{x}) - \nabla\RidgeForward(\boldsymbol{x})}_{F}^2\right],
\end{align*}
where the constants are given by
\begin{align*}
    C_2 &=\mathbb{E}_{\refCoord\sim\pi}\left[\norm{(\nabla\RidgeForward\circ\Map_{\MapParams})(\refCoord)\partial_{\MapParams} \Map_{\MapParams}(\refCoord)}^2\right],\\
    C_3 &=\mathbb{E}_{(\Obs,\refCoord)\sim\gamma\otimes\pi}\left[\norm{\nabla_{\MapParams} \Map_{\MapParams}(\refCoord)^{\transpose}\left((\RidgeForwardOpt\circ\Map_{\MapParams})(\refCoord)-\boldsymbol{V}^*\Obs\right)}^2\right].
\end{align*}
Since $\Forward,\FullRidgeForward\circ\Projector\in H^1_{\FullPrior}(\ParamSpace;\ObsCM)$, we have $\RidgeForwardOpt,\RidgeForward\in H^1_{\pi}(\R^{d_r};\R^{\ObsDim})$. Moreover, since $C_4=\text{ess sup}_{\refCoord\in\R^{d_r}}\|\nabla_{\MapParams}\Map(\refCoord)\|<\infty$, we have
\begin{equation*}
    C_2 \leq C_1C_4\mathbb{E}_{\refCoord\sim\pi}\left[\norm{\nabla\RidgeForward(\refCoord)}^2\right]<\infty,\quad
    C_3\leq C_1C_4\mathbb{E}_{(\Obs,\refCoord)\sim\gamma\otimes\pi}\left[\norm{\RidgeForwardOpt(\refCoord)-\boldsymbol{V}^*\Obs}^2\right]<\infty.
\end{equation*}
Therefore, we have
\begin{align*}
\left(\mathbb{E}_{\Obs\sim\gamma}\left[\norm{\nabla_{\MapParams}\mathcal{L}^{\Obs}(\MapParams) - \nabla_{\MapParams}\widetilde{\mathcal{L}}^{\Obs}(\MapParams)}\right]\right)^2&\leq \max\{C_2, C_3\}C_1C_4\\
&\quad\times\left(\mathbb{E}_{\boldsymbol{z}\sim\pi}\left[\norm{\RidgeForwardOpt(\boldsymbol{z}) - \RidgeForward(\boldsymbol{z})}^2 + \norm{\nabla\RidgeForwardOpt(\boldsymbol{z}) - \nabla\RidgeForward(\boldsymbol{z})}_{F}^2\right]\right).
\end{align*}
\end{proof}

\begin{proof}[Proof of \cref{corollary:optimality_gap}]
By the Polyak--Lojasiewicz inequality and the results in Part I, we have
\begin{align*}
    \mathbb{E}_{\Obs\sim\gamma}\left[\sqrt{\mathcal{L}^{\Obs}(\widetilde{\MapParams}^{\by,\dagger}) - \mathcal{L}^{\Obs}(\MapParams^{\by,\dagger})}\right]&\leq \mathbb{E}_{\Obs\sim\gamma}\left[\frac{1}{\sqrt{2\lambda^{\Obs}}}\norm{\nabla_{\MapParams}\mathcal{L}(\widetilde{\MapParams}^{\by,\dagger})}\right]\\
    &\leq \frac{1}{C_5}\max\{C_2, C_3\}C_1C_4  \\
    &\qquad\times\left(\mathbb{E}_{\boldsymbol{z}\sim\pi}\left[\norm{\RidgeForwardOpt(\boldsymbol{z}) - \RidgeForward(\boldsymbol{z})}^2 + \norm{\nabla\RidgeForwardOpt(\boldsymbol{z}) - \nabla\RidgeForward(\boldsymbol{z})}_{F}^2\right]\right)^{1/2},
\end{align*}
where $C_5 = \text{ess inf}_{\Obs\sim\gamma}\sqrt{2\lambda^{\Obs}}>0$.
\end{proof}

%% file: lazydino_detailed_appendix.tex
\section{Detailed Discussion of LazyDINO algorithm}
\label{appendix:detailed_lazydino_algorithm}

In this section, we describe algorithms for performing inference with a \texttt{LazyDINO} algorithm, comprised of an offline and two online phases. In~\ref{dino_construction}, we describe the \emph{offline} surrogate construction phase, that is, computing the encoder and decoder for the parameter latent space and then training \texttt{RB-DINO} in that latent space. Next, in~\ref{transport_map_training}, we describe the \emph{online} solving of a BIP for a given observed data vector $\Obs$.

\subsection{Offline phase: surrogate construction}\label{dino_construction}
\paragraph{Latent space: solving the generalized eigenvalue problem}\label{computing_ghep}
\SetKwComment{Comment}{$\triangleright$\ }{}

\begin{algorithm}[ht]
  \SetAlgoLined
  \DontPrintSemicolon
  \KwIn{
    \begin{itemize}
          \vspace{-0.25 cm}
      \item[(i)] prior distribution sampler: $\FullPrior$,  prior precision operator: $\PriorCov^{-1}$
      
      \vspace{-0.25 cm}
      \item[(ii)] noise precision matrix: $\NoiseCov^{-1}$, observable basis $\boldsymbol{V}$
          
          \vspace{-0.25 cm}
      \item[(iii)]  PtO map: $\FullParam \mapsto \Forward(\FullParam)$, Jacobian action: $D\Forward(\FullParam)$
          
          \vspace{-0.25 cm}
       \item[(iv)]   \# desired training dataset samples: $N \in \mathbb{N}$
          
          \vspace{-0.25 cm}
     \item[(v)] \# samples to compute encoder/decoder: $N_{L} \leq N$ 
     
     \vspace{-0.25 cm}
    \item[(vi)] embedding dimension: $\RedCoordDim$ or eigenvalue tail sum tolerance: $\epsilon_{L}$
          
          \vspace{-0.25 cm}
    \end{itemize}
  }
      \KwOut{
      \begin{itemize}
              
              \vspace{-0.25 cm}
      \item[(i)] encoder/decoder pair: $\Encoder, \Decoder$
              
              \vspace{-0.25 cm}
      \item[(ii)]  latent  space training dataset inputs: $\{\refCoord^{(j)}\}$, outputs: $\{\boldsymbol{g}^{(j)}\}$, $\{ \bJ_r^{(j)}\},\; j = 1,\ldots,N$
    
    \vspace{-0.25 cm}
    \end{itemize}
    }
  \Begin{
        \texttt{1.} $\FullParam^{(j)} \sim \FullPrior, \quad j = 1, \ldots, N$ \Comment*[r]{Sample prior}
         \texttt{2.} $\Forward(\FullParam^{(j)}), D_H \Forward(\FullParam^{(j)}),\quad j = 1, \ldots, N$ \Comment*[r]{Evaluate PtO map/Jacobian}
        \texttt{3.}  \texttt{Create encoder/decoder}:\\
         $\quad\;\;\{\psi_k\in\ParamSpace\}_{k=1}^{\RedCoordDim} \gets \texttt{eigenvalue\_problem} (\big\{D \Forward(\FullParam^{(j)})\big\}_{i=1}^{N_L},\NoiseCov^{-1}, \PriorCov^{-1},\epsilon_{L} \textrm{ or } d_r)$\Comment*[r]{~\cref{eq:derivative_gevp},~\cref{eq:approx_H}}
        $\quad\;\;\Decoder \refCoord:= \sum_{k=1}^{\RedCoordDim} \refCoord_k\psi_k,\; \Encoder := \Decoder^{\top}\PriorCov^{-1}$ \Comment*[r]{~\cref{eq:generic_encoder_decoder_definition},~\cref{eq:as_encoder_decoder_def}}
    \texttt{4.} \texttt{Embed dataset}: \Comment*[r]{~\cref{eq:latent_parameter},~\cref{eq:reduced_jacobian}}
    $\quad\;\;\refCoord^{(j)} \gets \Encoder \FullParam^{(j)},\;\RidgeForwardOptApprox^{(j)}\gets\boldsymbol{V}^{\top}\NoiseCov^{-1}\Forward(m^{(j)}),\;\bJ_r^{(j)} \gets \boldsymbol{V}^{\top}\NoiseCov^{-1}D \Forward(m^{(j)})\Decoder ,\quad j = 1,\ldots, N$
  }
    \caption{\texttt{LazyDINO}: define latent space and embed dataset}
    \label{alg:data_and_basis_generation_summary}
  \end{algorithm}

We describe here the computation of the \\\texttt{eigenvalue\_problem} in~\cref{alg:data_and_basis_generation_summary} to find the reduced basis $\Psi_{r}$ for our encoder and decoder. Motivated by \cref{theorem:kl_bound}, we take $\Psi_{r}$ to be composed of the leading $\ParamCM$-orthonormal eigenbasis functions of an MC approximation of the generalized eigenvalue problem defined in ~\cref{eq:derivative_gevp}:
\begin{align}\label{eq:approx_H}
    \PriorCov^{-1}\mathcal{H}_A \approx \frac{1}{N_{L}}\sum_{j=1}^{N_{L}} D \Forward(m^{(j)})^*\NoiseCov^{-1}D\Forward(m^{(j)}) ,\quad m^{(j)}\iid\FullPrior.
\end{align} 

Following~\cref{eq:as_ev_tail_sum}, the latent parameter space dimension $d_r \le \textrm{dim}(\ParamSpace)$ is chosen to capture the dominant information in $\mathcal{H_A}$; specifically, it is desirable to ensure that the eigenvalue tail sum is small. 
For many high- to infinite-dimensional BIPs, $d_r$ is expected to be small due to, e.g., Saint-Venant's principle for coercive elliptic PDEs, concentration of measure, and the low-rankness of sparse observations extracted from a PDE state. 
  
When the discretization dimension of the problem is large, the generalized eigenvalue decomposition must be computed matrix-free. To this end there are many suitable computational tools such as randomized methods \cite{halko2010findingstructurerandomnessprobabilistic, xiang2014randomizedalgorithmslargescaleinverse, saibaba2015randomizedalgorithmsgeneralizedhermitian} and Krylov subspace methods \cite{GolubGEP, SaadKrylov, sorensenQZGEP,doi:10.1137/S0895479802403459,CHOWDHURY1976439}.
In this case, since the computation of the eigenvalue tail is intractable, one can adaptively find a sufficiently large dimension, $d_r$, such that the eigenvalues have decayed sufficiently, i.e., $\lambda_{d_r}/ \lambda_1 $ is small.

The samples needed to compute the reduced basis can be reused as part of the training data set; far fewer samples are usually needed to compute the reduced basis than are needed to train \texttt{RB-DINO} to required low error tolerances, see e.g.~\cite{cao2024efficient}. 
In this case, the additional training sample latent Jacobians can be directly formed via its action or adjoint action, as in~\cref{eq:reduced_jacobian}. Specifically, once the PtO evaluation at $m^{(j)}$ is available, only $\text{min}(\ObsDim,\RedCoordDim)$ evaluations of the PtO map derivative or its adjoint action are needed for a latent Jacobian evaluation. The evaluation cost can often be reduced to a fraction of the PtO map evaluation cost for linear or highly nonlinear PDEs. For details on efficient means to form latent Jacobian matrices for PDE-constrained PtO maps, see \ref{section:bayes_pdes} and \cite[Section 4.3]{cao2024efficient}. Empirical evidence of the low relative cost of Jacobians can be seen in our numerical results in~\cref{amortized_actual_times} in~\cref{tab:alg_timing_1} and~\cref{tab:alg_timing_2}.

Using the encoder and decoder, defined previously in terms of the reduced basis $\Psi_{r}$ and prior precision $\PriorCov^{-1}$ in~\cref{eq:generic_encoder_decoder_definition} and \cref{eq:as_encoder_decoder_def}, we embed the training data into the latent space, resulting in the whitened latent inputs $\refCoord^{(j)}$, whitened PtO samples $\boldsymbol{g}^{(j)}$, and whitened latent Jacobian samples $\bJ_r^{(j)}$. A summary of these procedures is given in~\cref{alg:data_and_basis_generation_summary}.

\paragraph{\textit{RB-DINO} training}
Next, in~\cref{alg:dino_training}, we train \texttt{RB-DINO} using the embedded data set. This involves a straightforward empirical risk minimization arising from MC estimate of either~\cref{eq:l2_learning_latent} or~\cref{eq:h1_training_equations}. Any method for stochastic unconstrained optimization, e.g., stochastic gradient descent, Adam~\cite{kingma2017adammethodstochasticoptimization}, or second-order methods \cite{yao2021adahessian,o2024fast} can be used. If the conventional $L^2_\FullPrior$ empirical risk is employed, we refer to the surrogate as \texttt{RB-NO} (neural operator) instead of \texttt{RB-DINO} (derivative-informed neural operator).

\begin{algorithm}[h]
\DontPrintSemicolon
  \caption{\texttt{LazyDINO}: train reduced basis neural operator in latent space $\RedCoordParamSpace$}\label{alg:dino_training}
  \SetAlgoLined

  \KwIn{
    \begin{itemize}
       \vspace{-0.25 cm}\item[(i)] training dataset inputs: $\big\{\refCoord^{(j)}\big\}$, outputs:
     $\big\{\boldsymbol{g}^{(j)}\big\}, \big\{\bJ_r^{(j)} \big\},\; j = 1,\ldots,N$
    
    \vspace{-0.25 cm}\item[(ii)] untrained neural network: $\RidgeForward_w:\RedCoordParamSpace\times \mathbb{R}^{d_W} \rightarrow \ObsSpace$
    
    \vspace{-0.25 cm}\item[(iii)] choice of conventional $L^2_\FullPrior$ (\texttt{RB-NO}) or derivative-informed $H^1_{\FullPrior}$ (\texttt{RB-DINO}) objective
      \vspace{-0.25 cm}
    \end{itemize}
  }
  \KwOut{
  \begin{itemize}
          \vspace{-0.25 cm}
  \item[(i)] trained neural network: $\RidgeForward_{\boldsymbol{w}^*}$
  \vspace{-0.25 cm}
\end{itemize}
}
  \Begin{
    \texttt{1.} Train $\RidgeForward_{w}$ by minimizing an empirical risk:\\
        $\boldsymbol{w}^* = \underset{\boldsymbol{w} \in  \mathbb{R}^{d_W}}{\mathrm{argmin}} \frac{1}{N}\sum_{j=1}^N \bigg(\left\|\boldsymbol{g}^{(j)} - \RidgeForward_{\boldsymbol{w}} \left(\refCoord^{(j)}\right)\right\|^2
        + \; \underbrace{\left\|\boldsymbol{J}_r^{(j)} - \nabla_{\refCoord}\RidgeForward_{\boldsymbol{w}} \left(\refCoord^{(j)}\right)\right\|^2_F}_{\text{include for }H^1_{\FullPrior} \enskip(\texttt{RB-DINO})\text{ objective}}\bigg)$ 
  }
\end{algorithm}
 Equipped with a sufficiently accurate neural network approximation to the optimal latent PtO map $\RidgeForwardOpt$, including accurate approximations of its derivatives, one can proceed to transport map variational inference.
 
\subsection{Online phase: lazy map variational inference with surrogate latent objective function}~\label{transport_map_training}

We perform LMVI using a stochastic approximation of the surrogate rKL objective and its gradient defined in~\cref{eq:surrogate_transport_loss}. \cref{alg:LazyDINO_training} summarizes the \texttt{LazyDINO} training procedure when using first order methods.

\begin{algorithm}[h]
\DontPrintSemicolon
  \caption{\texttt{LazyDINO}: train transport map w/surrogate objective function in latent space $\RedCoordParamSpace$}
  \label{alg:LazyDINO_training}
  \SetAlgoLined

  \KwIn{
    \begin{itemize}
       \vspace{-0.25 cm}
      \item [(i)] whitened latent prior sampler: $\pi$
          \vspace{-0.25 cm}
      \item [(ii)] single-sample surrogate latent space rKL objective for observation $\Obs$: $\widetilde{\mathcal{L}}^{\Obs}_{1,r}(\cdot, \cdot\;; \boldsymbol{w}^*)$
      
            \vspace{-0.25 cm}
      \item [(iii)] untrained transport map with random initial weights: $\Map_\MapParams: \RedCoordParamSpace\to\RedCoordParamSpace$, $\MapParams^0$
        \vspace{-0.25 cm}
      \item [(iv)] $J$ batch sizes, learning rates, \# iterations: $(B_j, a_j, I_j), j=1,\ldots,J$
      \vspace{-0.25 cm}
    \end{itemize}
  }
  \KwOut{
  \begin{itemize}
          \vspace{-0.25 cm}
  \item
    [(i)] trained transport map with pushforward density: $\Map_{\MapParams^*}$, $(\Map_{\MapParams^*})_{\sharp}\RedCoordPrior$
  \vspace{-0.25 cm}
\end{itemize}
}
  \Begin{
    \For{$j = 1,\ldots, J$}{
    \For{$i = 1,\ldots, I_j$}{
            $\left\{\refCoord^{(k)}\right\}_{k=1}^{B_j} \sim \pi$ \Comment*[r]{Sample a new stochastic batch}
           $\Delta{\MapParams}^i\gets\frac{1}{B_j}\sum_{k=1}^{B_j}\nabla_{\MapParams}\widetilde{\mathcal{L}}^{\Obs}_{1,r}(\refCoord^{(k)}, {\MapParams}^{i-1};\boldsymbol{w}^*)$\Comment*[r]{Estimate gradient of objective function}
         ${\MapParams}^i\gets$\texttt{stochastic\_gradient\_based\_iteration}$(\Delta{\MapParams}^i, ({\MapParams}^0,\ldots,{\MapParams}^{i-1}), \alpha_{j})$ \Comment*[r]{e.g., Adamax}
         }
    ${\MapParams}^0 \gets {\MapParams}^{I_j}$
    }
    $\MapParams^* \gets\MapParams^{I_{J}}$  \Comment*[r]{Last parameter is approximately optimal}
  }
\end{algorithm}

Since our latent PtO surrogate is a neural network, an MC estimate of the gradient of the rKL objective with respect to transport map parameters, i.e., $\mathbb{E}_{\refCoord \sim \RedCoordPrior}\left[\widetilde{\mathcal{L}}^{\Obs}_{1,r}(\refCoord, \MapParams;\boldsymbol{w}^*)\right]$, can be computed rapidly on GPUs, especially when the surrogate objective function gradient is a compiled batch-vectorized expression. In practice, for each example we study (in~\cref{example1},~\cref{example2}), the evaluation time of the surrogate objective function gradient is orders of magnitude less than the evaluation time of the original PtO map-dependent \texttt{LazyMap} objective gradient.

Since iterations can be performed rapidly, we use rounds of stochastic approximation-based (SA) optimization, increasing the batch sample size~\cite{bollapragada2018adaptive} and decreasing the learning rate each time for a number of iterations that is computationally tractable. We found that using such a strategy was more successful than using either small or large batch sizes alone. This strategy is similar to \emph{retrospective approximation} (RA)~\cite{newton2024retrospectiveapproximationapproachsmooth}.

%% file: appendix_pde_bayes.tex
\section{Forming the reduced Jacobian through direct and adjoint sensitivities}\label{section:bayes_pdes}

While our discussion in the main body targets general BIPs, an important class is BIPs constrained by PDE models. This section mentions a few points regarding this class of problems. We consider a nonlinear variational residual problem involving an additional state variable $u \in \scrU$. The PtO map can then be written abstractly as
\begin{equation}
    \Forward(\FullParam):\FullParam \mapsto u\mapsto \boldsymbol{\calO}(u) \qquad \text{such that}\qquad \calR(u,\FullParam)= 0 \in \scrU',
\end{equation}
where $\boldsymbol{\calO}:\scrU\to\R^{\ObsDim}$ is an observation operator, $\calR:\scrU\times\scrM\to\scrU'$ is the residual of the PDE model, and $\scrU'$ is the topological dual of $\scrU$. 

Derivatives of $\Forward$ with respect to $\FullParam$ require implicit differentiation through the residual equation. They are well-defined if the conditions of the implicit function theorem are met (e.g., isolated solution, regular branch of solutions, stability, and sufficient resolution of the discretization of the PDEs).

At a given sample point, $\RedCoordDim$ derivative action or $\ObsDim$ derivative adjoint actions are required to form the latent Jacobian in \cref{eq:reduced_jacobian}. In particular, the derivative can be expressed as
\begin{equation}
    D \Forward(m) = - D\boldsymbol{\mathcal{O}}(u)\left[\partial_u \calR(u,m)\right]^{-1}\partial_m\calR(u,m).
\end{equation}
The dominant cost is in the inverse actions of $\partial_u \calR(u,m)$ or $\partial_u \calR(u,m)^*$ on the parameter or observable basis, where each action requires solving a linear PDE. This cost can be considerably reduced when sparse direct solvers are used, as one can amortize the factorization costs associated with these actions on all parameter or observable bases. The cost reduction is significant for large-scale linear PDE models and highly nonlinear PDE models.

%% file: appendix_additional_numerics.tex
\section{Defining the Laplace approximation baseline and its computation}\label{LA_background}
     The Laplace approximation has a long history rooted in the work of Laplace (1774) \cite{laplace1986memoir} and is especially important for approximate Bayesian inversion in high dimensions. Efficient estimates of the Laplace Approximation can be computed for problems that exhibit informativeness only in a parameter subspace. 

    We define the Laplace Approximation as the Gaussian distribution centered at the unique strong minimizer, the \emph{maximum-a-posteriori} (MAP) estimate $ \FullParam_{\textrm{MAP}}$, of the Onsager-Machlup functional $I_{\FullPost}:\ParamSpace  \to  \mathbb{R}^+$ of $\FullPost$, assuming it exists (see ~\cite{Kretschmann_2023}), with covariance defined as the inverse of the Hessian operator of the functional at the MAP estimate. For the posteriors considered in this work, the Onsager-Machlup functional is the sum of the potential function $\Misfit$ and the Onsager-Machlup functional of the Gaussian prior distribution,  
i.e.
    \begin{equation}
 \mu_{\text{LA}}=\mathcal{N}(\FullParam_{\text{MAP}}, \mathcal{C}_\text{LA}), \quad  \begin{cases}
        \FullParam_{\textrm{MAP}} = \operatorname*{argmin}\limits_{\FullParam \in \ParamSpace} I_{\FullPost}(\FullParam),\\I_{\FullPost}(\FullParam)=\Misfit(\FullParam) + \frac{1}{2}\|\FullParam \|_{\mathcal{C}^{-1}}^2,
        \\
        \mathcal{C}_\text{LA} = (D^2I(\FullParam_{\textrm{MAP}}))^{-1},
        \end{cases}
 \end{equation}
so long as the potential function is Lipschitz continuous. 
The Onsager-Machlup functional $I_{\FullPost}$ generalizes the commonly known \emph{negative log-posterior density} with respect to Lebesgue measure, $\log \pi^{\boldsymbol{y}}$, to posterior probability distributions, see e.g. \cite{Kretschmann_2023} for more.
 
In our numerical examples, we use an efficient Inexact Newton--Conjugate Gradients numerical optimization algorithm~\cite{Inexact_newton, eisenstat_walker} to find the MAP estimate, $\FullParam_{\text{MAP}}$ which converged usually within $O(10)-O(100)$ inexact Newton iterations. This is conservatively estimated to be equivalent in cost to 100 evaluations of the PtO map in the results section.

%% file: map-diagnostics.tex
\section{On estimating density-based diagnostics}\label{derivation_of_KL_divergences}
The key to computing density-based diagnostics is to evaluate the Radon--Nikodym derivative between the posterior approximation of interest and the prior, i.e., the approximate likelihood evaluations. Here, we provide the formula for this Radon--Nikodym derivative for Laplace approximation and transport map pushforward distributions.
\subsection{Laplace approximation formulae}
We consider the following decomposition of the LA covariance
\begin{equation*}
           \calC_{\text{LA}}=\mathcal{C} - \mathcal{D}_{\text{LA}}\left(\frac{\lambda_j}{\lambda_j+1}\delta_{jk}\right)\mathcal{E}_{\text{LA}}\mathcal{C},\quad
           \calC_{\text{LA}}^{-1} = \mathcal{C}^{-1} + \mathcal{C}^{-1} \mathcal{D}_{\text{LA}}(\lambda_j\delta_{jk})\mathcal{E}_{{\text{LA}}},
\end{equation*}
where $\mathcal{D}_{\text{LA}}$ and $\mathcal{E}_{\text{LA}}$ are the linear encoder and decoder based on the eigendecomposition of the prior-preconditioned Hessian of the potential at the MAP point $m_{\text{MAP}}$, and $\boldsymbol{\Lambda}_{\text{LA}} = \lambda_j\delta_{ij}$ is a diagonal matrix consisting of eigenvalues. The Radon--Nikodym derivative between  $\mu_{\text{LA}}$ and the prior $\mu$ is given by
\begin{equation}\label{eq:la_density}
\begin{aligned}
    \frac{\meas\mu_{\text{LA}}}{\meas\mu}(\FullParam) &= \frac{\meas\mu_{\text{LA}}}{\meas\mathcal{N}(0, \calC_{\text{LA}})}\times \frac{\meas\mathcal{N}(0, \calC_{\text{LA}})}{\meas\mu}\\
    &= \exp\Big(- \frac{1}{2} \norm{ m_{\text{MAP}}}_{\PriorCov^{-1}}^2 - \frac{1}{2}\norm{\mathcal{E}_{\text{LA}} m_{\text{MAP}}}^2_{\boldsymbol{\Lambda}_{\text{LA}}} + (\mathcal{E}_{\text{LA}}m_{\text{MAP}})^{\transpose}\boldsymbol{\Lambda}_{\text{LA}}(\mathcal{E}_{\text{LA}}m)  \\
    &\quad  + \left\langle m_{\text{MAP}}, m\right\rangle_{\PriorCov^{-1}} + \frac{1}{2}\sum_{j}\log(1+\lambda_j) - \frac{1}{2}\norm{\mathcal{E}_{\text{LA}}m}^2_{\boldsymbol{\Lambda}_{\text{LA}}} \Big).
\end{aligned}
\end{equation}
For example, the rKL between the Laplace approximation and the true posterior is given by 
\begin{align*}
    \Dkl(\mu_{\text{LA}}||\mu^{\by}) &= \mathbb{E}_{{\FullParam}\sim\mu_{\text{LA}}}\left[\log\left(\frac{\meas\mu_{\text{LA}}}{\meas\mu}(\FullParam)\frac{\meas\mu}{\meas\mu^{\by}}(\FullParam)\right)\right]\\
    &= \mathbb{E}_{{\FullParam}\sim\mu_{\text{LA}}}\left[\Phi^{\by}(\FullParam) + \log\left(\frac{\meas\mu_{\text{LA}}}{\meas\mu}(\FullParam)\right)\right] + \log \normalization^{\Obs},
\end{align*}
where the Radon--Nikodym derivative at parameters samples can be computed using \cref{eq:la_density}.

\subsection{lazy map pushforward posterior formulae}

Let $\FullMap$ be a lazy map, then we have
\begin{equation}\label{eq:lazy_map_density}
    \frac{\meas{\FullMap_{\sharp}\mu}}{\meas\mu}(\FullParam) = \frac{\Map_{\sharp}\pi(\Encoder m)}{\pi(\Encoder m)} = \frac{(\RedCoordPrior \circ \Map^{-1})(\Encoder m)  |\det\nabla \Map^{-1}(\Encoder m)|}{\pi(\Encoder m)},
\end{equation}
where $\Map$ is the latent space transport map. For example, the rKL between the lazy map pushforward and the posterior is given by:
\begin{align*}
    D_{\text{KL}}({\FullMap_{\sharp}\mu}||\mu^{\by}) &= \mathbb{E}_{{\FullParam}\sim{\FullMap_{\sharp}\mu}}\left[\log\left(\frac{\meas{\FullMap_{\sharp}\mu}}{\meas\mu}(\FullParam)\frac{\meas\mu}{\meas\mu^{\by}}(\FullParam)\right)\right]\\
    &= \mathbb{E}_{{\FullParam}\sim{\FullMap_{\sharp}\mu}}\left[\Phi^{\by}(\FullParam) + \log\left(\frac{\meas{\FullMap_{\sharp}\mu}}{\meas\mu}(\FullParam)\right)\right] + \log \normalization^{\by},
\end{align*}
where the Radon--Nikodym derivative at parameters samples can be computed using \cref{eq:lazy_map_density}.

%% file: appendix_numerical_results.tex
\section{Conventional Neural Operator training details}\label{convention_no_training_details}
For Example I, to train the neural operator using the conventional $L^2_\mu$ objective, the learning rate $\alpha^j$ and epoch number $E^j$ for each training data size $N^j$, reported here as ($\alpha^j$, $E^j$, $N^j$) are $\Big\{(2\times10^{-4}, 1500, 125),(2\times10^{-4}, 1500, 250),\ldots(2\times10^{-4}, 1500, 8k), (2\times10^{-4}, 500, 16k)\Big\}$. 
For Example II, we used $\Big\{((1\times10^{-4}, 3000, 125),(1\times10^{-4}, 3000, 250),\ldots(1\times10^{-4}, 3000, 4k), (1\times10^{-4}, 5000, 8k),(5\times10^{-5}, 5000, 16k)\Big\}.$ We used cross-validation on the test set to ensure training with these parameters led to similar training and generalization errors.

\section{Additional numerical results}\label{app:additional_results}
\begin{figure}[H]
    \centering
    \begin{tabular}{c c}
        \hspace{0.05\textwidth}Example I & \hspace{0.05\textwidth}Example II \\
        \includegraphics[width=0.4\linewidth]{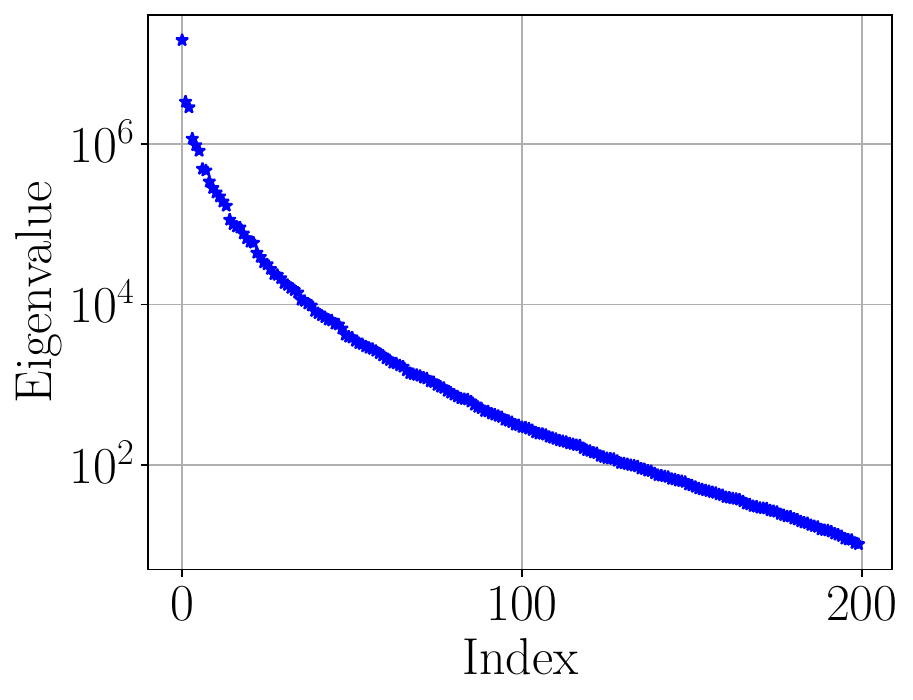}  & \includegraphics[width=0.4\linewidth]{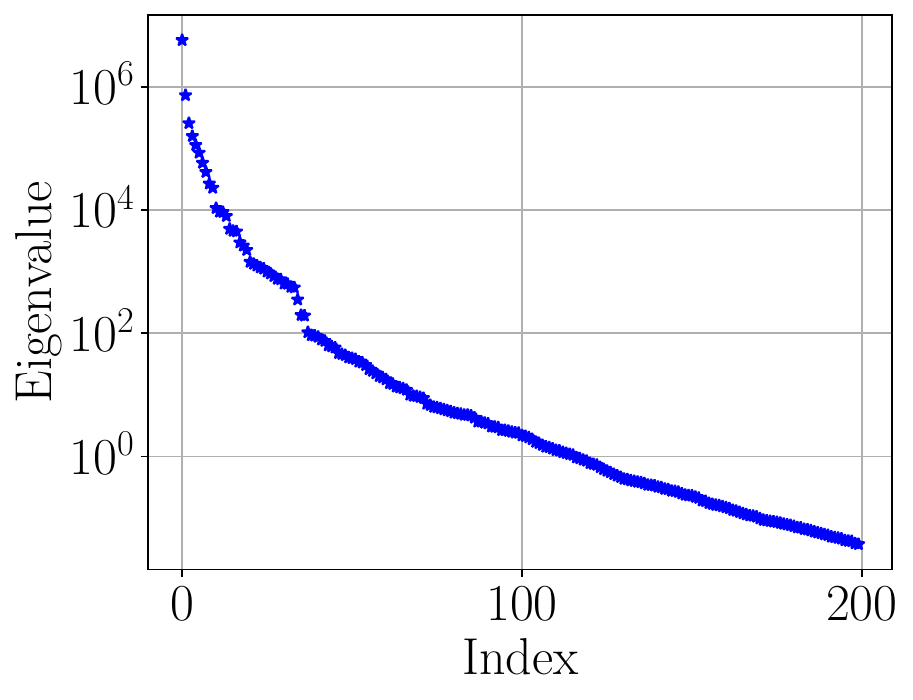} 
    \end{tabular}
    \caption{Visualization of eigenvalue decay for the generalized eigenvalue problem \cref{eq:derivative_gevp} for subspace identification in the two numerical examples.}
    \label{fig:enter-label}
\end{figure}
\begin{figure}[H]
    \centering
    \begin{tabular}{|c| c | c | c | c|} \hline
      \bf  &\bf \#1  &\bf \#5 &\bf \#25 &\bf \#125 \\\hline
       \bf \makecell{Decoder\\rows} & \raisebox{-.5\height}{\includegraphics[width = 0.15\linewidth]{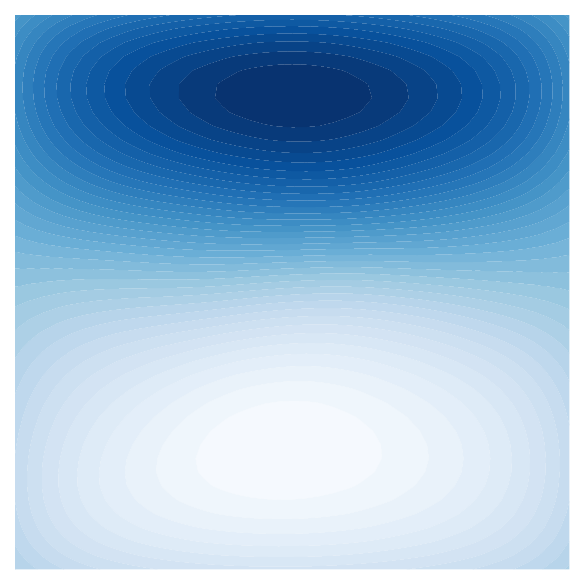}} & \raisebox{-.5\height}{\includegraphics[width = 0.15\linewidth]{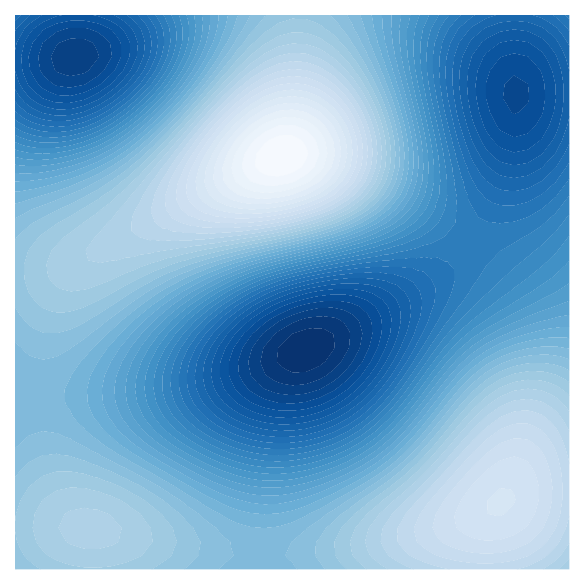}} & \raisebox{-.5\height}{\includegraphics[width = 0.15\linewidth]{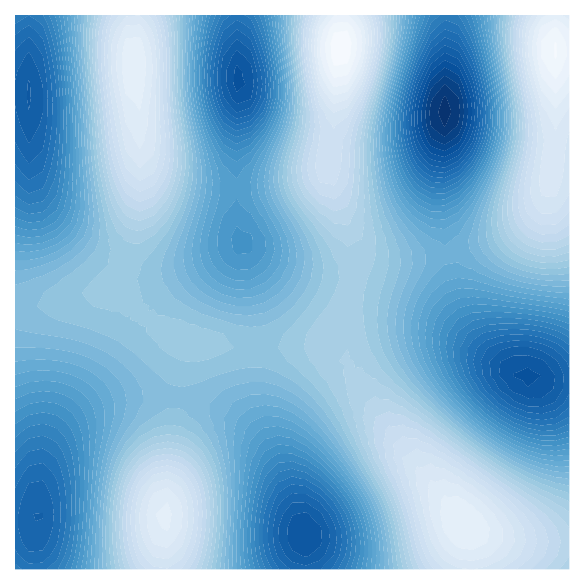}} & \raisebox{-.5\height}{\includegraphics[width = 0.15\linewidth]{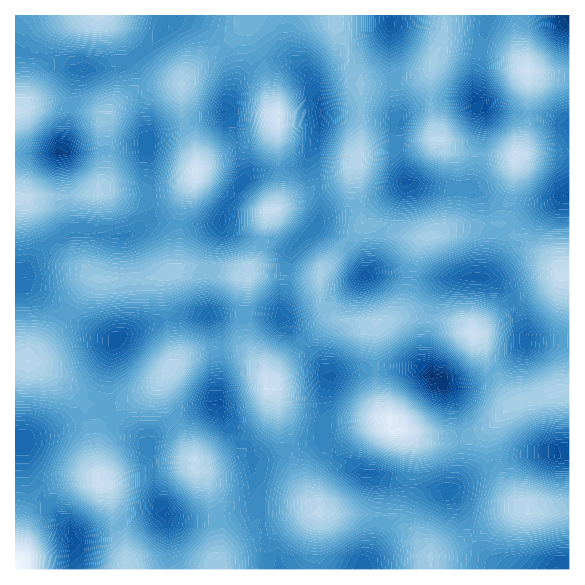}} \\\hline
       \bf \makecell{Encoder\\columns}  & \raisebox{-.5\height}{\includegraphics[width = 0.15\linewidth]{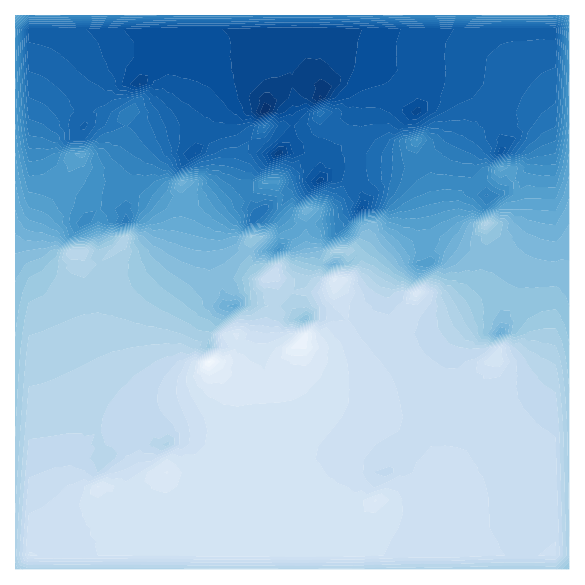}} & \raisebox{-.5\height}{\includegraphics[width = 0.15\linewidth]{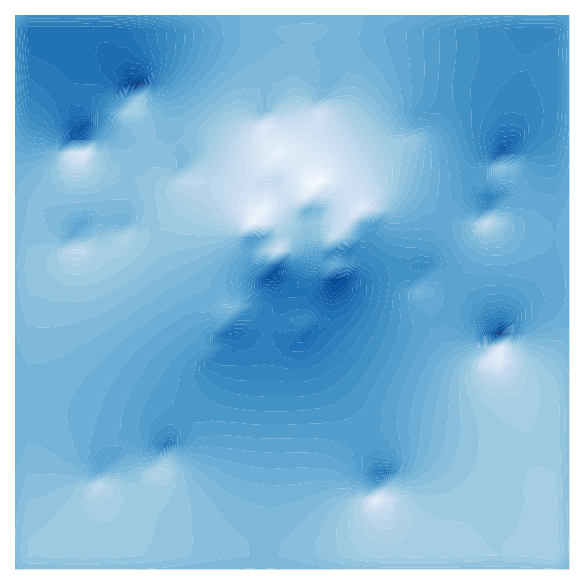}} & \raisebox{-.5\height}{\includegraphics[width = 0.15\linewidth]{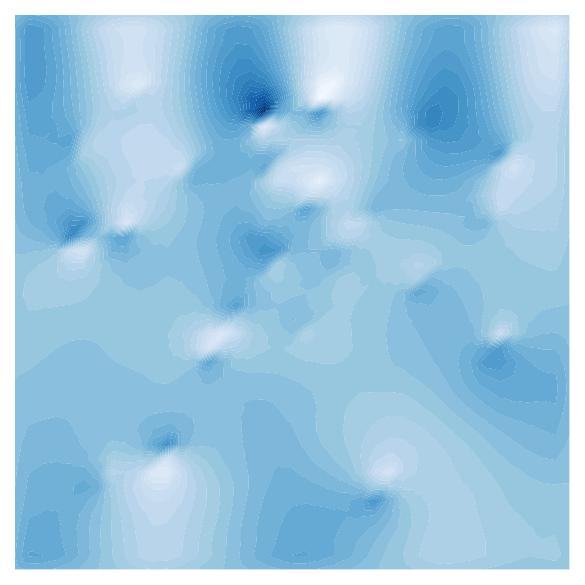}} & \raisebox{-.5\height}{\includegraphics[width = 0.15\linewidth]{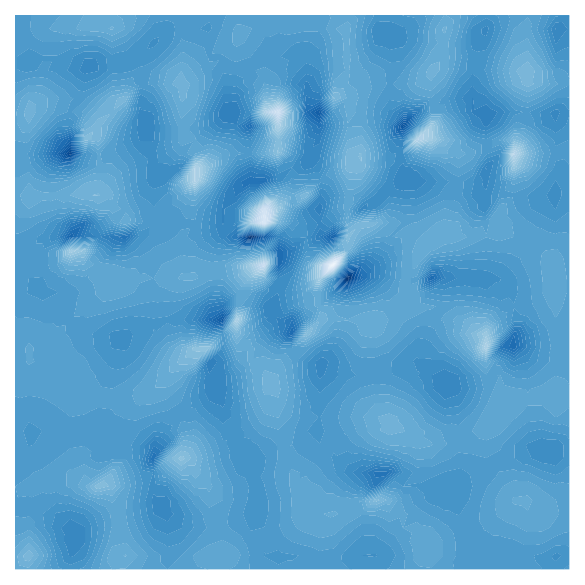}}  \\\hline
    \end{tabular}
    \caption{\textbf{Example I.} Visualization of selected decoder rows (basis functions) and encoder columns. We note that the encoder columns are computed using the action of the prior precision operator on the decoder rows. The encoder action on input is given by the vector space inner product of the encoder columns (discretized) on the input (discretized)}
    \label{fig:hyper_basis}
\end{figure}

\begin{figure}[H]
    \centering
    \begin{tabular}{|c| c | c | c | c|} \hline
      \bf  &\bf \#1  &\bf \#5 &\bf \#25 &\bf \#125 \\\hline
       \bf \makecell{Decoder\\rows} & \raisebox{-.5\height}{\includegraphics[width = 0.18\linewidth]{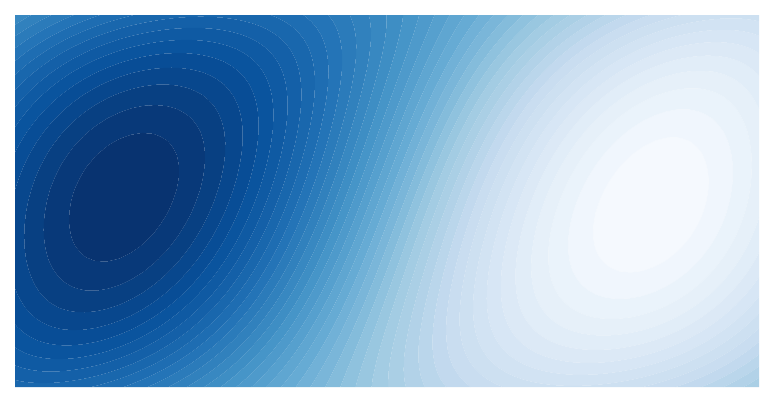}} & \raisebox{-.5\height}{\includegraphics[width = 0.18\linewidth]{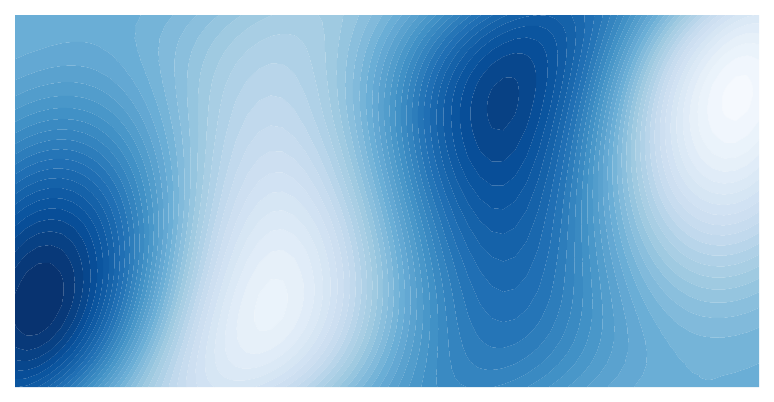}} & \raisebox{-.5\height}{\includegraphics[width = 0.18\linewidth]{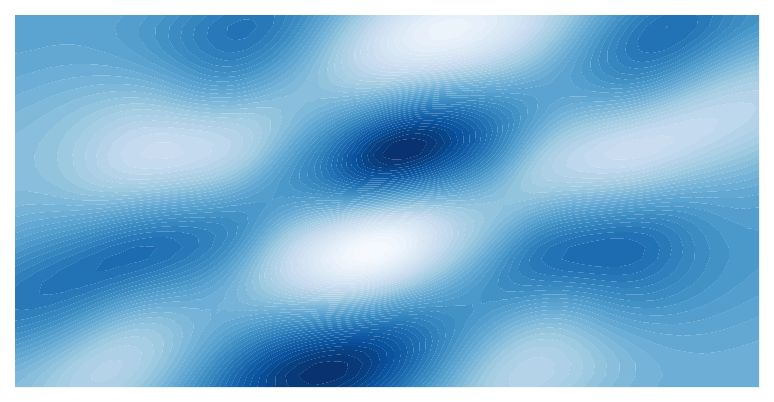}} & \raisebox{-.5\height}{\includegraphics[width = 0.18\linewidth]{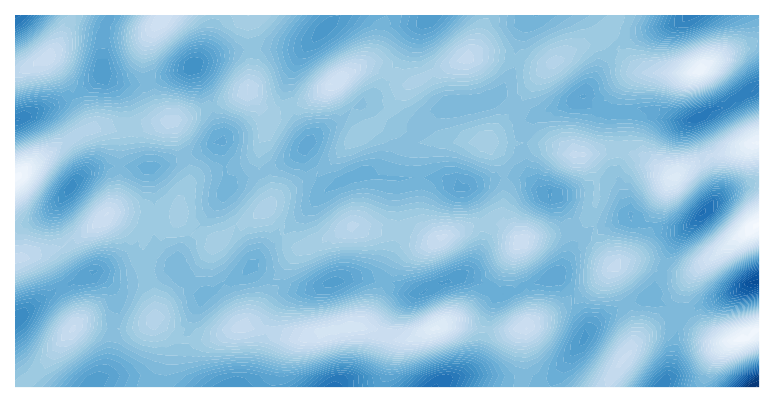}} \\\hline
       \bf \makecell{Encoder\\columns}  & \raisebox{-.5\height}{\includegraphics[width = 0.18\linewidth]{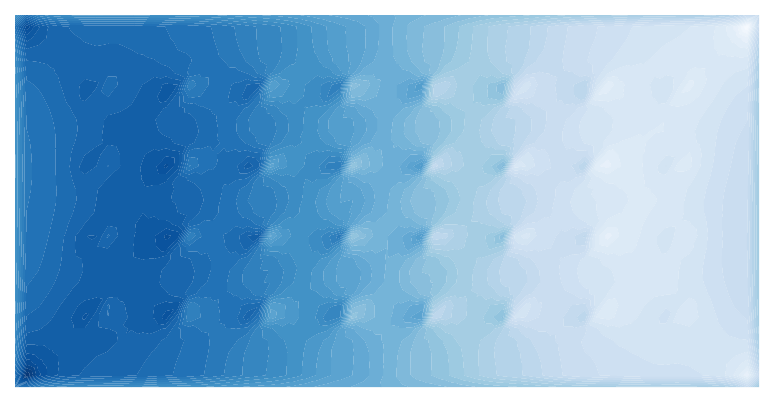}} & \raisebox{-.5\height}{\includegraphics[width = 0.18\linewidth]{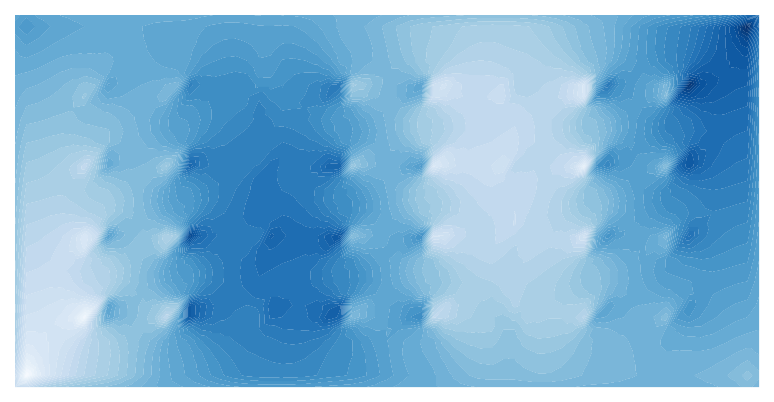}} & \raisebox{-.5\height}{\includegraphics[width = 0.18\linewidth]{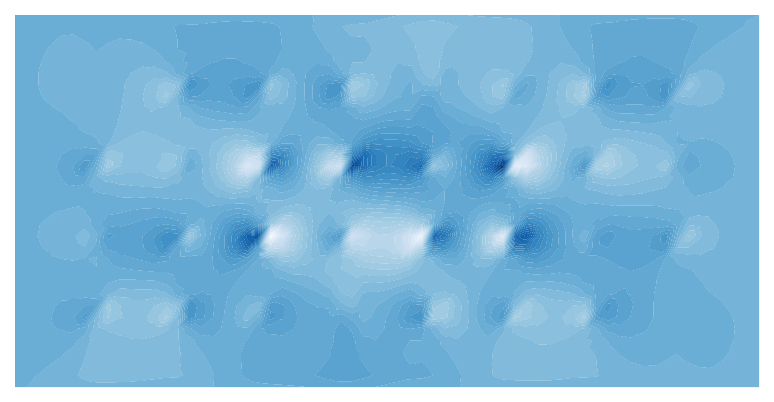}} & \raisebox{-.5\height}{\includegraphics[width = 0.18\linewidth]{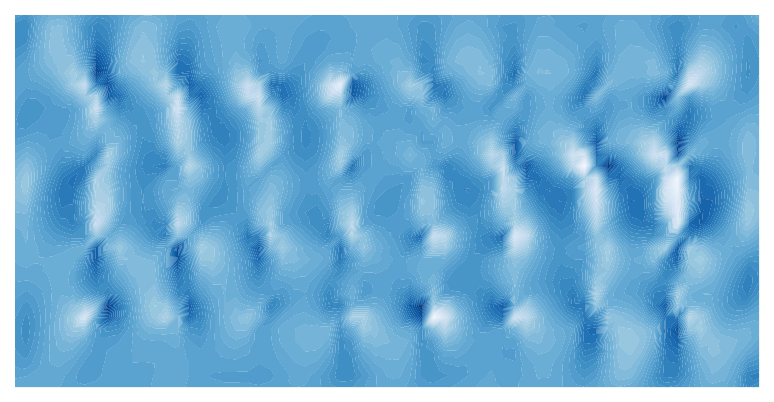}}  \\\hline
    \end{tabular}
    \caption{\textbf{Example II.} Visualization of selected decoder rows (basis functions) and encoder columns. We note that the encoder columns are computed using the action of the prior precision operator on the decoder rows. The encoder action on an input in $\scrM$ is given by the vector space inner product of the encoder columns (discretized) on the input (discretized)}
    \label{fig:hyper_basis}
\end{figure}

\begin{figure}[H]
    \centering
        {\begin{tabular}{l l l l} \includegraphics[width=0.04\textwidth]{figures/legend_dino_line.pdf}& \texttt{LazyDINO} &\includegraphics[width=0.04\textwidth]{figures/legend_mcmc_line.pdf} & True posterior via MCMC
    \end{tabular}}  
    \includegraphics[width=\linewidth]{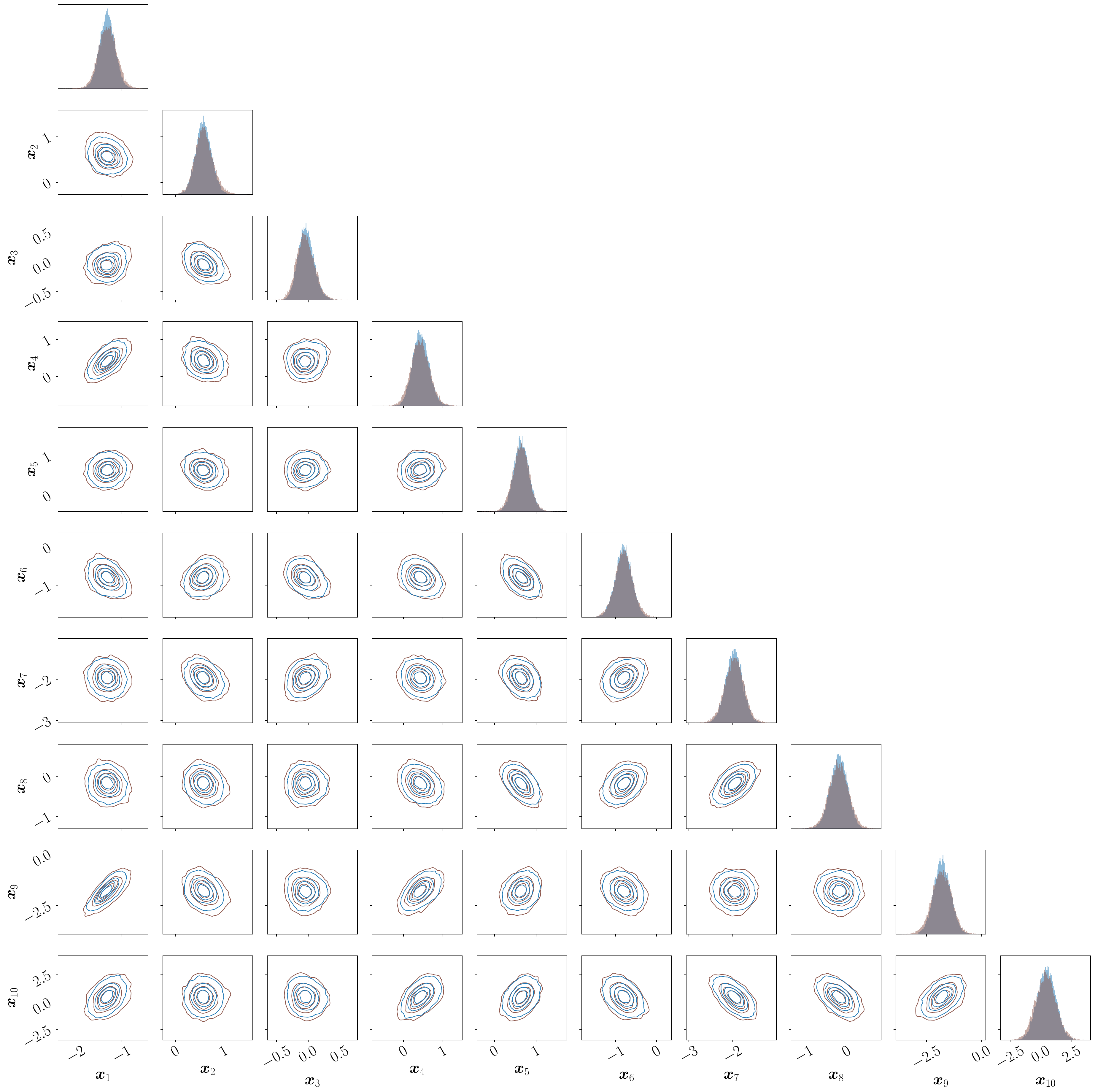}
    \caption{{\bf Example I BIP \#1} \texttt{LazyDINO} v.s. true posterior marginals at 16k \texttt{DINO} training samples in the ten leading dimensions of the latent space.}
    \label{fig:enter-label}
\end{figure}

\begin{figure}[H]
    \centering
        {\begin{tabular}{l l l l} \includegraphics[width=0.04\textwidth]{figures/legend_dino_line.pdf}& \texttt{LazyDINO} &\includegraphics[width=0.04\textwidth]{figures/legend_mcmc_line.pdf} & True posterior via MCMC
    \end{tabular}}  
    \includegraphics[width=\linewidth]{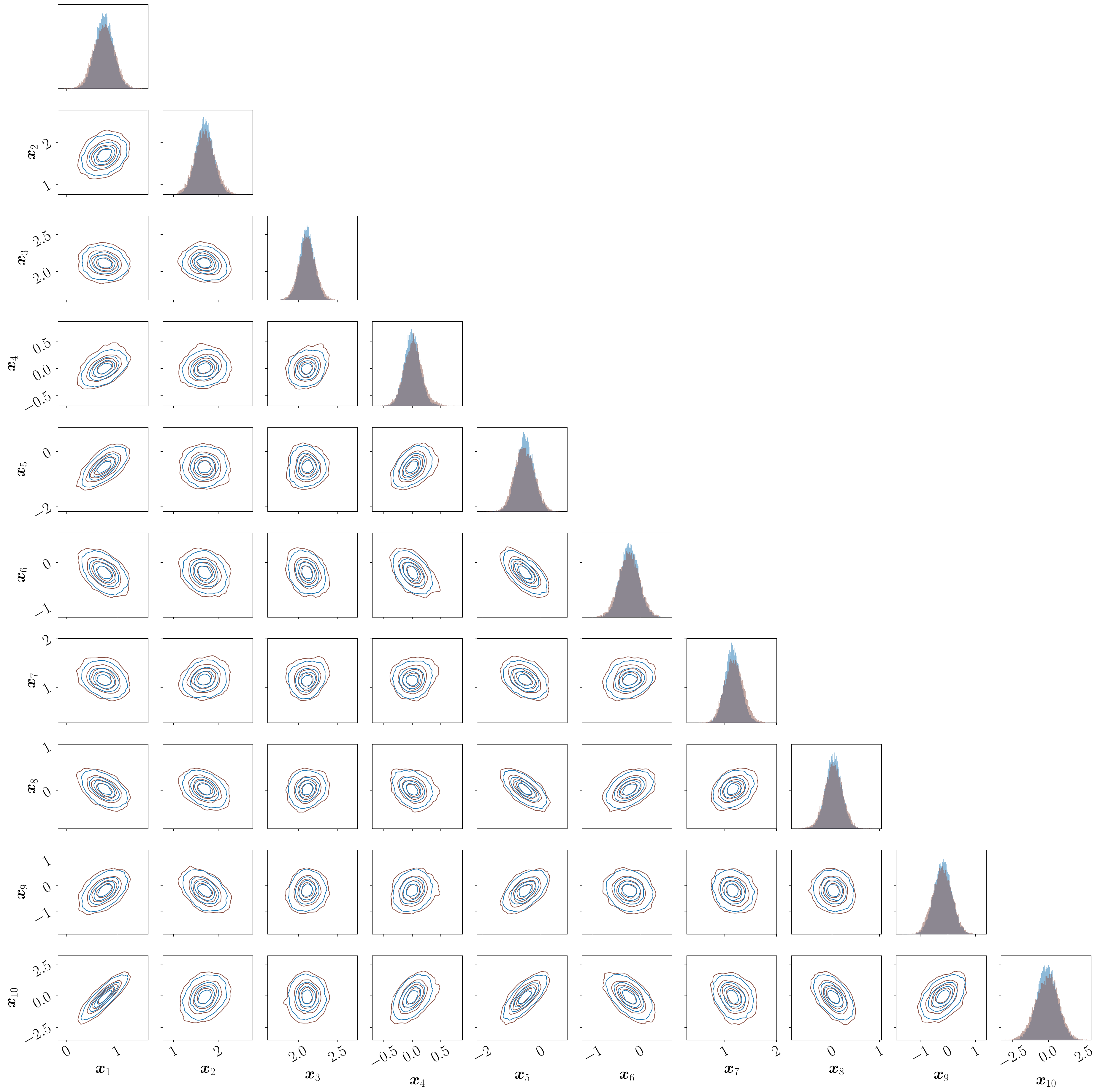}
    \caption{{\bf Example I BIP \#2} \texttt{LazyDINO} v.s. true posterior marginals at 16k \texttt{DINO} training samples in the ten leading dimensions of the latent space.}
    \label{fig:enter-label}
\end{figure}

\begin{figure}[H]
    \centering
        {\begin{tabular}{l l l l} \includegraphics[width=0.04\textwidth]{figures/legend_dino_line.pdf}& \texttt{LazyDINO} &\includegraphics[width=0.04\textwidth]{figures/legend_mcmc_line.pdf} & True posterior via MCMC
    \end{tabular}}  
    \includegraphics[width=\linewidth]{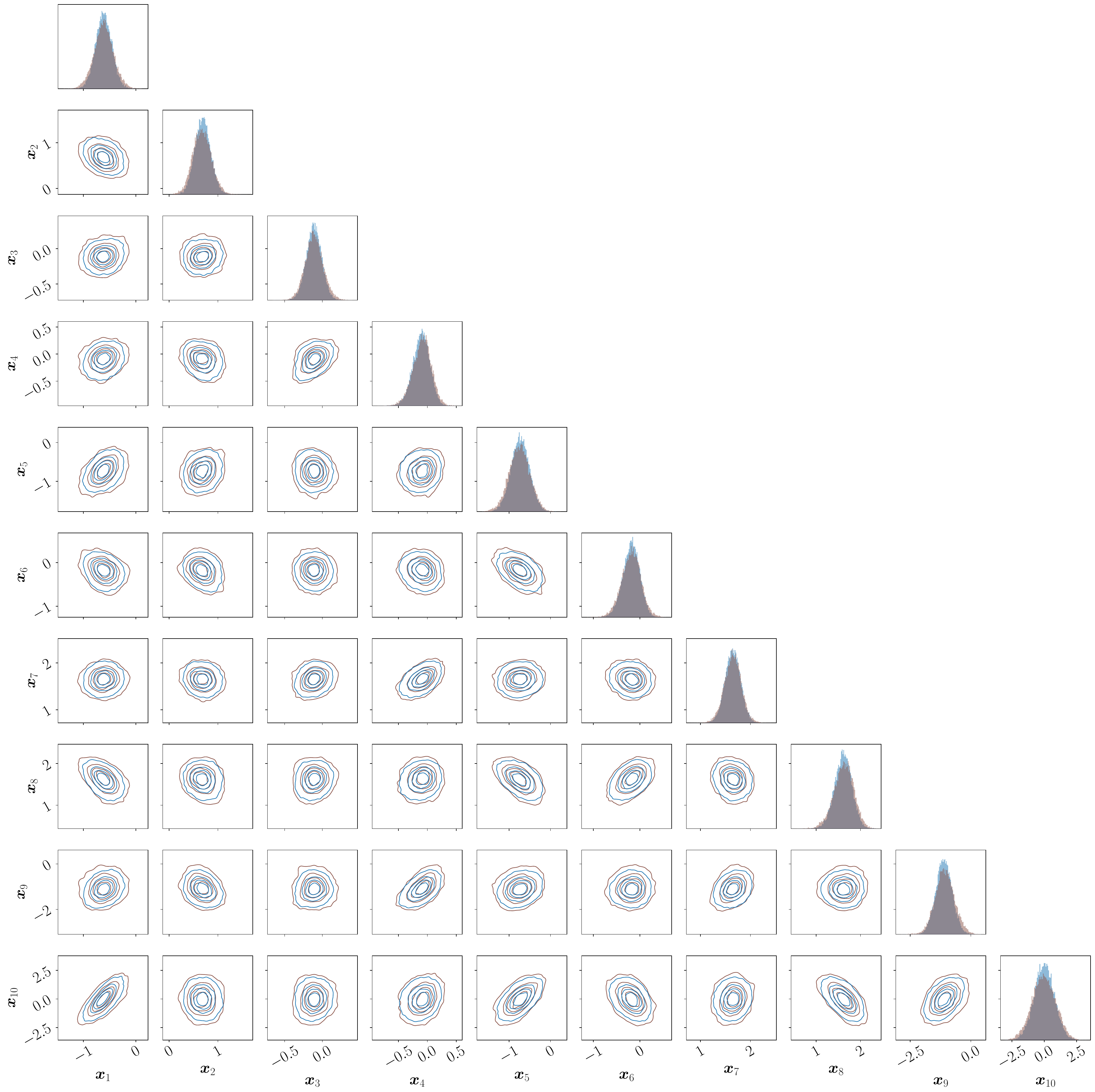}
    \caption{{\bf Example I BIP \#3} \texttt{LazyDINO} v.s. true posterior marginals at 16k \texttt{DINO} training samples in the ten leading dimensions of the latent space.}
    \label{fig:enter-label}
\end{figure}

\begin{figure}[H]
    \centering
        {\begin{tabular}{l l l l} \includegraphics[width=0.04\textwidth]{figures/legend_dino_line.pdf}& \texttt{LazyDINO} &\includegraphics[width=0.04\textwidth]{figures/legend_mcmc_line.pdf} & True posterior via MCMC
    \end{tabular}}  
    \includegraphics[width=\linewidth]{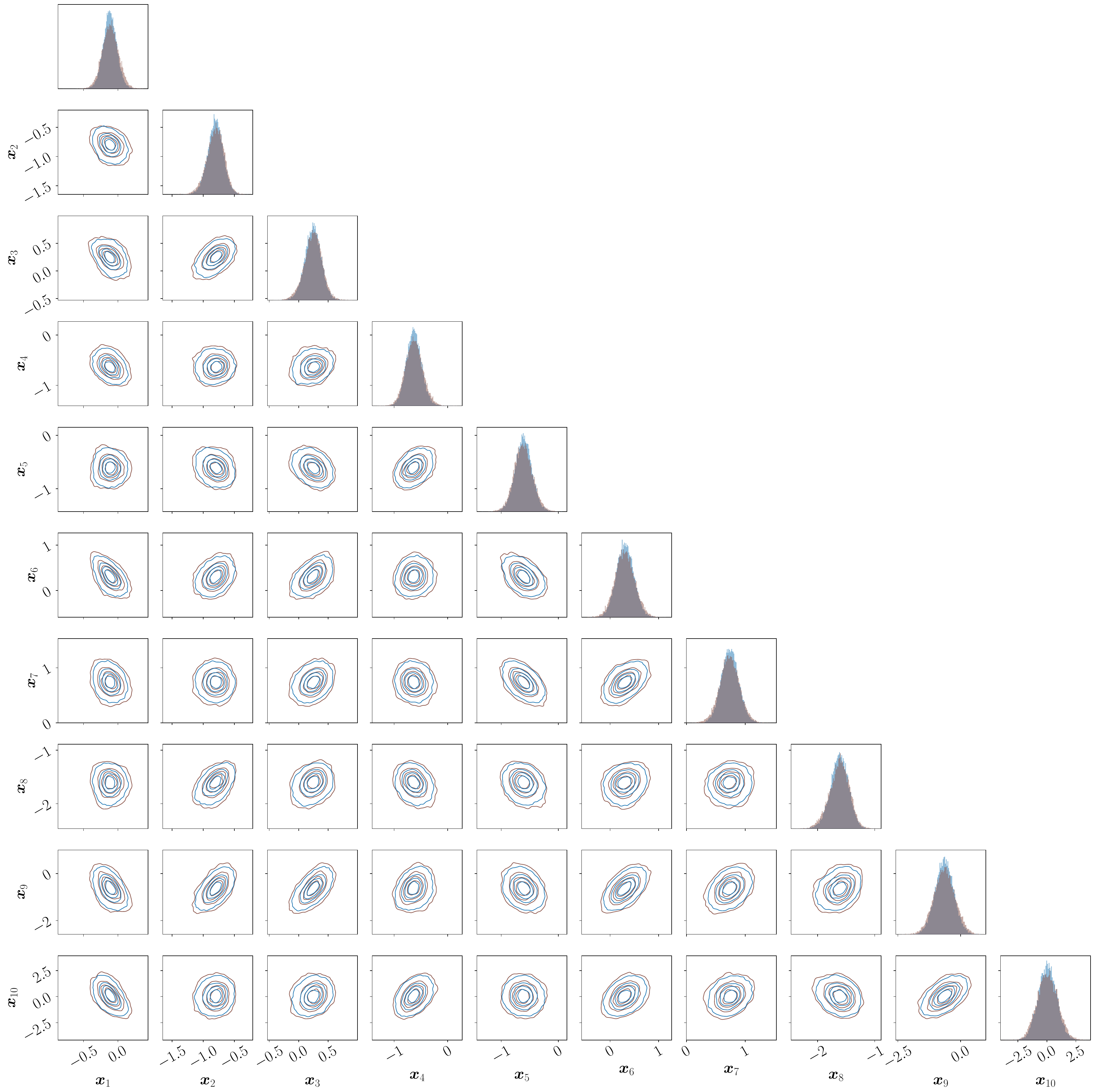}
    \caption{{\bf Example I BIP \#4} \texttt{LazyDINO} v.s. true posterior marginals at 16k \texttt{DINO} training samples in the ten leading dimensions of the latent space.}
    \label{fig:enter-label}
\end{figure}

\begin{figure}[H]
    \centering
    {\begin{tabular}{l l l l} \includegraphics[width=0.04\textwidth]{figures/legend_dino_line.pdf}& \texttt{LazyDINO} &\includegraphics[width=0.04\textwidth]{figures/legend_mcmc_line.pdf} & True posterior via MCMC
    \end{tabular}} 
    \includegraphics[width=\linewidth]{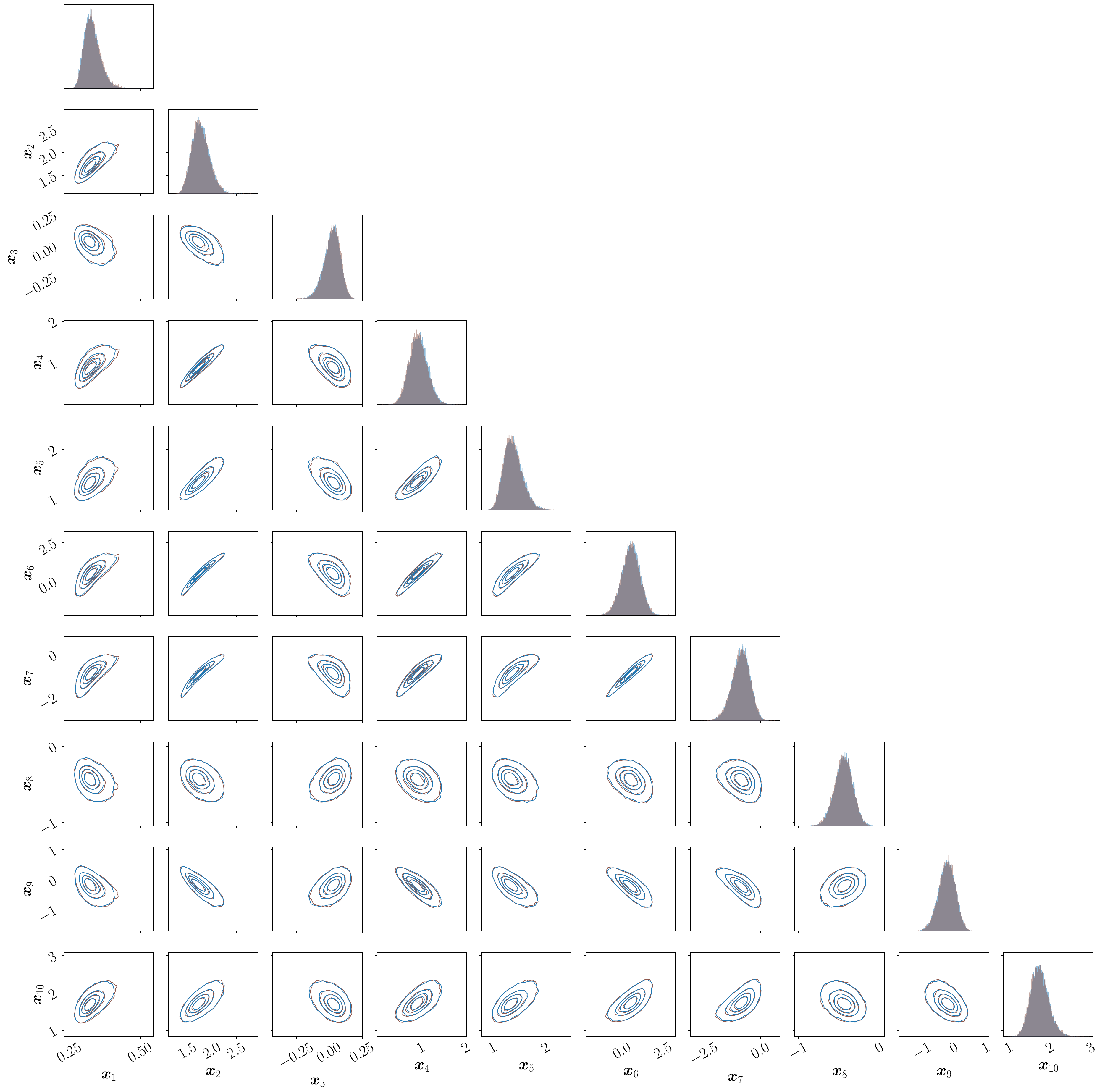}
    \caption{{\bf Example II BIP \#1 } \texttt{LazyDINO} v.s. true posterior marginals at 16k \texttt{DINO} training samples in the ten leading dimensions of the latent space.}
    \label{fig:enter-label}
\end{figure}

\begin{figure}[H]
    \centering
    {\begin{tabular}{l l l l} \includegraphics[width=0.04\textwidth]{figures/legend_dino_line.pdf}& \texttt{LazyDINO} &\includegraphics[width=0.04\textwidth]{figures/legend_mcmc_line.pdf} & True posterior via MCMC
    \end{tabular}} 
    \includegraphics[width=\linewidth]{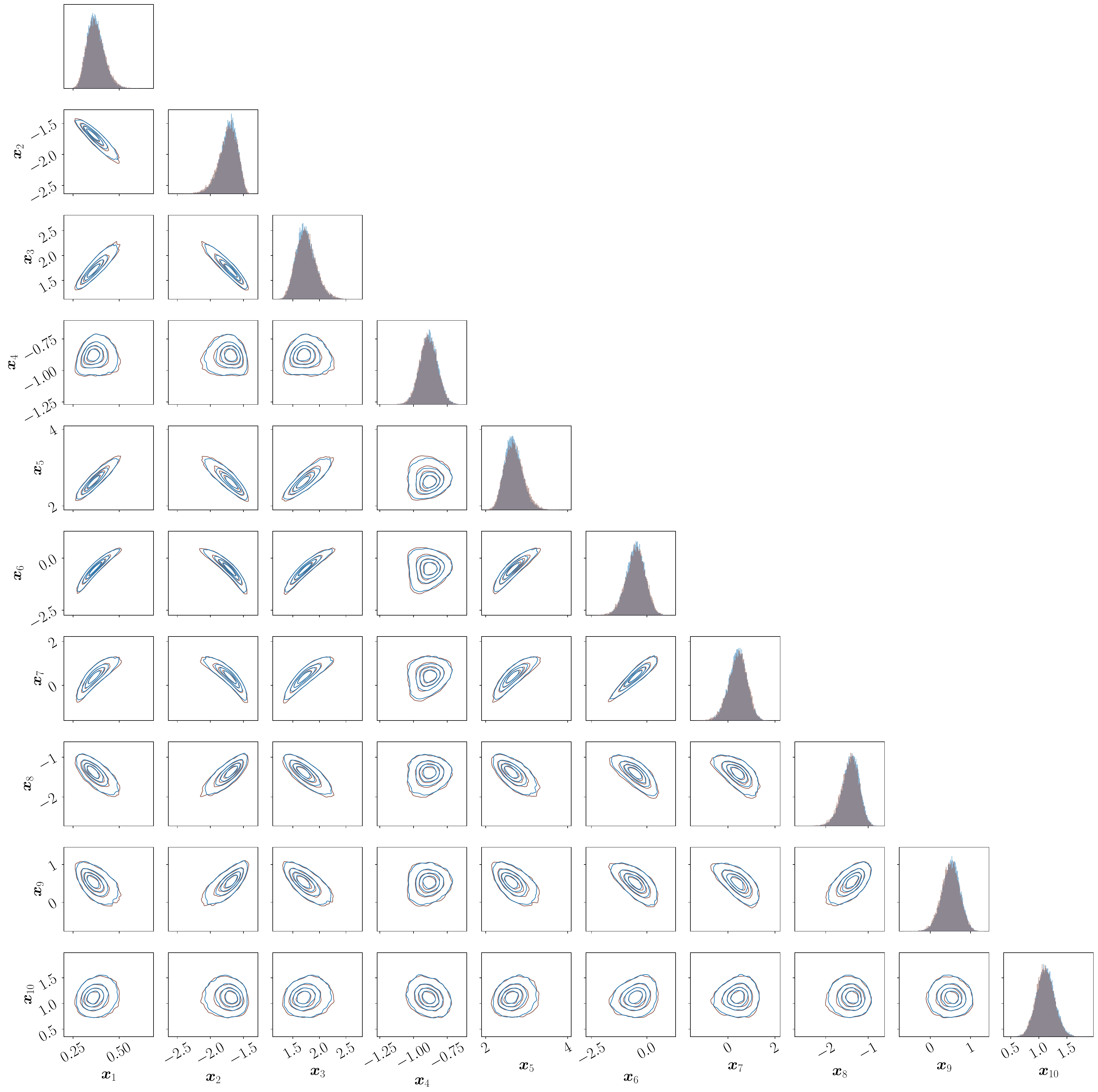}
    \caption{{\bf Example II BIP \#2 } \texttt{LazyDINO} v.s. true posterior marginals at 16k \texttt{DINO} training samples in the ten leading dimensions of the latent space in the ten leading dimensions of the latent space.}
    \label{fig:enter-label}
\end{figure}

\begin{figure}[H]
    \centering
    {\begin{tabular}{l l l l} \includegraphics[width=0.04\textwidth]{figures/legend_dino_line.pdf}& \texttt{LazyDINO} &\includegraphics[width=0.04\textwidth]{figures/legend_mcmc_line.pdf} & True posterior via MCMC
    \end{tabular}} 
    \includegraphics[width=\linewidth]{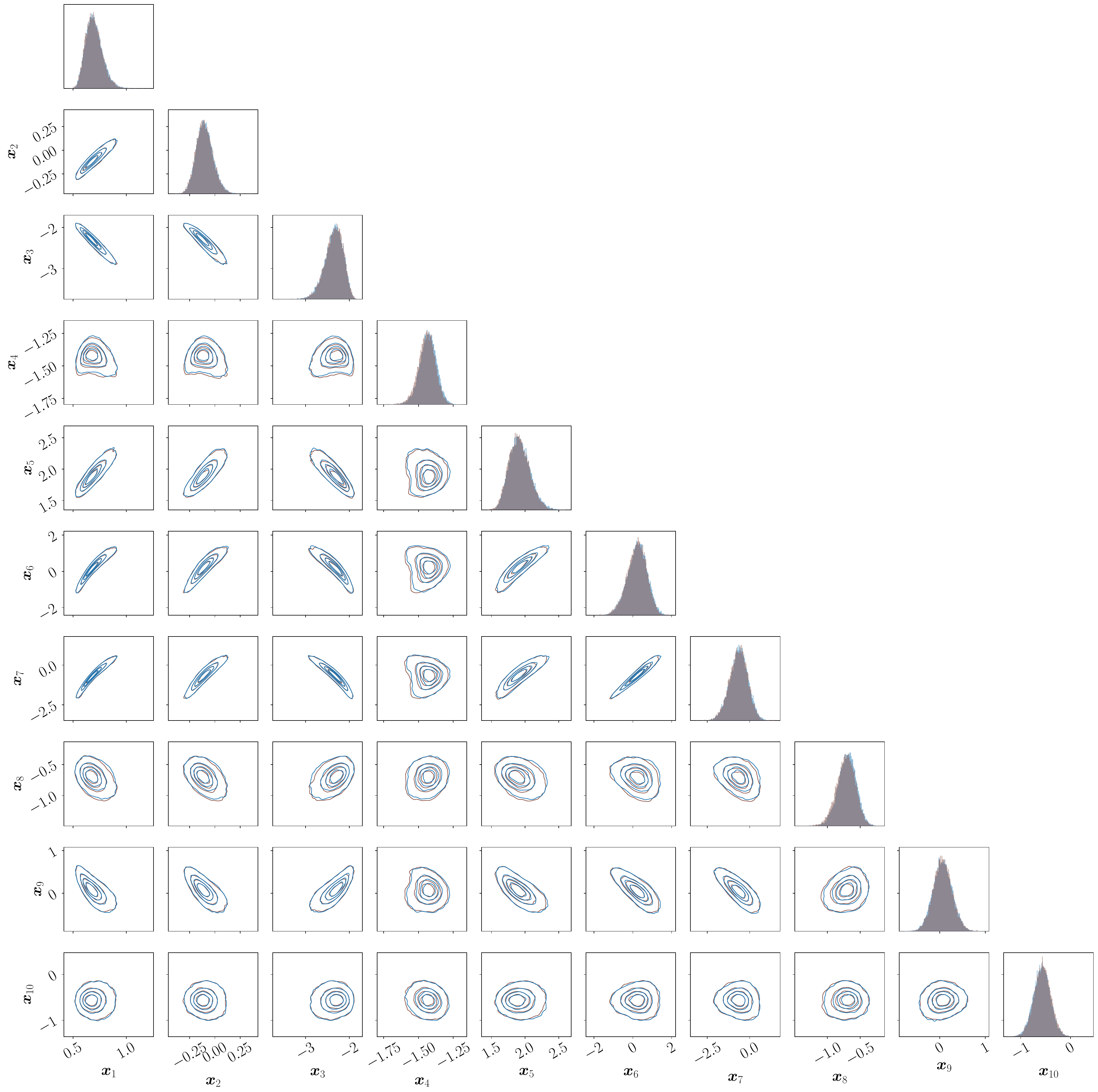}
    \caption{{\bf Example II BIP \#3} \texttt{LazyDINO} v.s. true posterior marginals at 16k \texttt{DINO} training samples in the ten leading dimensions of the latent space.}
    \label{fig:enter-label}
\end{figure}

\begin{figure}[H]
    \centering
    {\begin{tabular}{l l l l} \includegraphics[width=0.04\textwidth]{figures/legend_dino_line.pdf}& \texttt{LazyDINO} &\includegraphics[width=0.04\textwidth]{figures/legend_mcmc_line.pdf} & True posterior via MCMC
    \end{tabular}} 
    \includegraphics[width=\linewidth]{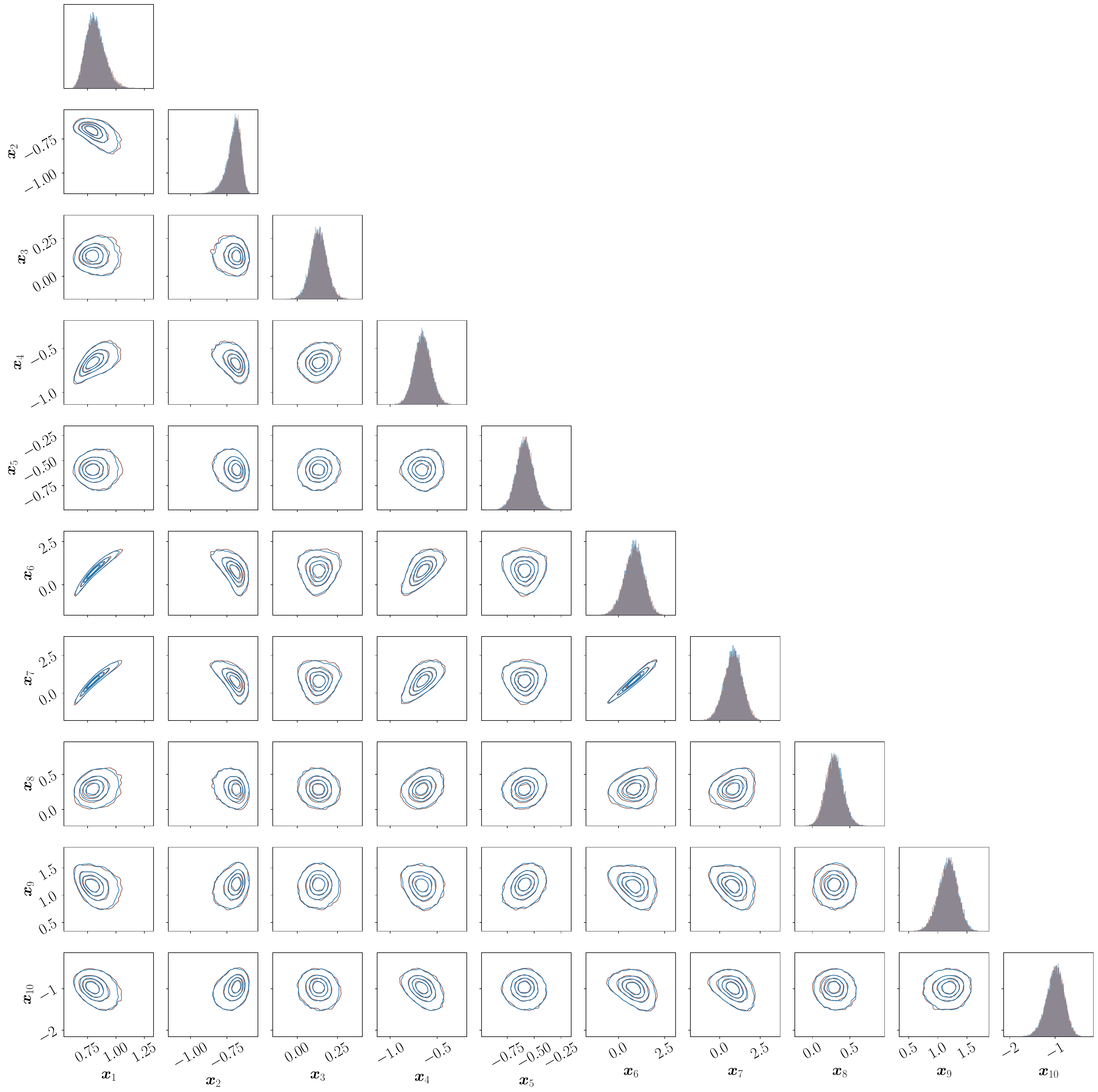}
    \caption{{\bf Example II BIP \#4} \texttt{LazyDINO} v.s. true posterior marginals at 16k \texttt{DINO} training samples in the ten leading dimensions of the latent space.}
    \label{fig:enter-label}
\end{figure}

\begin{figure}[H]
    \centering
    {\begin{tabular}{l l l l} \includegraphics[width=0.04\textwidth]{figures/legend_dino_line.pdf}& \texttt{LazyDINO} &\includegraphics[width=0.04\textwidth]{figures/legend_mcmc_line.pdf} & True posterior via MCMC
    \end{tabular}} 
    \includegraphics[width=\linewidth]{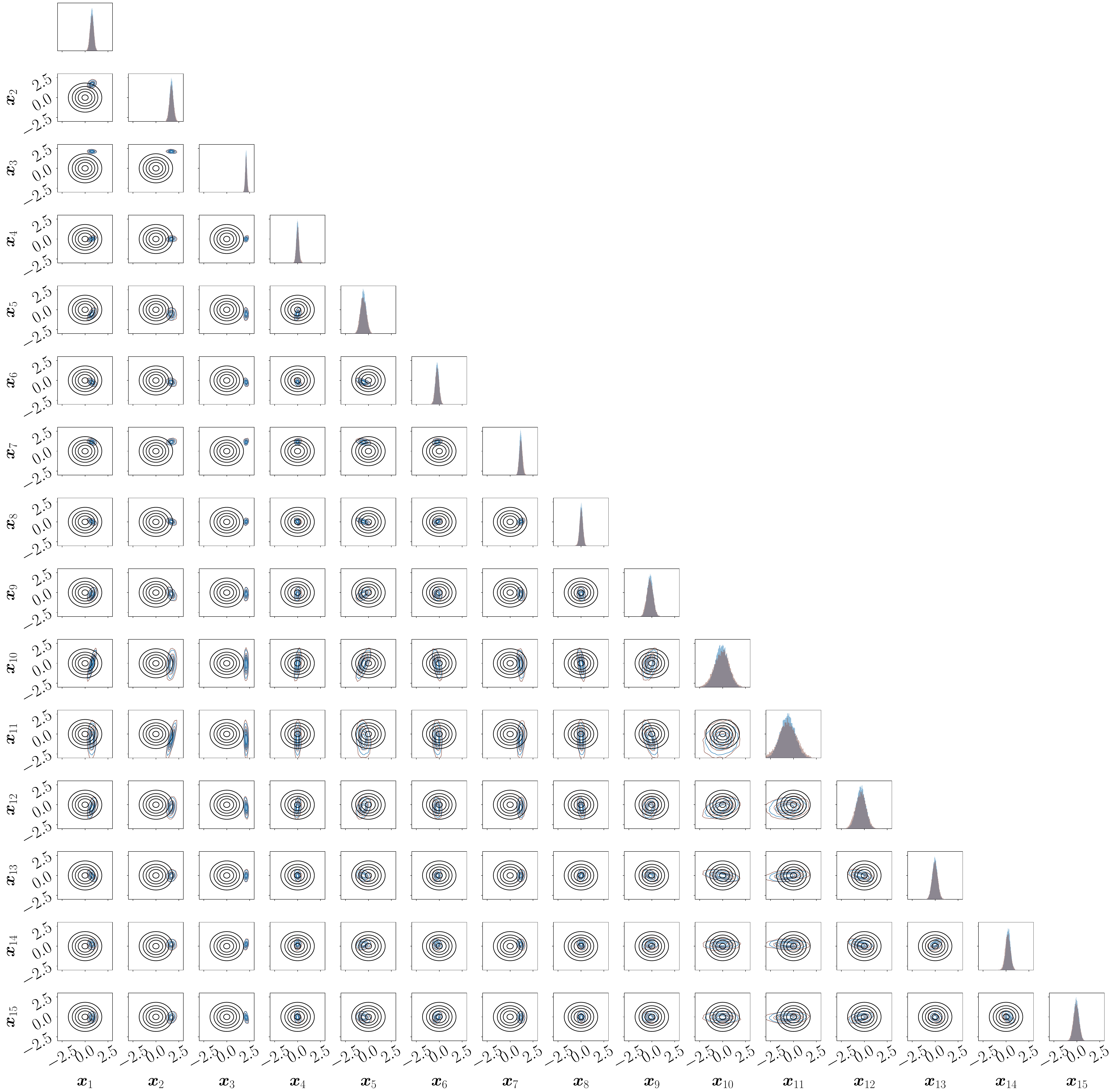}
    \caption{{\bf Example I BIP \#2} \texttt{LazyDINO} v.s. true posterior marginals compared to prior marginals (black contour lines) at 16k \texttt{DINO} training samples in the fifteen leading dimensions of the latent space. The posterior exhibits strong concentration relative to the prior. Impressively, the offline \texttt{DINO} surrogate training strategy is able to effectively equip Bayesian inversion even with such highly concentrated posteriors. }
    \label{fig:enter-label}
\end{figure}

\begin{figure}[H]
    \centering
    {\begin{tabular}{l l l l} \includegraphics[width=0.04\textwidth]{figures/legend_dino_line.pdf}& \texttt{LazyDINO} &\includegraphics[width=0.04\textwidth]{figures/legend_mcmc_line.pdf} & True posterior via MCMC
    \end{tabular}} 
    \includegraphics[width=\linewidth]{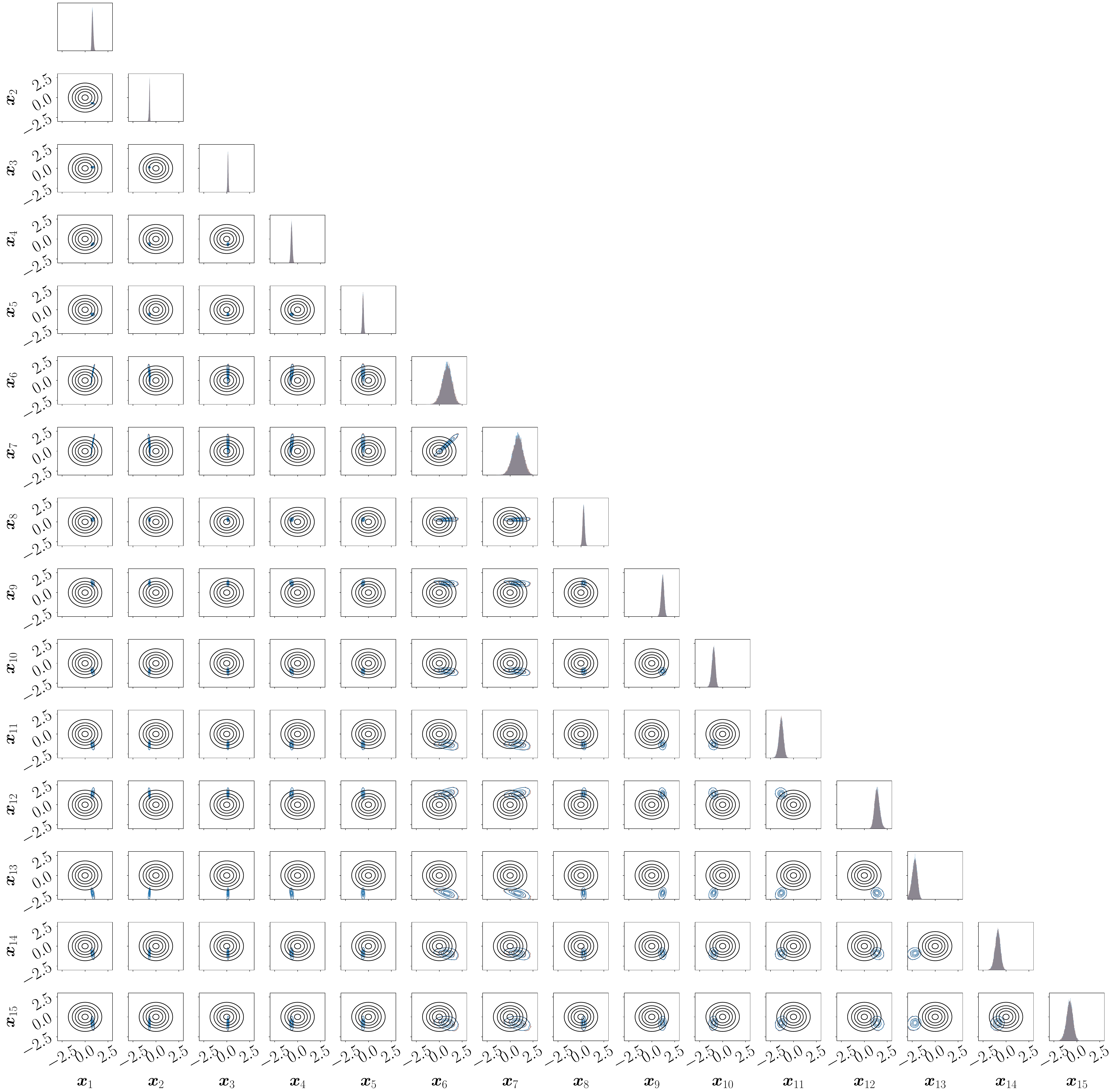}
    \caption{{\bf Example II BIP \#4} \texttt{LazyDINO} v.s. true posterior marginals compared to prior marginals at 16k \texttt{DINO} training samples in the fifteen leading dimensions of the latent space.}
    \label{fig:enter-label}
\end{figure}